\documentclass[letterpaper,11pt,twosides,onecolumn,final]{report}

\usepackage{multirow}

\usepackage[table]{xcolor}

\usepackage{stmaryrd}
\usepackage{amsfonts,amsmath,amssymb}

\usepackage{color}
\usepackage[table]{xcolor}

\usepackage{rotating}
\usepackage{mathabx}

\usepackage[sort&compress,square,comma,authoryear]{natbib}

\usepackage{relsize}

\usepackage[vlined,ruled,commentsnumbered]{algorithm2e}

\definecolor{aqua}{rgb}{0 , 0.4118, 0.5686}
\definecolor{violet}{rgb}{0.4706 ,   0.3137 ,   0.6667}
\definecolor{violetdark}{rgb}{0.3686  ,  0.2353   , 0.5216}
\definecolor{milkblue}{rgb}{0.2,0.2,0.7}
\definecolor{redaccent}{rgb}{0.3882  ,  0.1412  ,  0.1373}
\definecolor{bluedark}{rgb}{0.2,0.3,0.6}
\definecolor{QuartzPurple}{rgb}{0.302,	0.235,	0.573}
\definecolor{SapphireBlue}{rgb}{0.000,	0.388,	0.671}
\definecolor{DarkBlue}{rgb}{0.090,	0.212,	0.365}
\definecolor{BrightBlue}{rgb}{0.9137 ,   0.9686,    0.9922}

\colorlet{algshading_color}{BrightBlue}
\renewcommand\algcolor  {\color{black}}
\renewcommand\algcaptioncolor {\color{QuartzPurple}}

\usepackage{algorithmic}
\usepackage{booktabs}

\usepackage{url}

\usepackage{colortbl}

\definecolor{shadingcolor}{rgb} {.87,    0.92,    .98}
\newcommand{\shadingbox}[1]{   
    \fboxsep 0pt
    \colorbox{shadingcolor}{
        {\hskip -2pt #1}\hskip -2.5pt
    }
}

\newcommand{\minitab}[2][l]{\begin{tabular}{@{}#1}#2\end{tabular}}

\input{def1b.set}

\usepackage{mathpazo}

\usepackage[mathscr]{eucal}

\usepackage[normal]{subfigure}
\usepackage{blkarray}

\newtheorem{remark}{Remark}
\newtheorem{definition}{Definition}

\newcounter{example}
\newenvironment{example}
{\refstepcounter{example}\vspace{10pt}\par\noindent
\textbf{Example \theexample\ }
}
{}%

\newcommand{\condp}[2]{p\left(#1 \middle\vert #2 \right)}

\usepackage[titletoc,title]{appendix}

\usepackage[british]{babel}

\title{Tensor Networks  for  Dimensionality  Reduction and
 Large-Scale Optimizations\\
Part 2  Applications and Future Perspectives
{\footnote{Copyright  A.Cichocki {\it et al.} Please make reference to: \bf { A. Cichocki,  A.-H. Phan, Q. Zhao, N. Lee, I. Oseledets, and D.P. Mandic (2017), ``Tensor Networks for Dimensionality Reduction and Large-scale Optimization: Part 2 Applications and Future perspectives'', Foundations and Trends in Machine Learning: Vol. 9: No. 6, pp 431-673. }}} \\
\vspace{0.9cm}
A. Cichocki, \hspace{0.9cm} A-H. Phan, \\ Q. Zhao,  \hspace{0.6cm}  N. Lee, \\ I.V. Oseledets,  \hspace{0.9cm} M. Sugiyama, \\  D. Mandic}

\author{Andrzej CICHOCKI\\
Riken BSI, Japan, \\
Skolkovo Institute of Science and Technology (Skoltech), Russia  \\
Systems Research Institute, Polish Academy of Science, Poland\\
cia@brain.riken.jp\\
\\
Anh-Huy PHAN \\
Riken BSI, Japan \\
phan@brain.riken.jp
\\
\\
Qibin ZHAO \\
Riken AIP, Japan \\
zhao@brain.riken.jp\\\\
Namgil LEE \\
Kangwon National University\\
 Korea,
namgil.lee@riken.jp \\
\\
Ivan OSELEDETS\\
Skolkovo Institute of Science and Technology (SKOLTECH), and\\
Institute of Numerical Mathematics of Russian Academy of Sciences,\\
 Russia\\
i.oseledets@skolkovotech.ru\\
\\
Masashi SUGIYAMA \\
Riken AIP,  Japan\\
University of Tokyo \\
sugi@k.u-tokyo.ac.jp
\\
\\
Danilo P. MANDIC\\
Imperial College, UK\\
d.mandic@imperial.ac.uk
}

\date{}

\emergencystretch=\maxdimen
\usepackage{hyperref}

\usepackage{atbegshi}
\AtBeginDocument{\AtBeginShipoutNext{\AtBeginShipoutDiscard}}

\begin{document}

\sloppy

\pagenumbering{gobble}
\begingroup
\makeatletter
\let\ps@empty\ps@plain 
\maketitle
\endgroup


\pagenumbering{arabic}
\begin{abstract}
\thispagestyle{plain}

Part 2 of this monograph builds on the introduction to tensor networks and their operations presented in Part 1. It focuses on tensor network models for super-compressed higher-order representation of data/parameters and related cost functions, while providing an outline of their applications in machine learning and data analytics.

A particular emphasis is on  the tensor train (TT) and Hierarchical Tucker (HT) decompositions, and their physically meaningful interpretations which reflect the scalability of the tensor network approach. Through a graphical approach, we also elucidate how, by virtue of the underlying low-rank tensor approximations and sophisticated contractions of core tensors, tensor networks have the ability to perform distributed computations on otherwise prohibitively large volumes of data/parameters, thereby alleviating or even eliminating the curse of dimensionality.

The usefulness of this concept is illustrated over a number of applied areas, including generalized regression and classification (support tensor machines, canonical correlation analysis, higher order partial least squares), generalized eigenvalue decomposition, Riemannian optimization, and in the optimization of deep neural networks.

Part 1 and Part 2 of this work can be used either as stand-alone separate texts, or indeed as a conjoint comprehensive review of the exciting field of low-rank tensor networks and tensor decompositions.

\end{abstract}
\setcounter{page}{2}

\chapter{Tensorization  and Structured Tensors}
\label{Chapter1}
\sectionmark{Tensorization}

The concept of {\em tensorization} refers to the generation of higher-order structured tensors from the lower-order  data formats (e.g., vectors, matrices or even low-order tensors), or the representation of very large scale system parameters in low-rank tensor formats.
This is an essential step prior to multiway data analysis, unless the data itself is already collected in a multiway format; examples include  color image sequences where the R, G and B frames are stacked into a 3rd-order tensor, or multichannel EEG signals combined into a tensor with modes, e.g.,  channel $\times$ time $\times$ epoch.
For any given original  data format, the tensorization procedure may affect the choice and performance of a tensor decomposition in the next stage.

Entries of the so constructed tensor can be obtained through: $i$) a particular rearrangement, e.g., reshaping  of the original data to a tensor, $ii$) alignment of data blocks or epochs, e.g., slices of a third-order tensor are epochs of multi-channel EEG signals, or $iii$) data augmentation through, e.g., Toeplitz and Hankel matrices/tensors.
In addition, tensorization of fibers of a lower-order tensor will yield a tensor of higher order. A tensor can also be generated using transform-domain methods, for example, by a time-frequency transformation via the short time Fourier transform or wavelet transform. The latter procedure is most common for multichannel data, such as EEG, where, e.g., $S$ channels of EEG are recorded over $T$ time samples, to produce $S$ matrices of $F \times T$ dimensional time-frequency spectrograms stacked together into  an  $F \times T \times S$ dimensional third-order tensor. A tensor can also represent the data at multi-scale and orientation levels by using, e.g., the Gabor, countourlet, or pyramid steerable transformations.
When exploiting statistical independence of latent variables, tensors can be generated by means of higher-order statistics (cumulants) or by partial derivatives of the Generalised Characteristic Functions (GCF) of the observations. Such tensors are usually partially or fully  symmetric, and their entries represent mutual interaction between latent variables. This kind of tensorization is commonly used in ICA, BSS and blind identification of a mixing matrix. In a similar way, a symmetric tensor can be generated through measures of distances between observed entities, or their information exchange. For example, a third-order tensor, created to analyse common structures spread over EEG channels, can comprise distance matrices of pair-wise correlation or other metrics, such as causality over trials. A symmetric third-order tensor can involve three-way similarities. For such a tensorization, symmetric tensor decompositions with nonnegativity constraints are particularly well-suited.

Tensorization can also be performed through a suitable representation of the estimated parameters in some low-rank tensor network formats.
This method is often used when the number of estimated parameters is huge, e.g., in modelling system response in a  nonlinear system, in learning weights in a deep learning network. In this way, computation on the parameters, e.g., multiplication, convolution, inner product, Fourier transform, can be performed through core tensors of smaller scale.

One of the main motivations to develop various types of tensorization is to take advantage of data super-compression inherent in tensor network formats, especially in quantized tensor train (QTT) formats.
In general, the type of tensorization depends on a specific task in hand and the structure presented in data. The next sections introduce some common tensorization methods employed in blind source separation, harmonic retrieval, system identification, multivariate polynomial regression, and nonlinear feature extraction.

\section{Reshaping or Folding}

The simplest way of  tensorization is through the reshaping or folding operations, also known as segmentation \citep{Debals2015stoch,Bousse2015atensor}. This type of tensorization preserves the number of original data entries and their sequential ordering, as it only rearranges a vector to a matrix or tensor. Hence, folding does not require additional memory space.

{\noindent \bf Folding.} A tensor $\tY$ of size $I_1 \times I_2 \times \cdots \times I_N$ is considered a folding of a vector $\by$ of length $I_1 I_2 \cdots I_N$, if
\be
\tY(i_1,i_2,\ldots,i_N) = \by(i) \, ,
\ee
for all $1\le i_n \le I_n$, where  $i = 1 + \sum_{n = 1}^{N}  (i_n -1) \prod_{k = 1}^{n-1}I_k$ is a linear index of ($i_1, i_{2}, \ldots, i_{2}$).

In other words,  the vector $\by$ is vectorization of the tensor $\tY$, while $\tY$ is  a tensorization of $\by$.

As an  example,  the arrangement of elements  in a matrix of size $I \times L/I$, which is folded from a vector $\by$ of length $L$ is given by
\be
\bY &=& \left[
\begin{array}{cccc}
y(1) & y(I+1) & \cdots & y(L-I+1) \\
y(2) & y(I+2) &  \cdots  & y(L-I+2)\\
\vdots & \vdots  & \ddots & \vdots\\
y(I)  & y(2I) & \cdots & y(L)
\end{array}
\right]   \, .
\ee
Higher-order folding/reshaping refers to the application of the folding procedure several times, whereby a vector ${\bf y} \in \mathbb{R}^{I_1I_2\cdots I_N}$ is converted  into an $N$th-order tensor of size  $I_1\times I_2\times \cdots \times I_N$.

{\noindent \bf Application to BSS.}
It is important to notice that a higher-order folding (quantization) of a vector of length $q^N$ $(q=2,3,\ldots)$, sampled from an exponential function $y_k = az^{k-1}$, yields an $N$th-order tensor of rank 1.
 Moreover, wide classes of functions formed by products and/or sums of trigonometric, polynomial and rational functions can be quantized in this way to yield (approximate) low-rank tensor train (TT) network formats \citep{Khoromskij-SC,Khoromskij-TT,Ose2013constr}.
Exploitation of such low-rank representations allows us to separate the signals from a single or a few mixtures, as outlined below.

Consider a single mixture, $y(t)$, which is composed of $J$ component signals, $x_j(t)$, $j = 1, \ldots, J$, and corrupted by additive Gaussian noise, $n(t)$, to give
\be
	y(t) =  a_1  x_1(t) + a_2  x_2(t) + \cdots + a_J  x_J(t) + n(t).  \label{eq_mixing_model}
\ee
The aim is to extract the unknown sources (components) $x_j(t)$ from the observed signal $y(t)$.
Assume that higher-order foldings, $\tX_j$, of the component signals, $x_j(t)$, have  low-rank representations in, e.g., the  CP or Tucker format, given by
\be
	\tX_j =  \llbracket \tG_j; \bU_{j}^{(1)} , \bU_{j}^{(2)}, \ldots, \bU_{j}^{(N)} \rrbracket  \,, \notag
\ee
or in the TT format
\be
	\tX_j =  \llangle \tG_{j}^{(1)},\tG_{j}^{(2)}, \ldots, \tG_{j}^{(N)} \rrangle   \notag ,
\ee
or in any other tensor network format.
Because of the multi-linearity of this tensorization, the following relation between the tensorization of the mixture, $\tY$, and the tensorization of the hidden components, $\tX_j$, holds
\be
\tY &=& a_1  \tX_1 + a_2  \tX_2 + \cdots + a_J \tX_J + \tN \,,  \label{eq_mixing_folding}
\label{equ_btd3}
\ee
where $\tN$ is the tensorization of the noise $n(t)$.

\markright{\thesection.\quad Reshaping or Folding}

Now, by a decomposition of $\tY$ into $J$ blocks of tensor networks, each corresponding to a tensor network (TN) representation of a hidden component signal, we can find approximations of $\tX_j$ and the separate component signals up to a scaling ambiguity.
The separation method can be used in conjunction with the Toeplitz and Hankel foldings. Example~\ref{ex_TT_separation_exp} illustrates the separation of damped sinusoid signals.

\section{Tensorization through a Toeplitz/Hankel Tensor}
\sectionmark{Toeplitz and Hankel Tensors}
\subsection{Toeplitz Folding}

The Toeplitz matrix is a structured matrix with constant entries in each diagonal. Toeplitz matrices appear in many signal processing applications, e.g.,
through covariance matrices in prediction, estimation, detection, classification, regression, harmonic analysis, speech enhancement,
interference cancellation, image restoration, adaptive filtering, blind deconvolution and blind
equalization \citep{Bini1995,CIT-006}.

Before introducing a generalization  of a  Toeplitz matrix to a Toeplitz tensor, we shall first consider the discrete convolution between two vectors $\bx$ and $\by$ of respective lengths $I$ and $L > I$, given by
\be
	 \bz = \bx 	\ast \by  \label{eq_conv_ab} \,.
\ee
Now, we can write the entries $\bz_{I:L}  =  \left[z(I),  z(I+1),  \ldots, z(L)\right]^{\text{T}}$ in a linear algebraic form as
\begin{align}
	\bz_{I:L}
& =
\left[\begin{array}{cccccc}
	y(I)& y(I-1)& y(I-2) & \cdots & y(1) \\
	y(I+1)& y(I)& y(I-1) & \cdots & y(2) \\
	y(I+2)& y(I+1)& y(I) & \cdots & y(3) \\ \vdots & \vdots & \vdots &\ddots & \vdots \\
	y(L)&y(L-1) & y(L-2)&  \cdots& y(J)
	\end{array}
\right] \,
\left[\begin{array}{c} 	
x(1) \\ x(2)  \\ x(3) \\ \vdots \\ x(I)\end{array}
\right]
 \notag \\
&= \bY^{\text{T}}
 \bx   = \bY \bar{\times}_1  \bx \notag ,
\end{align}
where $J = L-I+1$. With this representation, the convolution can be computed through a linear matrix operator, $\bY$, which is called the {\emph{Toeplitz matrix}} of the generating vector $\by$.

{\noindent \bf Toeplitz matrix. }
A Toeplitz matrix of size $I \times J$, which is constructed from a vector $\by$ of length $L = I+J-1$, is defined as
\be
\bY = \calT_{I,J}(\by)  =   \left[
\begin{array}{cccc}
y(I) & y(I+1) & \cdots & y(L) \\
y(I-1) & y(I) &  \cdots  & y(L-1)\\
\vdots & \vdots  & \ddots & \vdots\\
y(1)  & y(2) & \cdots & y(L-I+1)
\end{array}
\right]\, .  \; \; \label{equ_Toeplitzation}
\ee
The first column and first row of the Toeplitz matrix represent its entire generating vector.

Indeed, all $(L+I-1)$ entries of $\by$ in the above convolution (\ref{eq_conv_ab}) can be expressed either by: (i) using a Toeplitz matrix formed from a zero-padded generating vector $[\0_{I-1}^{\text{T}}, \by^{\text{T}}, \0_{I-1}^{\text{T}}]^{\text{T}}$, with $[\by^{\text{T}}, \0_{I-1}^{\text{T}}]$ being the first row of this Toeplitz matrix, to give
\be
	\bz = \calT_{I,L+I-1}([\0_{I-1}^{\text{T}}, \by^{\text{T}}, \0_{I-1}^{\text{T}}]^{\text{T}})^{\text{T}} \, \bx \,,
\ee
or (ii) through a Toeplitz matrix of the generating vector $[\0_{L-1}^{\text{T}}, \bx^{\text{T}}, \0_{L-1}^{\text{T}}]^{\text{T}}$, to yield
\be
\bz =  \calT_{L,L+I-1}([\0_{L-1}^{\text{T}}, \bx^{\text{T}}, \0_{L-1}^{\text{T}}]^{\text{T}})^{\text{T}} \, \by \,.
\ee
The so expanded Toeplitz matrix is a circulant matrix of $[\by^{\text{T}}, \0_{I-1}^{\text{T}}]^{\text{T}}$.

Consider now a convolution of three vectors, $\bx_1$, $\bx_2$ and $\by$ of respective lengths $I_1$, $I_2$ and $(L \geq I_1 + I_2)$, given by
\be
	\bz = \bx_1 \ast \bx_2  \ast \by \, . \notag
\ee
For its implementation, we first construct a Toeplitz matrix, $\bY$, of size $I_1 \times (L-I_1 +1)$ from the generating vector $\by$. Then, we use the rows $\bY(k,:)$ to generate Toeplitz matrices, $\bY_k$ of size $I_2 \times I_3$. Finally, all $I_1$ Toeplitz matrices, $\bY_1$, \ldots, $\bY_{I_1}$, are stacked as horizontal slices of a third-order tensor $\tY$, i.e., $\tY(k,:,:) = \bY_k$, $k = 1, \ldots, I_1$. It can be verified that 
entries $[z(I_1+I_2-1), \ldots, z(L)]^T$ can be computed as
\be
\left[\begin{array}{cccc}
	z(I_1+I_2-1)  \\ \vdots \\ z(L)
	\end{array}
\right]  =  [\bx_1 \ast \bx_2 \ast \by]_{I_1+I_2-1:L} &=&  \tY \, \bar{\times}_1   \, \bx_1 \,\bar{\times}_2  \, \bx_2 \notag .
\ee
The tensor $\tY$ is referred to as the Toeplitz tensor of the generating vector $\by$.

\begin{figure}[t]
\centering
\subfigure[]{
\includegraphics[width=.6\linewidth, trim = 0.0cm 1cm 6cm 0cm,clip=true]{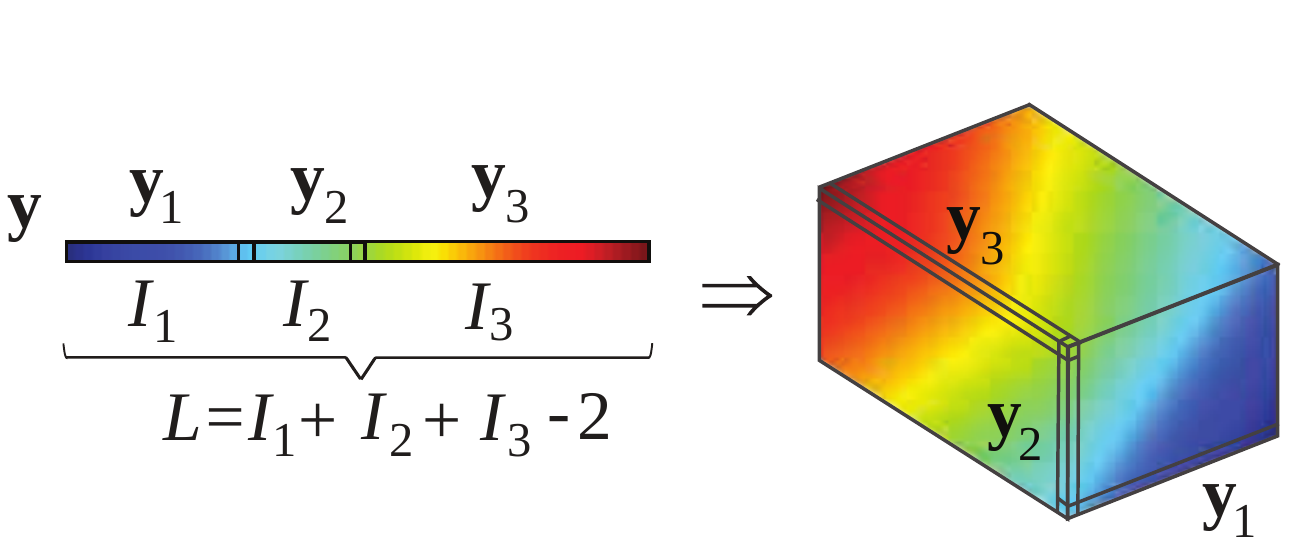}}
\subfigure[]{
\includegraphics[width=.29\linewidth, trim = 13.0cm 1.cm 0cm -.2cm,clip=true]{Toeplitzten_f2m}}
\caption{Illustration of a 3rd-order Toeplitz tensor of size $I_1 \times I_2 \times I_3$, generated from a vector $\by$ of length $L = I_1 + I_2 + I_3 -2$. (a) The highlighted fibers of the Toeplitz tensor form the generating vector $\by$. (b) The entries in every shaded diagonal intersection are identical and represent one element of $\by$.}
\label{fig_Toeplitz3}
\end{figure}

{\noindent \bf Toeplitz tensor.} An $N$th-order Toeplitz tensor of size $I_1 \times I_2 \times \cdots \times I_N$, which is represented by $\tY = \calT_{I_1, \ldots, I_N}(\by)$, is constructed from a generating vector $\by$ of length $L = I_1 + I_2 + \cdots + I_N - N+1$, such that its entries are defined as
\be
\tY(i_1, \ldots, i_{N-1}, i_N) = y(\bar{i}_1 + \cdots +  \bar{i}_{N-1} +  i_N) \, ,
\ee
where $\bar{i}_n = I_n-i_n$.
An example of the Toeplitz tensor is illustrated in Figure~\ref{fig_Toeplitz3}.
\begin{example}
Given a $3 \times 3 \times 3$ dimensional Toeplitz tensor of a sequence $1,2,\ldots,7$, the horizontal slices are Toeplitz matrices of sizes $3 \times 3$ given by
\begin{align}
\calT_{3,3,3}(1,\ldots,7)  &=
\left[
\begin{array}{cccc}
\calT_{3,3}(3,\ldots,7) \\
\calT_{3,3}(2,\ldots,6) \\
\calT_{3,3}(1,\ldots,5)
\end{array}
\right] \notag
= \left[ \begin{array}{c}
\left[
\begin{array}{cccc}
 5 &    6&     7\\
     4     &5 &    6\\
     3     &4 &    5
\end{array}
\right] \\
\left[
\begin{array}{cccc}
4 & 5 & 6\\
 3 &    4&     5\\
    2     &3 &    4\\
\end{array}
\right] \\
\left[
\begin{array}{cccc}
3 & 4 & 5\\
2 & 3 & 4\\
 1 &    2&     3\\
\end{array}
\right]
\end{array}
\right]   \, . \notag
\end{align}
\end{example}

{\noindent \bf Recursive generation.} An $N$th-order Toeplitz tensor of a generating vector $\by$ is of size $I_1 \times I_2 \times \cdots \times I_N$, can be constructed from an $(N-1)$th-order Toeplitz tensor of size $I_1 \times I_2 \times \cdots \times (I_{N-1} + I_N - 1)$ of the same generating vector, by a conversion of mode-$(N-1)$ fibers to Toeplitz matrices of size $I_{N-1} \times I_N$.

Following the definition of the Toeplitz tensor, the convolution of $(N-1)$ vectors, $\bx_n$ of respective lengths $I_n$, and a vector $\by$ of length $L$,
 can be represented as a tensor-vector product of an $N$th-order Toeplitz tensor  and vectors $\bx_n$, that is
\be
	[\bx_1 \ast \bx_2 \ast \cdots \ast \bx_{N-1} \ast \by]_{J:L} = \tY \,  \bar{\times}_1 \, \bx_1 \,\bar{\times}_2  \,\bx_2 \cdots \bar{\times}_{N-1} \, \bx_{N-1} \,, \notag
\ee
where $\tY = \calT_{I_1, \ldots, I_{N-1}, L-J}(\by)$ is a Toeplitz tensor of size $I_1 \times \cdots \times I_{N-1} \times (L-J)$ generated from $\by$, and $J = \sum_{n = 1}^{N-1} I_n - N+1$, or
\be
	\bx_1 \ast \bx_2 \ast \cdots \ast \bx_{N-1} \ast \by = \widetilde{\tY} \,  \bar{\times}_1 \, \bx_1 \,\bar{\times}_2  \,\bx_2 \cdots \bar{\times}_{N-1} \, \bx_{N-1} \,, \notag
\ee
where $\widetilde{\tY} = \calT_{I_1, \ldots,I_{N-1},L + J}([\0_{J}^{\text{T}}, \by^{\text{T}}, \0_{J}^{\text{T}}]^{\text{T}})$ is a Toeplitz tensor, of the zero-padded vector of $\by$, is of size $I_1 \times \cdots \times I_{N-1} \times (L+J)$.

\subsection{Hankel Folding}


The Hankel matrix and Hankel tensor have similar structures to the Toeplitz matrix and tensor and can also be used as linear operators in the convolution.

{\noindent \bf Hankel matrix. }
An $I \times J$ Hankel matrix of a vector $\by$, of length $L = I+J-1$, is defined as
\be
\bY = \calH_{I,J}(\by)  &=& \left[
\begin{array}{cccc}
y(1) & y(2) & \cdots & y(J) \\
y(2) & y(3) &  \cdots  & y(J+1)\\
\vdots & \vdots & \ddots & \vdots\\
y(I)  & y(I+1) & \cdots & y(L)
\end{array}
\right]   \, .
\ee

{\noindent {\bf  Hankel tensor.} \citep{NLA:NLA453}}
An $N$th-order Hankel tensor of size $I_1 \times I_2 \times \cdots \times I_N$, which is represented by $\tY = \calH_{I_1, \ldots, I_N}(\by)$,  is constructed from a generating vector $\by$ of length $L = \sum_{n} I_n - N + 1$, such that its  entries are defined as
\be
\tY(i_1, i_2, \ldots, i_N) = y(i_1 + i_2 + \cdots + i_N - N+1) \, .
\ee

\begin{remark} (Properties of a Hankel tensor)
\begin{itemize}
\item The generating vector $\by$ can be reconstructed by a concatenation of  fibers of the Hankel tensor $\tY(I_1, \ldots, I_{n-1},:, 1,\ldots, 1)$, where $ n = 1, \ldots, N-1$, and
\be
\by = \left[\begin{array}{c}
\tY(1:I_1-1, 1,\ldots, 1)\\
\vdots \\
\tY(I_1, \ldots, I_{n-1},1:I_n-1, 1,\ldots, 1)\\
\vdots \\
\tY(I_1, \ldots, I_{N-1},1:I_N)\\
\end{array}
\right]\,. \label{eq_hankel_fiber}
\ee
\item
Slices of a Hankel tensor $\tY$, i.e., any subset of the tensor produced by fixing $(N-2)$ indices of its entries and varying the two remaining indices, are also Hankel matrices.
\item
An $N$th-order Hankel tensor, $\calH_{I_1, \ldots, I_{N-1},I_N}(\by)$,  can be constructed from an $(N-1)$th-order Hankel tensor $\calH_{I_1, \ldots, I_{N-2},I_{N-1}+I_N-1}(\by)$ of size $I_1 \times   \cdots \times I_{N-2} \times (I_{N-1}+I_N-1)$ by converting its mode-$(N-1)$ fibers to Hankel matrices of size $I_{N-1} \times I_N$.
\item
Similarly to the Toeplitz tensor, the convolution of $(N-1)$ vectors, $\bx_n$ of lengths $I_n$, and a vector $\by$ of length $L$, can be represented as
\be
	[\bx_1 \ast \bx_2 \ast \cdots \ast \bx_{N-1} \ast \by]_{J:L} = \tY \,  \bar{\times}_1  \, \tilde{\bx}_1 \,\bar{\times}_2 \, \tilde{\bx}_2 \cdots \bar{\times}_{N-1} \, \tilde{\bx}_{N-1} ,\notag
\ee
or
\be
\bx_1 \ast \bx_2 \ast \cdots \ast \bx_{N-1} \ast \by  = \widetilde{\tY} \,  \bar{\times}_1  \, \tilde{\bx}_1 \,\bar{\times}_2 \, \tilde{\bx}_2 \cdots \bar{\times}_{N-1} \, \tilde{\bx}_{N-1} , \notag
\ee
where $\tilde{\bx}_n = [x_{n}(I_n), \ldots, x_n(2), x_n(1)]$, $J = \sum_n I_n - N+1$,
$\tY = \calH_{I_1, \ldots, I_{N-1}, L-J}(\by)$ is the $N$th-order Hankel tensor of $\by$, whereas $\widetilde{\tY} = \calH_{I_1, \ldots, I_{N-1},L+J}([\0_{J}^{\text{T}}, \by^{\text{T}}, \0_{J}^{\text{T}}]^{\text{T}})$ is the Hankel tensor of a zero-padded version of $\by$.

\item A Hankel tensor with identical dimensions $I_n = I$, for all $n$, is a symmetric tensor.
\end{itemize}
\end{remark}

\begin{example}
A $3 \times 3 \times 3$ -- dimensional Hankel tensor of a sequence $1,2,\ldots,7$ is a symmetric tensor, and is given by
\be
\calH_{3,3,3}(1:7)  &=& \left[ \left[
\begin{array}{cccc}
 1     &2&     3\\
     2&     3&     4\\
     3     &4&     5
\end{array}
\right], \left[
\begin{array}{cccc}
  2    & 3&     4\\
     3     &4 &    5\\
     4     &5&     6
\end{array}
\right],
\left[
\begin{array}{cccc}
     3     &4 &    5\\
     4     &5&     6\\
          5     &6&     7
\end{array}
\right]
\right]   \, . \notag
\ee
\end{example}

\subsection{\noindent \bf Quantized Tensorization}

It is important to notice that the tensorizations into the Toeplitz and Hankel tensors typically enlarge the number of data samples (in the sense that the number of entries of the corresponding  tensor is larger than the number of original samples). For example, when the dimensions $I_n = 2$ for all $n$, the so generated tensor to be a quantized tensor of order  $(L-1)$, while the number of entries of a such tensor increases from the original size $L$ to $2^{L-1}$. Therefore, quantized tensorizations are suited to analyse  signals of short-length, especially in multivariate autoregressive modelling.

\subsection{\noindent \bf Convolution Tensor}

\begin{figure}[t]
\centering
\includegraphics[width=.55\linewidth, trim = 0.0cm .0cm 0cm 0cm,clip=true]{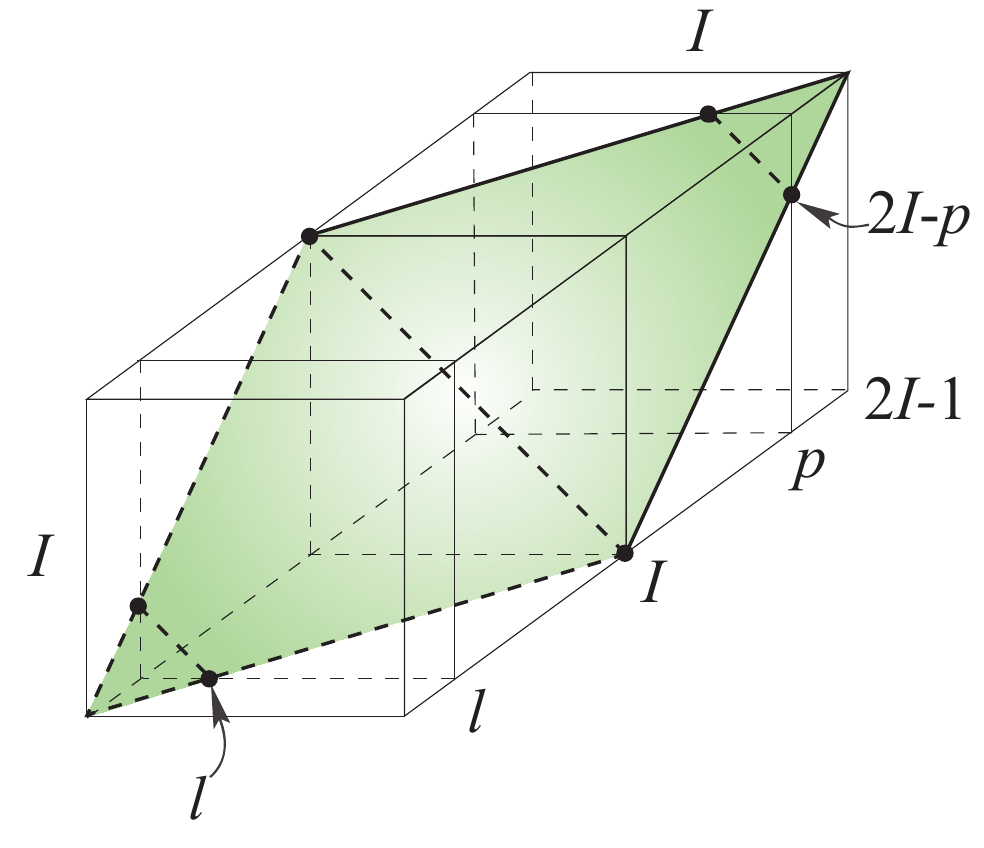}
\caption{Visualization of a convolution tensor of size $I \times I \times (2I-1)$. Unit entries are located on the shaded parallelogram.}
\label{fig_conv34}
\end{figure}

Consider again the convolution $\bx \ast \by$  of two vectors of respective lengths $I$ and $L$. We can then rewrite the expression for the entries-$(I, I+1, \ldots, L)$ as
\be
  [\bx \ast \by]_{I:L} = \tC \, {\bar{\times}}_1 \, \bx \,  {\bar{\times}}_3 \,\by \, , \notag
\ee
where $\tC$ is a third-order tensor of size $I \times J \times L$, $J = L-I+1$, for which the $(l-I)$-th diagonal elements of $l$-th slices are ones, and the remaining entries are zeros, for $l = 1, 2, \ldots, L$.  For example, the slices $\tC(:,:,l)$, for $l\le I$, are given by
\begin{equation}
\tC(:,:,l) = \begin{blockarray}{cccccc}
\begin{block}{c[ccccc]}
  & 0   & & & & 0\\
 \\
 & 1 &    &   \ddots   \\
& & \ddots &    & \\
&0 &  &  1 & & 0  \\
\end{block}
& & &  {l}
\end{blockarray}\notag
\,.
\end{equation}
The tensor $\tC$ is called the {\emph{convolution tensor}}. Illustration of a convolution tensor of size $I \times I \times (2I-1)$ is given in Figure~\ref{fig_conv34}.

Note that a product of this tensor with the vector $\by$ yields the Toeplitz matrix of the generating vector $\by$, which is of size $I\times J$, in the form
\be
\tC \, \bar{\times}_3\, \by = \calT_{I,J}(\by) \, ,\notag
\ee
while the tensor-vector product $\tC  \bar{\times}_1 \bx$ yields a Toeplitz matrix of the generating vector $[\0_{L-I}^{\text{T}}, \bx^{\text{T}}, \0_{J-1}^{\text{T}}]^{\text{T}}$, or a circulant matrix of $[\0_{L-I}^{\text{T}}, \bx^{\text{T}}]^{\text{T}}$
\be
\tC  \, \bar{\times}_1 \, \bx = \calT_{L,J}([\0_{L-I}^{\text{T}}, \bx^{\text{T}}, \0_{J-1}^{\text{T}}]^{\text{T}}) \,. \notag
\ee

%
In general, for a convolution of $(N-1)$ vectors, $\bx_1$, \ldots, $\bx_{N-1}$, of respective lengths $I_1, \ldots, I_{N-1}$ and a vector $\by$ of length $L$
\be
\bz = \bx_1 \ast \bx_2 \ast \cdots \ast \bx_{N-1} \ast \by    \, ,
\ee
the entries of $\bz$ can be expressed through a multilinear product of a convolution tensor, $\tC$, of $(N+1)$th-order and size  $I_1 \times I_2 \times \cdots \times I_N \times L$, $I_{N} = L - \sum_{n = 1}^{N-1} I_n + N-1$, and the $N$ input vectors
\be
\bz_{L-I_N+1:L} = \tC  \, {\bar\times}_1 \, \bx_1\, {\bar\times}_2 \, \bx_2 \cdots \, {\bar\times}_{N-1} \, \bx_{N-1} \, {\bar\times}_{N+1} \, \by  \, .
\ee
Most  entries of $\tC$ are zeros, except for those located at $(i_1,i_2,\ldots, i_{N+1})$, such that
\be
\sum_{n = 1}^{N-1} \bar{i}_n + i_{N}  - i_{N+1} = 0 \, ,
\ee
where $\bar{i}_n = I_n - i_n$, $i_n = 1, 2 , \ldots , I_n$.\pagebreak

 The tensor product $\tC {\bar{\times}}_{N+1}  \, \by$ yields the Toeplitz tensor of the generating vector $\by$, shown below
 \be
 	\tC \, {\bar\times}_{N+1}  \, \by = \calT_{I_1, \ldots, I_{N}}(\by).
 \ee

\subsection{\noindent \bf QTT Representation of the Convolution Tensor}

An important property of the convolution tensor is that it has a QTT representation with rank no larger than the number of inputs vectors, $N$.
To illustrate this property, for simplicity, we consider an $N$th-order Toeplitz tensor of size $I\times I \times \cdots \times I$ generated from a vector of length $(N \, I - N+1)$, where $I = 2^D$.
The convolution tensor of this Toeplitz tensor is of $(N+1)$th-order and of size $I\times I \times \cdots \times I \times (N \, I - N+1)$.

{\noindent \bf Zero-padded convolution tensor.}
By appending $(N-1)$ zero tensors of size $I\times I \times \cdots \times I$ before the convolution tensor, we obtain an $(N+1)$th-order convolution tensor, $\tC, $ of size $I\times I \times \cdots \times I \times I N$.

{\noindent \bf QTT representation.}
The zero-padded convolution tensor can be represented in the following QTT format
\be
\tC = \widetilde{\tC}^{(1)} \,  \skron \,\widetilde{\tC}^{(2)} \, \skron \cdots \skron \,\widetilde{\tC}^{(D)} \,\skron \, \widetilde{\tC}^{(D+1)}  \, ,\label{eq_qtt_convtesor}
\ee
where ``$\skron$'' represents the strong Kronecker product between block tensors{\footnote{A ``block tensor'' represents a multilevel matrix, the entries of which are matrices or tensors.}} $\widetilde{\tC}^{(n)} = [\widetilde{\tC}^{(n)}_{r,s} ]$ defined from the $(N+3)$th-order core tensors $\tC^{(n)}$ as $\widetilde{\tC}^{(n)}_{r,s} = \tC^{(n)}(r,:,\ldots,:,s)$.

The last core tensor $\tC^{(D+1)}$ represents an exchange (backward identity) matrix of size $N \times N$ which can represented as an $(N+3)$th-order tensor of size $N \times 1 \times \cdots  \times 1 \times N \times 1$.
 The first $D$ core tensors $\tC^{(1)}$, $\tC^{(2)}$, \ldots, $\tC^{(D)}$ are expressed based on the so-called elementary core tensor $\tS$ of size $N \times \underbrace{2 \times 2 \times \cdots \times 2}_{\text{$(N+1)$  dimensions}} \times N$, as
\be
\tC^{(1)} =  \tS(1,:,\ldots,:), \quad
\tC^{(2)} = \cdots = \tC^{(D)} = \tS \, .
\ee
The rigorous definition of the elementary core tensor is provided in Appendix \ref{sec:elementary_core_tensor}.

Table~\ref{tbl_ranksconvtt}  provides ranks of the QTT representation for various order of convolution tensors.
The elementary core tensor $\tS$ can  be further re-expressed  in a (tensor train) TT-format with $(N+1)$ sparse  TT cores, as
\be
\tS =  \llangle \tG^{(1)} , \tG^{(2)} ,  \ldots , \tG^{(N+1)}  \rrangle\, ,\notag
\ee
where $\tG^{(k)}$ is of size $(N+k-1) \times 2 \times (N+k)$, for $k = 1, \ldots, N$, and the last core tensor
  $\tG^{(N+1)}$ is of size $2N \times 2 \times N$. 


\begin{table}[t]\caption{Rank of QTT representations of convolution tensors of $(N+1)$th-order for $N = 2, \ldots, 17$.}
\centering
\shadingbox{
\begin{tabular}{cc|cc}
$N$  & QTT rank  & $N$  & QTT rank \\\hline
2       &  2, 2, 2, \ldots, 2	 & 10     & 6, 8, 9, \ldots, 9\\	
3       & 2, 3, 3,  \ldots, 3	& 11     & 6, 9, 10, \ldots, 10\\
4       & 3, 4, 4, \ldots, 4	& 12     & 7, 10, 11, \ldots, 11\\
5       & 3, 4, 5, \ldots, 5	& 13     & 7, 10, 12, \ldots, 12\\
6       & 4, 5, 6, \ldots, 6	& 14      & 8, 11, 13, \ldots, 13\\
7       & 4, 6, 7, \ldots, 7	& 15     & 8, 12, 14, \ldots, 14\\
8       & 5, 7, 8, \ldots, 8	& 16     & 9, 13, 15, \ldots, 15\\
9       & 5, 7, 8, \ldots, 8	& 17     & 9, 13, 15, \ldots, 15\\
\end{tabular}
}
\label{tbl_ranksconvtt}
\end{table}

\begin{figure}[t]
\centering
\includegraphics[width=.95\linewidth, trim = 0.0cm .0cm 0cm 0cm,clip=true]{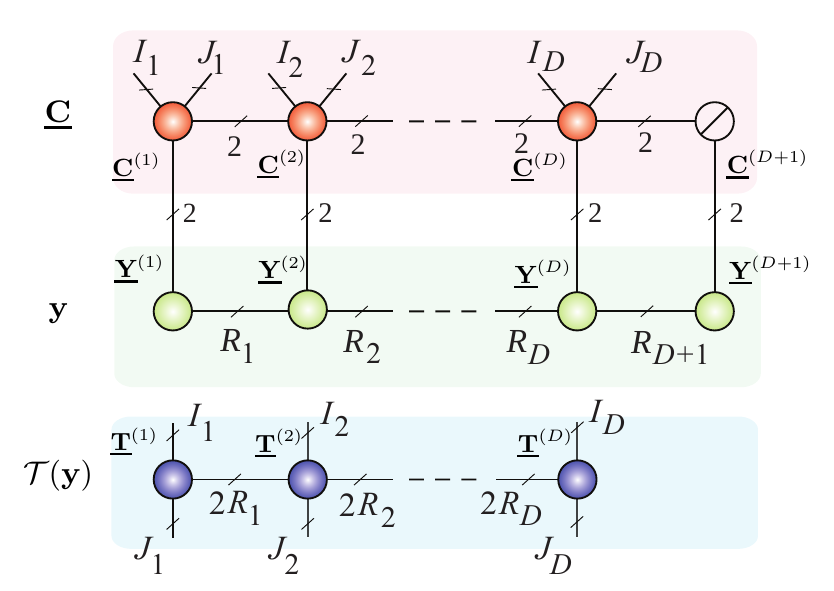}	%
\caption{Representation of the convolution tensor in QTT format. (Top) Distributed representation  of a convolution tensor $\tC$ of size $I \times J \times 2I$ in a QTT format, where $I = J = 2^D$. The first core tensor $\tC^{(1)}$ is of size $1 \times 2 \times 2 \times 2 \times 2$, the last core tensor $\tC^{(D+1)}$ represents a backward identity matrix, and the remaining 5th-order core tensors of size $2\times 2 \times 2 \times 2 \times 2 $ are identical.  A vector $\by$ is of length $2^{D+1}$ in a QTT format. (Bottom) Generation of the Toeplitz matrix, $\calT(\by)$, of the vector $\by$ from the convolution tensor  and its representation in the QTT format, $I_d = J_d = 2$ for $d = 1, \ldots, D$.}
\label{fig_conv3}
\end{figure}

\begin{example}{\bf Convolution tensor of 3rd-order.}\label{ex_conv_3}

For the vectors $\bx$ of length $2^D$ and $\by$ of length $(2^{D+1}-1)$, the expanded convolution tensor has size of  $2^D \times 2^D \times 2^{D+1}$.
The elementary core tensor $\tS$ is then of size $2\times 2 \times 2 \times 2 \times 2$ and its sub-tensors, $\tS(i, :,:,:,:)$,  are given in a $2\times 2$ block form of the last two indices  through four matrices, $\bS_1$, $\bS_2$, $\bS_3$ and $\bS_4$, of size $2 \times 2$, that is
\begin{align}
\tS(1,:,:,:,:) &= \left[\begin{array}{cc}
\bS_1 & \bS_3 \\ \bS_2 & \bS_4
\end{array}
\right]\notag \,, \quad
\tS(2,:,:,:,:)  =  \left[\begin{array}{cc}
\bS_2 & \bS_4\\
 \bS_3 & \bS_1
\end{array}
\right]\notag \,,
\end{align}
where
\be
	\bS_1 = \left[\begin{array}{cc}
	1 & 0 \\  0 & 1
	\end{array}
\right],\,
 	\bS_2 = \left[\begin{array}{cc}
	0 &1 \\  0 & 0
	\end{array}
\right]  , \,
	\bS_3 = \left[\begin{array}{cc}
	0 & 0 \\  0 & 0
	\end{array}
\right],\,
	\bS_4 = \left[\begin{array}{cc}
	0 & 0 \\  1 & 0
	\end{array}
\right]  .\,\, \notag
\ee

The convolution tensor can then be represented in a  QTT format of rank-2 \citep{KazeevToeplitz13}
with core tensors $\tC^{(2)} = \cdots = \tC^{(D)} = \tS$, $\tC^{(1)} = \tS(1,:,:,:,:)$,
and the last core tensor $\tC^{(D+1)} = \left[\begin{array}{cc} 0 &     1\\
     1    &  0\end{array}\right]$ which is of size $2 \times 1 \times 1 \times 2\times 1 $.
This QTT representation is useful to generate a Toeplitz matrix when
its generating vector is given in the QTT format.
     An illustration of the convolution tensor $\tC$ is provided in Figure~\ref{fig_conv3}.
\end{example}

\begin{example}{\bf Convolution tensor of fourth-order.}

For the convolution tensor of fourth order, i.e., Toeplitz order $N = 3$,
the elementary core tensor $\tS$ is of size $3 \times 2 \times 2 \times 2 \times2 \times 3$, and is given in a $2\times 3$ block form of the last two indices as
%
\begin{align}
\tS(1,:,\ldots,:) &= \left[\begin{array}{ccc}
\tS_1 & \tS_3 & \tS_5\\
\tS_2 & \tS_4 &  \tS_6
\end{array}
\right]\notag \,,
\quad
\tS(2,:,\ldots,:)  =  \left[\begin{array}{ccc}
\tS_2 & \tS_4 & \tS_6\\
\tS_5 & \tS_1 &  \tS_3
\end{array}
\right]\notag \,,\\
\tS(3,:,\ldots,:)  &=  \left[\begin{array}{ccc}
\tS_5 & \tS_1  & \tS_3 \\
\tS_6 & \tS_2 & \tS_4
\end{array}
\right]\notag \,.
\end{align}
where $\tS_n$ are of size $2 \times 2 \times 2$, $\tS_5$, $\tS_6$ are zero tensors, and
   \begin{align}
	\tS_1 &= \left[ \left[\begin{array}{cc}
	1 & 0    \\
	 0 & 0
	\end{array}
\right]  \left[\begin{array}{cc}
	 0 & 1 \\
	 1 & 0
	\end{array}
\right] \right]\;,\quad
	\tS_2 =  \left[ \left[\begin{array}{cc}
	 0 &  0  \\
	 0 & 0
	\end{array}
\right] \left[\begin{array}{cc}
	  1 & 0 \\
          0 & 0
	\end{array}
\right]\right]\;, \notag  \\
	\tS_3 &= \left[\left[\begin{array}{cc}
	 0 & 0    \\
	 0 & 1
	\end{array}
\right]  \left[\begin{array}{cc}
	 0 & 0 \\
	 0 & 0
	\end{array}
\right]\right]\;,\quad
	\tS_4 = \left[\left[\begin{array}{cc}
	 0 & 1    \\
	 1 & 0
	\end{array}
\right] \left[\begin{array}{cc}
 	 0 & 0 \\
	0 & 1
	\end{array}
\right]\right]\; . \notag
\end{align}
Finally, the zero-padded convolution tensor of size $2^D \times 2^D \times 2^D \times 3 \cdot 2^D$ has a QTT representation in (\ref{eq_qtt_convtesor}) with $\tC^{(1)} = \tS(1,:,:,:,:,[1, 2])$, $\tC^{(2)} = \tS([1, 2],:,:,:,:,:)$, $\tC^{(3)} = \cdots = \tC^{(D)} = \tS$,
and the last core tensor $\tC_{D+1} = \left[\begin{array}{ccc} 0 & 0 &    1\\ 0 & 1 & 0\\
     1    &  0 & 0\end{array}\right]$ which is of size $3 \times 1 \times 1 \times 3\times 1 $.
\end{example}

\subsection{\noindent \bf Low-rank Representation of Hankel and Toeplitz Matrices/Tensors}

The Hankel and Toeplitz foldings are multilinear tensorizations, and can be applied to the BSS problem, as in (\ref{equ_btd3}).
When the Hankel and Toeplitz tensors of the hidden sources are of low-rank in some tensor
 network representation,
the tensor of the mixture is expressed as a sum of low rank tensor terms.

For example, the Hankel and Toeplitz matrices/tensors of an exponential function, $v_k = az^{k-1}$, are rank-1 matrices/tensors,
 and consequently Hankel matrices/tensors of sums and/or products of exponentials, sinusoids, and polynomials will also be  of low-rank, which is equal to the degree of the function being considered.\\

{\noindent \bf Hadamard Product.} More importantly, when Hankel/Toeplitz tensors of two vectors $\bu$ and $\bv$ have low-rank CP/TT representations, the Hankel/\-Toeplitz tensor of their element-wise product, $\bw = \bu \circledast \bv$, can also be represented in the same CP/TT tensor format
\be
	\calH(\bu) \circledast \calH(\bv) &=& \calH(\bu \circledast \bv) \notag \\
	\calT(\bu) \circledast \calT(\bv) &=& \calT(\bu \circledast \bv) \notag \,.
\ee
The  CP/TT rank of $\calH(\bu \circledast \bv)$ or $\calT(\bu \circledast \bv)$ is not larger than the product of the CP/TT ranks of the tensors of $\bu$ and $\bv$.

\begin{example} \label{ex_hankel_tsin}

The third-order Hankel tensor of $u(t) = \sin(\omega t)$ is a rank-3 tensor, and the third-order Hankel tensor of $v(t) = t$ is of rank-2; hence the Hankel tensor of the $w(t) = t \, \sin(\omega t )$ has at most rank-6.
\end{example}

{\noindent \bf Symmetric CP and Vandermonde decompositions.}
It is important to notice that a Hankel tensor $\tY$ of size $I \times I \times \cdots \times I$ can always be represented by a symmetric CP decomposition
\be
	\tY = \tI \times_1 \bA  \times_2 \bA  \cdots \times_N \bA\,. \notag
\ee
Moreover, the tensor $\tY$ also admits a symmetric CP decomposition with Vandermonde structured factor matrix \citep{QiCMS15}
\be
\tY = \diag_N(\blambda) \times_1 \bV^{\text{T}}  \times_2 \bV^{\text{T}} \cdots \times_N \bV^{\text{T}}   \, ,\label{eq_hankel_vandermonde}
\ee
where $\blambda$ comprises $R$ non-zero coefficients, and $\bV$ is a  Vandermonde matrix generated from $R$ distinct values $\bv = [v_1, v_2, \ldots, v_R]$
\be
	\bV = \left[\begin{array}{ccccc}
	1 	&  v_1 & v_1^2 & \ldots & v_1^{I-1}\\
	1 	&  v_2 & v_2^2 & \ldots & v_2^{I-1}\\
	\vdots & \vdots & \vdots & \ddots & \vdots \\
	1 	&  v_R & v_R^2 & \ldots & v_R^{I-1}
	\end{array} \right] \, .
\ee
By writing the decomposition in (\ref{eq_hankel_vandermonde}) for the entries $\tY(I_1, \ldots, I_{n-1},:, 1,\ldots, 1)$ (see (\ref{eq_hankel_fiber})), the Vandermonde decomposition of the Hankel tensor $\tY$ becomes a Vandermonde factorization of $\by$ \citep{Chenthesis}, given by
\be
\by = \left[ \begin{array}{ccccc}
	1 	&  1 &   \ldots & 1\\
	v_1   & v_2  & \ldots & v_R \\
	v_1^2   & v_2^2  & \ldots & v_R^2 \\
	\vdots & \vdots &\ddots & \vdots \\
	v_1^{L-1}   & v_2^{L-1}  & \ldots & v_R^{L-1} \\
\end{array}\right] \blambda\, .	\notag
\ee
Observe that various Vandermonde decompositions of the Hankel tensors of the same vector $\by$, but of different tensor orders $N$, have the same generating Vandermonde vector $\bv$.
Moreover, the Vandemonde rank, i.e, the minimum of $R$ in the decomposition (\ref{eq_hankel_vandermonde}), therefore cannot exceed the length $L$ of the generating vector $\by$.\\

{\noindent \bf QTT representation of Toeplitz/Hankel tensor.}
As mentioned previously, the zero-padded convolution tensor of $(N+1)$th-order can be represented in a QTT format of rank of at most $N$.
Hence, if a vector $\by$ of length $2^D N$ has a QTT representation of rank-$(R_1, \ldots, R_D)$, given by
\be
\by = \widetilde{\bY}^{(1)} \skron \widetilde{\bY}^{(2)} \skron \cdots \skron  \widetilde{\bY}^{(D+1)} \, ,
\ee
where $\widetilde{\bY}^{(d)}$ is an $R_{d-1} \times R_d$ block matrix of the core tensor $\tY^{(d)}$ of size $R_{d-1} \times 2 \times R_d$, for $d = 1,\ldots, D$, or of $\tY^{(D+1)}$ of size $R_D \times N \times 1$, then following the relation between the convolution tensor and the Toeplitz tensor of the generating vector $\by$, we have
\be
	\calT(\by)  = \tC \, {\bar{\times}}_{N+1} \, \by  \,.\label{eq_qtt_Toeplitz1}
\ee
This $N$th-order Toeplitz tensor can also be represented by a QTT tensor with rank of at most $N(R_1, \ldots, R_D)$, as
\be
	\calT(\by)  = \widetilde{\tT}^{(1)} \skron \widetilde{\tT}^{(2)} \skron \cdots \skron  \widetilde{\tT}^{(D)}  \label{eq_qtt_Toeplitz} \,,
\ee
where $\widetilde{\tT}^{(d)}$ is a block tensor of the core tensor ${\tT}^{(d)}$.
The core $\tT^{(1)}$ is of size $1 \times 2\times  \cdots \times 2 \times N \, R_1$, and cores
$\tT^{(2)}$, \ldots, $\tT^{(D-1)}$ are of size $N R_{d-1} \times 2\times  \cdots \times 2 \times N R_{d}$, while the last core tensor $\tT^{(D)}$ is of size $NR_{D-1} \times 2 \times \cdots \times 2 \times 1$. These core tensors are core contractions between the two core tensors $\tC^{(d)}$ and $\tY^{(d)}$.
Figure~\ref{fig_conv3} illustrates the generation of a Toeplitz matrix  as a tensor-vector product of a third-order convolution tensor $\tC$ and a generating vector, $\bx$,  of length $2^{D+1}$, both in QTT-formats. The core tensors of $\tC$ are given in Example~\ref{ex_conv_3}.

{\bf Remarks:}
\begin{itemize}
\item
Because of zero-padding within the convolution tensor, the Toeplitz tensor of $\by$, generated in (\ref{eq_qtt_Toeplitz1}) and (\ref{eq_qtt_Toeplitz}), takes only entries
$\by(N),\by(N+1), \ldots, \by(2^D N)$,
i.e., it corresponds to the Toeplitz tensor of the generating vector
$\by(N),\by(N+1), \ldots, \by(2^D N)$.

\item The Hankel tensor also admits a QTT representation in the similar form to a Toeplitz tensor (cf. (\ref{eq_qtt_Toeplitz})).

\item  Low-rank TN representation of the Toeplitz and Hankel tensors has been exploited, e.g., in blind source separation and harmonic retrieval. By verifying a low-rank TN representation of the signal in hand, we can confirm the existence of a low-rank TN representation of Toeplitz/\-Hankel tensors of the signal.

\item
QTT rank of the Toeplitz tensor  in (\ref{eq_qtt_Toeplitz}) is at most $N$ times the QTT rank of the generating vector $\by$. The rank may not be minimal. For example, the sinusoid signal is of rank-2  in QTT format, and its Toeplitz tensor also has a rank-2 QTT representation.

\item \emph{Fast convolution of vectors in QTT formats.} A straightforward consequences is that when vectors $\bx_n$ are given in their QTT formats, their convolution $\bx_1 \ast \bx_2 \ast \cdots \ast \bx_N$ can be computed through core contractions between the core tensors of the convolution tensor and those of the vectors.

\end{itemize}


\section[Tensorization by Means of {L\"owner} Matrix]
{Tensorization by Means of {L\"owner} Matrix ({L\"owner}  Folding)}
\sectionmark{Tensorization by Means of {L\"owner} Matrix}

A {L\"owner} matrix of a vector ${\bf v} \in \mathbb{R}^{I+J}$ is formed from a function $f(t)$ sampled at $(I+J)$ distinct points $\{ x_1,\ldots, x_I, y_1,\ldots, y_J \}$, to give
	$$
	{\bf v} = \left[ f(x_1),\ldots,f(x_I),f(y_1),\ldots,f(y_J) \right]^{\text{T}} \in\mathbb{R}^{I+J},
	$$
so that the entries of ${\bf v}$ are partitioned into two disjoint sets, $\{f(x_i)\}_{i=1}^I$ and $\{f(y_j)\}_{j=1}^J$. The vector ${\bf v}$ is then converted into the L\"owner matrix, ${\bf L} \in \mathbb{R}^{I \times J}$, defined by
		$$
		{\bf L} = \left[ \frac{ f(x_i) - f(y_j) }{ x_i - y_j } \right]_{ij} \in \mathbb{R}^{I\times J}.
		$$

 L\"owner matrices appear as a powerful tool in fitting a model to data in the form of a rational (Pade form) approximation, that is $f(x) = A(x)/B(x)$. When considered as transfer functions, such type of approximations are much more powerful than the  polynomial approximations, as in this way it is also possible to model discontinuities and spiky data. The optimal order of such a rational approximation is given by the rank of the L\"owner matrix.
In the context of tensors, this allows us to construct a model of the original dataset which is amenable to higher-order tensor representation, has minimal computational complexity, and for which the accuracy is governed by the rank of the L\"owner matrix.
An example of L\"owner folding of a vector $[1/3, 1/4, 1/5,1/6,  1/8, 1/9, 1/10 ] $ is given below
\be
\begin{bmatrix}
\frac{1/3-1/8}{3-8} & \frac{1/3-1/9}{3-9}  & \frac{1/3-1/10}{3-10}  \\
\frac{1/4-1/8}{4-8} & \frac{1/4-1/9}{4-9}  & \frac{1/4-1/10}{4-10}  \\
\frac{1/5-1/8}{5-8} & \frac{1/5-1/9}{5-9}  & \frac{1/5-1/10}{5-10}  \\
\frac{1/6-1/8}{6-8} & \frac{1/6-1/9}{6-9}  & \frac{1/6-1/10}{6-10}  \\
\end{bmatrix}
=
- \begin{bmatrix}  1/3 \\ 1/4 \\ 1/5 \\ 1/6 \end{bmatrix}
\begin{bmatrix}  1/8 & 1/9 & 1/10 \end{bmatrix}
\notag.
\ee
More applications of this tensorization can be found in \citep{Debals2016loewner}.

\section[Tensorization based on Cumulants and Derivatives of GCFs]
{Tensorization based on Cumulant and Derivatives of \\the Generalised Characteristic Functions}
\sectionmark{Generalised Characteristic Functions}

The use of higher-order statistics (cumulants) or partial derivatives of the Generalised Characteristic Functions (GCF) as a means of tensorization is useful in the identification of a mixing matrix in a blind source separation.

\begin{figure}[t!]
\centering
\includegraphics[width=1\linewidth, trim = 0.0cm .0cm 0cm 0cm,clip=true]{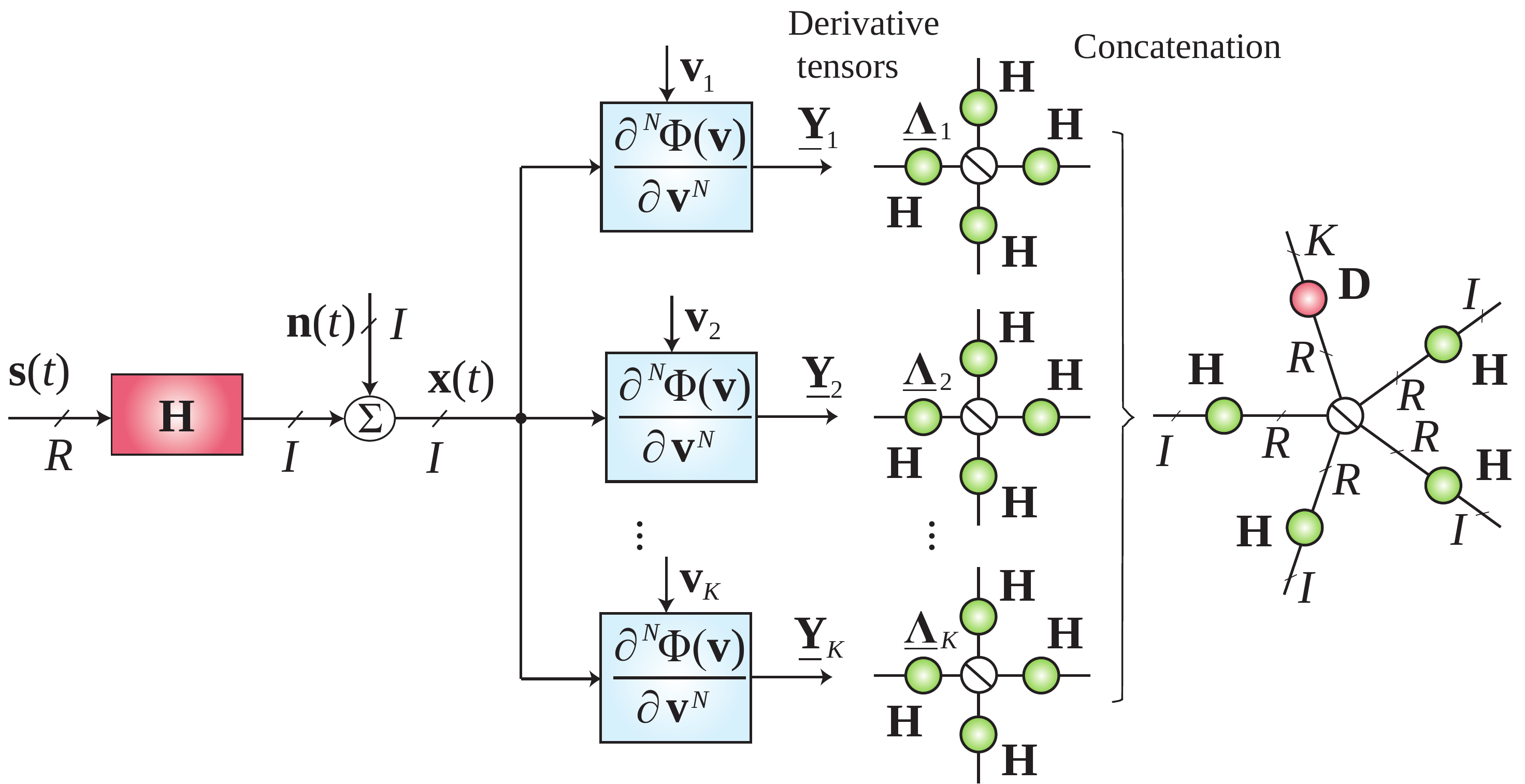}
\label{fig_bi_columnstyle1}
\caption{Tensorization  based on derivatives of the characteristic functions and tensor-based approach to blind identification. The task is to estimate the mixing matrix, $\bH$, from only the knowledge of the noisy  output observations $\bX=[\bx(1), \ldots, \bx(t), \ldots, \bx(T)] \in \Real^{I \times T}$, with $I < T$.
A high dimensional tensor $\tY$ is generated from the observations $\bX$ by means of higher-order statistics (cumulants) or partial derivatives of the second generalised characteristic functions of the observations. A CP decomposition of $\tY$  allows us to retrieve the mixing matrix $\bH$.}
\label{fig_bi_columnstyle}
\end{figure}

Consider linear mixtures of $R$ stationary sources, $\bS$, received by an array of $I$ sensors in the presence of additive noise, $\bN$ (see Figure~\ref{fig_bi_columnstyle} for a general principle). The task is to estimate a mixing matrix $\bH \in \Real^{I \times R}$ from only the knowledge of the noisy observations
\be
\bX=\bH \bS +\bN \, ,
\ee
 under some mild assumptions, i.e., the sources are statistically independent and non-Gaussian, their number is known, and the matrix $\bH$ has no pair-wise collinear columns (see also \citep{YeredorCaf,Comon20062271})
 %
%

A well-known approach to this problem is based on the decomposition of a high dimensional structured tensor, $\tY$, generated from the observations, $\bX$, by means
of partial derivatives of the second GCFs of the observations at multiple processing points. 

{\noindent \bf Derivatives of the GCFs.} More specifically, we next show how to generate the tensor $\tY$ from the observation, $\bX$. We shall denote the first and second GCFs of the observations evaluated at a vector $\bu$ of length $I$, respectively by
\begin{align}
\phi_{\bx}(\bu) = \mbox{E} \left[\exp(\bu^{\text{T}} \bx) \right], \quad
\Phi_{\bx}(\bu)  = \log{\phi_{\bx}(\bu)}  \, .
\end{align}
Similarly, $\phi_{\bs}(\bv)$ and $\Phi_{\bs}(\bv)$ designate the first and second GCFs of the sources, where $\bv$ is of length $R$.
Because the sources are statistically independent, the following holds
\be
\Phi_{\bs}(\bv) = \Phi_{s_1}(v_1)  + \Phi_{s_2}(v_2)  + \cdots + \Phi_{s_R}(v_R)  \,,
\ee
which implies that $N$th-order derivatives of $\Phi_{\bs}(\bv)$ with respect to $\bv$ result in $N$th-order diagonal tensors of size $R \times R \times \cdots \times R$, where $N = 2, 3, \ldots$, that is
\be
\tensor{\Psi}_{\bs}(\bv)  = \frac{\partial^N \Phi_{\bs}(\bv) }{\partial \bv^N}
 =    \diag_{N}\!\left\{\frac{d^N \Phi_{s_1}}{dv_1^N}, \frac{d^N \Phi_{s_2}}{ dv_2^N}, \ldots, \frac{d^N \Phi_{s_R}}{dv_R^N}\right\}\!.
\ee
In addition, for the noiseless case $\bx(t) = \bH \, \bs(t)$, and since $\Phi_{\bx}(\bu) = \Phi_{\bs}(\bH^{\text{T}} \bu)$, the $N$th-order derivative of $\Phi_{\bx}(\bu)$ with respect to $\bu$ yields a symmetric tensor of $N$th-order which admits a CP decomposition of rank-$R$ with $N$ identical factor matrices $\bH$, to give
\be
	\tensor{\Psi}_{\bx}(\bu) = \tensor{\Psi}_{\bs}(\bH^{\text{T}}\bu) \times_1 \bH  \times_2 \bH \cdots \times_N \bH \,. \label{eq_CP_gcf}
\ee
In order to improve the identification accuracy, the mixing matrix $\bH$ should be estimated as a joint factor matrix in decompositions of various derivative tensors, evaluated at distinct processing points $\bu_1$, $\bu_2$, \ldots, $\bu_K$. This is equivalent to a decomposition of an $(N+1)$th-order tensor $\tY$ of size $I\times I \times \cdots \times I \times K$ concatenated from the $K$ derivative tensors as
\be
	\tY(:,\ldots,:,k) = \tensor{\Psi}_{\bx}(\bu_k), \quad k = 1, 2, \ldots, K.
\ee
The CP decomposition of the tensor $\tY$ can be written in form of
\be
\tY =  \tI \times_1 \bH  \times_2 \bH \cdots \times_N \bH \times_{N+1} \bD\,,
\ee
where the last factor matrix $\bD$ is of size $K \times R$, and each row comprises the diagonal of the symmetric tensor $\tensor{\Psi}_{\bs}(\bH^{\text{T}}\bu_k)$.

In the presence of statistically independent, additive and stationary Gaussian noise,
we can eliminate the derivatives of the noise terms in the derivative tensor $\tensor{\Psi}_{\bx}(\bu)$ by subtracting any other derivative tensor $\tensor{\Psi}_{\bx}(\tilde{\bu})$, or by an average of derivative tensors.

{\noindent \bf Estimation of Derivatives of GCF.}
In practice, the GCF of the observation and its derivatives are unknown, but can be estimated from the sample first GCF \citep{YeredorCaf}.
Detailed expression and the approximation of the derivative tensor $\tensor{\Psi}_{\bx}(\bu)$ for some low orders $N = 2, 3,\ldots, 7$, are given in Appendix~\ref{sec:apx_gcf}.

{\noindent \bf Cumulants.} When the derivative is taken at the origin, $\bu = [0, \ldots, 0]^{\text{T}}$, the tensor $\calK_{\bx}^{(N)} = \tensor{\Psi}_{\bx}^{(N)}(\0)$ is known as the $N$th-order cumulant of $\bx$, and a joint diagonalization or the CP decomposition of
higher-order cumulants is a well-studied method for the estimation of the mixing matrix $\bH$.

For the sources with symmetric probabilistic distributions, their odd-order cumulants, $N = 3, 5, \ldots$, are zero, and the cumulants of the mixtures are only due to noise. Hence, a decomposition of such tensors is not able to retrieve the mixing matrix. However, the odd-order cumulant tensors can be used to subtract the noise term in the derivative tensors evaluated at other processing points.

\begin{example}{\bf Blind identification (BI) in a system of $2$ mixtures and $R$ binary signals.}
\label{ex_qam_bi_2xR}

 To illustrate the efficiency of higher-order derivatives of the second GCF in blind identification we  considered a system of two mixtures, $I = 2$, linearly composed by $R$ signals of length $T = 100  \times \, 2^R$, the entries of which can take the values 1 or $-1$, i.e., $s_{r,t}  =  1$ or $-1$. The mixing matrix $\bH$ of size $2 \times R$ was randomly generated, where $R = 4, 6, 8$.
The signal-to-noise ratio was SNR = 20 dB.
The main purpose of BI is to estimate the mixing matrix $\bH$.

We constructed 50 tensors $\tY_i$ $\;(i=1,\ldots, 50)$ of size $R\times \cdots \times R \times 3$, which comprise three derivative tensors evaluated at the two leading left singular vectors of $\bX$, and a unit-length processing point, generated such that its collinearity degree with the first singular vector  uniformly distributed over a range of $[-0.99, 0.99]$.
The average derivative tensor was used to eliminate the noise term in $\tY_i$.

%


{\emph{CP decomposition of derivative tensors.}}
The tensors $\tY_i$ were decomposed by CP decompositions of rank-$R$ to retrieve the mixing matrix $\bH$. 
The mean of Squared Angular Errors  $SAE(\bh_r, \hat{\bh}_r) = -20\log_{10}\arccos(\frac{\bh_r^{\text{T}} \hat{\bh}_r}{|\bh_r|_2 |\hat{\bh}_r|_2})$ over all columns $\bh_r$ was computed as a performance index for one estimation of the mixing matrix.

The averages of the mean and best MSAEs over 100 independent runs for the number of the unknown sources $R  = 4, 6, 8$ are plotted in Figure~\ref{fig_bi_2xR}. The results indicate that with a suitably chosen processing point, the decomposition of the derivative tensors yielded good estimation of the mixing matrix. Of more importance is that higher-order derivative tensors, e.g., 7th and 8th orders, yielded better performance than lower-order tensors, while the estimation accuracy deteriorated with the number of sources.

\begin{figure}[t!]
\centering
\includegraphics[width=.85\linewidth, trim = 0.0cm 0cm 0cm 0cm,clip=true]
{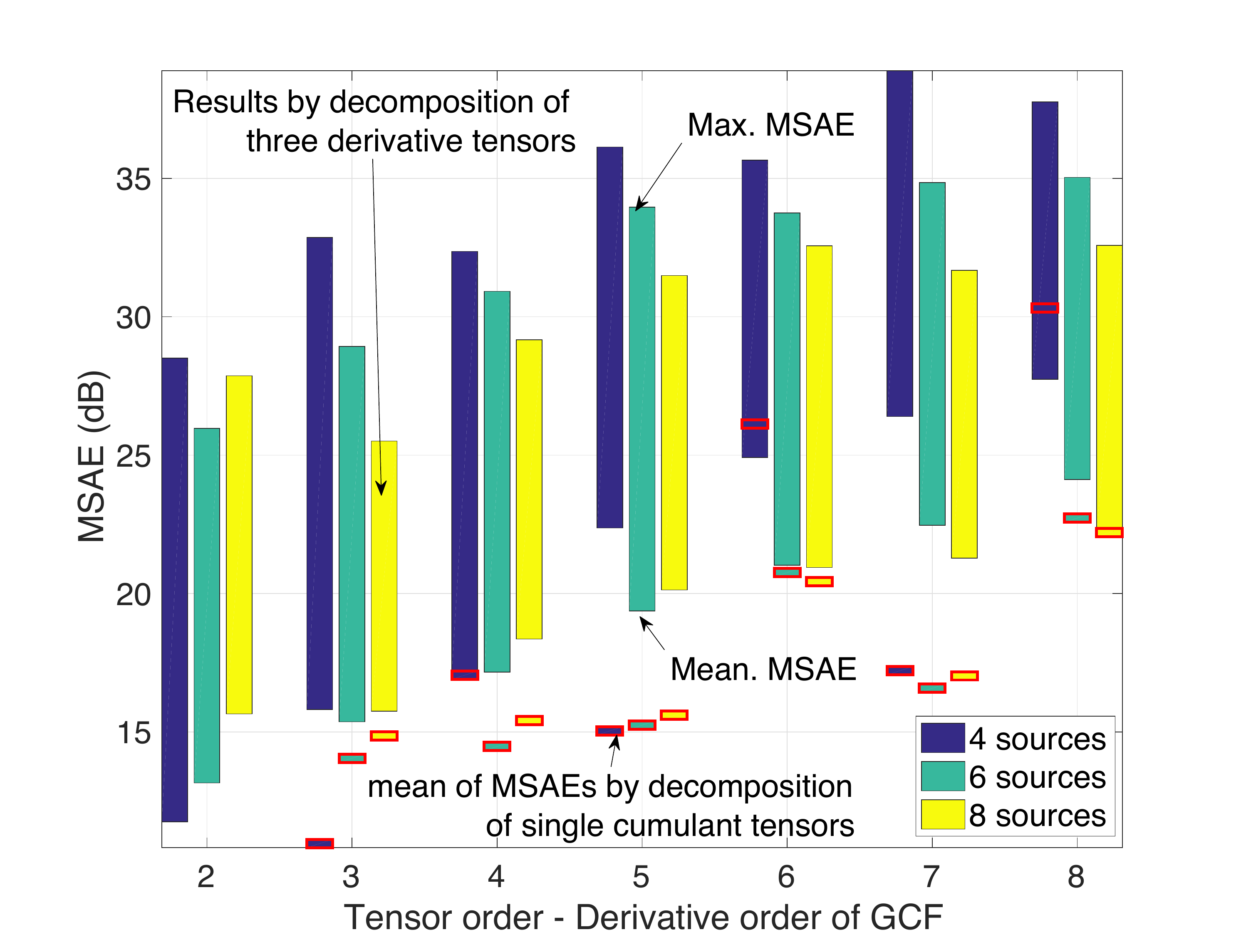}
\caption{Mean SAE (in dB) in the estimation of the mixing matrix $\bH$ from only two mixtures, achieved by CP decomposition of three $2\times 2 \times \cdots \times 2$ derivative tensors of the second GCFs. Small bars in red represent the mean of MSAEs, obtained by decomposition of single cumulant tensors.}
\label{fig_bi_2xR}
\end{figure}

{\emph{CP decomposition of cumulant tensors.}}
Because of symmetric pdfs, the odd order cumulants of the sources are zero.
Only decompositions of cumulants of order 6 or 8 were able to retrieve the mixing matrix $\bH$.
For all the test cases, better performances could be obtained by a decomposition of three derivative tensors.


{\emph{Tensor train decomposition of derivative tensors.}} The estimation of the mixing matrix $\bH$ can be performed in a two-stage decomposition
	\begin{itemize}
	\item A tensor train decomposition of high-order derivative tensors, e.g., tensor order exceeds 5.
	\item A CP decomposition of the tensor in TT-format, to retrieve the mixing matrix.
	\end{itemize}
Experimental results confirmed that the performances with prior TT-\-decomposition were more stable and yielded an approximately 2 dB higher mean SAE  than those using only CP decomposition for derivative tensors of orders 7 and 8 and a relatively high number of unknown sources.
\end{example}

\subsection{Tensor Structures in Constant Modulus Signal Separation}

Another method to generate tensors of relatively high order in BSS is through modelling modulus of the estimated signals as  roots of a polynomial.

Consider a linear mixing system $\bX = \bH \bS$ with $R$ sources of length $K$, and $I$ mixtures, where the modulus of the sources $\bS$ is drawn from a set of given moduli. For simplicity, we assume $I = R$.
For example, the binary phase-shift keying (BPSK) signal in telecommunication consists of a sequence of 1 and $-1$, hence, it has a constant modulus of unity. The quadrature phase shift keying (QPSK) signal takes one of the values $\pm 1 \pm 1 i$, i.e., it has a constant modulus $\sqrt{2}$. The 16-QAM signal has three squared moduli of 2, 10 and 18. For this BSS problem for single constant modulus signals, \citet{7080134} linked the problem to CP decomposition of a fourth-order tensor. For multi-constant modulus signals, \citet{DebalsCMA} established a link to a coupled CP decomposition.

A common method to extract the original sources $\bS$ is to use a demixing matrix $\bW$ of size $I \times R$ or a vector $\bw$ of length $I$ such that  $\by = \bw^{\text{T}} \bX$ is an estimate of one of the source signals.
The constant modulus constraints require that each entry, $|y_{k}|$, must be one of given moduli, $c_1, c_2, \ldots, c_M$. This means that for all entries of $\by$  the following holds
\be
	f(y_{k}) = \prod_{m = 1}^{M}  (|y_{k}|^2 - c_m) = 0\, .
\ee
In other words, $|y_{k}|^2$ are roots of an $M$th-degree polynomial, given by
\be
 	  p^M + \alpha_{m} p^{M-1} +  \cdots + \alpha_{2} p+ \alpha_{1}\, ,  \notag
\ee
with coefficients $\alpha_{M+1} = 1$, and $\alpha_1, \alpha_2, \ldots, \alpha_M$, given by
\be
	\alpha_m = (-1)^{m-1} \sum_{i_1, i_2, \ldots, i_m}  c_{i_1} c_{i_2} \cdots c_{i_m} \, .
\ee
By expressing $|y_{k}|^2 = (\bw \otimes \bw^*)^{\text{T}} (\bx_k \otimes \bx_k^*)$, and
\be
	|y_{k}|^{2m} = (\bw^{\otimes_m} \otimes (\bw^{\otimes_m})^*)^{\text{T}} (\bx_k^{\otimes_m} \otimes (\bx_k^{\otimes_m})^*)  \, , \notag
\ee
where the symbol ``$*$'' represents the complex conjugate, $\bx^{\otimes_m} = \bx \otimes \bx \otimes \cdots \otimes \bx$ denotes the Kronecker product of $m$ vectors $\bx$,
and bearing in mind that the rank-1 tensors $\bw^{\circ_m} = \bw \circ \bw \circ \cdots \circ \bw$ are symmetric, and in general have only $\frac{(R+m-1)!}{m! (R-1)!}$ distinct coefficients, the rank-1 tensors $\bw^{\circ_m} \circ (\bw^{\circ_m})^*$ have at least $\left(\frac{(R+m-1)!}{m! (R-1)!} \right)^2$ distinct entries.
We next introduce the operator $\calK$ which keeps only distinct entries of the symmetric tensor $\bw^{\circ_m} \circ (\bw^{\circ_m})^*$ or of the vector $\bw^{\otimes_m} \otimes (\bw^{\otimes_m})^*$.
The constant modulus constraint of $y_k$ can then be rewritten as
\begin{align}
f(y_k) &= \alpha_1 +  \sum_{m = 2}^{M+1}\!\!  \alpha_m (\bw^{\otimes_m} \otimes (\bw^{\otimes_m})^*)^{\text{T}}   (\bx_k^{\otimes_m} \otimes (\bx_k^{\otimes_m})^*) \notag \\
&=  \alpha_1 +  \sum_{m = 2}^{M+1}\!\!  \alpha_m  (\calK(\bw^{\otimes_m} \otimes (\bw^{\otimes_m})^*))^{\text{T}}  \diag(\bd_m)  \, \calK(\bx_k^{\otimes_m} \otimes (\bx_k^{\otimes_m})^*) \notag \\
&=  \alpha_1 + \left[\ldots,(\calK(\bw^{\otimes_m} \otimes (\bw^{\otimes_m})^*))^{\text{T}},\ldots\right]
\notag \\
&   \qquad\qquad [\ldots,    \calK(\bx_k^{\otimes_m} \otimes (\bx_k^{\otimes_m})^*)^{\text{T}} \diag(\alpha_m \bd_m), \ldots ]^{\text{T}}    \,,\notag
\end{align}
where $d_m(i)$ represents the number of occurrences of an entry of $\calK(\bx_k^{\otimes_m} \otimes (\bx_k^{\otimes_m})^*)$ in $\bx_k^{\otimes_m} \otimes (\bx_k^{\otimes_m})^*$.

The vector of the constant modulus constraints of $\by$ is now given by
\be
\vf &=& [\ldots, f(y_k), \ldots]^{\text{T}}  
= \alpha_1 {\bf 1} +\bQ  \bv  \, ,\label{equ_modulus}
\ee
where
\be
\bv = \left[\begin{array}{@{}c@{}} \vdots \\\calK(\bw^{\otimes_m} \otimes (\bw^{\otimes_m})^*)\\
\vdots
\end{array}
\right] , \;
\bQ = \left[\begin{array}{@{}c@{}} \vdots \\
\diag(\alpha_m \bd_m)   \calK(\bX^{\odot_m} \odot (\bX^{\odot_m})^*)
 \\
\vdots
\end{array}\right]^{\text{T}}
 \notag \! .
\ee
The constraint vector is zero for the exact case, and should be small for the noisy case.
For the exact case, from (\ref{equ_modulus}) and $f(y_{k+1}) - f(y_{k}) = 0$, this leads to
\be
	\bL \bQ \bv = \0	\, ,\notag
\ee
where $\bL$ is the first-order Laplacian implying that   the vector $\bv$ is in the null space of the matrix $\tilde{\bQ} = \bL \bQ$. The above condition holds for other demixing vectors $\bw$, i.e.,
$\tilde{\bQ} \bV = \0$, where $\bV = [\bv_1, \ldots, \bv_R]$, and each $\bv_r$ is constructed from a corresponding demixing vector $\bw_r$.

With the assumption $I = R$, and that the sources have complex values, and the mixing matrix does not have collinear columns, it can be shown that the kernel of the matrix $\tilde{\bQ}$ has the dimension of $R$ \citep{DebalsCMA}. Therefore, the basis vectors, $\bz_r$, $r = 1, \ldots, R$, of the kernel of $\tilde{\bQ}$ can be represented as linear combination of $\bV$, that is
\be
	\bz_r = \bV \blambda_r \,. \notag
\ee
Next we partition $\bz_r$ into $M$ parts,  $\bz_r = [\bz_{rm}]$, each of the length $\left(\frac{(R+m-1)!}{m! (R-1)!} \right)^2$, which can be expressed as
\be
\bz_{rm} =\sum_{s= 1}^R  \lambda_{rs}  \calK(\bw_s^{\otimes_m} \otimes (\bw_s^{\otimes_m})^* )  
=  	  \calK\left( \sum_{s = 1}^R  \lambda_{rs}   \bw_s^{\otimes_m} \otimes (\bw_s^{\otimes_m})^* \right),   \notag
\ee
thus implying that $\bW$ and $\bW^*$ are factor matrices of a symmetric tensor $\tZ_{rm}$ of $(2m)$th-order, constructed from the vector $\bz_{rm}$, i.e., $\calK(\vtr{\tZ_{rm}}) = \bz_{rm}$,
in the form
\be
	\tZ_{rm}  = \llbracket \diag_{2m}(\blambda_r); \underbrace{\bW, \ldots, \bW}_{\text{$m$ terms}}, \underbrace{\bW^*, \ldots, \bW^*}_{\text{$m$ terms}} \rrbracket \, .
\ee
By concatenating all $R$ tensors $\tZ_{1m}$, \ldots, $\tZ_{Rm}$ into one $(2m+1)$th-order tensor $\tZ_m$, the above $R$ CP decompositions become
\be
\tZ_{m} = \llbracket \tI;  \underbrace{\bW, \ldots, \bW}_{\text{$m$ terms}}, \underbrace{\bW^*, \ldots, \bW^*}_{\text{$m$ terms}} , \bLambda \rrbracket \,.
\ee
All together, the $M$ CP decompositions of $\tZ_1$, \ldots, $\tZ_M$ form a coupled CP tensor decomposition to find the two matrices $\bW$ and $\bLambda$.

\begin{example}[Separation of QAM signals.]

We performed the separation of two rectangular 32- or 64-QAM signals of length $1000$ from two mixture signals corrupted by additive Gaussian noise with SNR = 15 dB. Columns of the real-valued mixing matrix had unit-length, and a pair-wise collinearity of 0.4. The 32-QAM signal had $M = 5$ constant moduli of  2, 10, 18, 26 and 34, whereas the 64-QAM signal had $M = 9$ squared constant moduli of  2, 10, 18, 26, 34, 50, 58, 74 and 98.
Therefore, for the first case (32-QAM), the demixing matrix was estimated from 5 tensors of size $2\times 2 \times \cdots  \times 2$ and of respective orders 3, 5, 7, 9 and 11, while for the later case (64-QAM), we decomposed 9 quantized tensors of orders 3, 5, \ldots, 19. The estimated QAM signals for the two cases were perfectly reconstructed with zero bit error rates. Scatter plots of the recovered signals are shown in Figure~\ref{fig_qam32}.

\begin{figure}[t!]
\centering
\subfigure[32-QAM]{
\includegraphics[width=.45\linewidth, trim = 0.0cm 0cm 0cm 0cm,clip=true]{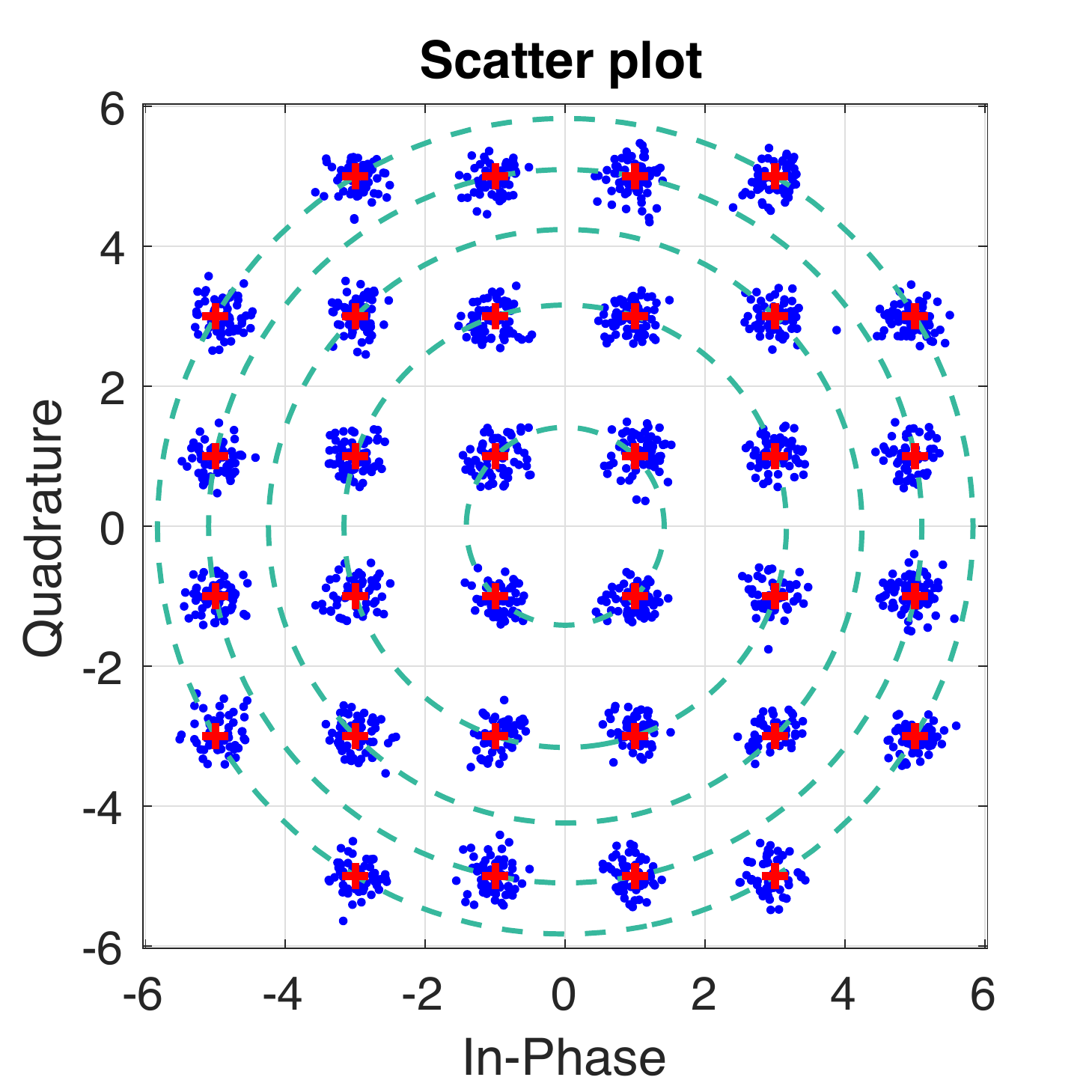}}
\hfill
\subfigure[64-QAM]{\includegraphics[width=.45\linewidth, trim = 0.0cm 0cm 0cm 0cm,clip=true]{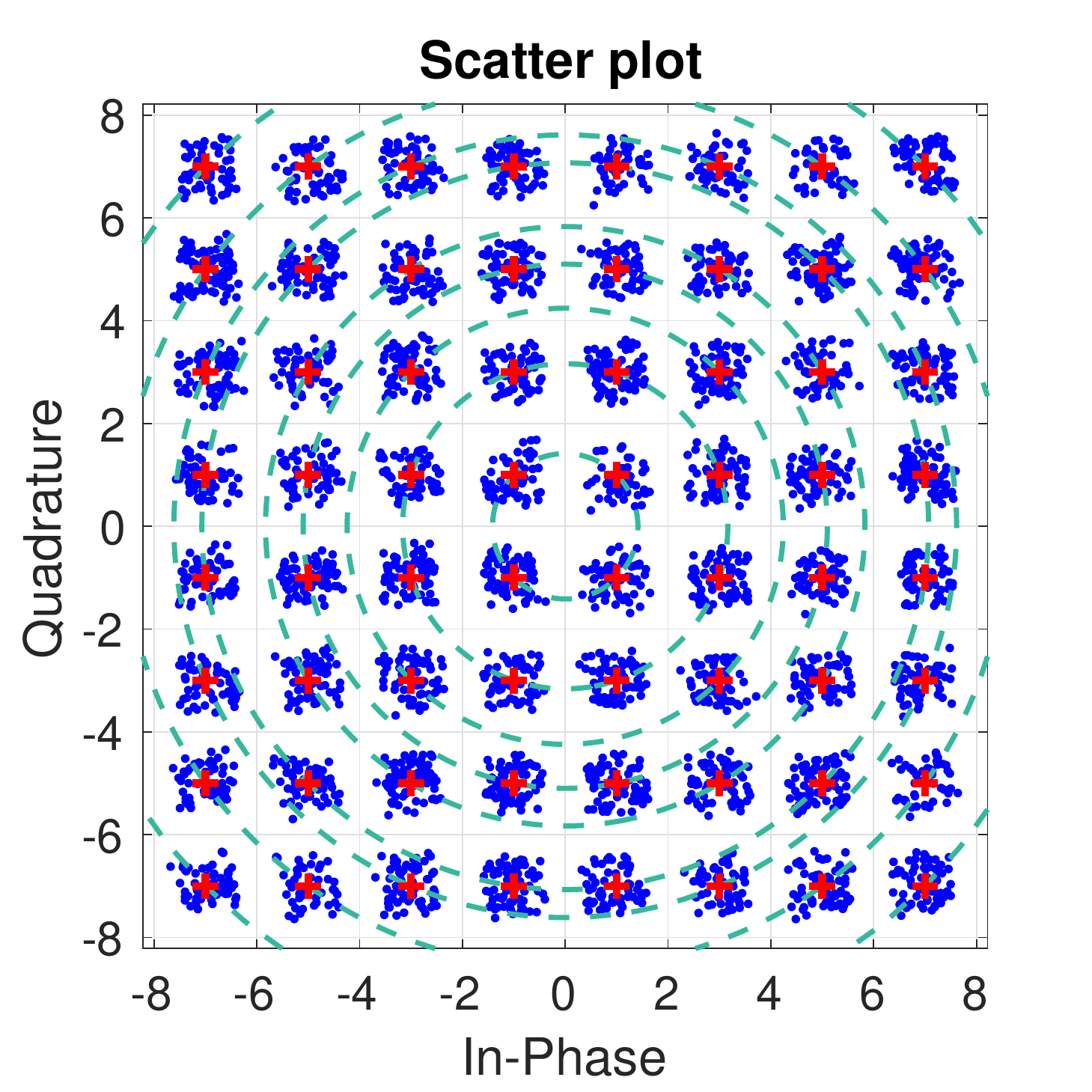}}
\caption{Scatter plots of the estimated sources (blue dots). Red dots indicate the ideal signal constellation, the values of which are located on one of dashed circles.}\label{fig_qam32}
\end{figure}

\end{example}

\section{Tensorization by Learning Local Structures}

Different from the previous tensorizations, this tensorization approach generates tensors from local blocks (patches) which are similar or closely related.
For the example of an image, given that the intensities of pixels in a small window are highly correlated, hidden structures which represent relations between small patches of pixels can be learnt in local areas. These structures can then be used to reconstruct the image as a whole in, e.g., an application of image denoising \citep{Phan_TT_part1}.

\begin{figure}[t!]
\centering
\includegraphics[width=.8\linewidth, trim = 0.0cm 0cm 0cm 0cm,clip=true]{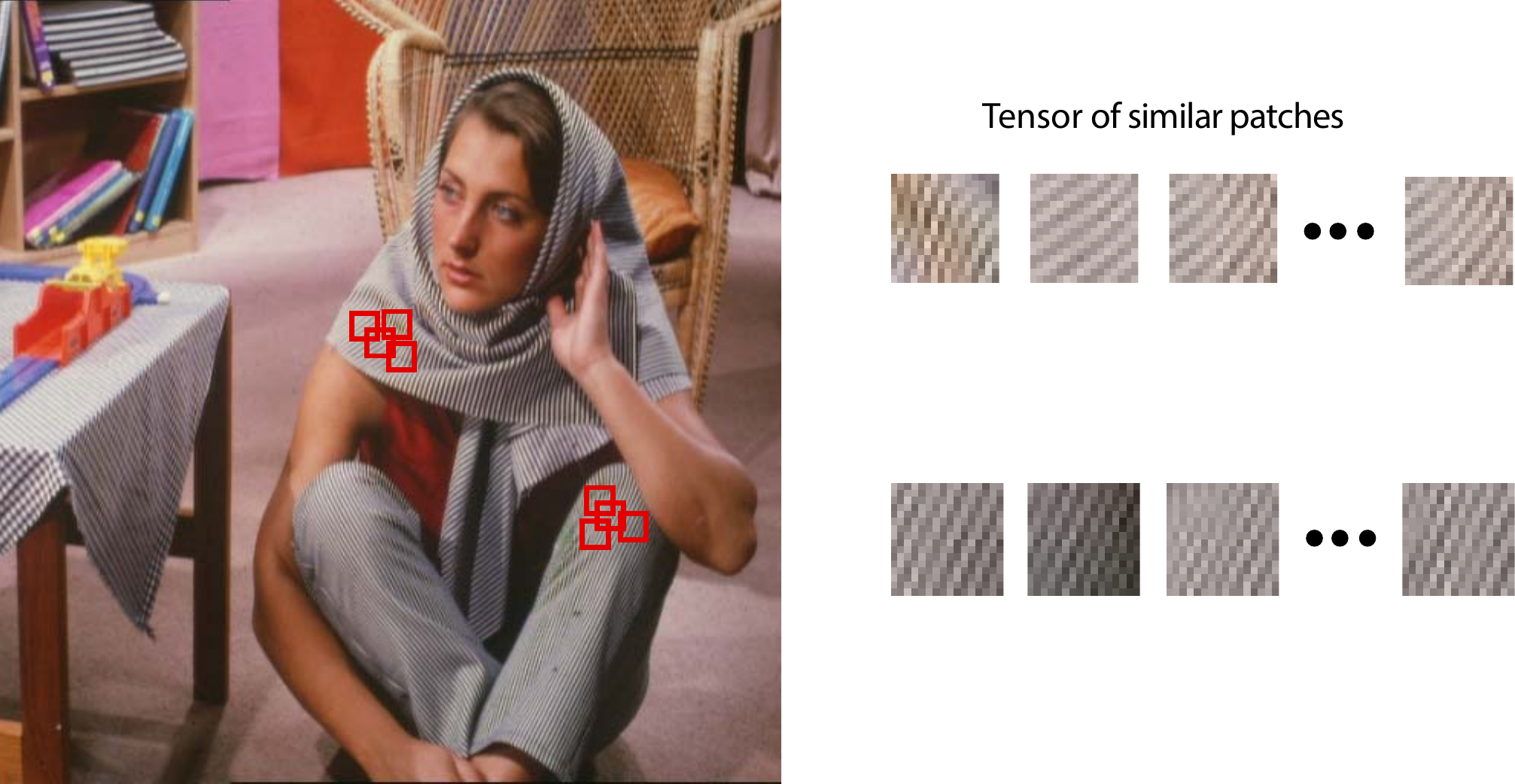}
\caption{A ``local-structure'' tensorization method generates 5th-order tensors of size $h \times w \times 3 \times (2d+1) \times (2d+1)$ from similar image patches, or patches in close spatial proximity.}\label{fig_localpatches}
\end{figure}


For a color RGB image $\tY$ of size $I  \times J \times 3$,  each block of pixels of size $h \times w  \times 3$ is denoted as
\be
	\tY_{r,c} = \tY(r:r+h-1,c:c+w-1,:)		\notag .
\ee
A small tensor, $\tZ_{r,c}$, of size $h \times w \times 3 \times (2d+1) \times (2d+1)$, comprising $(2d+1)^2$ blocks centered around $\tY_{r,c}$,  with $d$ denoting the neighbourhood width, can be constructed in the form
\be
	\tZ_{r,c}(:,:,:,d+1+i,d+1+j) =   \tY_{r+i,c+j}, 		\notag
\ee
where $i, j = -d, \ldots, 0 , \ldots, d$, as illustrated in Figure~\ref{fig_localpatches}.
Every $(r,c)$-th block $\tZ_{r,c}$ is then approximated through a constrained tensor decomposition
\be
	\|\tZ_{r,c} - \hat{\tZ}_{r,c}\|_F^2 \le \varepsilon^2 \,, \label{equ_block_approx}
\ee
where the noise level $\varepsilon^2$ can be determined by inspecting the coefficients of the image in the high-frequency bands.
A pixel is then reconstructed as the average of all its approximations which cover that pixel.

\begin{figure}[t!]
\centering
\subfigure[Noisy image]{\includegraphics[width=.315\linewidth, trim = 0cm 0cm 0cm 0cm,clip=true]{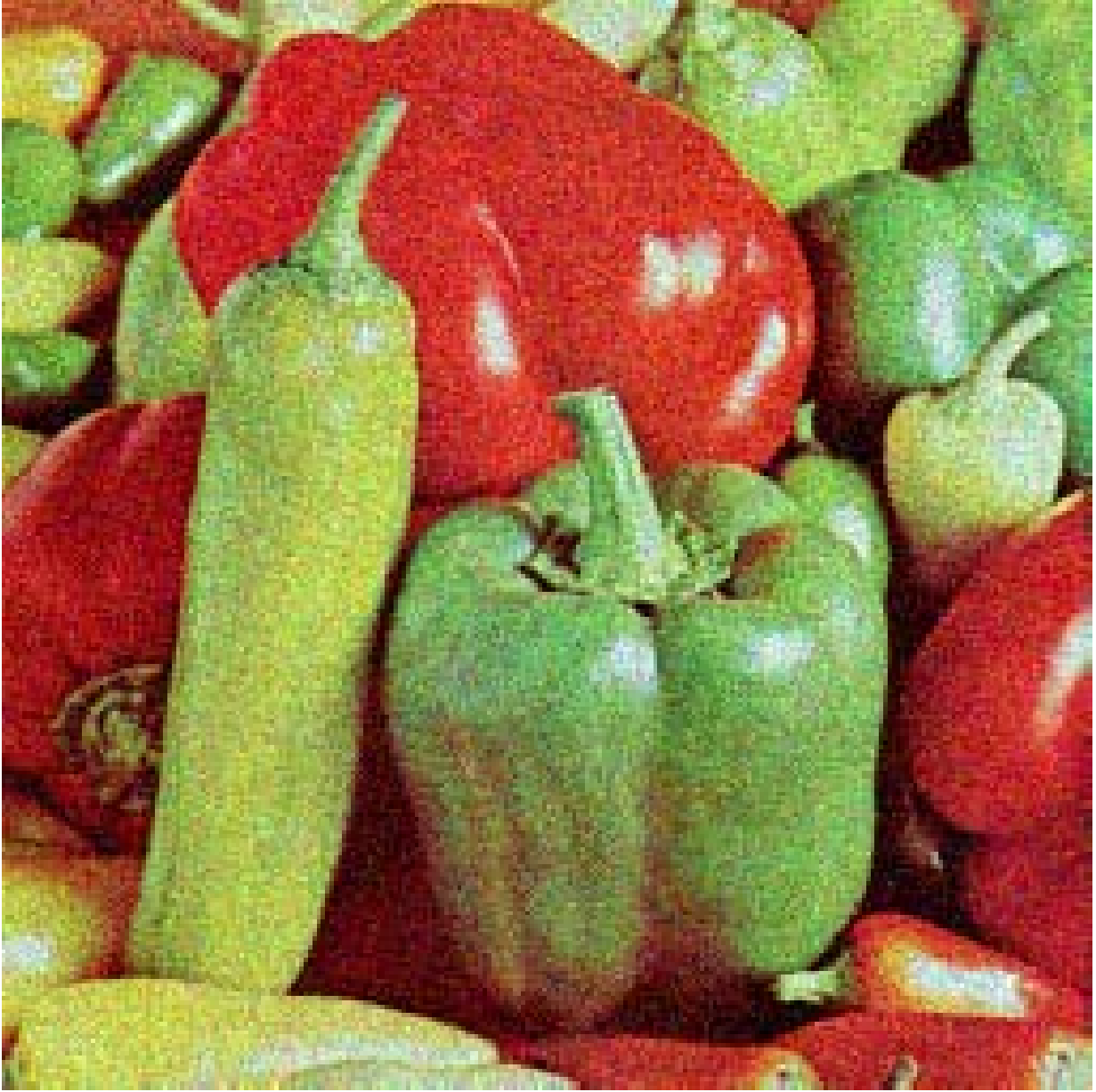}}
\subfigure[TT, PSNR = {31.64} dB]{\includegraphics[width=.315\linewidth, trim =  0cm 0cm 0cm 0cm,clip=true]{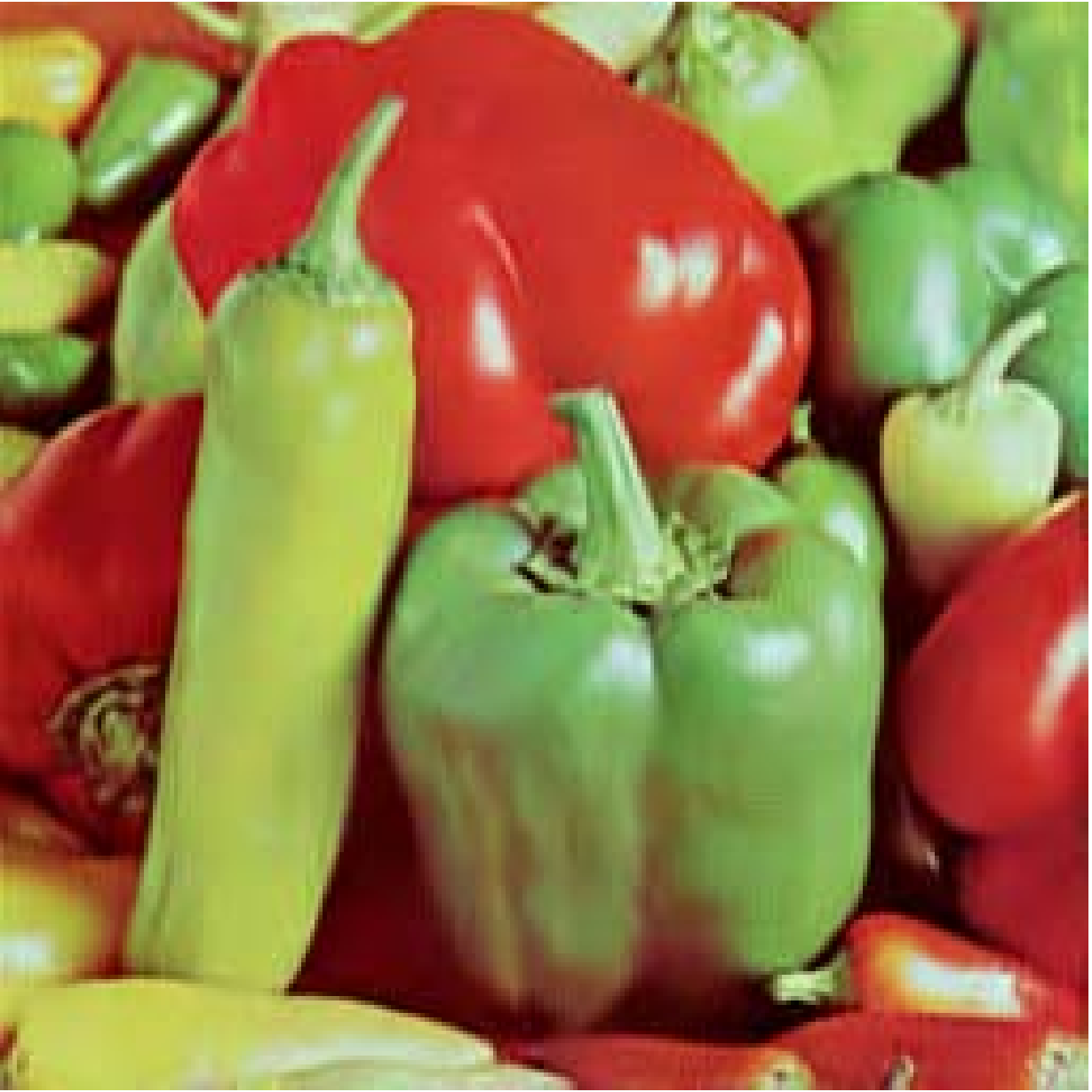}
\label{fig_pepper_256_10dB_tt_ascu}}
\subfigure[CP, PSNR = 28.90 dB]{\includegraphics[width=.315\linewidth, trim =  0cm 0cm 0cm 0cm,clip=true]{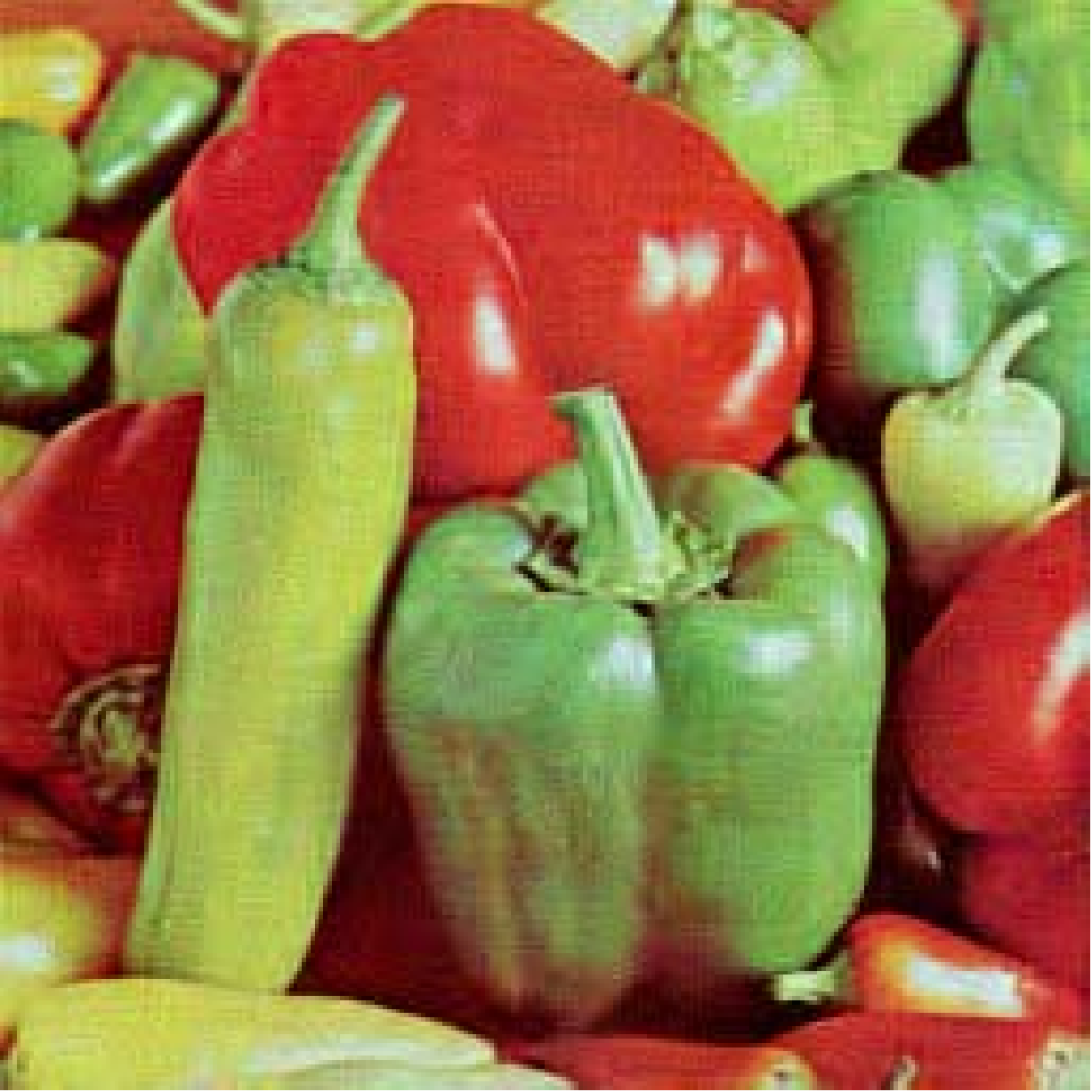}
\label{fig_pepper_256_10dB_brtf}}
%
\caption{Tensor based image reconstruction in Example~\ref{ex_TT_denoising_image}. The Pepper image with added noise at 10 dB SNR (left), and the images reconstructed using the TT (middle) and CP (right) decompositions.}
\label{fig_pepper_10dB}
\end{figure}

\begin{example}{\bf Image denoising.}\label{ex_TT_denoising_image}
The principle of tensorization from learning the local structures is next demonstrated in an image denoising application
for  the benchmark ``peppers'' color image of size $256 \times 256 \times 3$, which was corrupted by white Gaussian noise at SNR = 10 dB.
Latent structures were learnt for patches of sizes $8 \times 8 \times 3$ (i.e., $h = w = 8$) in the search area of width $d = 3$.
To the noisy image, we applied the DCT spatial filtering before their block reconstruction.
The results are shown in Figure~\ref{fig_pepper_10dB}, and illustrate the advantage of the tensor network approach over a CP decomposition approach.
\end{example}

\section{Tensorization based on Divergences, Similarities or Information Exchange}
\sectionmark{Tensorization based on Distance Measure or Similarity}

For a set of $I$ data points $\bx_i$, $i = 1, 2, \ldots, I$, this type of tensorization generates an $N$th-order nonnegative symmetric tensor of size $I \times I \times \cdots \times I$, the entries of which represent $N$-way similarities or dissimilarities between $\bx_{i_1}$, $\bx_{i_2}$, \ldots, $\bx_{i_N}$, where $i_n = 1, \ldots, I$, so that
\be
	\tY(i_1, i_2, \ldots, i_N) = d(\bx_{i_1}, \bx_{i_2}, \ldots, \bx_{i_N}) \, .
\ee
Such metric function can express pair-wise distances between the two observations $\bx_i$ and $\bx_j$. In a general case, $d(\bx_{i_1}, \bx_{i_2}, \ldots, \bx_{i_N})$ can compute the volume of a convex hull formed by $N$ data points.

The so generated tensor can be expanded to $(N+1)$th-order tensor, where the last mode expresses the change of data points over e.g., time or trials.
Tensorizations based on divergences  and similarities are useful for the analysis of interaction between observed entities, and for their clustering or classification.

\section{Tensor Structures in Multivariate Polynomial Regression}
\label{sec:mpr}

The Multivariate Polynomial Regression (MPR) is an extension of the linear and multilinear regressions which allows us to model nonlinear interaction between independent variables \citep{doi:10.1080/00207178908559683,billings2013nonlinear,Taylorfit}.
For illustration, consider a simple example of fitting a curve to data with two independent variables $x_1$ and $x_2$, in the form
\be
	y &=& w_0 + w_1 x_1 + w_2 x_2 + w_{12} x_1 x_2\,
	\label{eq_mpr2_4} .
\ee
The term $w_{12}$ then quantifies the strength of interaction  between the two independent variables in the data, $x_1$ and $x_2$. Observe that the model is still linear with respect to the variables $x_1$ and $x_2$, while involving the cross-term $w_{12} x_1 x_2$.
The above model can also have more terms, e.g., $x_1^2, x_1 x_2^2$, to describe more complex functional behaviours. For example, the full quadratic polynomial regression for two independent variables, $x_1$ and $x_2$, can have up to 9 terms, given by
\be
y &=&  w_0 + w_1 x_1 + w_2 x_2 + w_{12} x_1 x_2 \notag \\ &&+ w_{11} x_1^2 + w_{22} x_2^2 + w_{112} x_1^2 x_2 + w_{122} x_1 x_2^2 + w_{1122} x_1^2 x_2^2 \,.\;\;\;\;\label{eq_mpr2_9}
\ee




{\noindent \bf Tensor representation of the system weights.}
The simple model for two independent variables in  (\ref{eq_mpr2_4}) can be rewritten in a bilinear form as
\be
y =   \left[\begin{array}{cc}1   &  x_1 \end{array}\right] \, \left[\begin{array}{cc} w_0 & w_{2} \\ w_{1} & w_{12} \end{array}\right]
\left[\begin{array}{c} 1 \\ x_2\end{array}\right]  \,, \notag
\ee
whereas the full model in (\ref{eq_mpr2_9}) has an equivalent bilinear  expression
\be
y &=&  \left[\begin{array}{ccc}1   &  x_1 & x_1^2 \end{array}\right] \, \left[\begin{array}{ccc} w_0 & w_{2} & w_{22} \\ w_{1} & w_{12} & w_{122} \\ w_{11} & w_{112} & w_{1122} \end{array}\right]
\left[\begin{array}{c} 1 \\ x_2  \\ x_2^2 \end{array}\right]  \,, \notag
\ee
or a tensor-vector product representation
\be
y = \tW\, {\bar{\times}_1} \,   \left[\begin{array}{c}1   \\  x_1 \end{array}\right]   \, {\bar{\times}_2} \,   \left[\begin{array}{c}1   \\  x_1 \end{array}\right] \,
{\bar{\times}_3} \,   \left[\begin{array}{c}1   \\  x_2 \end{array}\right]\,
{\bar{\times}_4} \,   \left[\begin{array}{c}1   \\  x_2 \end{array}\right] \,  , \label{eq_mpr_2x4}
\ee
where the 4th-order weight tensor $\tW$ is of size $2\times 2 \times 2 \times 2$, and is given by
\begin{align}
\tW(:,:,1,1) &=
\left[\begin{array}{cc} w_0 & \frac{1}{2}w_1 \\ \frac{1}{2} w_1  & w_{11}\end{array}\right], \;\;
\tW(:,:,2,2) = \left[\begin{array}{cc} w_{22} &   \frac{1}{2} w_{122} \\ \frac{1}{2}  w_{122}  & w_{1122}\end{array}\right]  \, ,\notag \\
\tW(:,:,1,2) &= \tW(:,:,2,1)  =
\frac{1}{2} \left[\begin{array}{cc}   w_2 & \frac{1}{2} w_{12}  \\ \frac{1}{2} w_{12}   &  w_{112} \end{array}\right]
 \,.\;\; \notag
\end{align}


%

\begin{figure}[t]
(a) \hspace{6cm} (b)
\centering
\includegraphics[width=.995\linewidth, trim = 0.0cm .0cm 0cm 0cm,clip=true]{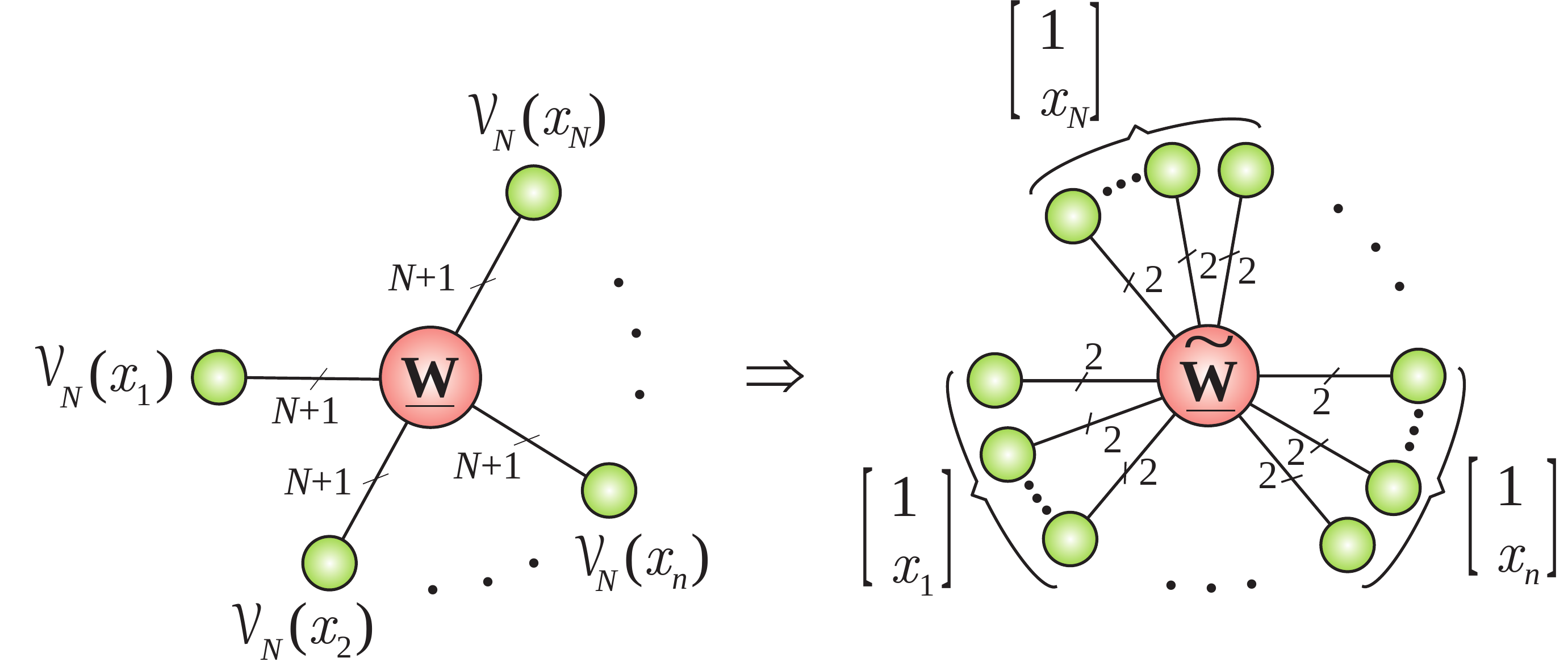}	 
\caption{Graphical illustration of  Multivariate Polynomial Regression (MPR). (a) The MPR for multiple input units $x_1$, \ldots, $x_N$, where the nonlinear function $h(x_1, \ldots, x_N)$ is expressed as a multilinear tensor-vector product of an $N$th-order  tensor, $\tW$, of size $(N+1) \times (N+1) \times \cdots \times (N+1)$, and Vandermonde vectors $\calV_N(x_n)$ of length $(N+1)$. (b) An equivalent MPR model but with
quantized $N^2$th-order  tensor $\widetilde{\tW}$  of size $2 \times 2 \times \cdots \times 2$.}
\label{Fig:MPR1}
\end{figure}

It is now obvious that for a generalised system with $N$ independent variables, $x_1$, \ldots, $x_N$, the MPR can be written as a tensor-vector product as \citep{doi:10.1080/00207178908559683}
\be
	y &=& \sum_{i_1 = 0}^{N}\sum_{i_2 = 0}^{N} \cdots\sum_{i_N = 0}^{N}  w_{i_1, i_2, \ldots, i_N} \, x_1^{i_1} \, x_2^{i_2} \cdots \, x_N^{i_N} \,  \notag \\
	&=& \tW \,  \bar{\times}_1 \, \calV_N(x_1)\,\bar{\times}_2 \, \calV_N(x_2) \, \cdots \, \bar{\times}_N \, \calV_N(x_N) \,  , \label{eq_mpr_N}
\ee
where $\tW$ is an $N$th-order tensor of size $(N+1) \times (N+1) \times \cdots \times (N+1)$, and $\calV_N(x)$ is the length-$(N+1)$ Vandermonde vector of $x$, given by
\be
\calV_N(x) = \left[ \begin{array}{ccccc} 1& x & x^2 & \ldots & x^N \end{array} \right]^{\text{T}}\, .
\ee
Similarly to the representation in (\ref{eq_mpr_2x4}), the MPR model in (\ref{eq_mpr_N}) can be equivalently expressed as a product of a tensor of $N^2$th-order and size $2\times 2 \times \cdots \times 2$ with $N$ vectors of length-2, to give
\be
y = \widetilde{\tW} \,  \bar{\times}_{1:N} \left[\begin{array}{c} 1\\ x_1 \end{array}\right]  \, \bar{\times}_{N+1:2N} \left[\begin{array}{c} 1\\ x_2 \end{array} \right] \, \cdots \bar{\times}_{N(N-1)+1:N^2} \left[\begin{array}{c} 1 \\ x_N \end{array}\right] \, .\;\;\;\label{eq_mpr_Nsquared}
\ee
An illustration of the MPR is given in Figure~\ref{Fig:MPR1}, where the input units are scalars.

The MPR has found numerous applications, owing to its ability to model any smooth, continuous nonlinear input-output system , see e.g. \citep{Taylorfit}. However, since the number of parameters in the model in (\ref{eq_mpr_N}) grows exponentially with the number of variables, $N$, the MPR demands a huge amount of data in order to yield a good model, and therefore,  it is computationally intensive in a raw tensor format, and thus not suitable for very high-dimensional data. To this end, low-rank tensor network representation emerges as a viable approach to accomplishing MPR.
For example, the weight tensor $\tW$ can be constrained to be in low rank TT-format \citep{DBLP:journals/corr/ChenBSW16}.
An alternative approach would be to consider a truncated model which takes only two entries along each mode of $\tW$ in (\ref{eq_mpr_N}). In other words, this truncated model becomes linear with respect to each variable $x_n$ \citep{novikov2016exponential}, leading to
\be
	y 
	&=& \tW_{t} \,  \bar{\times}_1 \,\left[\begin{array}{c} 1\\ x_1 \end{array}\right]  \, \bar{\times}_2 \,\left[\begin{array}{c} 1\\ x_2 \end{array} \right] \, \cdots \bar{\times}_N \,\left[\begin{array}{c} 1 \\ x_N \end{array}\right] \, ,\label{eq_mpr_Nlin}
\ee
where $\tW_{t}$ is a tensor of size $2 \times 2 \times \cdots \times 2$ in the QTT-format.
Both (\ref{eq_mpr_Nsquared}) and (\ref{eq_mpr_Nlin}) represent the weight tensors in the QTT-format, however, the tensor $\widetilde{\tW}$ in (\ref{eq_mpr_Nsquared}) has $N^2$ core tensors of the full MPR, whereas $\tW_{t}$ in (\ref{eq_mpr_Nlin}) has $N$ core tensors for the truncated model.

\section{Tensor Structures in Vector-variate Regression}
\label{sec:tensorregression}

The MPR in (\ref{eq_mpr_N}) is formulated for  scalar data. When the  observations are vectors or tensors, the model can be extended straightforwardly. For illustration, consider a simple case of  two independent vector inputs $\bx_1$ and $\bx_2$. Then, the nonlinear function which maps the input to the output $y = h(\bx_1, \bx_2)$ can be approximated in a linear form as
\begin{align}
	y = h(\bx_1, \bx_2) &= w_{0} + \bw_1^\tp \bx_1  + \bw_2^\tp \bx_2  + \bx_1^\tp \bW_{12} \bx_2    \, \label{eq_tenreg_twounits}  \\
	&= [1,  \;\bx_1^\tp] \, \left[\begin{array}{cc} w_0 & \bw_2^{\tp} \\ \bw_1 & \bW_{12}\end{array}\right]    \left[\begin{array}{cc} 1 \\ \bx_2 \end{array}\right] \notag ,
\end{align}
or in a quadratic with 9 terms, including one bias, two vectors, three matrices, two third-order tensors and one fourth-order tensor,
given by
\begin{align}
	h(\bx_1, \bx_2) &= w_{0} + \bw_1^\tp \bx_1  + \bw_2^\tp \bx_2  + \bx_1^\tp \bW_{12} \bx_2    + \bx_1^\tp \, \bW_{11} \bx_1  + \bx_2^\tp \, \bW_{22} \bx_2	\notag \\
	&	+  \tW_{112}  \,\bar{\times}_1 \, \bx_1 \, \bar{\times}_2  \, \bx_1 \, \bar{\times}_3 \, \bx_2
	+  \tW_{122} \,  \bar{\times}_1  \, \bx_1 \, \bar{\times}_2 \,  \bx_2 \, \bar{\times}_3  \,\bx_2 \notag \\
	& +  \tW_{1122}  \,  \bar{\times}_1 \,  \bx_1 \, \bar{\times}_2  \, \bx_1 \, \bar{\times}_3   \, \bx_2 \, \bar{\times}_4  \, \bx_2  \notag \\
	&= [1,  \;\bx_1^\tp, \; (\bx_1 \otimes \bx_1)^T] \; \bW\;
	\left[\begin{array}{cc} 1 \\ \bx_2 \\ \bx_2 \otimes \bx_2 \end{array}\right] \notag  ,
\end{align}
where the matrix $\bW$ is given
\be
\bW = 	\left[\begin{array}{ccc}
		 w_0 & \bw_2^{\tp} & \vtr{\bW_{22}}^\tp \\
		 \bw_1 & \bW_{12}  & [\tW_{122}]_{(1)}  \\
		 \vtr{\bW_{11}} & [\tW_{112}]_{(1,2)}  & [\tW_{1122}]_{(1,2)}  		
	\end{array}\right]   \, .
\ee
and $[\tW_{112}]_{(1,2)}$ represents the mode-(1,2) unfolding of the tensor $\tW_{112}$.
Similarly to (\ref{eq_mpr_2x4}), the above model has an equivalent expression of through the tensor-vector product of a fourth-order tensor $\tW$, in the form
\be
y = \tW \, {\bar{\times}_1} \,   \left[\begin{array}{c}1   \\  \bx_1 \end{array}\right]   \, {\bar{\times}_2} \,   \left[\begin{array}{c}1   \\  \bx_1 \end{array}\right] \,
{\bar{\times}_3} \,   \left[\begin{array}{c}1   \\  \bx_2 \end{array}\right]\,
{\bar{\times}_4} \,   \left[\begin{array}{c}1   \\  \bx_2 \end{array}\right] \,.
\ee

\begin{figure}[t]
(a) \hspace{6cm} (b)
\centering
\includegraphics[width=.995\linewidth, trim = 0.0cm .0cm 0cm 0cm,clip=true]{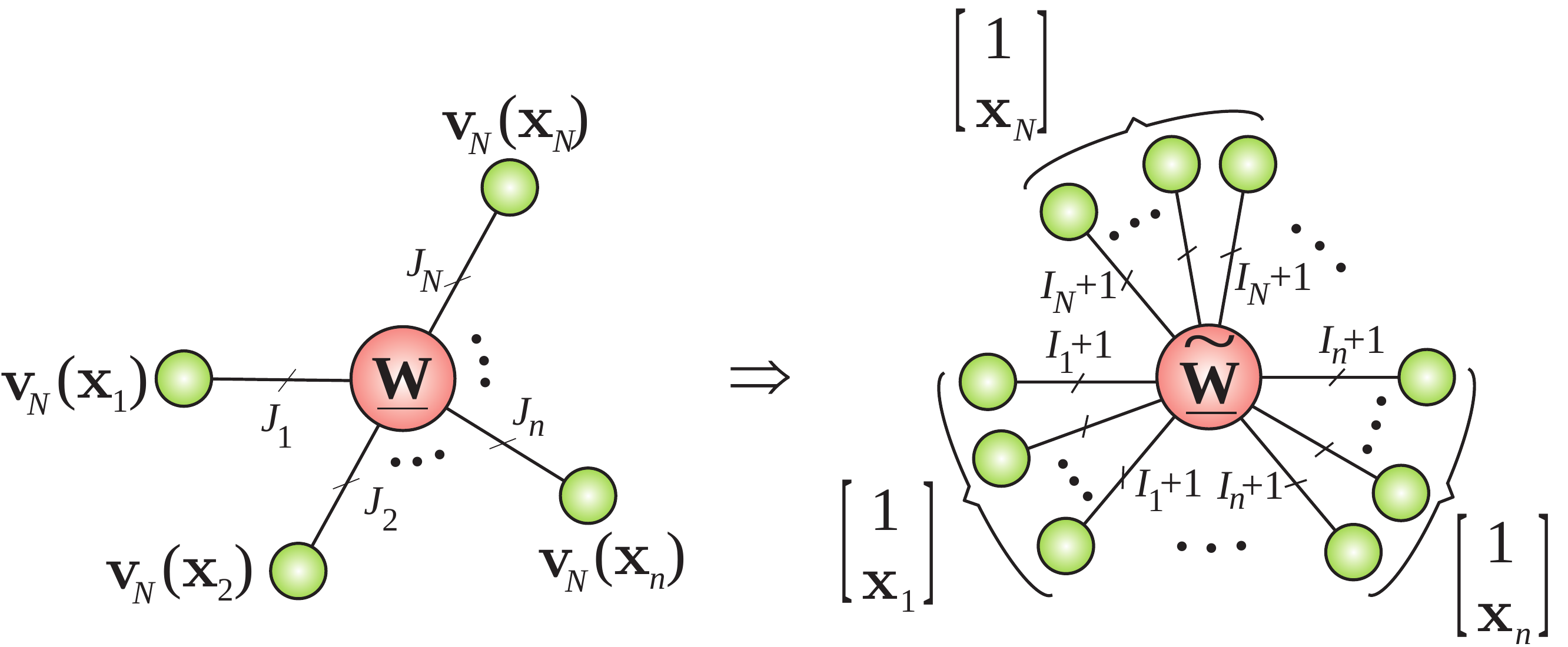}	
\caption{Graphical illustration of the vector-variate regression. (a) The vector-variate regression for multiple input units $\bx_1$, \ldots, $\bx_N$, where the nonlinear function $h(\bx_1, \ldots, \bx_N)$ is expressed as a tensor-vector product of an $N$th-order core tensor, $\tW$, of size $J_1 \times J_2 \times \cdots \times J_N$, and Vandermonde-like vectors $\bv_N(\bx_n)$ of length $J_n$, where $J_n = \frac{I_n^{N+1}-1}{I_n - 1}$. (b) An equivalent regression model but with an $N^2$th-order tensor of   size $(I_1+1)\times \cdots \times (I_1+1) \times (I_2+1)\times  \cdots \times (I_N+1) \times \cdots \times (I_N+1)$.
When the input units are scalars, the tensor $\widetilde{\tW}$ is of size $2\times 2 \times \cdots \times 2$.}
\label{fig_boltzmann}
\end{figure}

In general, the regression for a system with $N$ input vectors, $\bx_n$ of lengths $I_n$, can be written as
\begin{align}
	h(\bx_1, \ldots, \bx_N) = w_0 +
	 \sum_{d = 1}^{N^2} \sum_{i_1, i_2, \ldots, i_d = 1}^{N}   \tW_{i_1, i_2, \ldots, i_d}  \bar{\times}  \,  (\bx_{i_1} \circ  \bx_{i_2} \, \circ  \cdots \circ   \bx_{i_d} )  \, ,\label{eq_tenreg_N}
\end{align}
where $\bar{\times}$ represents the inner product between two tensors, and the tensors $\tW_{i_1, \ldots, i_d}$ are of $d$-th order, and of size $I_{i_1} \times I_{i_2} \times \cdots \times I_{i_d}$,  $d = 1, \ldots, N^2$.
The representation of the generalised model as a tensor-vector product of an $N$th-order tensor of size $J_1 \times J_2 \times \cdots \times J_N$, where $J_n = \frac{I_n^{N+1} -1}{I_n-1}$, comprising all the weights, is given by
\be
	h(\bx_1, \ldots, \bx_N)  = \tW \, \bar{\times}_1 \, \bv_N(\bx_1)  \, \bar{\times}_2 \, \bv_N(\bx_2)  \, \cdots \, \bar{\times}_N \, \bv_N(\bx_N),  \label{eq_tenreg_N2}
\ee
where
\be
 \bv_N(\bx) = \left[ \begin{array}{ccccc} 1& \bx^\tp & (\bx \otimes \bx)^T & \ldots & (\bx \otimes \cdots \otimes \bx)^T \end{array} \right]^{\text{T}}\, ,
\ee
or,  in a more compact form, with a very high-order tensor $\widetilde{\tW}$ of $N^2$th-order and of size $(I_1+1)\times \cdots \times (I_1+1) \times (I_2+1) \times \cdots \times (I_N+1) \times \cdots \times (I_N+1)$, as
\be
h(\bx_1, \ldots, \bx_N)  = \widetilde{\tW} \, \bar{\times}_{1:N} \, \left[\begin{array}{c}1   \\  \bx_1 \end{array}\right]   \,  \cdots \bar{\times}_{N(N-1)+1:N^2} \,   \left[\begin{array}{c}1   \\  \bx_N \end{array}\right]\, . \label{eq_tenreg_NN}
\ee
The illustration of this generalized model is given in Figure~\ref{fig_boltzmann}.

{\noindent \bf Tensor-variate  model.} When the observations are matrices, $\bX_n$, or higher-order tensors, $\tX_n$, the models in (\ref{eq_tenreg_N}), (\ref{eq_tenreg_N2}) and (\ref{eq_tenreg_NN}) are still applicable and operate by replacing the original vectors,  $\bx_n$, by the vectorization of the higher-order  inputs. This is because the inner product between two tensors can be expressed as a product of their two vectorizations.

{\noindent \bf Separable representation of the weights.} Similar to the MPR, the challenge in the generalised tensor-variate regression is the curse of dimensionality of the weight tensor $\tW$ in (\ref{eq_tenreg_N2}), or of the tensor $\widetilde{\tW}$ in (\ref{eq_tenreg_NN}).

A common method to deal with the problem is to restrict the model to some low order, i.e., to the first order. 
The weight tensor is now only of size $(I_1+1) \times (I_2+1) \times \cdots \times (I_N+1)$. The large weight tensor can then be represented in the canonical form \citep{nguyen2015tensor,qi2016matrix},
 the TT/MPS tensor format \citep{Stoudenmire_optim}, or the hierarchical Tucker tensor format  \citep{Cohen_IB}.

\section{Tensor Structure in Volterra Models of Nonlinear Systems}
\sectionmark{Tensor Structure in Volterra Models}

\subsection{Discrete Volterra Model}

System identification is a paradigm which aims to provide a mathematical description of a system from the observed system inputs and outputs  \citep{billings2013nonlinear}.
In practice, tensors are inherently present in {\emph{Volterra operators}} which model the system response of a nonlinear system which maps an input signal $x(t)$ to an output signal $y(t)$ in the form
\be
y(t)  = V(x(t)) =  h_0 + H_1(x(t)) + H_2(x(t)) + \cdots + H_n(x(t)) + \cdots    \notag
\ee
where $h_0$ is a constant and $H_n(x(t))$ is the $n$th-order Volterra operator, defined as a generalised convolution of the integral {\emph{Volterra kernels}} $h^{(n)}(\tau_1, \ldots, \tau_n)$ and the input signal, that is
\be
	H_n(x(t)) = \int h^{(n)}(\tau_1, \ldots, \tau_n) x(t-\tau_1) \cdots x(t-\tau_n) d\tau_1 \cdots d\tau_n \,. \quad \ee
The system, which is assumed to be time-invariant and continuous, is treated as a black box, and needs to be represented by appropriate Volterra operators.

\begin{figure}[t!]
\centering
\includegraphics[width=.75\textwidth, trim = 0.0cm 0cm 0cm 0cm,clip=true]{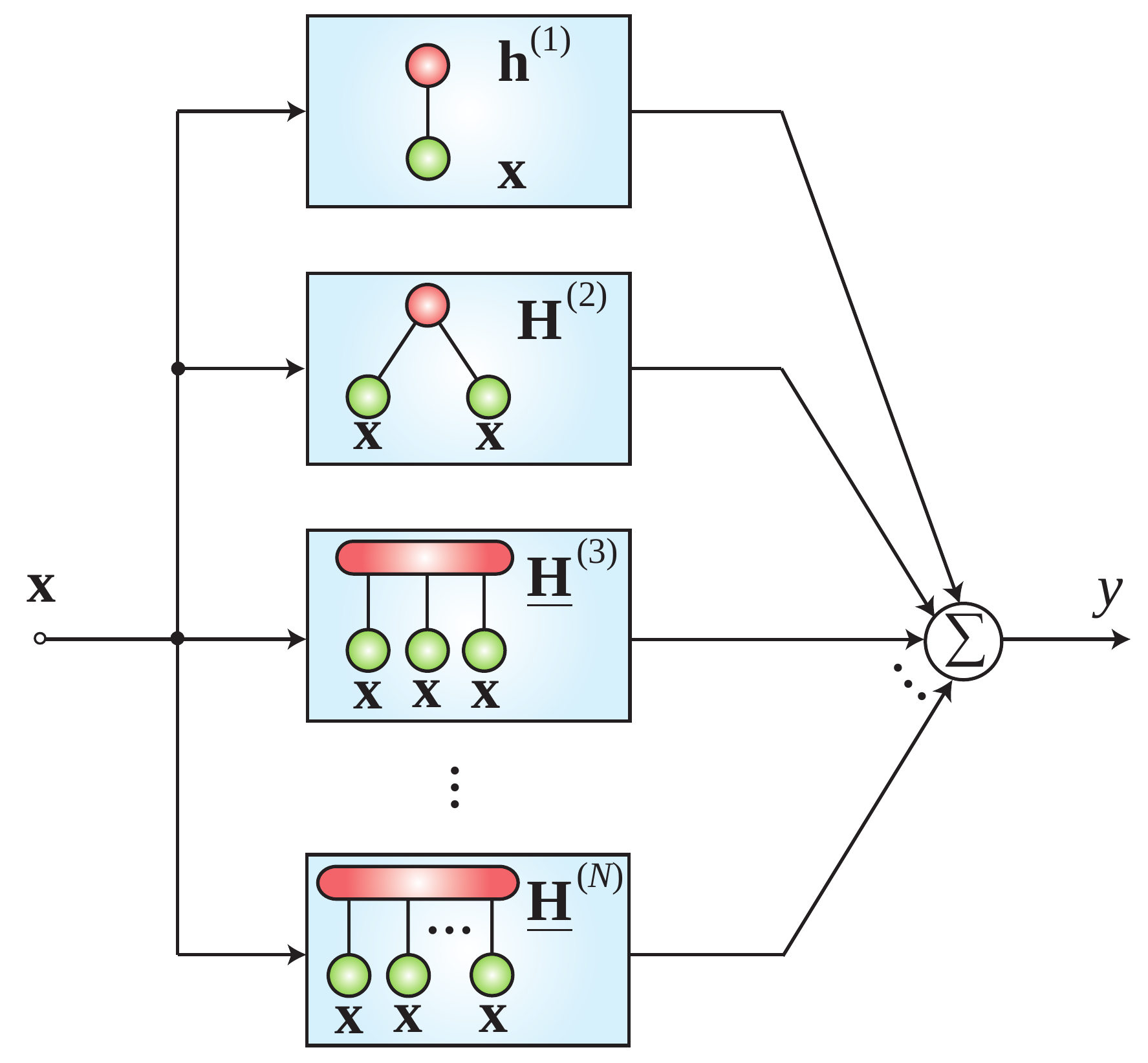}\\
\vspace{0.4cm}
\caption{A Volterra model of a nonlinear system with  memory of length $M$. Each block computes the tensor product between an $n$th-order Volterra kernel, $\tH^{(n)}$, and the vector $\bx$ of length $M$, which comprises $M$ samples of the input signal. The system identification task amounts to estimating the Volterra kernels, $\tH^{(n)}$, directly or in  suitable tensor network formats.}
\label{fig_Volterra}
\end{figure}

In practice, for a finite duration sample input data, $\bx$, the discrete system can be modelled using truncated Volterra kernels of size $M \times M \times \cdots \times M$, given by
\be
	H_n(\bx) &=& \sum_{i_1 = 1}^{I}  \cdots \sum_{i_n = 1}^{I} h^{(n)}_{i_1, \ldots, i_n} x_{i_1} \ldots x_{i_2}  \notag\\
	&=& \tH^{(n)}  \, \bar{\times}_1  \bx  \,  \bar{\times}_2  \bx   \, \cdots \, \bar{\times}_n  \bx. 
\ee
For simplicity, the Volterra kernels $\tH^{(n)} = [h^{(n)}_{i_1, \ldots, i_n}]$ are assumed to have the same size $M$ in each mode, and, therefore, to yield a symmetric tensor. Otherwise, they can be symmetrized. 

{\noindent \bf Curse of dimensionality.}
The output which corresponds to the input $\bx$ is written as a sum of $N$ tensor products (see in Figure~\ref{fig_Volterra}), given by
\be
	y = h_0  +  \sum_{n = 1}^{N}   \tH^{(n)}  \, \bar{\times}_1  \,\bx \, \bar{\times}_2 \, \bx \, \cdots \, \bar{\times}_n \, \bx    .
\ee
Despite the symmetry of the Volterra kernels, $\tH^{(n)}$, the number of actual coefficients of the $n$th-order kernel to be estimated is still huge, especially for higher-order kernels, and is given by
$
\frac{(M+n-1)!}{n! (M-1)!}$.
As a consequence, the estimation requires a large number of measures (data samples), so that the method  for a raw tensor format is only feasible for systems with a relatively small memory and low-dimensional input signals.

\subsection{Separable Representation of Volterra Kernel}

In order to deal with the curse of dimensionality in Volterra kernels, we  consider the kernel $\tH^{(n)}$ to be separable, i.e., it can be expressed in some low rank tensor format, e.g., as a  CP tensor or in any other suitable tensor network format (for the concept of general separability of variables, see  Part 1).

{\noindent \bf Volterra-CP model.} The first and simplest separable Volterra model, proposed in \citep{ACS:ACS1272}, represents the kernels by symmetric tensors of rank $R_n$ in the CP format, that is
\be
	\tH^{(n)} = \tI \times_1 \bA_n \times_2 \bA_n \cdots \times_n \bA_n \,.
\ee
For this tensor representation, the identification problem simplifies into the estimation of $N$
factor matrices, $\bA_n$, of size $M \times R_n$ and an offset, $h_0$, so that the number of parameters reduces to $M\sum_{n} R_n+1$ (note that $R_1 = 1$).
Moreover, the implementation of the Volterra model becomes
\be
	y_k = h_0 + \sum_{n = 1}^{N}   (\bx_k^{\text{T}} \, \bA_n)^{\bigcdot n} \, \1_{R_n} \,,
\ee
where $\bx_k = [x_{k-M+1},  \ldots, x_{k-1}, x_k]^{\text{T}}$ comprises $M$ samples of the input signal, and $(\bigcdot)^{\bigcdot n}$ represents the element-wise power operator.
The entire output vector $\by$ can be computed in a simpler way through the convolution of the input vector $\bx$ and the factor matrices $\bA_n$, as \citep{7827594}
\be
\by =  h_0 + \sum_{n=1}^{N} (\bx \ast \bA_n)^{\bigcdot n} \,  \, \1_{R_n} \, .
\ee

{\noindent \bf Volterra-TT model.}   Alternatively, the Volterra kernels, $\tH^{(n)}$, can be represented in the TT-format, as
\be
\tH^{(n)} = \llangle \tG_{n}^{(1)}  ,   \tG_{n}^{(2)} , \ldots ,     \tG_{n}^{(n)} \rrangle\,.
\ee
By exploiting the fast contraction over all modes between a TT-tensor and $\bx_k$, we have
\be
 \tH^{(n)} \bar{\times} \bx_k  = (\tG_{n}^{(1)}   \bar{\times}_2 \bx_k) (\tG_{n}^{(2)}  \bar{\times}_2 \bx_k) \cdots   (\tG_{n}^{(n)} \bar{\times}_2 \bx_k) \, . \notag
\ee
 The output signal, can be then computed through the convolution of the core tensors and the input vector, as
\be
y_k = h_0 + \sum_{n = 1}^{N}  \tZ_{n,1}(1,k,:) \, \tZ_{n,2}(:,k,:)  \, \cdots \, \tZ_{n,n-1}(:,k,:) \tZ_{n,n}(:,k) \,, \notag
\ee
where $\tZ_{n,m} = \tG_{n}^{(m)} \ast_2 \bx$ is a mode-2 partial convolution of the input signal $\bx$ and the core tensor $\tG_{n}^{(m)}$, for $m = 1, \ldots, n$.
A similar method, but with only one TT-tensor, is considered in   \citep{BatselierCW16}.

\subsection{Volterra-based Tensorization for Nonlinear Feature Extraction}
\subsectionmark{Volterra-based Tensorization}

Consider  nonlinear feature extraction in a supervised learning system, such that the extracted features maximize the Fisher score \citep{KumarVolterraface}. In other words, for a data sample $\bx_k$, which can be a recorded signal in one trial or a vectorization of an image, a feature extracted from $\bx_k$ by a nonlinear process is denoted by $y_k = f(\bx_k)$. Such constrained (discriminant) feature extraction can be treated as a maximization of the Fisher score
\be
	\max \frac{\sum_{c} (\bar{y}_{c} - \bar{y})^2}{\sum_{k} (y_{k} - \bar{y}_{c_k})^2}  \,, \label{eq_disc_gen}
\ee
where $\bar{y}_{c_k}$ is the mean feature of the samples in class-$k$, and $\bar{y}$ the mean feature of all the samples.

Next, we model the nonlinear system $f(\bx)$ by a truncated Volterra series representation
\be
y_k &=& \sum_{n = 1}^{N} \tH^{(n)}   \bar{\times} (\bx_k \circ \bx_k \circ \cdots \circ \bx_k)   = \bh^{\text{T}} \, \bz_k   \,,\label{eq_volterra_feature}
\ee
where $\bh$ and $\bx_k$ are vectors comprising all coefficients of the Volterra kernels and
\be
	\bh &=& [\vtr{\bH^{(1)}}^{\text{T}}, \vtr{\bH^{(2)}}^{\text{T}}, \ldots, \vtr{\bH^{(N)}}^{\text{T}}]^{\text{T}} \notag \, , \\
	{\bz}_k &=& [\bx_k^{\text{T}}, (\bx_k^{\otimes_{\,2}})^{\text{T}}, \ldots,  (\bx_k^{\otimes_{\,N}})^{\text{T}}]^{\text{T}} \notag\,.
\ee
The shorthand $\bx^{\otimes_{\,n}} = \bx \otimes \bx \otimes \cdots \otimes \bx$ represents the Kronecker product of $n$ vectors $\bx$.
The offset coefficient, $h_0$, is omitted in the above Volterra model because it will be eliminated in the objective function (\ref{eq_disc_gen}). The vector $\bh$ can be shortened by keeping only distinct coefficients, due to symmetry of the Volterra kernels. The augmented sample $\bz_k$ needs a similar adjustment but multiplied with the number of occurrences.

 Observe that the nonlinear feature extraction, $f(\bx_k)$, becomes a linear mapping, as in (\ref{eq_volterra_feature}) after $\bx_k$ is tensorized into ${\bz}_k$.
Hence, the nonlinear discriminant in (\ref{eq_disc_gen}) can be rewritten in the form of a standard linear discriminant analysis
 \be
	\max \frac{\bh^{\text{T}}  \bS_b \bh}{\bh^{\text{T}} \bS_w \bh}  \, ,
\label{eq_disc_gen2}
 \ee
 where $\bS_b = \sum_c (\bar{\bz}_c  -  \bar{\bz}) (\bar{\bz}_c  -  \bar{\bz})^{\text{T}}$ and
 $\bS_w = \sum_k ({\bz}_k  -  \bar{\bz}_{c_k}) ({\bz}_k  -  \bar{\bz}_{c_k})^{\text{T}}$ are respectively between- and within-scattering matrices of $\bz_k$. The problem then boils down to finding generalised principal eigenvectors of $\bS_b$ and $\bS_w$.

 {\noindent \bf Efficient implementation.} The problem with the above analysis is that the length of eigenvectors, $\bh$, in (\ref{eq_disc_gen2})  grows exponentially with the data size, especially for higher-order Volterra kernels. To this end, \citet{KumarVolterraface} suggested  to split the data into small patches. Alternatively, we can impose low rank-tensor structures, e.g., the CP or TT format, onto the Volterra kernels, $\bH^{(n)}$, or the entire vector $\bh$.

 \section[{Tensor Representations of Sinusoid Signals and BSS}]{Low-rank Tensor Representations of Sinusoid Signals and their Applications to BSS and Harmonic Retrieval}
  \sectionmark{Tensor Representations of Sinusoid Signals and BSS}
  \label{sec:sinusoid}

Harmonic signals are fundamental in many practical applications. This section addresses low-rank structures of sinusoid signals under several tensorization methods. These properties can then be exploited in the blind separation of sinusoid signals or their modulated variants, e.g., the exponentially decaying  signals, the examples of which are
\begin{align}
	  x(t) &=   \sin(\omega \,  t + \phi)   \,,        &x(t) &= t \sin(\omega \,  t + \phi)\, \label{eq_t_sin} \, , \\
       x(t) &= \exp(-\gamma t)  \sin(\omega \,  t + \phi)\,,      &x(t)  &=  t \exp(-\gamma t)  \label{eq_t_exp}\,,
\end{align}
for $t = 1, 2, \ldots, L$, $\omega \neq 0$.

\subsection{Folding - Reshaping of  Sinusoid}

{\noindent \bf Harmonic matrix.} The harmonic matrix $\bU_{\omega, I}$ is a matrix of size $I \times 2$ defined over the two variables, the angular frequency $\omega$ and the folding size $I$, as
\be
\bU_{\omega, I} &=&  \left[
	\begin{array}{cc}
	1 & 0 \\
	\vdots & \vdots \\
	\cos(k\omega) & \sin(k\omega)\\
	\vdots & \vdots \\
	\cos((I-1)\omega) & \sin((I-1)\omega) \\
	\end{array}
	\right] \, . \label{eq_UwI}
\ee

{\noindent \bf Two-way folding.}
A matrix of size $I \times J$, folded from a sinusoid signal $x(t)$ of length $L = IJ$, is of rank-2, and can be decomposed as
\be
\bY &=&   \bU_{\omega,I} \, \bS \, \bU_{\omega I, J}^{\text{T}}  \, ,\label{lem_mat_sin}
\ee
where $\bS$ is invariant to the folding size $I$, depends only on the phase $\phi$, and takes the form
\be
	\bS   = \left[\begin{array}{cc}
		\sin(\phi)  & \cos(\phi) \\
		\cos(\phi) & -\sin(\phi)
	\end{array}	\right] . \label{eq_S}
\ee
{\noindent \bf Three-way folding.}
A third-order tensor of size $I\times J\times K$, where $I,J,K>2$, reshaped from a sinusoid signal of length $L$, can take the form of a multilinear rank-(2,2,2) or rank-3 tensor
\be
\tY =  \llbracket \tH;  \bU_{\omega, I}  ,  \bU_{\omega I, J}  , \bU_{\omega IJ,K}	 \rrbracket \, , 	 \label{eq_Y_rnk222}\label{lem_sinu_T_tensor3}
\ee
where $\tH = \tG \times_3 \bS$ is a small-scale tensor of size $2 \times 2 \times 2$, and
\be
\tG(:,:,1) &=& \left[
	\begin{array}{cc}
		1  & 0 \\
		0 & -1
	\end{array}		
	\right],   \qquad
\tG(:,:,2)  =  	\left[
	\begin{array}{cc}
		0 & 1 \\
		1 & 0
	\end{array}		
	\right] 	\,  . \label{eq_G}
\ee	
The above expression can be derived by folding the signal $y(t)$ two times.
%
We can prove by contradiction that the so-created core tensor $\tG$ does not have rank-2, but  has the following rank-3 tensor representation
\be
\tG = \frac{1}{2}
\left[
	\begin{array}{@{}c@{}}
		1  \\
		1
	\end{array}		
	\right]  \circ
\left[
	\begin{array}{@{}c@{}}
		1  \\
		1
	\end{array}		
	\right]
\circ
\left[
	\begin{array}{@{}c@{}}
		1  \\
		1
	\end{array}		
	\right] 	 	
	 -
	 \frac{1}{2}
\left[	 \begin{array}{@{}c@{}}
		-1  \\
		1
	\end{array}		
	\right]  \circ
\left[
	\begin{array}{@{}c@{}}
		-1  \\
		1
	\end{array}		
	\right]
\circ
\left[
	\begin{array}{@{}c@{}}
		-1  \\
		1
	\end{array}		
	\right]
	+
	2
	\left[\begin{array}{@{}c@{}}
		0  \\
		1
	\end{array}		
	\right]  \circ
\left[
	\begin{array}{@{}c@{}}
		0  \\
		1
	\end{array}		
	\right]
\circ
\left[
	\begin{array}{@{}c@{}}
		-1  \\
		0
	\end{array}		
	\right] \, .\notag
\ee
Hence, $\tY$ is also a rank-3 tensor. Note that $\tY$ does not have a unique rank-3 decomposition.

\begin{remark}\label{rem_sin_folding3}
The Tucker-3 decomposition in (\ref{eq_Y_rnk222}) has a fixed core tensor $\tG$, while the factor matrices are identical for signals of the same frequency.
\end{remark}

{\noindent \bf {Higher-order folding - TT-representation.}}
An $N$th-order tensor of size $I_1 \times I_2 \times \cdots \times I_N$, where $I_n \ge 2$, which is reshaped from a sinusoid signal, can be represented by a multilinear rank-(2,2,\ldots, 2) tensor
\be
\tY = \llbracket \tH ;  \bU_{\omega, I_1} ,  \bU_{\omega J_1,  I_2}  , \ldots,  \bU_{\omega J_{N-1}, I_N}	 \rrbracket \,, 	 \label{eq_Y_rnk2_2}
\ee
where $\tH =  \llangle \underbrace{\tG,\tG, \ldots, \tG}_{(N-2) \text{terms}},\bS \rrangle $ is an $N$th-order tensor of size $2 \times 2 \times \cdots \times 2$, and $J_n = \prod_{k = 1}^{n} I_k$.
%

\begin{remark}[TT-representation]
Since the tensor $\tH$ has  TT-rank of (2,2,\ldots,2), the folding tensor $\tY$ is also a tensor in TT-format of rank-(2,2,\ldots,2), that is
\be
\tY =  \llangle \tA_1, \tA_2, \,  \ldots ,    \tA_N \rrangle \, ,
\ee
where $\tA_1 = \bU_{\omega, I_1}$, $\tA_N = \bS \bU_{\omega J_{N-1}, I_N}^{\text{T}}$ and  $\tA_n =  \tG \times_2 \, \bU_{\omega J_{n-1},I_n}$ for $n = 2, \ldots, N-1$.
\end{remark}

\begin{remark}[QTT-Tucker representation]
When the folding sizes $I_n = 2$, for $n = 1, \ldots, N$, the representation of the folding tensor $\tY$ in (\ref{eq_Y_rnk2_2}) is also known as the \emph{QTT-Tucker format}, given by
\be
\tY = \llbracket \tH;  \bA_1 , \ldots , \bA_{N-1}, \bA_{N} \rrbracket ,
\ee
where $\bA_n = \left[
	\begin{array}{@{}cc@{}}
		1 & 0 \\
		\cos(2^{n-1} \omega) & \sin(2^{n-1} \omega)
	\end{array}		
	\right]$.
\end{remark}

\begin{example}{\bf Separation of damped sinusoid signals.}\label{ex_TT_separation_exp}

This example demonstrates the use of multiway folding in a single channel separation of damped sinusoids.
We considered a vector composed of $P$ damped sinusoids,
\be
y(t) &=& \sum_{p = 1}^{P}  a_p  \, x_p(t) + n(t) \, ,
\ee
where
\be
x_p(t) = \exp({\frac{-5t}{L p}}) \, \sin(\frac{2\pi f_p}{f_s} t  + \frac{(p-1)\pi}{P}) \, ,\notag
\ee
with frequencies $f_p$ = 10, 12 and 14 Hz, and the sampling frequency  $f_s = 10 f_P$. Additive Gaussian noise, $n(t)$, was  generated at a specific signal-noise-ratio   (SNR). The weights, $a_p$, were set such that the component sources were equally contributing to the mixture, i.e., $a_1 \|\bx_1\| = \cdots = a_P \|\bx_P\|$, and
the signal length was $L =  2^{d} \, P^2$.
\begin{figure}[t!]
\centering
{\includegraphics[width=.65\textwidth, trim = 0.0cm 0cm 0cm 0cm,clip=true]{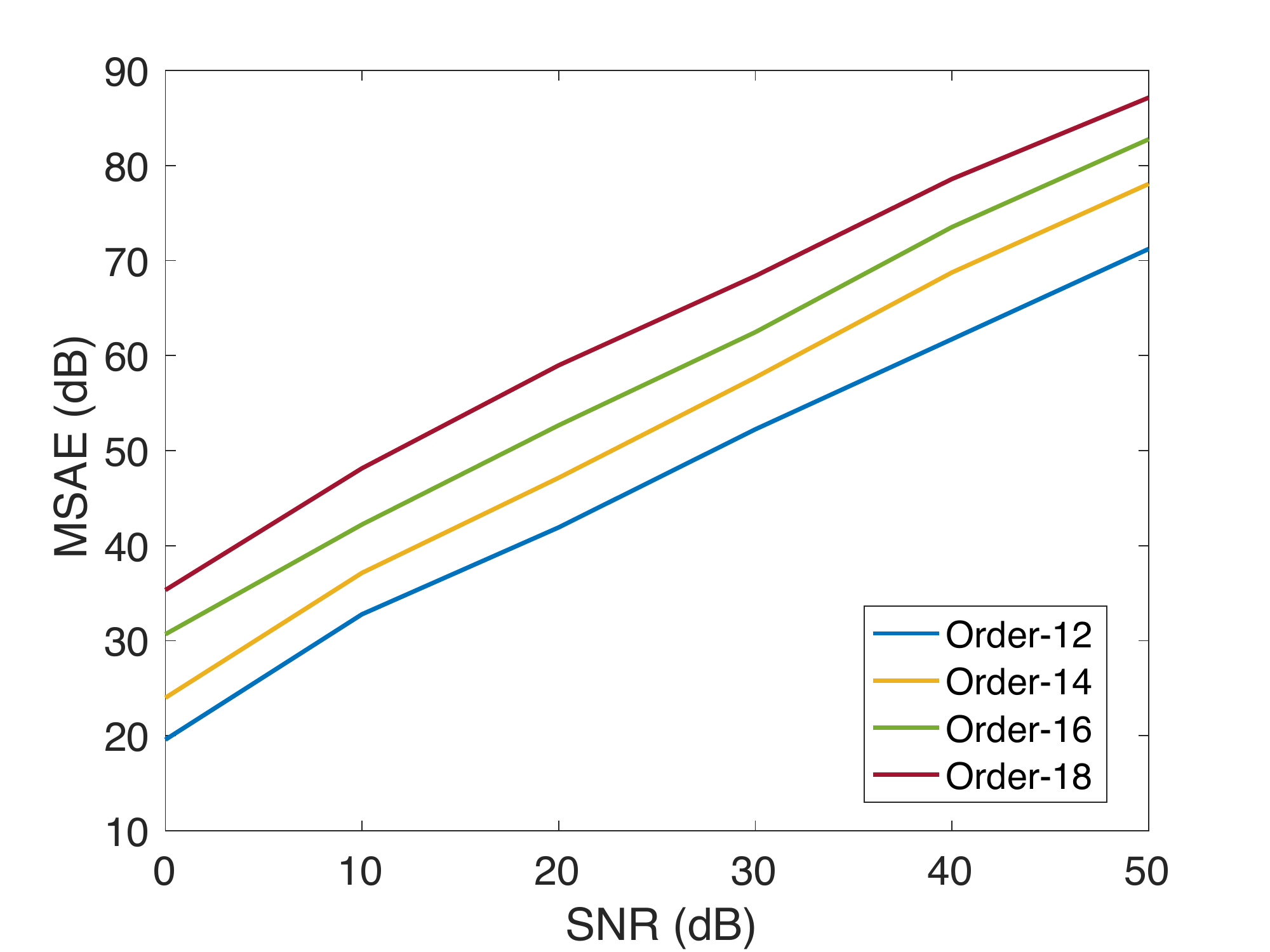}
\label{fig_TTseparation_R3_exps_SAE}}
\hfill
\caption{Comparison of the mean SAEs for various noise levels SNR, signal lengths, and tensor orders.}
\label{fig_ex_TT_shortsignal_exps}
\end{figure}

In order to separate the three signals $x_p(t)$ from the mixture $y(t)$,  we tensorized the mixture to a $d$th-order tensor of size $2R \times 2 \times \cdots \times 2 \times 2R$. Under this tensorization, the exponentially decaying signals $\exp(\gamma t)$ yielded rank-1 tensors, while  according to (\ref{eq_Y_rnk2_2}) the sinusoids have TT-representations of rank-$(2,2,\ldots, 2)$.
Hence, the tensors of $x(t)$ can also be represented by tensors in the TT-format of rank-$(2,2,\ldots, 2)$.
We were, therefore, able to approximate $\tY$ as a sum of $P$ TT-tensors $\tX_r$ of rank-$(2,2,\ldots, 2)$, that is, through the minimization \citep{Phan_TT_part1}
\be
\min  \quad \|\tY - \tX_1 - \tX_2 - \cdots - \tX_P\|_F^2 \, .\label{eq_fit_RtermsTT}
\ee
For this purpose, a tensor  $\tX_p$ in a TT-format was fitted sequentially to the residual $\tY_p =  \tY - \sum_{s \neq p} \tX_s$, calculated by the difference between the data tensor $\tY$ and its approximation by the other TT-tensors $\tX_s$ where $s \neq p$, that is,
\be
\argmin_{\tX_p} \|\tY_{p} - \tX_p \|_F^2	\, ,\label{eq_fit_Yr}
\ee
for $p = 1, \ldots, P$.
Figure~\ref{fig_ex_TT_shortsignal_exps} illustrates the mean SAEs (MSAE) of the estimated signals for various noise levels SNR = 0, 10, \ldots, 50 dB, and different signal lengths $K = 9 \times 2^d$, where $d = 12, 14, 16, 18$.

On average, \emph{an improvement of 2 dB SAE is achieved if the signal is two times longer.}
If the signal has less than $L = 9 \times 2^6 = 576$ samples, the estimation quality will deteriorate by about 12 dB compared to the case when signal length of $L = 9 \times 2^{12}$.
For such cases, we suggest to augment the signals using other tensorizations before performing  the source extraction, e.g., by construction of multiway Toeplitz or Hankel tensors. Example~\ref{ex_TT_separation_short_exp} further illustrates the separation of short length signals.
\end{example}

\subsection{Toeplitz Matrix and Toeplitz Tensors of Sinusoidal Signals}
\subsectionmark{Toeplitz Tensors of Sinusoidal Signals}

{\noindent \bf Toeplitz matrix of sinusoid.}
The Toeplitz matrix, $\bY$, of a sinusoid signal, $y(t)  = \sin(\omega \,  t + \phi)$,
is of rank-2 and can be decomposed as
\be
\bY
=
\left[
\begin{array}{cccc}
y(1) & y(2) \\
y(2) & y(3)\\
\vdots &   \vdots\\
y(I)  & y(I+1)
\end{array}
\right]
\bQ_T
\left[
\begin{array}{cccc}
y(I) & \cdots & y(L) \\
y(I-1) & \cdots & y(L-1) \\
\end{array}
\right] \,,\;\label{lem_Toeplitz_sin}
\ee
where $\bQ_T$ is invariant to the selection of folding length $I$, and has the form
\be
\bQ_T = \frac{1}{\sin^2(\omega)} \,\left[
\begin{array}{cc}
-y(3) & y(2)\\
y(2)  &   -y(1)
\end{array}
\right] \,\, . \label{equ_toeplitz2_Qt}
\ee
The above expression follows from the fact that
\be
\left[\begin{array}{@{\hspace{.1ex}}c@{\hspace{1ex}}c@{\hspace{.1ex}}}
y(i) & y(i+1)
\end{array}
\right] \,\left[
\begin{array}{@{\hspace{.1ex}}c@{\hspace{1ex}}c@{\hspace{.1ex}}}
-y(3) & y(2)\\
y(2)  &   -y(1)
\end{array}
\right] \,
\left[\begin{array}{@{\hspace{1ex}}c@{\hspace{1ex}}c@{\hspace{1ex}}cc}
y(j) \\ y(j-1)]
\end{array}
\right]
= \sin^2(\omega)  \, y(j-i+1)
\notag
\, .
\ee

{\noindent \bf Toeplitz tensor of sinusoid.}
An $N$th-order Toeplitz tensor, tensorized from a sinusoidal signal, has a TT-Tucker representation
\be
\tY =  \llbracket  \tG;  \bU_1, \ldots,   \bU_{N-1} ,  \bU_N \rrbracket \label{eq_toep_N_tt_tucker}
\ee
where  the factor matrices $\bU_n$ are given by
\begin{align}
\bU_1 &=  \left[
	\begin{array}{@{}cc@{}}
	y(1) & y(2) \\
	\vdots & \vdots \\
	y(J_1) &  	y(J_1+1)
	\end{array}		
	\right]
\, , \quad
\bU_N =  \left[
	\begin{array}{@{}cc@{}}
	y(J_{N-1}-1) &  y(J_{N-1}-2)\\
	\vdots & \vdots \\	
	y(L) 	     & y(L-1)
	\end{array}
	\right] \,  ,   \notag \\
\bU_{n} &=  \left[
	\begin{array}{@{}cc@{}}
	y(J_{n-1}) & y(J_{n-1}+1)\\
	\vdots & \vdots \\
	y(J_{n}-1) &  y(J_{n})
		\end{array}
	\right]\, ,  \quad n = 2, \ldots, N-1 \,,\label{eq_toeplitzN_factor2}
\end{align}
in which $J_n = I_1 + I_2 + \cdots + I_n$.
The core tensor $\tG$ is an $N$th-order tensor of size $2 \times 2 \times \cdots \times 2$, in a TT-format, given by
\be
\tG = \llangle \tG^{(1)},  \tG^{(2)},   \ldots,   \tG^{(N-1)} \rrangle, \label{eq_core_toeplitzN}
\ee
where $\tG^{(1)} =   \bT(1)$ is a matrix of size $1 \times 2 \times 2$, while the core tensors $\tG^{(n)}$, for $n = 2,\ldots, N-1$, are of size $2\times 2 \times 2$ and have two horizontal slices, given by
\begin{align}
\tG^{(n)}(1,:,:)  &=  \bT(J_{n-1}-n+2)\, ,  \quad \tG^{(n)}(2,:,:)  =  \bT(J_{n-1}-n+1) \, 	  , \notag
\end{align}
with
\be
\bT(I) = \frac{1}{\sin^2(\omega)} \, \left[
	\begin{array}{@{}cc@{}}
		-y(I+2) &  y(I+1)\\
		y(I+1)  &  -y(I)
	\end{array}		
	\right] \label{eq_T}\,.
\ee
Following the two-stage Toeplitz tensorization, and upon applying (\ref{lem_Toeplitz_sin}), we can deduce the decomposition in (\ref{eq_toep_N_tt_tucker}) from that for the $(N-1)$th-order Toeplitz tensor.

\begin{remark}
For second-order tensorization, the core tensor $\tG$ in (\ref{eq_core_toeplitzN}) comprises only $\bG^{(1)}$, which is identical to the matrix $\bQ_T$ in (\ref{equ_toeplitz2_Qt}).
\end{remark}

{\noindent \bf Quantized Toeplitz tensor.} An $(L-1)$th-order Toeplitz tensor of a sinusoidal signal of length $L$ and  size $2 \times 2 \times \cdots \times 2$  has a TT-representation with $(L-3)$ identical core tensors $\tG$, in the form
\be
\tY = \llangle \tG, \tG , \ldots , \tG,    \left[
	\begin{array}{@{}cc@{}}
		y(L-1) &  y(L)\\
		y(L-2)  &  y(L-1)
	\end{array}		
	\right]  \rrangle \,,   \notag
\ee
where
\be
	\tG(1,:,:) =  \left[
	\begin{array}{@{}cc@{}}
		1 & 0\\
		0 & 1
	\end{array}		
	\right]
, \quad
	\tG(2,:,:) =  \left[
	\begin{array}{@{}cc@{}}
		0 & 1\\
		-1 & 2\cos(\omega)
	\end{array}		
	\right] \,. \notag
\ee

\begin{example}{\bf Separation of short-length damped sinusoid signals.}\label{ex_TT_separation_short_exp}

This example illustrates the use of  Toeplitz-based tensorization in the separation of damped sinusoid signals from a short-length observation.
We considered a single signal composed by $P = 3$ damped sinusoids of length $L =  66$, given by
\be
y(t) &=& \sum_{p = 1}^{P}  a_p  \, x_p(t) + n(t) \,,
\ee
where
\be
x(t) = \exp(\frac{-pt}{30}) \, \sin(\frac{2\pi f_p}{f_s} t  + \frac{p\pi}{7})\,
\ee
with frequencies $f_p$  = 10, 11 and 12 Hz, the sampling frequency  $f_s$ = 300 Hz, and the mixing factors $a_p = p$. Additive Gaussian noise $n(t)$ was generated at a specific signal-noise-ratio.

\begin{figure}[t!]
\centering
{\includegraphics[width=.65\textwidth, trim = 0.0cm 0cm 0cm 0cm,clip=true]{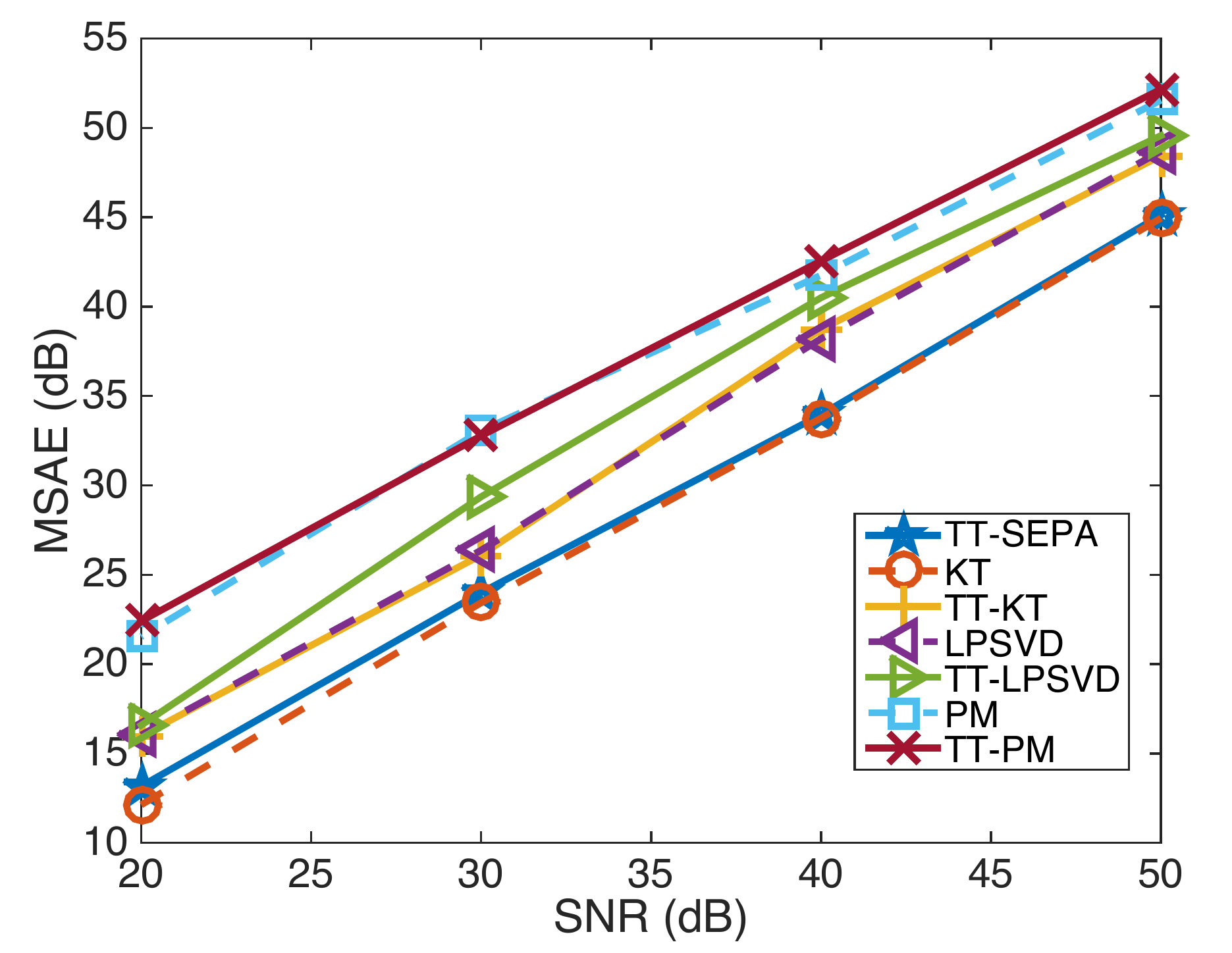}
\label{fig_TTseparation_R3_shortexps_SAE}}
\caption{
Mean SAEs (MSAE) of the estimated signals in Example~\ref{ex_TT_separation_short_exp}, for various noise levels SNR.}
\label{fig_ex_TT_shortsignal_exps2}
\end{figure}

In order to separate the three signals, $x_p(t)$, from the mixture $y(t)$,  we first tensorized the observed signal to a 7th-order Toeplitz  tensor of size $16 \times 8 \times 8 \times 8 \times 8 \times   8 \times 16$, then folded this tensor to a 23th-order tensor of  size $2 \times 2 \times \cdots  \times 2$.
With this tensorization, according to (\ref{eq_toep_N_tt_tucker}) and (\ref{eq_Y_rnk2_2}), each damped sinusoid $x_p(t)$ had a TT-representation of rank-$(2,2,\ldots, 2)$.
The result produced by minimizing the cost function (\ref{eq_fit_RtermsTT}), annotated by TT-SEPA, is shown  in Figure~\ref{fig_ex_TT_shortsignal_exps2} as a solid line with star marker.
The so obtained performance was much better than in Example~\ref{ex_TT_separation_exp}, even for the signal length of only 66 samples.

We note that the parameters of the damped signals can be estimated
using linear self-prediction (auto-regression) methods, e.g., singular value decomposition of the Hankel-type matrix as in the Kumaresan-Tufts (KT) method \citep{Kumaresan}. As shown in Figure~\ref{fig_ex_TT_shortsignal_exps2}, the obtained results based on the TT-decomposition were slightly better than those using the KT method.
For this particular problem, the estimation performance can even be higher when applying self-prediction algorithms, which exploit the low-rank structure of damped signals, e.g., TT-KT, and TT-linear prediction methods  based on SVD.
For a detailed derivation of these algorithms, see \citep{Phan_TT_part2}.
\end{example}

\subsection{Hankel Matrix and Hankel Tensor of Sinusoidal Signal}
\subsectionmark{Hankel Tensor of Sinusoidal Signal}

{\noindent \bf Hankel tensor of sinusoid.}
The Hankel tensor of a sinusoid signal
 $y(t)$ is a TT-Tucker tensor,
 \be
 	\tY = \llbracket \tG; \bU_1,  \bU_2, \ldots,  \bU_N \rrbracket \, ,\label{lem_hankel_sin_N}
 \ee
for which the factor matrices are defined in (\ref{eq_toeplitzN_factor2}).
The core tensor $\tG$ is an $N$th-order tensor of size $2 \times 2 \times \cdots \times 2$, in the TT-format, given by
\be
\tG =  \llangle \tG^{(1)},  \tG^{(2)},  \ldots,  \tG^{(N-1)} \rrangle \, ,
\ee
where $\tG^{(1)}  =   \bH(J_1)$ is a matrix of size $1 \times 2 \times 2$, while the core tensors $\tG^{(n)}$, for $n = 2,\ldots, N-1$, are of size $2\times 2 \times 2$ and have two horizontal slices, given by
\begin{align}
\tG^{(n)}(1,:,:)  &=  \bH(J_{n}-n+1)\, ,  \quad \tG^{(n)}(2,:,:)  =  \bH(J_{n}-n+2) \, 	  , \notag
\end{align}
with
\be
\bH(I) = \frac{1}{\sin^2(\omega)} \, \left[
	\begin{array}{@{}cc@{}}
		y(I) &  -y(I+1)\\
		-y(I-1)  &  y(I)
	\end{array}		
	\right] \label{eq_H}\,.
\ee

 \begin{remark}
The two TT-Tucker representations of the Toeplitz and Hankel tensors of the same sinusoid have similar factor matrices $\bU_n$, but their core tensors are different.
 \end{remark}

\begin{figure}[t]
\centering
\shadingbox{
\begin{tabular}{@{}l@{\hspace{.5ex}}l@{}l@{}}
\multicolumn{2}{@{}l}{1. Folded tensor} \\
\multirow{2}{*}{
\includegraphics[width=.42\linewidth, trim = 0.0cm .0cm 0cm 0cm,clip=true]{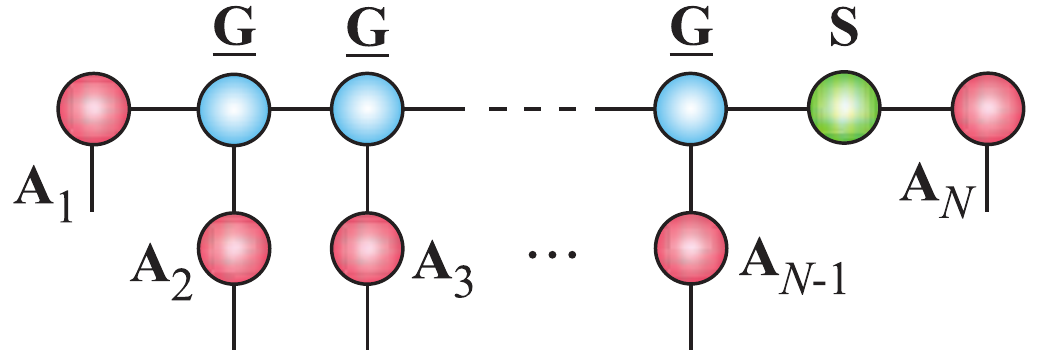}}
& \multicolumn{2}{@{}l}{
	$	\tG = \left[ \left[
	\begin{array}{@{}cc@{}}
		1 & 0\\
		0 & 1
	\end{array}		
	\right],  \left[
	\begin{array}{@{}cc@{}}
		0 & 1\\
		-1 & 0
	\end{array}		
	\right]   \right]$
}
\\[1.em]
& 	\multicolumn{2}{@{}l}{$
	\bS   = \left[\begin{array}{@{}cc@{}}
		\sin(\phi)  & \cos(\phi) \\
		\cos(\phi) & -\sin(\phi)
	\end{array}	\right] $,  $\bA_n = \bU_{\omega, 2^{n-1}}$}
\\[1.em]\hline
\multicolumn{2}{@{}l}{2. Toeplitz tensor}\\
\multirow{2}{*}{\includegraphics[width=.42\linewidth, trim = 0.0cm .0cm 0cm 0cm,clip=true]{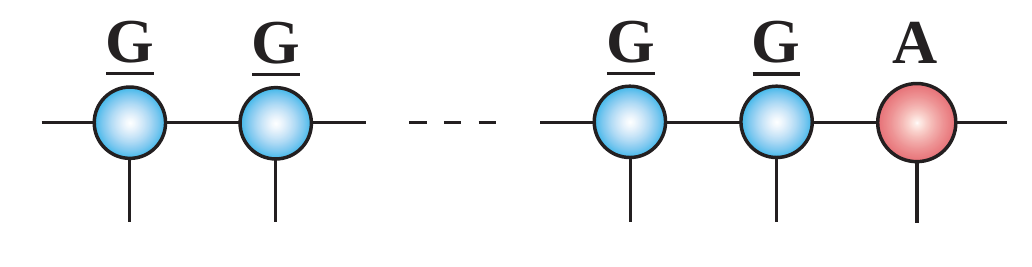}}
& \multicolumn{2}{@{}l}{
	$	\tG =  \left[  \left[
	\begin{array}{@{}cc@{}}
		1 & 0\\
		0 & 1
	\end{array}		
	\right] ,    \left[
	\begin{array}{@{}c@{\hspace{1ex}}c@{}}
		0 & 1\\
		-1 & 2\cos(\omega)
	\end{array}		
	\right] 	\right]   $}\\
& $\bA =  \left[
	\begin{array}{@{}cc@{}}
		y{(L-1)} &  y{(L)}\\
		y{(L-2)}  &  y{(L-1)}
	\end{array}		
	\right] $
\\[1.5em]\hline
\multicolumn{2}{@{}l}{3. Hankel tensor}\\
\multirow{2}{*}{\includegraphics[width=.42\linewidth, trim = 0.0cm .0cm 0cm 0cm,clip=true]{HanT}}
& \multicolumn{2}{@{}l}{
	$\tG=  \left[ \left[
	\begin{array}{@{}c@{\hspace{1ex}}c@{}}
		2\cos(\omega) & -1 \\
		1 &  0
	\end{array}		
	\right]  ,     \left[
	\begin{array}{@{}cc@{}}
		1 & 0\\
		0 & 1
	\end{array}		
	\right]	\right]$ } \\
& 	$\bA =  \left[
	\begin{array}{@{}cc@{}}
		y{(L-2)} &  y{(L-1)}\\
		y{(L-1)}  &  y{(L)}
	\end{array}		
	\right] $ \\
\end{tabular}
}
\caption{Representations of a sinusoid signal in different quantized tensor formats of size $2\times 2 \times \cdots \times 2$.}
\label{fig_quantum_sinu}
\end{figure}

{\noindent \bf Quantized Hankel tensor.}
An $(L-1)$th-order Hankel tensor of size $2 \times 2 \times \cdots \times 2$ of a sinusoid signal of length $L$ has a TT-representation with $(N-2)$ identical core tensors $\tG$, in the form
\be
\tY = \llangle \tG , \tG,  \ldots,  \tG,   \left[
	\begin{array}{@{}cc@{}}
		y(L-2) &  y(L-1)\\
		y(L-1)  &  y(L)
	\end{array}		
	\right]  \rrangle    \, ,\notag
\ee
where
\be
	\tG(1,:,:) =  \left[
	\begin{array}{@{}cc@{}}
		 2\cos(\omega)  & -1\\
		1  & 0
	\end{array}		
	\right] \,, \quad \tG(2,:,:) =  \left[
	\begin{array}{@{}cc@{}}
		1 & 0\\
		0 & 1
	\end{array}		
	\right]. \notag
\ee

Finally, representations of the sinusoid signal in various tensor format of  size are summarised in Figure~\ref{fig_quantum_sinu}.

\section{Summary}
\label{sec:conclusion}

This chapter has introduced several common tensorization methods, together with their properties
and illustrative applications in blind source separation, blind identification, denoising, and  harmonic retrieval.
The main criterion for choosing a suitable tensorization is that the tensor generated from lower-order original data must reveal the underlying low-rank tensor structure in some tensor format.
For example, the folded tensors of mixtures of damped sinusoid signals have low-rank QTT representation, while the derivative tensors in blind identification admit the CP decomposition.
The Toeplitz and Hankel tensor foldings augment the number of signal entries, through the replication of signal segments (redundancy), and in this way become suited to modeling of  signals of short length.
A property crucial to the solution via the tensor networks shown in this chapter, is that the tensors can be generated in the TT/QTT format, if the generating vector admits a low-rank QTT representation.

In  modern data analytics problems,  such as regression and deep learning, the number of model parameters can be huge,
which renders the model intractable. Tensorization can then serve as a remedy, by representing the parameters in some low-rank tensor format.
For further discussion on tensor representation of parameters in tensor regression,  we refer to Chapter \ref{Chapter2}.
A  wide class of optimization problems including  of solving linear systems, eigenvalue decomposition, singular value decomposition, Canonical Correlation Analysis (CCA) are addressed in Chapter \ref{Chapter3}. The tensor structures for Boltzmann machines  and convolutional  deep neural networks (CNN) are provided in Chapter \ref{Chapter4}.


\chapter{Supervised Learning with Tensors}
\label{Chapter2}

\vspace{-1.2cm}

\begin{figure}[t]
\centering
  \includegraphics[width=0.65\textwidth]{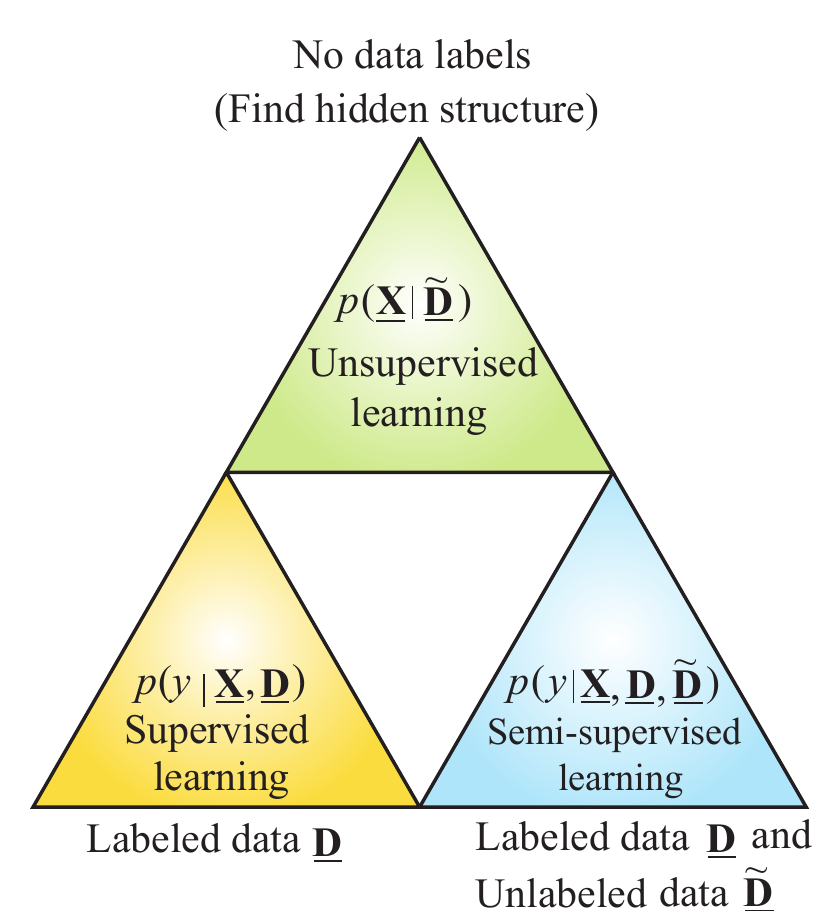}
  \caption{Graphical illustration of three fundamental learning approaches: Supervised, unsupervised and  semi-supervised learning.}
  \label{Fig:learning}
\end{figure}

Learning a statistical model that formulates a hypothesis for the data distribution merely from multiple input data samples, $\bx$, without knowing the corresponding values of the response variable, $y$, is refereed to as  \textit{unsupervised learning}. In contrast, \textit{supervised learning} methods, when seen from a probabilistic perspective,  model either the joint distribution $p(\bx, y)$ or the conditional distribution $p(y|\bx)$, for given training data pair $\{\mathbf x, y\}$. Supervised learning  can be categorized into regression, if $y$  is continuous, or classification, if $y$ is categorical (see also Figure~\ref{Fig:learning}).

Regression models can  be categorized
into linear regression and nonlinear regression.
In particular, multiple linear regression is associated with multiple smaller-order predictors, while multivariate regression corresponds to a single linear regression model but with multiple predictors and multiple responses. Normally,  multivariate regression tasks are encountered when the predictors are
arranged as vectors, matrices or tensors of variables. A basic linear regression model in the vector form is defined as
\be
y= f(\bx;\bw,b) = \langle \bx, \bw \rangle  +b =\bw^{\text{T}} \bx +b,
\ee
where $\bx \in \Real^I$ is the input vector of independent variables,
 $\bw \in \Real^I$  the vector of regression coefficients or weights, $b$ the bias, and $y$  the
 regression output or a dependent/target variable.

Such simple linear models can be applied not only for  regression but also for feature selection and classification. In all the cases, those models approximate the target variable $y$ by a weighted linear combination of input variables, $\bw^{\text{T}} \bx + b$.

Tensor representations are often very useful in mitigating the small sample size problem in discriminative subspace selection, because the information about the structure of objects is inherent in tensors and is a natural constraint which helps reduce the number of unknown parameters in the description of a learning model. In other words, when the number of training measurements is limited, tensor-based learning machines are expected to perform better than the corresponding vector-based learning machines,
as  vector representations are associated with several problems, such as loss of information for structured data and over-fitting for high-dimensional data.

\section{Tensor Regression}

Regression is at the very core of  signal processing and machine learning, whereby the output  is typically estimated based on a linear combination of regression coefficients  and the input regressor, which can be a vector, matrix, or tensor.  In this way, regression analysis can be used to predict  \textit{dependent variables} (responses,
outputs, estimates),  from a set of \textit{independent variables} (predictors, inputs, regressors),
by exploring the correlations
among these variables as well as explaining the inherent factors behind the observed patterns. It is also often convenient, especially regarding ill-posed cases of matrix inversion, which is inherent to regression to jointly perform regression and dimensionality reduction through, e.g., principal component regression (PCR) \citep{jolliffe1982note}, whereby regression is performed on a well-posed low-dimensional subspace defined through most significant principal components. With tensors being  a natural generalization of vectors and matrices, tensor  regression can be defined in an analogous way.

A well established and important supervised learning technique is  linear or nonlinear Support Vector Regression (SVR) \citep{smola1997support}, which allows for the modeling of streaming data and is quite closely related to Support Vector Machines (SVM) \citep{cortes1995support}. The model produced by SVR only depends on a subset of training data (support vectors), because the cost function for building the model ignores any training data that is close (within a threshold $\varepsilon$) to the model prediction.

Standard support vector regression techniques  have been naturally   extended to
Tensor Regression (TR) or Support Tensor Machine (STM) methods \citep{tao2005supervised}.
In the (raw) tensor format, the TR/STM can be formulated as
\be
\label{eq:tensorregression}
y= f(\tensor{X};\tensor{W}, b) = \langle \tensor{X}, \tensor{W} \rangle  +b,
\ee
where $\tensor{X} \in \Real^{I_1 \times \cdots \times I_N}$ is the input tensor regressor,
 $\tensor{W} \in \Real^{I_1 \times \cdots \times I_N}$   the tensor of weights (also called regression tensor or model tensor), $b$ the bias, and $y$  the regression output, $\langle \tX,\tW \rangle = \vtr{\tX}^{\text{T}} \vtr{\tW}$ is the inner product of two tensors.

We shall denote input samples of multiple predictor variables (or features) by $\underline \bX_{1}, \ldots, \underline \bX_M$ (tensors), and the actual continuous or categorical response variables by $y_1, y_2, \ldots, y_M$ (usually scalars).
The training process, that is the estimation of the weight tensor $\tensor{W}$
and bias $b$, is carried out based on the set of available training samples $\{\tensor{X}_m, y_m\}$ for $m=1,\ldots,M$.
Upon arrival of a new training
sample, the TR model is used to make predictions for that sample.

The problem is usually formulated as a minimization of the following \textit{squared cost function}, given by
\begin{equation}
J(\tensor{X}, y \; | \; \tensor{W}, b) = \sum_{m=1}^{M}\Big(y_m - \big(\langle\tensor{W},\tensor{X}_m\rangle + b\big)\Big)^2
\end{equation}
or the logistic loss function (usually employed in classification problems), given by
\begin{equation}
J(\tensor{X}, y  \;| \; \tensor{W},b) = \sum_{m=1}^{M} \log \left(\frac{1}{1+ e^{-y_m(\langle\tensor{W},\tensor{X}_m\rangle+b)}}\right).
\end{equation}

In practice, for very large scale problems, tensors are expressed approximately in tensor network formats, especially using Canonical Polyadic (CP), Tucker or Tensor Train (TT)/Hierarchical Tucker (HT) models \citep{OseledetsTT11,grasedyck2010hierarchical}. In this case, a suitable representation of the weight tensor, $\tensor{W}$, plays a key role in the model performance. For example, the  CP representation of the weight tensor, in the form
\be
\tensor{W} &=& \sum_{r=1}^R \mathbf{u}_r^{(1)}\circ \mathbf{u}_r^{(2)}\circ \cdots \circ\mathbf{u}_r^{(N)} \notag \\
&=& \tI \times_1\mathbf{U}^{(1)}\times_2 \mathbf{U}^{(2)} \cdots\times_N\mathbf{U}^{(N)},
\ee
where ``$\circ$'' denotes the outer product of vectors, leads to a generalized linear model (GLM), called the CP tensor regression \citep{zhou2013tensor}.

Analogously, upon the application of Tucker multilinear rank tensor representation
\be
\tensor{W} = \tensor{G}\times_1\mathbf{U}^{(1)}\times_2 \mathbf{U}^{(2)} \cdots\times_N\mathbf{U}^{(N)},
\ee
we obtain Tucker tensor regression \citep{Hoff_MTR,li2013tucker,yulearning}.

An alternative form of the multilinear Tucker regression model, proposed by \citet{Hoff_MTR}, assumes that the replicated observations $\{\tensor{X}_m, \tensor{Y}_m\}_{m=1}^M$  are stacked in concatenated
tensors $\tensor{X} \in \Real^{I_1 \times \cdots \times I_N \times M}$ and  $\tensor{Y} \in \Real^{J_1 \times \cdots \times J_N \times M}$, which admit  the following model
\be
\tensor{Y} = \tensor{X} \times_1 \bW_1 \times_2 \bW_2 \cdots \times_N \bW_N \times_{N+1} \bD_M + \tensor{E},
\ee
where $\bD_M$ is an $M \times M$ diagonal matrix, $\bW_n \in \Real^{J_n \times I_n}$ are the weight matrices within the Tucker model, and $\tensor{E}$ is a zero-mean residual tensor of the same order as
$\tensor{Y}$.  The unknown regression coefficient matrices, $\bW_n$ can be found using the procedure  outlined in Algorithm~\ref{alg:MTR}.

It is important to highlight that the Tucker regression model offers several benefits over the CP regression model, which include:
\begin{enumerate}
  \item A more parsimonious modeling capability and  a more compact model, especially for a limited number of available samples;
  \item Ability to fully exploit multi-linear ranks, through the freedom to choose a different rank for each mode, which is essentially useful when data is skewed in dimensions (different sizes in modes);
  \item Tucker decomposition explicitly models the interaction between  factor matrices in different modes, thus allowing for a finer grid search over a larger modeling space.
\end{enumerate}

\begin{algorithm}[t!]
\caption{\textbf{Multilinear Tucker Regression}}
\label{alg:MTR}
{\small
\begin{algorithmic}[1]
\REQUIRE $\tensor{X}\in\mathbb R^{I_1\times I_2\times \cdots\times I_N\times M}$ and $\tensor Y\in \mathbb R^{J_1\times I_2\times \cdots\times J_N\times M}$.
\ENSURE $\{\mathbf W_n\}, n=1,\ldots, N$.
\STATE Initialize randomly $\bW_n$ for $n=1, \ldots, N$.
\WHILE {not converged or iteration limit is not reached}
\FOR {$n=1$ to $N$}
\STATE $\underline{\bX}^{(n)}= \underline{\bX} \times_1 \bW_1  \cdots \times_{n-1} \bW_{n-1}\times_{n+1} \bW_{n+1} \cdots \times_N \bW_N$
\STATE Matricize tensors $\underline \bX^{(n)}$ and $\underline \bY$ into their respective unfolded matrices \\$\bX^{(n)}_{(n)}$ and $\bY_{(n)}$
\STATE Compute $\bW_n= \bY_{(n)} (\bX_{(n)}^{(n)})^{\text{T}} \left(\bX_{(n)}^{(n)} \; (\bX_{(n)}^{(n)})^{\text{T}}\right)^{-1}$
\ENDFOR
\ENDWHILE
\end{algorithmic}
}
\end{algorithm}

 Both the CP and Tucker tensor regression models can be solved by alternating least squares (ALS) algorithms which sequentially estimate one factor matrix at a time while keeping other factor matrices fixed. To deal with the curse of dimensionality, various  tensor  network decompositions can be applied,  such as the  TT/HT decomposition for very high-order tensors \citep{OseledetsTT11,grasedyck2010hierarchical}.
When the weight tensor $\tensor{W}$ in (\ref{eq:tensorregression}) is represented by a low-rank HT decomposition, this is referred to as the H-Tucker tensor regression \citep{hou2016phd}.

{\bf Remark 1.} In some applications, the weight tensor, $\underline {\bW}$, is of a higher-order than input tensors, $\underline \bX_m$, this yields a  more general  tensor regression model
\be
\tensor{Y}_m =   \langle \tensor{X}_m | \tensor{W} \rangle  + \tensor{E}_m, \qquad m=1,\ldots,M,
\label{eq:GTReg}
\ee
where $\langle \tensor{X}_m | \tensor{W} \rangle$ denotes a tensor contraction along the  first $N$ modes of an $N$th-order input (covariate) tensor, $\tensor{X}_m \in \Real^{I_1 \times \cdots \times I_N}$, and a $P$th-order weight tensor,
 $\tensor{W} \in \Real^{I_1 \times \cdots \times I_P}$,  with $P >N$, while $\tensor{E} \in \Real^{I_{P+1} \times \cdots \times I_P} $ is the residual tensor and $\tensor{Y} \in \Real^{I_{P+1} \times \cdots \times I_P}$ the response  tensor.

Observe that the tensor inner product  is equivalent to a contraction of two tensors of the same order (which is a scalar)  while a contraction of two tensors of different orders, $\tensor{X} \in \Real^{I_1 \times \cdots \times I_N}$  and   $\tensor{W} \in \Real^{I_1 \times \cdots \times I_P}$,  with $P > N$,  is defined as a tensor $\underline \bY \in \Real^{I_{N+1} \times \cdots \times I_P}$ of $(P-N)$th-order  with entries
\be
\langle \tensor{X} | \tensor{W} \rangle_{i_{N+1},\ldots,i_P} = \sum_{i_1=1}^{I_1} \cdots \sum_{i_N=1}^{I_N} x_{i_1,\ldots,i_N} \; w_{i_1,\ldots,i_N,i_{N+1},\ldots,i_P}.
\ee
Many  regression problems are special cases of the general tensor regression model in (\ref{eq:GTReg}), including the multi-response regression, vector autoregressive model and pair-wise interaction tensor model (see \citep{Raskutti15} and references therein).\\

In summary, the aim of tensor regression is to estimate the entries of a weight tensor, $\tensor{W}$, based on available input-output observations $\{\tensor{X}_m, \tensor{Y}_m\}$.
In a more general scenario, support tensor  regression (STR) aims to identify a  nonlinear function, $f: \mathbb{R}^{I_1\times I_2\times\cdots\times I_N} \rightarrow \mathbb{R}$, from a collection of observed input-output data pairs $\{\tensor{X}_m, y_m\}_{m=1}^M$ generated from the model
\begin{equation}
y_m = f(\tensor{X}_m) + \varepsilon,
\end{equation}
where the input $\tensor{X}_m$ has a tensor structure,  the output $y_m$ is a scalar, and $\varepsilon$ is also a scalar which represents the error.

Unlike linear regression models, nonlinear regression models have the ability to characterize complex
 nonlinear dependencies
in data, in which the responses are modeled
through  a combination of nonlinear parameters and predictors, with the nonlinear parameters
usually taking the form of an exponential function, trigonometric function, power function,
etc.  The nonparametric counterparts, so called nonparametric regression models, frequently appear in  machine
learning, such as  Gaussian Processes (GP) (see \citep{zhao2013tensor, hou2016phd}  and references therein).

\section{Regularized Tensor Models}

Regularized  tensor models aims to  reduce the complexity of tensor regression models through  constraints (restrictions) on the model  parameters, $\underline \bW$. This is particularly advantageous for problems  with a large number of features but a small number of data samples.

Regularized linear tensor models can be generally formulated as
 \be
\min_{\underline  \bW, b} \, f (\underline  \bX,y)= J(\underline \bX, y \;| \;  \underline \bW, b)+ \lambda R(\underline \bW),
\ee
where $J(\underline \bX, y \;| \; \underline \bW, b)$ denotes a loss (error) function, $R(\underline \bW)$ is a regularization term, while the parameter $\lambda >0$ controls the trade-off between the contributions of the original  loss function and regularization term.

One such well-known regularized linear model is the Frobenius-norm regularized model for support vector regression, called Grouped LASSO \citep{Lasso1996}, which employs  logistic loss and $\ell_1$-norm regularization for simultaneous classification (or regression) and feature selection. Another very popular model of this kind is the support vector machines (SVM) \citep{cortes1995support} which uses a hinge loss in the form $l(\hat y) = \text{max}(0,1- y \hat y)$, where $\hat y $ is the prediction and $y $  the true label.

In regularized tensor regression, when the regularization is performed via the Frobenius norm, $\|\underline \bW\|^2_F=\langle\underline \bW,\underline \bW\rangle$, we arrive at the standard Tikhonov regularization, while using the  $\ell_1$ norm, $\|\cdot\|_{\ell_1}$, we impose classical sparsity constraints. The advantage of tensors over vectors and matrices is that  we can exploit the flexibility in the choice of sparsity profiles. For example, instead of imposing global sparsity for the whole tensor,
we can impose sparsity for   slices or fibers if there is a need  to enforce  some fibers or slices to have most of their entries zero. In such a case,  similarly to group LASSO, we can apply group-based $\ell_1$-norm  regularizers,  such as the  $\ell_{2,1}$ norm, i.e., $\|\mathbf X\|_{2,1} = \sum_{j=1}^J \| \mathbf x_{j}\|_2$.

Similarly to matrices, the various rank properties can also be employed for tensors, and are much richer and more complex  due to  the multidimensional structure of tensors. In addition to sparsity constraints, a low-rank of a tensor can be exploited as a regularizer, such as the canonical rank or multilinear (Tucker)  rank,  together with various more sophisticated  tensor  norms like the nuclear norm or latent norm of the weight tensor $\underline \bW$.

The low-rank regularization problem can be formulated as
\begin{equation}
\begin{split}
\min_{\tensor{W}} \; J(\tensor{X}, y \;| \; \tensor{W}, b ), \quad \mbox{s.t.}
\quad \text{multilinear rank}(\tensor{W}) \leq R.
\end{split}
\end{equation}
Such a formulation  allows us to estimate a weight tensor, $\tensor{W}\in\mathbb{R}^{I_1\times\cdots\times I_N}$, that minimizes the empirical loss $J$, subject to the constraint that the multilinear rank of $\tensor{W}$ is at most $R$. Equivalently, this implies that the weight tensor, $\tensor{W}$, admits a low-dimensional factorization in the form $\tensor{W} = \tensor{G}\times_1\mathbf{U}^{(1)}\times_2 \mathbf{U}^{(2)}\cdots\times_N \mathbf{U}^{(N)}$.

Low rank structure of tensor data has been successfully utilized in  applications including  missing data imputation \citep{liu2013tensor,zhao2015bayesian}, multilinear robust principal component analysis \citep{inoue2009robust}, and subspace clustering \citep{zhang2015low}. Instead of low-rank properties of data itself, low-rank regularization can be also applied to learning coefficients in regression and classification \citep{yulearning}.

Similarly, for very large-scale problems, we can also apply the tensor train (TT) decomposition as the constraint, to yield
\begin{equation}
\begin{split}
\min_{\tensor{W}} \, J(\underline \bX, y \;| \;  \underline \bW, b), \quad
\mbox{s.t.}  \quad \tensor{W} = \llangle \tensor{G}^{(1)}, \ldots, \tensor{G}^{(N)} \rrangle.
\end{split}
\end{equation}
In this way, the weight tensor $\tensor{W}$ is approximated by a low-TT rank tensor of the  TT-rank $(R_1, \ldots, R_{N-1})$ of at most $R$.

The low-rank constraint for the tensor $\tensor{W}$ can also be formulated through  the
tensor  norm, in the form \citep{wimalawarne2016theoretical}
\begin{equation}
\min_{\tensor{W}, b} \, (J(\underline \bX, y \;| \;  \underline \bW, b) + \lambda \|\tensor{W}\|),
\end{equation}
where $\|\cdot\|$ is a suitably chosen tensor/matrix norm.


 One of the important and useful  tensor norms is the tensor nuclear norm \citep{liu2013tensor} or the (overlapped) trace norm \citep{wimalawarne2014multitask}, which can be defined  for a tensor $\tensor{W}\in\mathbb{R}^{I_1\times \cdots\times I_N}$ as
\begin{equation}
\label{eq:tensornorm}
\| \tensor{W}\|_* = \sum_{n=1}^N \|\mathbf{W}_{(n)}\|_*,
\end{equation}
where  $\|\mathbf W_{(n)}\|_* = \sum_{i} \sigma_i$, $\sigma_i$ denotes the $i$th  singular value of $\mathbf W_{(n)}$.
The overlapped tensor trace norm can be viewed as a direct extension of the matrix trace norm, since it uses unfolding matrices of a tensor in all of its modes,  and then computes the sums of trace norms of those unfolding matrices. Regularization based on the overlapped trace norm can also be viewed as an overlapped group regularization, due to the fact that the same tensor is unfolded over all modes and regularized using the trace norm.

Recently  \citet{tomioka2013convex} proposed the {\it latent trace norm} of a tensor,  which takes a mixture of $N$ latent tensors, $\tensor{W}^{(n)}$,  and regularizes each of them separately, as in (\ref{eq:tensornorm}), to give
\begin{equation}
\|\tensor{W}\|_{latent} = \inf_{ \tensor{W}} \; \sum_{n=1}^N \|\mathbf{W}_{(n)}^{(n)}\|_*,
\end{equation}
where $\tensor{W}=\tensor{W}^{(1)} + \tensor{W}^{(2)}+ \cdots+ \tensor{W}^{(N)}$ and $\mathbf{W}^{(n)}_{(n)}$ denotes the unfolding of $\tensor{W}^{(n)}$ in its $n$th mode.
In contrast to the nuclear norm, the latent tensor trace norm effectively regularizes different latent tensors in each unfolded mode, which makes it possible for  the latent tensor trace norm to identify (in some cases) a latent tensor with the lowest possible rank. In general,   the content of each latent tensor depends on the rank of its unfolding matrices.

{\bf Remark 2.} A major drawback of the latent trace norm approach is its inability to identify the rank of a tensor,  when the
rank  value is relatively close  to its dimension (size). In other words, if a tensor has a mode with a  dimension (size)  much smaller than  other modes,  the latent trace norm  may  be  incorrectly estimated.
To deal with this problem, the scaled latent norm  was proposed by  \citet{wimalawarne2014multitask} which is defined as
\begin{equation}
\label{eq:scaletensornorm}
\|\tensor{W}\|_{scaled} = \inf_{ \tensor{W}} \; \sum_{n=1}^N \frac{1}{\sqrt{I_n}} \|\mathbf{W}^{(n)}_{(n)}\|_*,
\end{equation}
where $ \tensor{W}= \tensor{W}^{(1)} + \tensor{W}^{(2)}+ \cdots + \tensor{W}^{(N)}$. Owing to the normalization by the mode size, in the form of (\ref{eq:scaletensornorm}),  the scaled latent trace norm scales with mode dimension, and thus estimates more reliably the rank,  even when the sizes of  various modes are quite different.

The state-of-the-art tensor regression models therefore employ the scaled  latent trace norm,  and are solved by the following  optimization problem
%
\begin{equation}
P(\tensor{W}, b ) = \min_{ \tensor{W},b} \, \bigg(\sum_{m=1}^{M}\Big(y_m - \big(\langle\tensor{W},\tensor{X}_m\rangle + b\big)\Big)^2 + \sum_{n=1}^N \lambda_n\|\mathbf{W}^{(n)}_{(n)}\|_*\bigg),
\end{equation}
where $\tensor{W}= \tensor{W}^{(1)} + \tensor{W}^{(2)}+ \cdots + \tensor{W}^{(N)}$, for $n=1, \ldots, N$, and for any given regularization parameter $\lambda$, where $\lambda_n=\lambda$ in the case of latent trace norm and $\lambda_n =\frac{\lambda}{\sqrt{I_n}}$ in the case of the scaled latent trace norm, while  $\mathbf{W}^{(n)}_{(n)}$ denotes the unfolding of $\tensor{W}^{(n)}$ in its $n$th mode.

{\bf Remark 3.} Note that the above optimization problems involving the latent and scaled latent trace norms require a large number
($N \prod_n I_n$) of variables in $N$ latent tensors.

 Moreover, the existing methods based on SVD are infeasible for very  large-scale applications. Therefore,  for very large scale high-order tensors, we need to
use tensor network decompositions and perform all operations in tensor networks formats, taking advantage of their lower-order core tensors.

\section{Higher-Order Low-Rank Regression (HOLRR)}

Consider a full multivariate regression task in which  the response has a tensor structure. Let $f:\mathbb{R}^{I_0} \rightarrow \mathbb{R}^{I_1\times I_2\times\cdots\times I_N}$ be a function we desire to learn from  input-output data, $\{\mathbf{x}_m, \tensor{Y}_m\}_{m=1}^M$, drawn from the model $\tensor{Y} =\tensor{W} \; \bar \times_1 \; \mathbf{x} + \tensor{E}$, where $\tensor{E}$ is an error tensor, $\bx \in \Real^{I_0}$, $\underline \bY \in \Real^{I_1 \times \cdots \times I_N}$, and
$\tensor{W}\in\mathbb{R}^{I_0\times I_1\times \cdots \times I_N}$ is a tensor of regression coefficients.
The extension of  low-rank regression methods to tensor-structured responses can be achieved by enforcing low multilinear rank of the regression tensor, $\tensor{W}$, to yield the so-called higher-order low rank regression (HOLRR) \citep{rabusseau2016higher,sun2016sparse}. The aim is to find a low-rank regression tensor, $\tensor{W}$, which minimizes the loss function based on the training data.

To avoid numerical instabilities and to prevent overfitting, it is convenient to employ  a ridge regularization of the objective function, leading to the following minimization problem
\begin{equation}
\label{eq:lrtr}
\begin{split}
\min_{\tensor{W}} & \sum_{m=1}^M \|\tensor{W} \; \bar \times_1 \; \mathbf{x}_m- \tensor{Y}_m\|_F^2 + \gamma \|\tensor{W}\|_F^2 \\
\mbox{s.t.}  \quad & \mbox{multilinear rank}(\tensor{W}) \leq (R_0, R_1,\ldots,R_N),
\end{split}
\end{equation}
for some given integers $R_0, R_1,\ldots,R_N$.


  Taking into account the fact that input   vectors, $\bx_m$, can be concatenated into an input
   matrix,  $\bX =[\bx_1, \ldots, \bx_M]^{\text{T}} \in \Real^{M \times I_0}$,  the above optimization problem  can be reformulated in the following form
\begin{equation}
\begin{split}
& \min_{\tensor{W}} \, \|\tensor{W}\times_1 \mathbf{X}- \tensor{Y}\|_F^2 + \gamma \|\tensor{W}\|_F^2 \\
\mbox{s.t.} \quad & \tensor{W} = \tensor{G}\times_1\mathbf{U}^{(0)}\times_2\mathbf{U}^{(1)} \cdots\times_{N+1}\mathbf{U}^{(N)}, \\
& \mathbf{U}^{(n)\,\text{T}}\mathbf{U}^{(n)}=\mathbf{I}  \quad \text{for}  \quad n=1, \ldots, N,
\end{split}
\end{equation}
where
the output tensor $\tY$ is obtained by stacking the output tensors $\tY_{m}$ along the first mode  $\tY(m,:,\ldots,:) = \tY_m$.
%
The regression function can  be rewritten as
\begin{equation}
f: \mathbf{x} \mapsto \tensor{G} \times_1 \mathbf{x}^{\text{T}}\mathbf{U}^{(0)}\times_2 \mathbf{U}^{(1)} \cdots\times_{N+1}\mathbf{U}^{(N)}.
\end{equation}
This minimization problem can be reduced to finding ($N+1$) projection matrices onto subspaces of dimensions, $R_0, R_1,\ldots, R_N$, where  the
core tensor $\tensor{G}$ is given by
\begin{equation}
\tensor{G} = {\tensor{Y}} \times_1 (\mathbf{U}^{(0)\,\text{T}} ({\mathbf{X}}^{\text{T}} {\mathbf{X}} + \gamma\mathbf{I}) \mathbf{U}^{(0)})^{-1} \mathbf{U}^{(0)\,\text{T}} {\mathbf{X}}^{\text{T}} \times_2 \mathbf{U}^{(1)\,\text{T}} \cdots\times_{N+1}\mathbf{U}^{(N)\,\text{T}},
\end{equation}
and  the orthogonal matrices $\mathbf{U}^{(n)}$ \;  for $ n=0,1,\ldots,N$ can be computed by the eigenvectors of
\begin{equation}
\Bigg\{ \begin{array}{ll}
          (\mathbf{X}^{\text{T}}\mathbf{X} +\gamma\mathbf{I})^{-1}\mathbf{X}^{\text{T}}\mathbf{Y}_{(1)}\mathbf{Y}_{(1)}^{\text{T}}\mathbf{X}, &  \quad n=0, \\
          \mathbf{Y}_{(n)}\mathbf{Y}_{(n)}^{\text{T}}, & otherwise .
        \end{array}
\end{equation}
In other words, by the computation of $R_n$ largest eigenvalues and the corresponding eigenvectors. The  HOLRR procedure is outlined  in Algorithm~\ref{alg:HOLRR}.

\begin{algorithm}[t!]
\caption{\textbf{Higher-Order Low-Rank Regression (HOLRR)}}
\label{alg:HOLRR}
{\small
\begin{algorithmic}[1]
\REQUIRE $\mathbf{X}\in\mathbb R^{M\times I_0}$, $\tensor{Y} \in\mathbb R^{M\times I_1\times\cdots\times I_N}$, rank $(R_0,R_1,\ldots,R_N)$ \\and a regularization parameter $\gamma$.
\ENSURE $\tensor{W} = \tensor{G}\times_1\mathbf{U}^{(0)}\times_2\mathbf{U}^{(1)} \cdots\times_{N+1}\mathbf{U}^{(N)}$
\STATE $\mathbf U^{(0)} \leftarrow$ top $R_0$ eigenvectors of $(\mathbf{X}^{\text{T}}\mathbf{X} +\gamma\mathbf{I})^{-1}\mathbf{X}^{\text{T}}\mathbf{Y}_{(1)}\mathbf{Y}_{(1)}^{\text{T}}\mathbf{X}$
\FOR {$n=1$ to $N$}
\STATE $\mathbf U^{(n)} \leftarrow $ top $R_n$ eigenvectors of $\mathbf Y_{(n)}\mathbf Y_{(n)}^{\text{T}}$
\ENDFOR
\STATE $\mathbf{T} =\left(\mathbf{U}^{(0)\,\text{T}} ({\mathbf{X}}^{\text{T}} {\mathbf{X}} + \gamma\mathbf{I}) \mathbf{U}^{(0)}\right)^{-1} \mathbf{U}^{(0)\,\text{T}} {\mathbf{X}}^{\text{T}}$
\STATE $\tensor{G} \leftarrow {\tensor{Y}} \times_1 \mathbf T \times_2 \mathbf{U}^{(1)\,\text{T}} \cdots\times_{N+1}\mathbf{U}^{(N)\,\text{T}}$
\end{algorithmic}
}
\end{algorithm}

\section{Kernelized HOLRR}

The HOLRR can be  extended to its kernelized version, the kernelized \-HOLRR (KHOLRR).  Let $\phi: \mathbb{R}^{I_0} \rightarrow \mathbb{R}^L$ be a feature map and let $\mbi \Phi\in\mathbb{R}^{M\times L}$ be a  matrix with row vectors $\phi(\mathbf{x}_m^{\text{T}})$ for $m\in\{1,\ldots, M\}$. The HOLRR problem in the feature space boils down to the ridge regularized  minimization problem
\begin{equation}
\begin{split}
& \min_{\tensor{W}} \,   \|\tensor{W}\times_1 \mbi \Phi- \tensor{Y}\|_F^2 + \gamma \|\tensor{W}\|_F^2 \\
\mbox{s.t.} \quad & \tensor{W} = \tensor{G}\times_1\mathbf{U}^{(0)}\times_2\mathbf{U}^{(1)} \cdots\times_{N+1}\mathbf{U}^{(N)}, \\
& \text{and}\quad  \mathbf{U}^{(n)\,\text{T}}\mathbf{U}^{(n)}=\mathbf{I} \quad \text{for} \quad n=0, \ldots, N.
\end{split}
\end{equation}
The tensor $\tensor W$ is represented in a Tucker format.
Then, the core tensor $\tensor{G}$  is given by
\begin{equation}
\tensor{G} = {\tensor{Y}} \times_1 (\mathbf{U}^{(0)\,\text{T}} {\mathbf{K}} ({\mathbf{K}} + \gamma\mathbf{I}) \mathbf{U}^{(0)})^{-1} \mathbf{U}^{(0)\,\text{T}} {\mathbf{K}} \times_2 \mathbf{U}^{(1)\,\text{T}}  \cdots\times_{N+1}\mathbf{U}^{(N)\,\text{T}},
\end{equation}
where $\mathbf K$ is the kernel matrix and  the orthogonal matrices $\mathbf{U}^{(n)}, n=0,1,\ldots,N$, can be computed via the eigenvectors of
\begin{equation}
\Bigg\{ \begin{array}{ll}
          (\mathbf{K} +\gamma\mathbf{I})^{-1}\mathbf{Y}_{(1)}\mathbf{Y}_{(1)}^{\text{T}}\mathbf{K}, &  \quad n=0, \\
          \mathbf{Y}_{(n)}\mathbf{Y}_{(n)}^{\text{T}}, & otherwise,
        \end{array}
\end{equation}
which  corresponds to the computation of $R_n$ largest eigenvalues and the associated eigenvectors, where  $\mathbf K = \mbi \Phi\mbi\Phi^{\text{T}}$ is the kernel matrix.  The  KHOLRR procedure is outlined in  Algorithm~\ref{alg:KHOLRR}.

\begin{algorithm}[t!]
\caption{\textbf{Kernelized Higher-Order Low-Rank Regression (KHOLRR)}}
\label{alg:KHOLRR}
{\small
\begin{algorithmic}[1]
\REQUIRE Gram matrix $\mathbf{K}\in\mathbb R^{M\times M}$, $\tensor{Y} \in\mathbb R^{M\times I_1\times\cdots\times I_N}$,\\ rank $(R_0,R_1,\ldots,R_N)$ and a regularization parameter $\gamma$.
\ENSURE $\tensor{W} = \tensor{G}\times_1\mathbf{U}^{(0)}\times_2\mathbf{U}^{(1)} \cdots\times_{N+1}\mathbf{U}^{(N)}$
\STATE $\mathbf U^{(0)} \leftarrow$ top $R_0$ eigenvectors of $(\mathbf{K} +\gamma\mathbf{I})^{-1}\mathbf{Y}_{(1)}\mathbf{Y}_{(1)}^{\text{T}}\mathbf{K}$
\FOR {$n=1$ to $N$}
\STATE $\mathbf U^{(n)} \leftarrow $ top $R_n$ eigenvectors of $\mathbf Y_{(n)}\mathbf Y_{(n)}^{\text{T}}$
\ENDFOR
\STATE $\mathbf{T} =\left(\mathbf{U}^{(0)\,\text{T}} {\mathbf{K}} ({\mathbf{K}} + \gamma\mathbf{I}) \mathbf{U}^{(0)}\right)^{-1} \mathbf{U}^{(0)\,\text{T}} {\mathbf{K}} $
\STATE $\tensor{G} \leftarrow {\tensor{Y}} \times_1 \mathbf T \times_2 \mathbf{U}^{(1)\,\text{T}} \cdots\times_{N+1}\mathbf{U}^{(N)\,\text{T}}$
\end{algorithmic}
}
\end{algorithm}

Note that the kernel  HOLRR  returns the weight tensor, $\tensor{W}\in \mathbb{R}^{M\times I_1\times\cdots\times I_N}$, which is used to define the regression function
\begin{equation}
f: \mathbf{x} \mapsto \tensor{W} \bar \times_1 \mathbf{k}_x = \sum_{m=1}^M k(\mathbf{x}, \mathbf{x}_m) \tensor{W}_m,
\end{equation}
where $\tensor{W}_m = \tensor{W}(m,:,\ldots,:) \in\mathbb{R}^{I_1\times \cdots\times I_N}$ and the $m$th component of $\mathbf k_x$ is $k(\mathbf{x}, \mathbf{x}_m)$.

\section{Higher-Order Partial Least Squares (HOPLS)}
\sectionmark{Higher-Order Partial Least Squares}

In this section, Partial Least Squares (PLS) method is briefly presented followed by its generalizations  to tensors.

\subsection{Standard Partial Least Squares }
\begin{figure}[t]
\centering
\includegraphics[width=0.99 \textwidth]{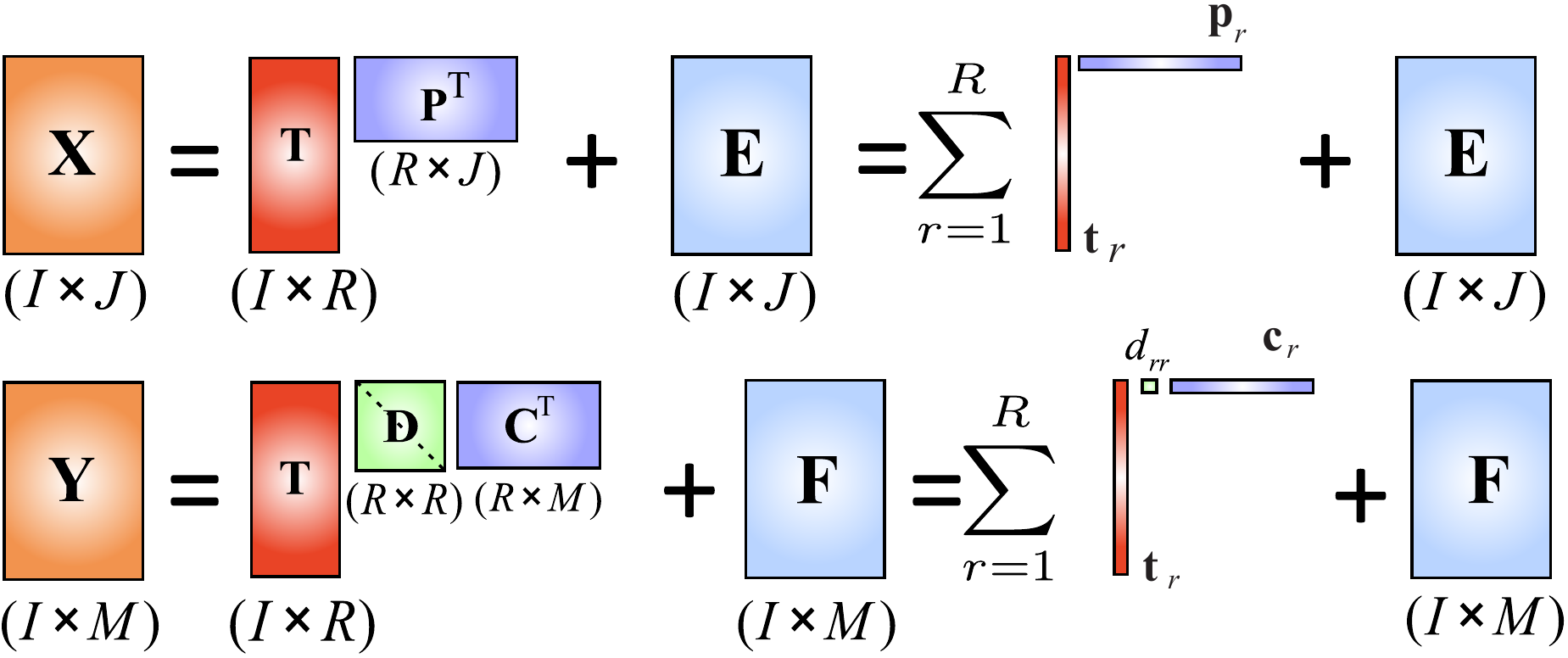}
\caption{The standard PLS model which performs data decomposition as a sum of rank-one matrices.}
\label{fig:PLS}
\end{figure}

The principle behind the PLS method is to search for a common set of latent vectors in the independent variable $\mathbf X\in\mathbb{R}^{I\times J}$ and the dependent variable $\mathbf Y\in\mathbb{R}^{I\times M}$ by performing their simultaneous decomposition,  with the constraint that the components obtained through such a decomposition explain as much as possible of the covariance between $\mathbf{X}$ and $\mathbf{Y}$. This problem can be formulated as (see also Figure~\ref{fig:PLS})
\begin{eqnarray}\label{Eq:PLS-X}
 \mathbf{X} &=& \mathbf{T}\mathbf{P}^{\text{T}} + \mathbf{E} = \sum_{r=1}^{R}{\mathbf{t}_{r}\mathbf{p}_{r}^{\text{T}}} + \mathbf{E},\\
\label{Eq:PLS-Y}
 \mathbf{Y} &=& \mathbf{T} \mathbf{D} \mathbf{C}^{\text{T}} + \mathbf{F} = \sum_{r=1}^{R} d_{rr} {\mathbf{t}_{r}\mathbf{c}_{r}^{\text{T}}} + \mathbf{F},
\end{eqnarray}
where $\mathbf{T}=[\mathbf{t}_1, \mathbf{t}_2, \ldots , \mathbf{t}_R]\in\mathbb{R}^{I\times R}$ consists of $R$ orthonormal latent variables from $\mathbf{X}$, and a matrix $\mathbf{U} = \mathbf{T} \mathbf{D} =[\mathbf{u}_1, \mathbf{u}_2,\ldots ,\mathbf{u}_R]\in\mathbb{R}^{I\times R}$ represents  latent variables from $\mathbf{Y}$ which have maximum covariance with the latent variables, $\bT$, in $\mathbf X$.
 The matrices $\mathbf{P}$ and $\mathbf{C}$ represent loadings (PLS regression coefficients), and $\mathbf{E}, \mathbf{F}$ are respectively the residuals for $\mathbf{X}$ and $\mathbf{Y}$, while $\mathbf{D}$ is a scaling diagonal matrix.

The PLS procedure is recursive, so that in order to obtain  the  set of first latent components in $\mathbf{T}$, the standard PLS algorithm  finds  the two sets of weight vectors, $\mathbf{w}$ and  $\mathbf{c}$, through the following optimization problem
\begin{equation}\label{Eq:PLSoptimization}
\max_{\{\mathbf{w}, \mathbf{c}\}} \; \left( \mathbf{w}^{\text{T}}\mathbf{X}^{\text{T}}\mathbf{Y}\mathbf{c}\right)^2, \quad \mbox{s.t. }\quad \mathbf{w}^{\text{T}}\mathbf{w} = 1, \quad \mathbf{c}^{\text{T}}\mathbf{c} = 1.
\end{equation}
The  obtained latent variables are  given by $\mathbf{t}={\mathbf{X}\mathbf{w}}/{\|\mathbf{X}\mathbf{w}\|_2^2}$ and $\mathbf{u}=\mathbf{Yc}$. In doing so, we have made the following two assumptions:  i) the latent variables $\{\mathbf{t}_r\}_{r=1}^{R}$ are good predictors of $\mathbf{Y}$; ii) the linear (inner) relation between the latent variables $\mathbf{t}$ and $\mathbf{u}$ does exist, that is $\mathbf{U} = \mathbf{TD} + \mathbf{Z}$, where $\mathbf{Z}$ denotes the matrix of Gaussian \emph{i.i.d.} residuals. Upon combining  with the decomposition of $\mathbf{Y}$, in (\ref{Eq:PLS-Y}), we have
\begin{equation}\label{Eq:PLS-YY}
\mathbf{Y} =\mathbf{TDC}^{\text{T}} + (\mathbf{ZC}^{\text{T}} + \mathbf{F}) = \mathbf{T}\mathbf{D}\mathbf{C}^{\text{T}} + \mathbf{F}^*,
\end{equation}
where $\mathbf{F}^* $ is the residual matrix. Observe from (\ref{Eq:PLS-YY})  that the problem boils down to finding common latent variables, $\mathbf{T}$, that best explain the variance in both $\mathbf{X}$ and $\mathbf{Y}$. The prediction of new dataset $\mathbf{X}^*$ can then be performed by
$\mathbf Y^* \approx \mathbf X^* \mathbf W \mathbf D \mathbf C^{\text{T}}$.

\subsection{The $N$-way PLS Method}

\begin{figure}[t]
\centering
\includegraphics[width=0.85\textwidth]{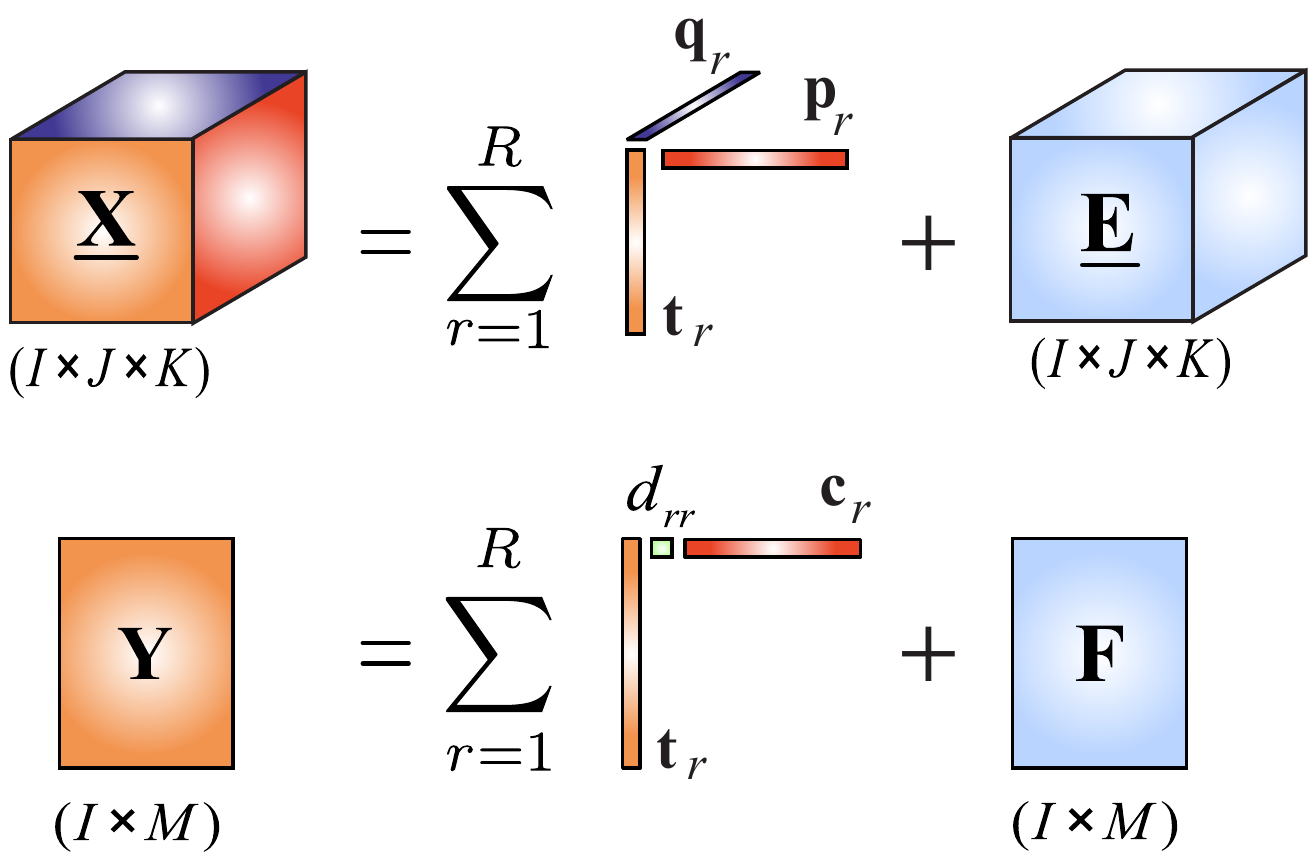}
\caption{The $N$-way PLS model performs a joint decomposition of data tensors/matrices as a sum of rank-one tensors through standard CP decomposition for the independent variables, $\underline{\mathbf{X}}$, and a sum of rank-one matrices for the responses, $\mathbf{Y}$. }
\label{fig:NPLS}
\end{figure}

The multi-way PLS (called $N$-way PLS) proposed by  \citet{bro1996multiway}  is a simple extension
of the standard PLS. The  method factorizes  an $N$th-order tensor, $\underline{\mathbf{X}}$, based on the CP decomposition, to predict response variables represented by  $\mathbf{Y}$, as shown in Figure~\ref{fig:NPLS}. For a 3rd-order tensor, $\underline{\mathbf{X}}\in\mathbb{R}^{I\times J\times K}$, and a multivariate response matrix, $\mathbf{Y}\in\mathbb{R}^{I\times M}$, with the respective elements $x_{ijk}$ and $y_{im}$,  the tensor of independent variables, $\underline{\mathbf{X}}$, is decomposed into one latent vector $\mathbf{t}\in\mathbb{R}^{I\times 1}$ and two loading vectors, $\mathbf{p} \in\mathbb{R}^{J\times 1}$ and $\mathbf{q}\in\mathbb{R}^{K\times 1}$, i.e., one loading vector per mode. As shown in Figure~\ref{fig:NPLS}, the $3$-way PLS (the $N$-way PLS for $N=3$) performs the following simultaneous tensor and matrix decompositions
\begin{equation}
\label{Eq:NPLS}
 \underline{\mathbf{X}}= \sum_{r=1}^{R}\mathbf{t}_{r}\circ\mathbf{p}_{r}\circ\mathbf{q}_{r} + \underline{\mathbf{E}}, \qquad \mathbf{Y} = \sum_{r=1}^{R} d_{rr} \, \mathbf{t}_r \mathbf{c}_r^{\text{T}} + \mathbf{F}.
\end{equation}

\noindent
The  objective is to find the vectors $\mathbf{p}_r,\mathbf{q}_r$ and $\mathbf{c}_r$,  which are the solutions of the following optimization problem
\begin{eqnarray*}
\label{Eq:NPLS-costfun}
\{\mathbf{p}_r,\mathbf{q}_r, \mathbf{c}_r\} &=& \arg\max_{\mathbf{p}_r,\mathbf{q}_r,\mathbf{c}_r}\ \mbox{cov}(\mathbf{t}_r,\mathbf{u}_r),\\
\mbox{s.t.} \quad \mathbf{t}_r = \tensor{X}\bar{\times}_1\mathbf{p}_r\bar{\times}_2\mathbf{q}_r, \qquad \mathbf{u}_r &=& \mathbf{Yc}_r, \quad \|\mathbf{p}_r\|^2_2 =\|\mathbf{q}_r\|^2_2=\|\mathbf{c}_r\|^2_2=1.
\end{eqnarray*}
Again, the problem boils down to finding the common latent variables, $\mathbf t_r$, in both $\tensor X$ and $\tensor Y$, that best explain the variance in $\tensor X$ and $\tensor Y$. The prediction of a new dataset, $\tensor X^*$, can then be performed by
$\mathbf Y^*_{(1)} \approx \mathbf X_{(1)}^* (\mathbf Q \bigodot \mathbf P)^{-1} \mathbf D \mathbf C^{\text{T}}$, where $\mathbf P = [\mathbf p_1,\ldots,\mathbf p_R]$, $\mathbf Q = [\mathbf q_1,\ldots,\mathbf q_R]$, and $\mathbf D = \text{diag}(d_{rr})$.

\subsection{HOPLS using Constrained Tucker Model}

\begin{figure}[htp]
\centering
\includegraphics[width=0.99\textwidth]{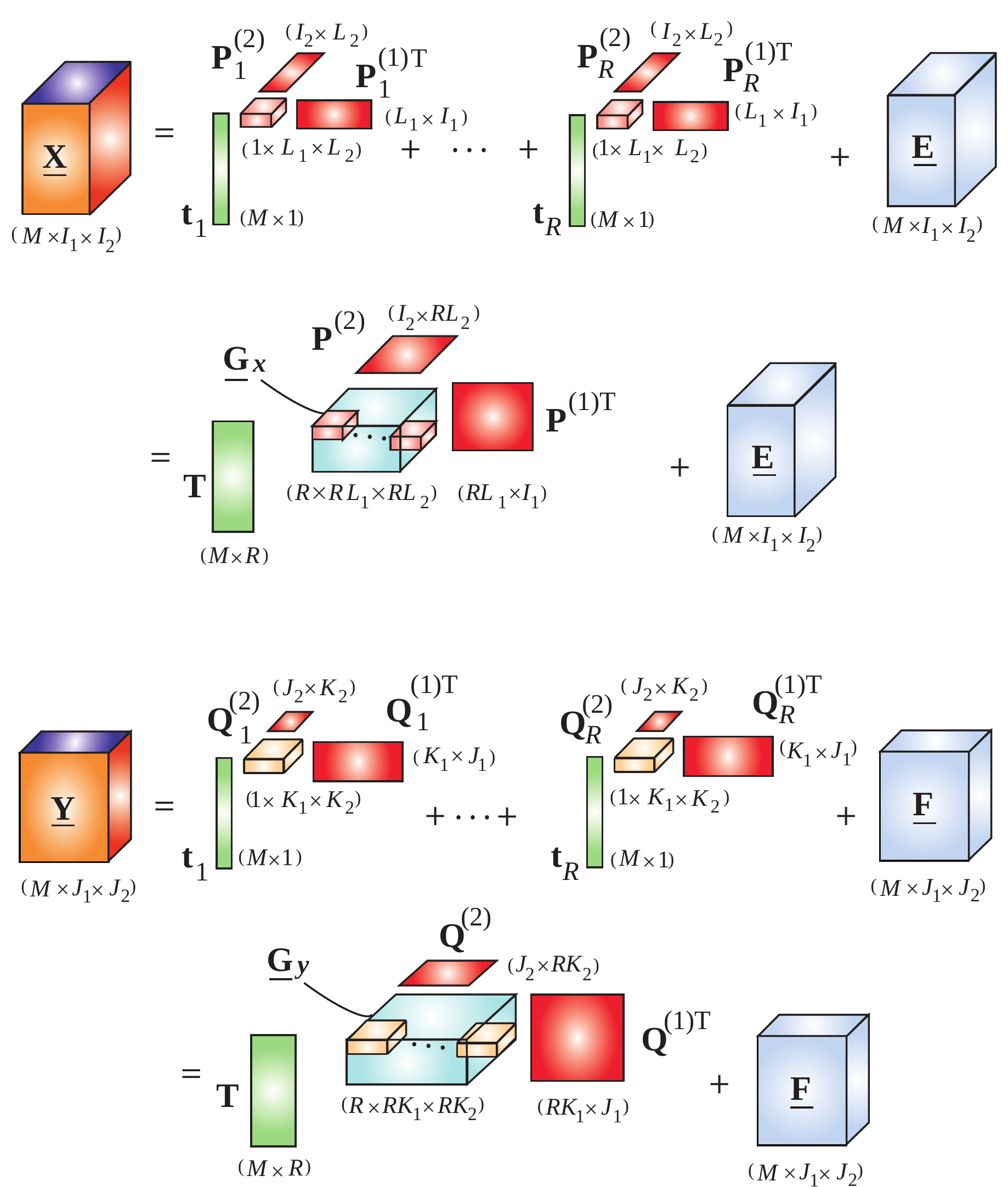}
\caption{The HOPLS model which approximates the independent variables, $\underline{\mathbf{X}}$, as a sum of rank-$(1, L_1, L_2)$ tensors. The approximation for the dependent variables, $\underline{\mathbf{Y}}$, follows a similar principle, whereby the common latent components, $\mathbf{T}$, are shared between $\tensor X$ and $\tensor Y$.}
\label{fig:HOPLS}
\end{figure}

An alternative, more flexible and general multilinear regression model, termed the higher-order partial least squares (HOPLS)
 \citep{NIPSQibin,Qibin-HOPLS}, performs simultaneously constrained Tucker decompositions for an $(N+1)$th-order independent tensor, $\underline{\mathbf{X}}\in\mathbb{R}^{M\times I_1\times\cdots\times I_N}$, and an $(N+1)$th-order dependent tensor, $\underline{\mathbf{Y}}\in\mathbb{R}^{M\times J_1\times\cdots\times J_N}$, which have the same size in the first mode, i.e., $M$ samples. Such a model allows us to find the optimal subspace approximation of $\underline{\mathbf{X}}$, in which the independent and dependent variables share a common set of latent vectors in one specific mode (i.e., samples mode). More specifically, we assume that $\underline{\mathbf{X}}$ is decomposed as a sum of rank-($1, L_1, \ldots, L_N$) Tucker blocks, while $\underline{\mathbf{Y}}$ is decomposed as a sum of rank-($1, K_1, \ldots, K_N$) Tucker blocks, which can be expressed as
\begin{equation}\label{Eq:HOPLS_T2T}
\begin{split}
\underline{\mathbf{X}} &= \sum_{r=1}^{R} \underline{\mathbf{G}}_{xr} \! \times_{1} \mathbf{t}_{r}\! \times_{2}\mathbf{P}_{r}^{(1)}\!  \cdots\! \times_{N+1}\mathbf{P}_{r}^{(N)} +\!  \underline{\mathbf{E}}_R,\\
\underline{\mathbf{Y}} &= \sum_{r=1}^{R}\underline{\mathbf{G}}_{yr}\!  \times_{1} \mathbf{t}_{r}\! \times_{2}\mathbf{Q}_{r}^{(1)}\!  \cdots\! \times_{N+1} \! \mathbf{Q}_{r}^{(N)} + \! \underline{\mathbf{F}}_R,
\end{split}
\end{equation}
where $R$ is the number of latent vectors, $\mathbf{t}_{r}\in\mathbb{R}^{M}$ is the $r$-th latent vector, $\left\{\mathbf{P}^{(n)}_{r}\right\}_{n=1}^{N}\in\mathbb{R}^{I_{n}\times L_{n}}$ and $\left\{\mathbf{Q}^{(n)}_{r}\right\}_{n=1}^{N}\in\mathbb{R}^{J_{n}\times K_{n}} $ are the loading matrices in mode-$n$, and $\underline{\mathbf{G}}_{xr}\in\mathbb{R}^{1\times L_1\times\cdots\times L_N}$ and $\underline{\mathbf{G}}_{yr}\in\mathbb{R}^{1\times K_1\times\cdots\times K_N}$ are core tensors.
By defining a latent matrix $\mathbf{T} = [\mathbf{t}_{1},\ldots,\mathbf{t}_{R}]$, \mbox{mode-$n$} loading matrix $\overline{\mathbf{P}}^{(n)} =[\mathbf{P}^{(n)}_{1},\ldots,\mathbf{P}^{(n)}_{R}]$, mode-$n$ loading matrix $\overline{\mathbf{Q}}^{(n)} =[\mathbf{Q}^{(n)}_{1},\ldots,\mathbf{Q}^{(n)}_{R}]$ and core tensors
\begin{equation}
\begin{split}
{\underline{\mathbf{G}}}_x &= \mbox{blockdiag}(\underline{\mathbf{G}}_{x1},\ldots,\underline{\mathbf{G}}_{xR})\in\mathbb{R}^{R\times RL_1\times\cdots\times RL_N},\\
{\underline{\mathbf{G}}}_y &= \mbox{blockdiag}(\underline{\mathbf{G}}_{y1},\ldots,\underline{\mathbf{G}}_{yR})\in\mathbb{R}^{R\times RK_1\times\cdots\times RK_N},
\end{split}
\end{equation}
the HOPLS model in (\ref{Eq:HOPLS_T2T}) can be rewritten as
\begin{equation}\label{Eq:HOPLS-model-T}
\begin{split}
   \underline{\mathbf{X}} &= {\underline{\mathbf{G}}}_x\times_{1}\mathbf{T}\times_{2}\overline{\mathbf{P}}^{(1)} \cdots\times_{N+1}\overline{\mathbf{P}}^{(N)}+\underline{\mathbf{E}}_R, \\
      \underline{\mathbf{Y}} &= {\underline{\mathbf{G}}}_y\times_{1}\mathbf{T}\times_{2}\overline{\mathbf{Q}}^{(1)} \cdots\times_{N+1}\overline{\mathbf{Q}}^{(N)}+\underline{\mathbf{F}}_R,
\end{split}
\end{equation}
where $\underline{\mathbf{E}}_R$ and $\underline{\mathbf{F}}_R$ are the residuals obtained after extracting $R$ latent components. The core tensors, ${\underline{\mathbf{G}}}_x$ and ${\underline{\mathbf{G}}}_y$, have a special block-diagonal structure (see Figure~\ref{fig:HOPLS}) and their elements indicate the level of local interactions between the corresponding latent vectors and loading matrices.

Benefiting from the advantages of Tucker decomposition over the CP model, HOPLS  generates  approximate latent components that better  model the data than $N$-way PLS. Specifically, the HOPLS differs substantially from the $N$-way PLS model in the sense that the sizes of loading matrices are controlled by a hyperparameter, providing a tradeoff between the model fit  and model complexity.
Note that HOPLS simplifies into $N$-way PLS if we impose the rank-1 constraints for the Tucker blocks, that is  $\forall n: \{L_n\}=1$ and $\forall m: \{K_m\}=1$.

The optimization of subspace transformation according to (\ref{Eq:HOPLS_T2T}) can be formulated as a problem of determining a set of orthogonal loadings and latent vectors. Since these terms can be optimized sequentially,  using the same criteria and based on deflation, in the following, we shall illustrate the procedure based on  finding the first latent vector, $\mathbf{t}$, and two sequences of loading matrices, $\mathbf{P}^{(n)}$ and $\mathbf{Q}^{(n)}$. Finally,  we shall denote
tensor contraction in mode one by $\langle\tensor{X},\tensor{Y}\rangle_{\{1;1\}}$, so that the cross-covariance tensor is defined by
\begin{equation}
\underline{\mathbf{C}}=\mbox{COV}_{\{1;1\}}(\underline{\mathbf{X}},\underline{\mathbf{Y}})\in\mathbb{R}^{I_1\times\cdots\times I_N \times J_1 \times \cdots\times J_N}.
\end{equation}
The  above optimization problem can now be  formulated as
\begin{align}
\label{Eq:HOPLS-finalobj}
\nonumber
   \!\! \!\!\max_{\left\{\mathbf{P^{(n)}},\mathbf{Q^{(n)}}\right\}} &\left \| [\![ \underline{\mathbf{C}}; \mathbf{P}^{(1)\,\text{T}},\!\ldots,\!\mathbf{P}^{(N)\,\text{T}},\mathbf{Q}^{(1)\,\text{T}},\!\ldots\!,\mathbf{Q}^{(N)\,\text{T}}   ]\!]\right\|_F^{2}  \\
 \mbox{s.t.} \quad & \mathbf{P}^{(n)\,\text{T}}\mathbf{P}^{(n)} =\mathbf{I}_{L_{n}}, \quad \mathbf{Q}^{(n)\,\text{T}}\mathbf{Q}^{(n)} = \mathbf{I}_{K_{n}},
\end{align}
where $[\![\ldots]\!]$ denotes the multilinear products between a tensor and a set of matrices,  $\mathbf{P}^{(n)}$ and $\mathbf{Q}^{(n)}, \; n=1,\ldots,N$, comprise the unknown parameters. This is equivalent to finding the best subspace approximation of $\tensor{C}$
which can be obtained by \mbox{rank-($L_1,\ldots,L_N,K_1,\ldots,K_N$)} HOSVD of tensor $\underline{\mathbf{C}}$. The higher-order orthogonal iteration (HOOI) algorithm, which is known to converge fast, can be employed to find the parameters $\mathbf{P}^{(n)}$ and $\mathbf{Q}^{(n)}$ by orthogonal Tucker decomposition of $\underline{\mathbf{C}}$. The detailed procedure of HOPLS is shown in Algorithm~\ref{alg:HOPLS}.

The prediction for a new sample set $\tensor X^*$ is then performed as $\mathbf Y_{(1)}^* \approx \mathbf{X}_{(1)}^* \mathbf W^x\mathbf W^y$, where
the $r$th column of $\mathbf W^x_r$ is given by $\mathbf w^x_r = (\mathbf{P}_r^{(N)}\otimes\cdots\otimes \mathbf P^{(1)}_r ) {\mathbf {G}_x}_{r(1)}^{\dagger}$ and
the $r$th row of $\mathbf W^y_r$ is given by $\mathbf w^y_r =  {\mathbf G_y}_{r(1)}^{\dagger} (\mathbf{Q}_r^{(N)}\otimes\cdots\otimes \mathbf Q^{(1)}_r )^{\text{T}}$.

\begin{algorithm}[t!]
\caption{\textbf{Higher-Order Partial Least Squares (HOPLS)}}
\label{alg:HOPLS}
{\small
\begin{algorithmic}[1]
\REQUIRE $\tensor{X}\in\mathbb R^{M\times I_1\times I_2\times \cdots\times I_N}$ and $\tensor Y\in \mathbb R^{M\times J_1\times J_2\times \cdots\times J_N}$.
\ENSURE $\{\mathbf P_r^{(n)}\}$, $\{\mathbf Q_r^{(n)}\}$, $\{{\tensor G_x}_r\}$, $\{{\tensor G_y}_r\}$, $\mathbf T$,  $n=1,\ldots, N, r=1,\ldots, R$.
\FOR {$r=1$ to $R$}
\STATE $\tensor C \leftarrow \text{COV}_{1;1}(\tensor X, \tensor Y) \in\mathbb R^{I_1\times \cdots \times I_N \times J_1\times \cdots \times J_N}$
\STATE Compute $\{\mathbf P_r^{(n)}\}$ and $\{\mathbf Q_r^{(n)}\}$ by HOOI decomposition of $\tensor C$.
\STATE $\mathbf t_r\leftarrow$ the principal eigenvector of $\tensor X \times_2 \mathbf P_r^{(1)} \cdots\times_{N+1} \mathbf P^{(N)}_{r}$
\STATE ${\tensor G_x}_r \leftarrow \tensor X \bar\times_1 \mathbf t_r\times_2 \mathbf P_r^{(1)\,\text{T}} \cdots\times_{N+1} \mathbf P^{(N)\,\text{T}}_{r}$
\STATE ${\tensor G_y}_r \leftarrow \tensor Y \bar\times_1 \mathbf t_r \times_2 \mathbf Q_r^{(1)\,\text{T}} \cdots\times_{N+1} \mathbf Q^{(N)\,\text{T}}_{r}$
\STATE Deflation: $\tensor X \leftarrow \tensor X - {\tensor G_x}_r\times_1 \mathbf t_r\times_2 \mathbf P_r^{(1)} \cdots\times_{N+1} \mathbf P^{(N)}_{r}$
\STATE Deflation: $\tensor Y \leftarrow \tensor Y - {\tensor G_y}_r\times_1 \mathbf t_r\times_2 \mathbf Q_r^{(1)} \cdots\times_{N+1} \mathbf Q^{(N)}_{r}$
\ENDFOR
\end{algorithmic}
}
\end{algorithm}

\section{Kernel HOPLS}

We next briefly introduce the concept of kernel-based tensor PLS (KTPLS or KHOPLS) \citep{zhao2013kernelization,hou_qibin2016aaai},  as a natural extension of
the HOPLS to possibly infinite dimensional and nonlinear kernel spaces. Consider $N$ pairs of tensor observations $ \{(\tensor{X}_{m},\tensor{Y}_{m})\}_{m=1}^M $, where $\tensor{X}_{m}$ denotes an $N$th-order independent tensor and $\tensor{Y}_{m}$  an $L$th-order dependent tensor.
Note that the data tensors, $\tensor X_m$ and $\tensor Y_m$,  can be concatenated to form an $(N+1)$th-order tensor $\tensor{X}\in\mathbb{R}^{M \times I_1\times\cdots\times I_N}$ and $(L+1)$th-order tensor $\tensor{Y}\in\mathbb{R}^{M\times J_1\times\cdots\times J_L}$.
 Within the kernel framework, the data tensors $\tensor{X}$ and $\tensor{Y}$ are mapped onto the Hilbert space by a reproducing kernel mapping $\phi: \quad \tensor{X}_{m}$ $\longmapsto \phi\left(\tensor{X}_{m}\right)$. For simplicity, we shall denote $\phi(\tensor{X})$ by  $\underline{\boldsymbol\Phi}$ and $\phi(\tensor{Y})$ by  $\underline{\boldsymbol\Psi}$.  The KTPLS then aims to perform the  tensor decompositions, corresponding to those in (\ref{Eq:HOPLS-model-T}), such that
\begin{equation}
\begin{split}
\underline{\boldsymbol\Phi} &= \tensor{G}_{\tensor{X}}\times_{1}\mathbf{T}\times_{2}{\mathbf{P}}^{(1)} \cdots\times_{N+1}{\mathbf{P}}^{(N)}+\tensor{E}_{\tensor{X}},\\
\underline{\boldsymbol\Psi}&= \tensor{G}_{\tensor{Y}}\times_{1}\mathbf{U}\times_{2}{\mathbf{Q}}^{(1)} \cdots\times_{L+1}{\mathbf{Q}}^{(L)}+\tensor{E}_{\tensor{Y}},\\
\mathbf{U} &= \mathbf{TD} + \mathbf{E}_{U}.
\end{split}
\end{equation}
Since within the kernel framework, the tensors $ \widetilde{\tensor{G}}_{\tensor{X}} = \tensor{G}_{\tensor{X}}\times_{2}{\mathbf{P}}^{(1)} \cdots\times_{N+1}{\mathbf{P}}^{(N)}$  and $ \widetilde{\tensor{G}}_{\tensor{Y}} =\tensor{G}_{\tensor{Y}}\times_{2}{\mathbf{Q}}^{(1)} \cdots\times_{L+1}{\mathbf{Q}}^{(L)}$  can be respectively represented as a linear combination of $\{\phi(\tensor{X}_{m})\}$ and $\{\phi(\tensor{Y}_{m})\}$, i.e., $\widetilde{\tensor{G}}_{\tensor{X}} = \boldsymbol\Phi\times_1 \mathbf{T}^{\text{T}}$ and $\widetilde{\tensor{G}}_{\tensor{Y}} = \boldsymbol\Psi\times_1 \mathbf{U}^{\text{T}}$, for the KTPLS solution we only need to   find  the latent vectors of $\mathbf{T} = [\mathbf{t}_1,\ldots,\mathbf{t}_R]$ and $\mathbf{U}=[\mathbf{u}_1,\ldots,\mathbf{u}_R]$, such that they exhibit  maximum pair-wise covariance, through solving sequentially the following optimization problems
\be
\label{eq:TKPLSobj1}
  & \max_{\{\mathbf{w}^{(n)}_r,\mathbf{v}^{(l)}_r\}}  \quad [\mbox{cov}(\mathbf{t}_r,\mathbf{u}_r)]^2,
  \ee
  where
\begin{equation}
\label{eq:TKPLSobj}
\begin{split}
 &\mathbf{t}_r = \underline{\boldsymbol\Phi}\bar{\times}_2 \mathbf{w}^{(1)}_r \cdots\bar{\times}_{N+1}\mathbf{w}^{(N)}_r,  \\  & \mathbf{u}_r =  \underline{\boldsymbol\Psi}\bar{\times}_2\mathbf{v}^{(1)}_r \cdots\bar{\times}_{L+1}\mathbf{v}^{(L)}_r.
\end{split}
\end{equation}
Upon rewriting (\ref{eq:TKPLSobj}) in a matrix form, this yields  $\mathbf{t}_r = \boldsymbol\Phi_{(1)}\widetilde{\mathbf{w}}_r$ where $\widetilde{\mathbf{w}}_r = \mathbf{w}_r^{(N)}\otimes\cdots\otimes\mathbf{w}_r^{(1)}$ and $\mathbf{u}_r=\boldsymbol\Psi_{(1)}\widetilde{\mathbf{v}}_r$ where $\widetilde{\mathbf{v}}_r = \mathbf{v}_r^{(N)}\otimes\cdots\otimes\mathbf{v}_r^{(1)}$, which can be solved by a kernelized version of the eigenvalue problem, i.e.,  $\boldsymbol\Phi_{{(1)}}\boldsymbol\Phi_{(1)}^{\text{T}}\boldsymbol\Psi_{(1)}\boldsymbol\Psi_{(1)}^{\text{T}} \mathbf{t}_r = \lambda \mathbf{t}_r$ and $\mathbf{u}_r = \boldsymbol\Psi_{(1)}\boldsymbol\Psi_{(1)}^{\text{T}} \mathbf{t}_r$ \citep{rosipal2002kernel}. Note that since the term $\boldsymbol\Phi_{{(1)}}\boldsymbol\Phi_{(1)}^{\text{T}}$ contains only the inner products between vectorized  input tensors, it can be replaced by an $M\times M$ kernel matrix $\mathbf{K}_{\tensor{X}}$, to give $\mathbf{K}_{\tensor{X}}\mathbf{K}_{\tensor{Y}}\mathbf{t}_r = \lambda \mathbf{t}_r$ and $\mathbf{u}_r =\mathbf{K}_{\tensor{Y}}\mathbf{t}_r$.

In order to exploit the rich multilinear structure of tensors,  kernel matrices should be computed using the kernel functions for tensors, i.e., $(\mathbf{K}_{\tensor{X}})_{mm'} = k\left(\tensor{X}_{m}, \tensor{X}_{m'}\right)$ and $(\mathbf{K}_{\tensor{Y}})_{mm'} = k\left(\tensor{Y}_{m}, \tensor{Y}_{m'}\right)$ (for more details see Section \ref{sec:tensorkernel}).

 Finally, for a new data sample, $\tensor{X}^{*}$, the prediction of $\mathbf y$, denoted by $\mathbf y^*$ is performed by computing the vector
\begin{equation}\label{Eq:ktplspredionfun}
\mathbf{y}^{*T} =   \mathbf{k}^{*T}\mathbf{U}(\mathbf{T}^{\text{T}} \mathbf{K}_{\tensor{X}} \mathbf{U})^{-1} \mathbf{T}^{\text{T}} \mathbf{Y}_{(1)},
\end{equation}
where $(\mathbf{k}^*)_m = k\left(\tensor{X}_{m},\tensor{X}^{*}\right)$ and the vector $\mathbf{y}^*$ can be reshaped to a tensor  $\tensor{Y}^*$. The detailed procedure of KHOPLS is shown in Algorithm~\ref{alg:KHOPLS}.

\begin{algorithm}[t!]
\caption{\textbf{Kernel Higher-Order Partial Least Squares (KHOPLS)}}
\label{alg:KHOPLS}
{\small
\begin{algorithmic}[1]
\REQUIRE $\tensor{X}\in\mathbb R^{M\times I_1\times I_2\times \cdots\times I_N}$ and $\tensor Y\in \mathbb R^{M\times J_1\times I_2\times \cdots\times J_L}$.
\ENSURE $\mathbf T, \mathbf U, \mathbf K_{\tensor X}$.
\STATE Compute kernel matrix $\mathbf K_{\tensor{X}}$ and $\mathbf K_{\tensor Y}$ by using tensor kernel functions
\STATE Compute $\mathbf t_r$ by solving $\mathbf K_{\tensor X}\mathbf K_{\tensor Y} \mathbf t_r = \lambda \mathbf t_r$, $r=1,\ldots, R$
\STATE $\mathbf u_r = \mathbf K_{\tensor{Y}}\mathbf t_r$, $r=1,\ldots, R$
\STATE $\mathbf T = [\mathbf t_1,\ldots, \mathbf t_R]$ and $\mathbf U = [\mathbf u_1,\ldots, \mathbf u_R]$
\end{algorithmic}
}
\end{algorithm}

The intuitive interpretation of (\ref{Eq:ktplspredionfun})
can take several forms: i) as a linear combination of $M$ observations $\{\tensor{Y}_{m}\}$ and the coefficients $\mathbf{k}^{*T}\mathbf{U}(\mathbf{T}^{\text{T}} \mathbf{K}_{\tensor{X}} \mathbf{U})^{-1} \mathbf{T}^{\text{T}}$; ii)  as a linear combination of $M$ kernels, each  centered around a training point, i.e., $y_{j}^* = \sum_{m=1}^M \alpha_m k\left(\tensor{X}_{m},\tensor{X}^{*}\right)$, where $\alpha_m = \left(\mathbf{U}(\mathbf{T}^{\text{T}} \mathbf{K}_{\tensor{X}} \mathbf{U})^{-1} \mathbf{T}^{\text{T}} \mathbf{Y}_{(1)}\right)_{mj}$; iii)  the vector  $\mathbf{t}^*$ is obtained by nonlinearly projecting the tensor $\tensor{X}^{*}$ onto the latent space spanned by $\mathbf{t}^{*T}=\mathbf{k}^{*T}\mathbf{U}(\mathbf{T}^{\text{T}} \mathbf{K}_{\tensor{X}} \mathbf{U})^{-1}$, then $\mathbf{y}^{*T}$ is predicted by a linear regression against $\mathbf{t}^*$, i.e., $\mathbf{y}^{*T} = \mathbf{t}^{*T} \mathbf{C}$, where the regression coefficients are given by the matrix $\mathbf{C} = \mathbf{T}^{\text{T}}\mathbf{Y}_{(1)}$.

Note that, in general, to ensure a strict linear relationship between the latent vectors and output in the original spaces, the kernel functions for data $\tensor{Y}$ are restricted to linear  kernels.

\section{Kernel Functions in  Tensor Learning}
\label{sec:tensorkernel}

\begin{figure}[t]
\centering
  \includegraphics[width=0.71\textwidth]{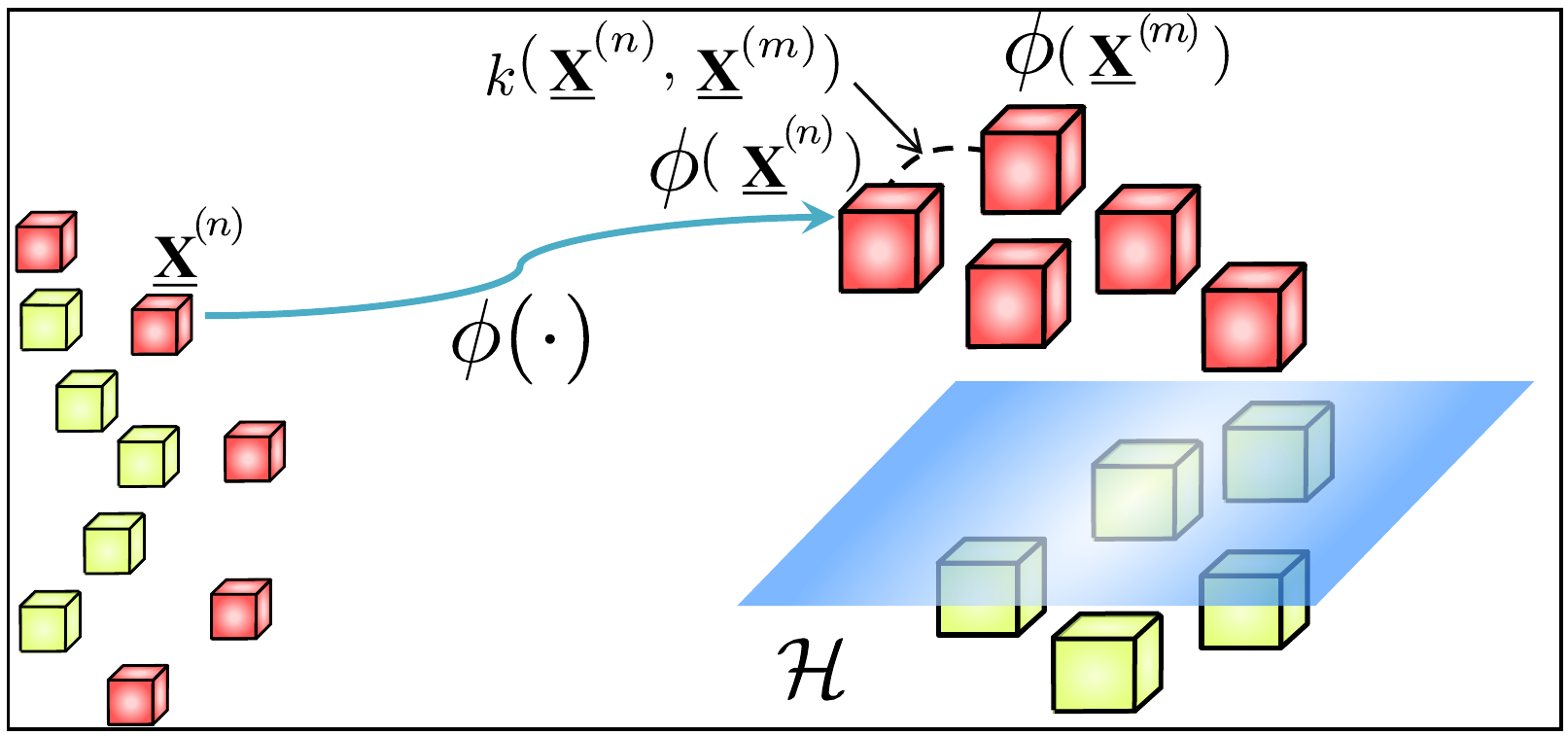}
  \caption{The concept of kernel learning for tensors. Tensor observations are mapped into the reproducing kernel Hilbert space $\mathcal{H}$ by a nonlinear mapping function $\phi(\cdot)$. The kernel function serves as a similarity measure between two tensors.}
  \label{fig:mappingtensor}
\end{figure}

Kernel functions can be considered as a means for defining a new topology which implies  {\em a priori} knowledge about the invariance in the input space \citep{scholkopf2002learning} (see Figure~\ref{fig:mappingtensor}). Kernel algorithms can also be used for estimation of vector-valued nonlinear and nonstationary signals \citep{tobar2014multikernel}. In this section, we discuss the kernels for tensor-valued inputs, which should exploit multi-way structures   while maintaining the notion of  similarity measures.

Most straightforward tensor-valued reproducing kernels are a direct generalization from vectors to $N$th-order tensors, such as the following kernel functions $k: \tensor{X} \times \tensor{X} \rightarrow \mathbb{R}$,  given by
\begin{equation}
\begin{split}
&\mbox{Linear kernel:} \quad  k(\tensor{X},\tensor{X}')  = \langle \mbox{vec}(\tensor{X}), \mbox{vec}(\tensor{X}') \rangle,  \\
&\mbox{Gaussian-RBF kernel:} \quad  k(\tensor{X},\tensor{X}') = \mbox{exp}\left( -\frac{1}{2\beta^2} \|\tensor{X}-\tensor{X}'\|_F^2  \right).
\end{split}
\end{equation}
In order to define a similarity measure that directly exploits  the multilinear algebraic structure of input tensors, \citet{signoretto2011kernel, signoretto2012classification} proposed a tensorial kernel which both exploits the algebraic geometry  of tensors spaces and provides a similarity measure between the different subspaces spanned by higher-order tensors. Another such kernel is the  product kernel which can be defined by $N$ factor kernels, as $k(\tensor{X},\tensor{X'}) = \prod_{n=1}^N k\left(\mathbf{X}_{(n)},\mathbf{X'}_{(n)}\right)$,
where each factor kernel represents a similarity measure between  \mbox{mode-$n$} matricizations of two tensors, $\tensor X$ and $\tensor{X'}$.

One  similarity measure between matrices is the so called Chordal distance, which is a projection of the Frobenius norm on a Grassmannian manifolds \citep{signoretto2011kernel}. For example, for   an $N$th-order tensor, $\tensor{X}$, upon the applications of SVD to its mode-$n$ unfolding, that is, $\mathbf{X}_{(n)} = \mathbf{U}_{\mathbf{X}}^{(n)} \mathbf{\Sigma}^{(n)}_{\mathbf{X}}\mathbf{V}^{(n)\,\text{T}}_{\mathbf{X}}$, the Chordal distance-based kernel for tensorial data is defined as
\begin{equation}\label{Eq:kernelfun-tensor}
k(\tensor{X},\tensor{X'}) = \prod_{n=1}^N \mbox{exp}\left(-\frac{1}{2\beta^2} \left\|\mathbf{V}^{(n)}_{\mathbf{X}}\mathbf{V}^{(n)\,\text{T}}_{\mathbf{X}} - \mathbf{V}^{(n)}_{\mathbf{X}'}\mathbf{V}^{(n)\,\text{T}}_{\mathbf{X}'} \right\|_F^2 \right),
\end{equation}
where $\beta$ is a kernel parameter.
It should be emphasized that such a tensor kernel  ensures (provides) rotation and reflection invariance for elements on the Grassmann manifold \citep{signoretto2011kernel}.

\citet{zhao2013kernelization}  proposed  a whole family of probabilistic product kernels  based on generative models. More specifically, an $N$th-order tensor-valued  observation is first mapped onto an $N$-dimensional model space, then information divergence is applied as a similarity measure in such a model space (see also \citep{Cichocki-Entropy11,cichocki2015log}).

The advantage of probabilistic tensor kernels is that they can deal with multiway data with missing values and with variable data length. Given that probabilistic kernels provide a way to model one tensor from $N$ different viewpoints which  correspond to different lower-dimensional vector spaces, this makes it possible for multiway relations to  be captured within a similarity measure.

A general similarity  measure between two tensors, $\tensor{X}$ and $\tensor{X}'$, in mode-$n$ can be defined as
\begin{equation}\label{eq:tensordivergence}
S_n(\tensor{X}\|\tensor{X}')= D\left(p(\mathbf{x}|\boldsymbol\Omega_n^{\tensor{X}})\| q(\mathbf{x}|\boldsymbol\Omega^{\tensor{X}'}_n)\right),
\end{equation}
where $p$ and  $q$ are distributions which respectively represent the probability density function for $\tensor{X}$ and $\tensor{X}'$,  $D(p\|q)$ is an information divergence between two distributions, while $\Omega_n^{\tensor{X}}$ denotes the parameters of a mode-$n$ distribution of $\tensor{X}$, which depends on the model assumption.

For simplicity, we usually assume Gaussian models for all modes of tensor $\tensor{X}$, so that  $\Omega$ can be expressed by the mean values and covariance matrices of that distribution.
One simple and very useful information divergence is the standard  \emph{symmetric Kullback-Leibler} (sKL) divergence \citep{moreno2003kullback}, expressed as
\begin{equation}
\begin{split}
D_{sKL}\left(p(\mathbf{x}|\boldsymbol \lambda)\|q(\mathbf{x}|\boldsymbol \lambda')\right) &= \frac{1}{2}  \int_{-\infty}^{+\infty} p(\mathbf{x}|\boldsymbol\lambda) \log  \frac{p(\mathbf{x}|\boldsymbol{\lambda})}{q(\mathbf{x}|\boldsymbol{\lambda}')}   d\mathbf{x}\\  &+  \frac{1}{2} \int_{-\infty}^{+\infty} q(\mathbf{x}|\boldsymbol{\lambda}') \log  \frac{q(\mathbf{x}|\boldsymbol{\lambda}')}{p(\mathbf{x}|\boldsymbol{\lambda})}  d\mathbf{x},
\end{split}
\end{equation}
where $\boldsymbol\lambda$ and $\boldsymbol\lambda'$ are  parameters of the probability distributions.
An alternative simple divergence is the symmetric Jensen-Shannon (JS) divergence
\citep{chan2004family, endres2003new,Cichocki-Entropy11}, given by
\begin{equation}\label{eq:JSdivergence}
D_{JS}(p\|q) = \frac{1}{2} \text{KL}(p\|r) + \frac{1}{2} \text{KL}(q\|r),
\end{equation}
where  $\text{KL}(\cdot\|\cdot)$ denotes the Kullback-Leibler divergence and $r(\mathbf{x}) = \frac{p(\mathbf{x}) + q(\mathbf{x})}{2}$.

\renewcommand{\arraystretch}{1.5}
 \begin{table*}[t]
\caption{Kernel functions for tensor}
\label{tab:tensorkernels}
\centerline{
\shadingbox{
{\begin{tabular}{@{}c@{}}
\hline \hline
$k(\tensor X, \tensor Y) = \Big\langle \Phi (\text{vec}(\tensor{X})), \Phi (\text{vec}(\tensor{Y})) \Big\rangle $ \\ \hline
$k(\tensor X, \tensor Y) = \Big\langle \Phi (\mathbf X_{(1)},\ldots, \mathbf X_{(N)} ), \Phi (\mathbf Y_{(1)},\ldots, \mathbf Y_{(N)} ) \Big\rangle $ \\ \hline
$k(\tensor X, \tensor Y) = \Big\langle \Phi (\tensor{X}^{'}), \Phi (\tensor{Y}^{'}) \Big\rangle $ , where $\tensor{X}^{'}$ is a CP decomposition of $\tensor X$\\ 
\hline
$k(\tensor X, \tensor Y) = \Big\langle \Phi\big(\tensor{G}_X^{(1)}, \ldots, \tensor{G}_X^{(N)}\big), \Phi\big(\tensor{G}_Y^{(1)}, \ldots, \tensor{G}_Y^{(N)}\big) \Big\rangle$\\
where  $\tensor{G}_X^{(1)}, \ldots, \tensor{G}_X^{(N)}$ are cores obtained by TT decomposition of $\tensor{X}$,\\
 and $\tensor{G}_Y^{(1)}, \ldots, \tensor{G}_Y^{(N)}$ are cores obtained by TT decomposition of $\tensor{Y}$.\\
\hline\hline
\end{tabular}}
}}
\end{table*}

In summary, a probabilistic product kernel for tensors can be defined as
\begin{equation}\label{Eq:tensorkernel-divergence}
k(\tensor{X},\tensor{X}') =\alpha^2  \prod_{n=1}^{N} \mbox{exp}\left(-\frac{1}{2\beta_n^2} S_n(\tensor{X}\|\tensor{X}')\right),
\end{equation}
where $S_n(\tensor{X}\|\tensor{X}')$ is a suitably chosen probabilistic divergence, $\alpha$ denotes a magnitude parameter and $[\beta_1, \ldots, \beta_N]$ are length-scale parameters.


An intuitive interpretation of the kernel function in (\ref{Eq:tensorkernel-divergence}) is that $N$th-order tensors are assumed to be generated from $N$ generative models, while the similarity of these models, measured by information divergence, is employed to provide a multiple kernel with well defined properties and conditions. This kernel can then effectively capture the statistical properties of tensors, which promises to be a powerful tool for multidimensional structured data analysis, such as video classification and multichannel  feature extraction from brain electrophysiological responses.

\begin{figure}[tp!]
\centering
  \includegraphics[width=0.999\textwidth]{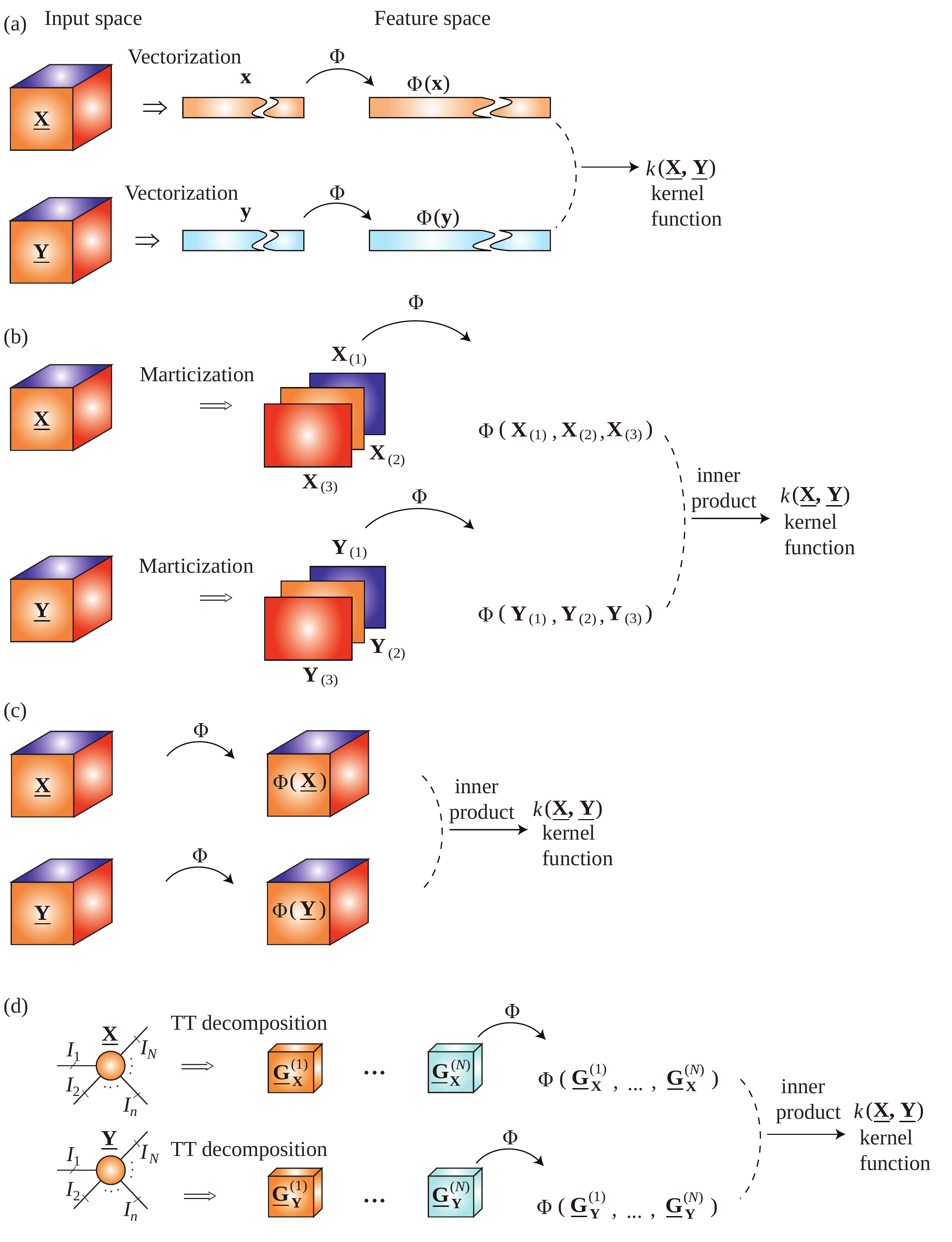}\\
  \caption{Approaches for the construction of tensor kernel functions. The kernel function can be defined based on (a) vectorization of data tensors, (b) matricization of data tensors, (c)  raw tensors format (direct approach), and (d)    core tensors from low-rank tensor decompositions.  }
  \label{fig:tensorkernel}
\end{figure}

There are many possibilities to define  kernel functions for tensors, as outlined in Table~\ref{tab:tensorkernels} and Figure~\ref{fig:tensorkernel}.

\section{Tensor Variate Gaussian Processes (TVGP)}

Gaussian processes (GP) can be considered as a class of probabilistic models which  specify a distribution over a function space, where the inference is performed directly  in the function space
 \citep{rasmussen2006gaussian}. The GP model for tensor-valued input spaces, called the Tensor-based Gaussian Processes (Tensor-GP), is designed so as to take into account the tensor structure of data
  \citep{zhao2013kernelization,zhao2014tensor,hou2015online}.

Given a paired dataset of $M$ observations $\mathcal{D}=\{(\tensor{X}_m, y_m)|m =1,\ldots, M\}$, the tensor inputs for all $M$ instances (cases) are aggregated into an $(N+1)$th-order {\it design concatenated tensor} $\tensor{X} \in\mathbb{R}^{M \times I_1\times \cdots\times I_N}$, while the targets are collected in the vector $\mathbf{y}=[y_1,\ldots,y_M]^{\text{T}}$. After observing the training data $\mathcal{D}=\{\tensor{X},\mathbf{y}\}$, the aim is to make inferences about the relationship between the inputs $\underline \bX$ and targets (output $\mathbf y$), i.e., to perform  estimation of the conditional distribution of the targets given the inputs, and in doing so to perform the  prediction based on a new input, $\tensor{X}_*$, which has not been seen in the training set.

The distribution of observations can be factored over cases in the training set by
$
\mathbf{y} \sim \prod_{m=1}^{M} \mathcal{N}(y_m| f_m, \sigma^2),
$
where $f_m$ denotes the latent function $f(\tensor{X}_m)$, and $\sigma^2$ denotes noise variance.

A Gaussian process prior can be placed over the latent function, which implies that any finite subset of latent variables has a multivariate Gaussian distribution, denoted by
\begin{equation}
f(\tensor{X})  \sim \mathcal{GP}(m(\tensor{X}), k(\tensor{X},\tensor{X'})\mid \boldsymbol{\theta}),
\end{equation}
where $m(\tensor{X})$ is the mean function which  is usually set to zero for  simplicity, and $k(\tensor{X},\tensor{X'})$ is the covariance function (e.g., kernel function) for tensor data, with a set of hyper-parameters denoted by $\boldsymbol\theta$.

The hyper-parameters from the observation model and the GP prior are collected in $\Theta =\{\sigma, \boldsymbol\theta\}$. The model is then hierarchically extended to the third level by  also giving priors over the hyperparameters in $\Theta$.

To incorporate the knowledge that the training data provides about the function $\mathbf f(\tensor{X})$, we can use the Bayes rule to infer the posterior of the latent function $\mathbf{f}(\tensor X) = [f(\tensor{X}_1), \ldots, f(\tensor{X}_M)]^{\text{T}} $ by
\begin{equation}\label{eq:posteriorf}
\condp{\mathbf{f}}{\mathcal{D}, \Theta} = \frac{\condp{\mathbf{y}}{\mathbf{f}, \sigma} \condp{\mathbf{f}}{\tensor{X},\theta}}{\int \condp{\mathbf{y}}{\mathbf{f},\sigma} \condp{\mathbf{f}}{\tensor{X}, \theta}d \mathbf{f}},
\end{equation}
where the denominator in (\ref{eq:posteriorf}) can be interpreted as the marginal likelihood obtained by the integration over $\mathbf{f}$, to yield
\begin{equation}
\mathbf{y}|\tensor{X}, \theta, \sigma^2 \sim \mathcal{N}(\mathbf{y}|\mathbf{0},\mathbf{K}+\sigma^2\mathbf{I}),
\end{equation}
where $(\mathbf{K})_{ij} = k(\tensor{X}_i, \tensor{X}_j)$ denotes the covariance matrix or \emph{kernel matrix}.

Note that, since the Gaussian observation model is analytically tractable and so it  avoids the approximate inference, the conditional posterior of the latent function $\mathbf{f}$ is Gaussian, and the posterior of $f_*$ is also Gaussian, together with the observation $y_*$.

Finally, the predictive distribution of $y_*$ corresponding to $\tensor{X}_*$ can be inferred as $y_*|\tensor{X}_*,\tensor{X},\mathbf{y}, \Theta \sim \mathcal{N}(\overline{y}_*, \mbox{cov}(y_*))$, where
\begin{equation}
\begin{split}
\overline{y}^* &= k(\tensor{X}_*,\tensor{X}) (k(\tensor{X},\tensor{X}) + \sigma_n^2 \mathbf{I})^{-1} \mathbf{y},\\
\mbox{cov}(y^*) &= k(\tensor{X}_*,\tensor{X}_*) - k(\tensor{X}_*,\tensor{X})(k(\tensor{X},\tensor{X})+ \sigma^2\mathbf{I})^{-1}k(\tensor{X},\tensor{X}_*).
\end{split}
\end{equation}
The classification problem consisting of $N$th-order tensors $\tensor{X}_m\in\mathbb{R}^{I_1\times\cdots\times I_N}$ which are associated with target classes $y_m \in\{1,\ldots,C\}$, where $C>2$ and  $m=1,\ldots,M$, was investigated by \citet{zhao2013tensor}.
All class labels are collected in the $M\times 1$ target vector $\mathbf{y}$, and all tensors are concatenated in an $(N+1)$th-order tensor $\tensor{X}$ of size $M\times I_1\times\cdots\times I_N$.
Given the latent function $\mathbf{f}_m =[{f}_m^1, {f}_m^2,\ldots,{f}_m^C]^{\text{T}}= \mathbf{f}(\tensor{X}_m)$ at the observed input location $\tensor{X}_m$, the class labels $y_m$ are assumed to be independent and identically distributed, as defined by a multinomial probit likelihood model $p(y_m|\mathbf{f})$. The latent vectors from all observations are denoted by $\mathbf{f} = [f_1^1, \ldots, f_M^1, f_1^2,\ldots, f_M^2, \ldots,$ $f_1^C, \ldots, f_M^C]^{\text{T}}$.

The objective of TVGP is to predict the class membership for a new input tensor, $\tensor{X}_*$, given the observed data, $\mathcal{D}=\{\tensor{X}, \mathbf{y}\}$. Gaussian process priors are placed on the latent function related to each class, which is the common assumption in multi-class GP classification (see
 \citep{rasmussen2006gaussian,riihimaki2012nested}).

  Such a specification results in the  zero-mean Gaussian prior for $\mathbf{f}$, given by
\begin{equation}\label{eq:gpPrior}
p(\mathbf{f}|\tensor{X}) = \mathcal{N}(\mathbf{0}, \mathbf{K}),
\end{equation}
where $\mathbf{K}$ is a $C M \times C M$ block-diagonal covariance matrix and every matrix $\mathbf{K}^1,\ldots, \mathbf{K}^C$ (of size $M \times M$) on its diagonal corresponds to each class. The element $K_{i,j}^c$ in a $c$th class covariance matrix defines the prior covariance between $f_i^c$ and $f_j^c$, which is governed by a kernel function $k(\tensor{X}_i, \tensor{X}_j)$, i.e., $K_{i,j}^c = k(\tensor{X}_i,\tensor{X}_j) = \text{Cov}(f_i^c,f_j^c)$ within the class $c$. Note that since the kernel function should be defined in tensor-variate input space, hence  commonly used kernel functions, such as Gaussian RBF, are infeasible. Therefore, the framework of probabilistic product kernel for tensors, described in Sec.\ref{sec:tensorkernel}, should be applied. In a kernel function, the hyperparameters are defined so as to control the smoothness properties and overall variance of latent functions, and  are usually combined  into one vector $\boldsymbol\theta$.
For simplicity, we use the same $\boldsymbol\theta$ for all classes. For the likelihood model, we consider the multinomial probit, which is a generalization of the probit model, given by
\begin{equation}\label{eq:probit}
p(y_m|\mathbf{f}_m) = \text{E}_{p(u_m)} \left\{ \prod_{c=1, c\neq y_m}^{C} \Phi(u_m + f_m^{y_m} - f_m^c)\right\},
\end{equation}
where $\Phi$ denotes the cumulative distribution function of the standard normal distribution, and the auxiliary variable $u_m$ is distributed as $p(u_m) = \mathcal{N}(0,1)$.

By applying the Bayes theorem, the posterior distribution of the latent function is given by
\begin{equation}
p(\mathbf{f}|\mathcal{D},\boldsymbol\theta) = \frac{1}{Z} p(\mathbf{f}|\tensor{X}, \boldsymbol\theta)\prod_{m=1}^M p(y_m|\mathbf{f}_m),
\end{equation}
where $Z = \int p(\mathbf{f}|\tensor{X},\boldsymbol\theta)\prod_{m=1}^M p(y_m|\mathbf{f}_m)d\mathbf{f}$ is known as the marginal likelihood.

The inference for a test input, $\tensor{X}_*$, is performed in two steps. First, the posterior distribution of the latent function, $\mathbf{f}_*$, is given as $p(\mathbf{f}_*|\mathcal{D},\tensor{X}_*,\boldsymbol\theta) =\int p(\mathbf{f}_*|\mathbf{f},\tensor{X}_*,
\boldsymbol\theta)p(\mathbf{f}|\mathcal{D},\boldsymbol\theta)d\mathbf{f}$. Then, we compute the posterior predictive probability of $\tensor{X}_*$, which is given by $p(y_*|\mathcal{D},\tensor{X}_*,\boldsymbol\theta) = \int p(y_*|\mathbf{f}_*)p(\mathbf{f}_*|\mathcal{D},\tensor{X}_*,\boldsymbol\theta)d\mathbf{f}_*$.

Since the non-Gaussian likelihood model results in an analytically intractable posterior distribution,  variational approximative methods can be used  for approximative inference. The additive multiplicative nonparametric regression (AMNR) model constructs $f$ as the sum of local functions which take  the components of a rank-one tensor as inputs \citep{imaizumi2016doubly}. In this approach, the function space and the input space are simultaneously decomposed.

For example, upon applying  the CP decomposition, to give $\tensor{X} = \sum_{r=1}^R \mathbf{u}_r^{(1)}\circ\cdots\circ\mathbf{u}^{(N)}_r $,
function $f$ is decomposed into a set of local functions as
\begin{equation}
f(\tensor{X}) = f(\mathbf{U}^{(1)}, \ldots,\mathbf{U}^{(N)}) = \sum_{r=1}^R \prod_{n=1}^N f^{(n)}_{r} (\mathbf{u}^{(n)}_{r}).
\end{equation}
For each local function, $f^{(n)}_r$, consider the GP prior $GP(f^{(n)}_r)$, which is represented as a multivariate Gaussian distribution $\mathcal{N}(0_{I_n}, \mathbf K^{(n)}_r)$.
The likelihood function is then $\prod_{m} \mathcal{N}(y_m | f(\tensor{X}_m),\sigma^2)$. By employing the Gaussian processes regression, we can obtain the posterior distribution $p(f|\tensor{X}, \mathbf y)$, and the predictive distribution of $p(f^*| \tensor{X}, \mathbf y, \tensor{X}^*)$.

\section{Support Tensor Machines}

In this section, Support Vector Machine (SVM) is briefly reviewed followed by the generalizations of SVM to matrices and tensors.

\subsection{Support Vector Machines}

\begin{figure}[t!]
\centering
  \includegraphics[width=0.47\textwidth]{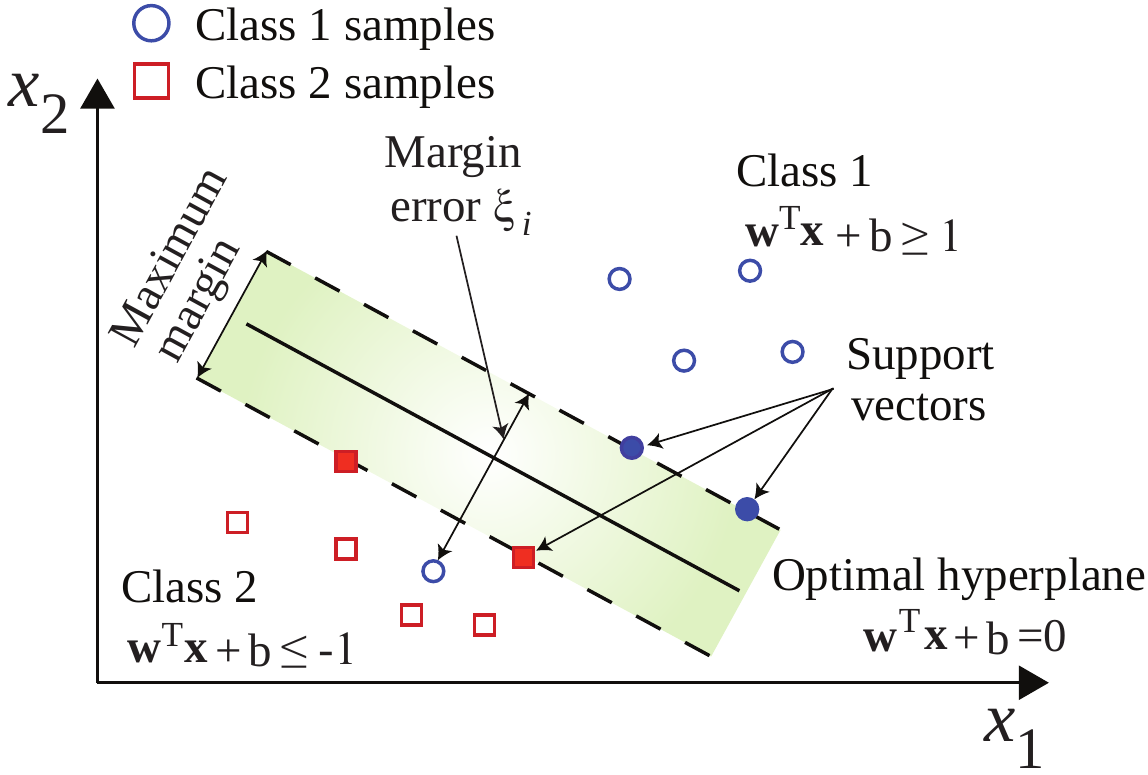}\hspace{0.01cm}
  \includegraphics[width=0.50\textwidth]{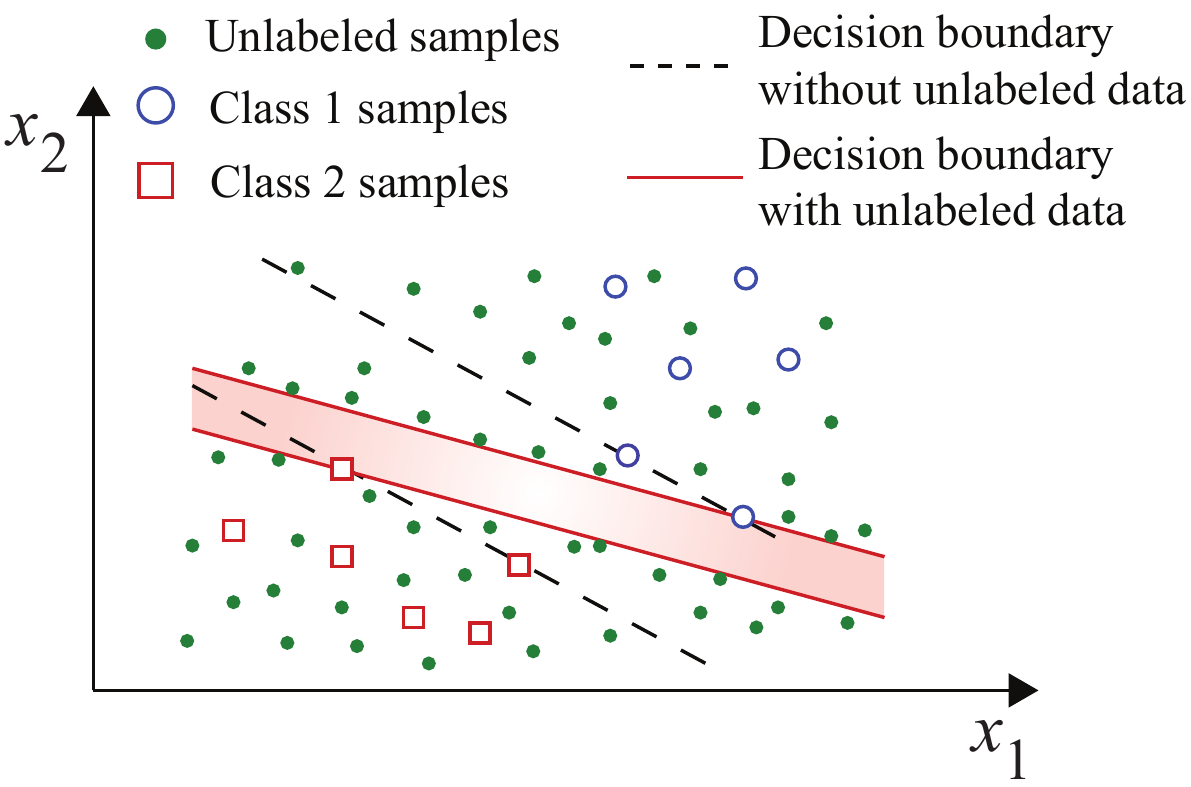}\\
  {\small(a)\hspace{12pc}(b)}\\
  \caption{Principle of SVM. (a) The support vector machine (SVM) classifier aims to maximize the margin between the Class 1 and Class 2 training measurements. (b) Semi-supervised learning -- exploits both labeled and unlabeled data.}
  \label{fig:svm}
\end{figure}

The classical SVM \citep{cortes1995support} aims to find a classification hyperplane which maximizes the margin between the `positive' measurements and the `negative' measurements, as illustrated in Figure~\ref{fig:svm}. Consider $M$ training measurements, $\mathbf{x}_m \in \mathbb{R}^{L}
(m=1,\ldots, M) $, associated with one of the  two class labels  of interest $y_m\in\{+1, -1\}$. The standard SVM, that is, the soft-margin SVM, finds a projection vector $\mathbf{w}\in\mathbb{R}^L$ and a bias $b\in\mathbb{R}$ through the minimization of the cost function
\begin{equation}
\label{eq:svm}
\begin{split}
\min_{\mathbf{w}, b, \boldsymbol\xi}\  & J(\mathbf{w}, b, \boldsymbol\xi) = \frac{1}{2}\|\mathbf{w}\|^2
 + C\sum_{m=1}^M \xi_m,\\
\mbox{s.t.} \quad & y_m(\mathbf{w}^{\text{T}}\mathbf{x}_m +b) \geq 1- \xi_m, \quad \boldsymbol\xi \geq 0,
 \quad m=1, \ldots, M,
 \end{split}
\end{equation}
where $\boldsymbol \xi = [\xi_1, \xi_2,\ldots,\xi_M]^{\text{T}}\in\mathbb{R}^M$ is the vector of all slack variables\footnote{In an optimization problem, a slack variable is a variable that is added to an inequality constraint to transform it to an equality.}  required to deal with this linearly nonseparable problem. The parameter $\xi_m,  \;\; m=1, \ldots, M$, is also called the marginal error for the $m$th training measurement, while the margin is $\frac{2}{\|\mathbf{w}\|^2}$. When the classification problem is linearly separable, we can set $\boldsymbol \xi =0$. The decision function for such classification is the binary $y(\mathbf{x}_m) = \text{sign} [\mathbf{w}^{\text{T}} \mathbf{x}_m + b]$.

The soft-margin SVM can be simplified into the least squares SVM
 \citep{suykens1999least,zhao2014least}, given by
\begin{equation}
\label{eq:lssvm}
\begin{split}
\min_{\mathbf{w},b,\boldsymbol\varepsilon}\ & J(\mathbf{w},b, \boldsymbol\varepsilon) = \frac{1}{2}\|\mathbf{w}\|^2 + \frac{\gamma}{2} \boldsymbol\varepsilon^{\text{T}}\boldsymbol\varepsilon,\\
\mbox{s.t.} \quad & y_m( \mathbf{w}^{\text{T}}\mathbf{x}_m+b) = 1- \varepsilon_m, \quad m=1, \ldots, M,
\end{split}
\end{equation}
where the penalty  coefficient $\gamma > 0 $.

 The two differences between the soft-margin SVM and the least squares SVM are: 1) the inequality constraints in (\ref{eq:svm}) are replaced by equality constraints in (\ref{eq:lssvm}); and 2) the loss $\sum_{m=1}^M \xi_m (\xi_m \geq 0)$ is replaced by a squared loss. These two modifications enable the solution of the least-squares SVM  to be more conveniently obtained, compared to the soft-margin SVM.

{\bf Remark 4}: According to the statistical learning theory \citep{vapnik1998statistical},  SVM-based learning  performs well when the number of  training measurements is larger than the complexity of the model. Moreover, the complexity of the model and the number of parameters to describe the model are always in a  direct proportion.

\subsection{Support Matrix Machines (SMM)}

Recall that the general SVM problem in (\ref{eq:svm}) deals with data expressed in a  vector form, $\mathbf{x}$. If, on the other hand, the data are collected as matrices, $\mathbf X$, these are typically  first vectorized by $\text{vec}(\mathbf{X})$, and then fed to the SVM. However, in many classification problems, such as EEG classification, data is naturally expressed by matrices, the structure  information  could be exploited for a better solution.  For matrices, we have $\langle \mathbf{W},\mathbf{W} \rangle=\text{tr}(\mathbf{W}^{\text{T}}\mathbf{W})$, and so intuitively we could consider the following formulation for the soft margin support matrix machine \citep{luo2015support}
\begin{equation}
\begin{aligned}
& \min_{\mathbf{W},b,\boldsymbol \xi} \, \frac{1}{2}   \text{tr}(\mathbf{W}^{\text{T}}\mathbf{W})  +C\sum_{m=1}^{M}\xi_m\\
\text{s.t.} \hspace{2mm} &  \hspace{3mm} y_m(  \text{tr}(\mathbf{W}^{\text{T}}\mathbf{X}_m)   +b)\geq 1-\xi_m,  \, \boldsymbol\xi \geq 0,
 \quad m=1, \ldots, M.
\end{aligned}
\end{equation}

However, observe that when $\mathbf{w}=\text{vec}(\mathbf{W}^{\text{T}})$ the above formulation is essentially equivalent to (\ref{eq:svm}) and does not exploit the correlation between the data channels inherent to a matrix structure, since

\begin{equation}
\begin{aligned}
\text{tr}(\mathbf{W}^{\text{T}}\mathbf{X}_m)&=\text{vec}(\mathbf{W}^{\text{T}})^{\text{T}}\text{vec}(\mathbf{X}_m^{\text{T}})=\mathbf{w}^{\text{T}}\mathbf{x}_m,\\
\text{tr}(\mathbf{W}^{\text{T}}\mathbf{W})&=\text{vec}(\mathbf{W}^{\text{T}})^{\text{T}}\text{vec}(\mathbf{W}^{\text{T}})=\mathbf{w}^{\text{T}}\mathbf{w},\\
\end{aligned}
\end{equation}

In order  to include correlations between the rows or columns of  data matrices, $\mathbf{X}_m$, the nuclear norm can be introduced so that  the problem becomes

\begin{equation}\label{smm}
\begin{aligned}
& \text{arg}\underset{\mathbf{W},b, \xi_m}{\text{min}}\ \frac{1}{2}   \text{tr}(\mathbf{W}^{\text{T}}\mathbf{W}) +\tau \|\mathbf{W}\|_*     +C\sum_{m=1}^{m}\xi_m\\
&\text{s.t.}  \hspace{3mm}y_m(  \text{tr}(\mathbf{W}^{\text{T}}\mathbf{X}_m)   +b)\geq 1-\xi_m.
\end{aligned}
\end{equation}
\noindent Note that this is a generalization of the standard  SVM, which is obtained for $\tau=0$. The solution to (\ref{smm}) is obtained as

\begin{equation}
\tilde{\mathbf{W}}=D_\tau\left(\sum_{m=1}^{M}\tilde{\beta}_m y_m\mathbf{X}_m\right),
\end{equation}
where $D_\tau(\mathbf{\cdot})$ is the singular value thresholding operator, which suppresses singular values below $\tau$ to zeros. Denote by $\mathbf{\Omega}=\sum_{m=1}^{M}\tilde{\beta}_m y_m\mathbf{X}_m$ a combination of $\mathbf{X}_m$ associated to non-zero ${\tilde{\beta}_m}$, which are  the so called support matrices \citep{luo2015support}.

The solution to (\ref{smm}) can  be obtained by rewriting the problem as
\begin{equation}
\text{arg}\underset{\mathbf{W},b}{\text{min}}\ \frac{1}{2}   \text{tr}(\mathbf{W}^{\text{T}}\mathbf{W}) +\tau \|\mathbf{W}\|_*  +C\sum_{m=1}^{M}   [  1-y_m  (\text{tr}  (\mathbf{W}^{\text{T}}\mathbf{X}_m)+b) ]_+,
\end{equation}
which is equivalent to
\be
\text{arg}\underset{\mathbf{W},b,\mathbf{S}}{\text{min}} \; H(\mathbf{W},b)+G(\mathbf{S}) \quad \text{s.t.} \hspace{4mm} \mathbf{S}-\mathbf{W}= \mathbf 0,
\ee
and
\begin{equation}
\begin{aligned}
H(\mathbf{W},b)=\frac{1}{2}   \text{tr}(\mathbf{W}^{\text{T}}&\mathbf{W})    +C\sum_{m=1}^{M}   [  1-y_m  (\text{tr}  (\mathbf{W}^{\text{T}}\mathbf{X}_m)+b) ]_+ \\
&G(\mathbf{S})=\tau\|\mathbf{S}\|_* ,
\end{aligned}
\end{equation}
where $[x]_+=\max\{x,0\}$.
The solution is obtained using an augmented Lagrangian form
\begin{equation}
L(\mathbf{W},b,\mathbf{S},\mathbf{\Lambda})=H(\mathbf{W},b)   +  G(\mathbf{S})    + \text{tr}(\mathbf{\Lambda}^{\text{T}}(\mathbf{S-W}))+\frac{\rho}{2}\|  \mathbf{S-W}  \|_F^2,
\end{equation}
where $\rho>0$ is a hyperparameter.

\subsection{Support Tensor Machines (STM)}

Consider now a typical problem in computer vision, where the objects are  represented by data tensors and the number of the training measurement is limited. This naturally leads to a tensor extension of SVM, called the \textit{support tensor machine} (STM). Consider a general supervised learning scenario with $M$ training measurements, $\{\tensor{X}_m, y_m\}, m=1,\ldots, M$, represented by $N$th-order tensors, $\tensor{X}_m\in\mathbb{R}^{I_1\times I_2\times\cdots\times I_N}$, which are associated with the scalar variable $y_m$.
There are two possible scenarios: 1) $y_m$ takes a continuous set of values, which leads to the \textit{tensor regression} problem, and 2)  $y_m\in\{+1,-1\}$ that is, it takes categorical values, which is a standard \textit{classification problem}.
For the classification case,  STM \citep{tao2005supervised,biswas2016linear,hao2013linear} can be formulated through
the following minimization problem
\begin{equation}
\begin{split}
&\min_{\mathbf{w}_n, b, \boldsymbol\xi} \; J(\mathbf{w}_n, b, \boldsymbol\xi) = \frac{1}{2}\left\|\otimes_{n=1}^{N}\mathbf{w}_n\right\|^2 + C\sum_{m=1}^M \xi_m\\
\mbox{s.t.} & \quad  y_m\left(\tensor{X}_m \bar\times_1 \mathbf{w}_1 \cdots\bar\times_N\mathbf{w}_N + b\right) \geq 1- \xi_m,
 \quad \boldsymbol\xi\geq 0, \\
 & \quad m=1, \ldots, M.
\end{split}
\end{equation}
Here, $\boldsymbol\xi = [\xi_1,\xi_2,\ldots,\xi_M]^{\text{T}}\in\mathbb{R}^M$ is the vector of all slack variables which helps to deal with linearly nonseparable problems.

The STM problem is therefore composed of $N$ quadratic programming (QP)\-sub-problems with inequality constraints, where the $n$th QP sub-problem can be formulated, as follows
 \citep{hao2013linear}:
\begin{equation}
\label{eq:stm}
\begin{split}
&\min_{\mathbf{w}_n, b, \boldsymbol\xi} \; \frac{1}{2} \; \|\mathbf{w}_n\|^2 \, \prod_{1\leq i\leq N}^{i\neq n} (\|\mathbf{w}_i\|^2)  + C \sum_{m=1}^M \xi_m\\
\mbox{s.t.}  &\quad y_m\left(\mathbf{w}_n^{\text{T}} \left(\tensor{X}_m {\bar\times}_{i\neq n}{\bw_i} \right) + b\right) \geq 1-\xi_m,  \quad \boldsymbol\xi\geq 0, \\
 &  m=1, \ldots, M.
\end{split}
\end{equation}
Based on the least squares SVM  in (\ref{eq:lssvm}), its tensor extension, referred to as the least squares STM (LS-STM), can be formulated as
\begin{equation}
\begin{split}
&\min_{\mathbf{w}_n, b, \boldsymbol\varepsilon}  \; J(\mathbf{w}_n, b, \boldsymbol\varepsilon) = \frac{1}{2}\left\|\bw_1\otimes \bw_2\otimes \cdots \otimes \bw_N \right\|^2 + \frac{\gamma}{2}\boldsymbol\varepsilon^{\text{T}} \boldsymbol\varepsilon\\
\mbox{s.t.} &\quad  y_m\left(\tensor{X}_m \bar \times_1 \mathbf{w}_1 \cdots \bar \times_N\mathbf{w}_N + b\right) = 1- \varepsilon_m,
  \quad m=1, \ldots, M.
\end{split}
\end{equation}
Then, the  LS-STM  solution can be obtained by the following alternating least squares optimization problem (see Algorithm~\ref{alg:LS-STM})
\begin{equation}
\begin{split}
&\min_{\mathbf{w}_n, b, \boldsymbol\varepsilon} \; \frac{1}{2} \; \|\mathbf{w}_n\|^2  \; \prod_{1\leq i\leq N}^{i\neq n} (\|\mathbf{w}_i\|^2 )\; +  \frac{\gamma}{2}\boldsymbol\varepsilon^{\text{T}} \boldsymbol\varepsilon\\
\mbox{s.t.} & \quad y_m\left(\mathbf{w}_n^{\text{T}} \left(\tensor{X}_m {\bar{\times}}_{i\neq n}\bw_i \right) + b\right) = 1- \varepsilon_m,
  \quad m=1, \ldots, M.
\end{split}
\end{equation}
Once the STM solution has been obtained, the class label of a test example, $\tensor{X}_*$, can be predicted by
a nonlinear transformation
\begin{equation}
y(\tensor{X}) = \text{sign}\left(\tensor{X}_*\bar\times_1\mathbf{w}_1 \cdots\bar\times_N\mathbf{w}_N + b\right).
\end{equation}

\begin{algorithm}[t!]
\caption{\textbf{Least Squares Support Tensor Machine (LS-STM)}}
\label{alg:LS-STM}
{\small
\begin{algorithmic}[1]
\REQUIRE $\{\tensor{X}_m, y_m\}, m=1,\ldots,M$ where $\tensor{X}_m\in\mathbb R^{I_1\times I_2\times \cdots\times I_N}$,   \\ $y_m\in\{+1,-1\}$ and the penalty coefficient $\gamma$.
\ENSURE $\{\mathbf w_n\}, n=1,\ldots, N$ and $b$.
\STATE Initialization: Random vectors $\{\mathbf w_n\}, n=1,\ldots, N$.
\WHILE {not converged or iteration limit is not reached}
\FOR {$n=1$ to $N$}
\STATE $\eta \leftarrow \prod_{1\leq i\leq N}^{i\neq n} \|\mathbf w_i\|^2$,\\ \vspace*{+3pt}
\STATE $\mathbf x_m = \tensor{X}_m {\bar\times}_{i \neq n} \bw_i$
\STATE Compute  $\mathbf w_n$  by solving the optimization problem
\be
&&\min_{\mathbf{w}_n, b, \boldsymbol\varepsilon} \, \frac{\eta}{2}  \|\mathbf{w}_n\|^2 +  \frac{\gamma}{2}\boldsymbol\varepsilon^{\text{T}} \boldsymbol\varepsilon \notag\\
\mbox{s.t.} && \quad y_m \left( \mathbf w_n^{\text{T}} \mathbf x_m + b\right) = 1- \varepsilon_m, \notag \\
 &&  \quad m=1, \ldots, M.\notag
\ee
\ENDFOR
\ENDWHILE
\end{algorithmic}
}
\end{algorithm}

In practice, it is often more convenient to solve the  optimization problem   (\ref{eq:stm}) by considering the \textit{dual problem}, given by
\begin{equation}
\label{eq:dualstm}
\begin{split}
&\max_{\{\alpha_i\}_{i=1}^{M}} \sum_{m=1}^M \alpha_m -\frac{1}{2}\sum_{i,j=1}^M \alpha_i\alpha_j y_i y_j\langle\tensor{X}_i,\tensor{X}_j\rangle,\\
\mbox{s.t.} & \quad \sum_{m=1}^M \alpha_m y_m =0, \quad 0\leq\alpha_m\leq C, \quad m=1,\ldots, M,
\end{split}
\end{equation}
where $\alpha_m$ are the Lagrange multipliers and $\langle\tensor{X}_i,\tensor{X}_j\rangle$  the inner products of $\tensor{X}_i$ and $\tensor{X}_j$.

It is obvious that when the input samples, $\tensor{X}_i$, are vectors, this optimization model simplifies into the standard vector SVM. Moreover, if the original input tensors are used to compute $\langle\tensor{X}_i, \tensor{X}_j\rangle$, then the optimal solutions are the same as those produced by the SVM. Consider now the rank-one decompositions of $\tensor{X}_i$ and $\tensor{X}_j$ in the form $\tensor{X}_i\approx \sum_{r=1}^R \mathbf{x}_{ir}^{(1)}\circ \mathbf{x}_{ir}^{(2)}\circ\cdots\circ\mathbf{x}_{ir}^{(N)}$ and $\tensor{X}_j\approx \sum_{r=1}^R \mathbf{x}_{jr}^{(1)}\circ \mathbf{x}_{jr}^{(2)}\circ\cdots\circ\mathbf{x}_{jr}^{(N)}$. Then, the inner product of $\tensor{X}_i$ and $\tensor{X}_j$ is given by \citet{hao2013linear}
\begin{equation}
\langle \tensor{X}_i, \tensor{X}_j\rangle \approx \sum_{p=1}^R\sum_{q=1}^R \left\langle\mathbf{x}_{ip}^{(1)},\mathbf{x}_{jq}^{(1)}\right\rangle \left\langle\mathbf{x}_{ip}^{(2)},\mathbf{x}_{jq}^{(2)}\right\rangle\cdots \left\langle\mathbf{x}_{ip}^{(N)},\mathbf{x}_{jq}^{(N)}\right\rangle,
\end{equation}
and (\ref{eq:dualstm}) can be solved by a sequential QP optimization algorithm. The class label of a test example $\tensor{X}_*$ is then predicted as
\begin{equation}
y(\tensor{X}_*) = \text{sign}\left(\sum_{m=1}^M\sum_{p=1}^R\sum_{q=1}^R \alpha_m y_m \prod_{n=1}^N \left\langle \mathbf{x}_{mp}^{(n)},\mathbf{x}_{*q}^{(n)} \right\rangle +b    \right).
\end{equation}

\section{Higher Rank Support Tensor Machines (HRSTM)}
\sectionmark{Higher Rank Support Tensor Machines}

Higher Rank STMs  (HRSTM) aim  to estimate a set of parameters in the form of  the sum of rank-one tensors \citep{kotsia2012higher}, which  defines a separating hyperplane between classes of data. The benefits of this scheme are twofold:
\begin{enumerate}
  \item The use of a direct CP  representation is intuitively closer to the idea of properly processing tensorial input data, as the data structure is more efficiently retained;
  \item The use of  simple CP decompositions allows for multiple projections of the input tensor along each mode, leading to considerable improvements in the discriminative ability of the resulting classifier.
\end{enumerate}

The corresponding optimization problems can be  solved in an iterative manner utilizing, e.g., the standard CP decomposition, where at each iteration the parameters corresponding to the projections along a single tensor mode are estimated by solving a typical STM-type optimization problem.

The aim of the  so formulated STM is therefore to learn a multilinear decision function, $g: \mathbf{R}^{I_1\times I_2\times\cdots\times I_N} \rightarrow [-1,1]$, which  classifies a test tensor $\tensor{X}\in\mathbf{R}^{I_1\times I_2\times\cdots\times I_N}$, by using the nonlinear transformation $g(\tensor{X}) = \text{sign}(\langle \tensor{X}, \tensor{W}\rangle + b)$.

The weight tensor, $\tensor{W}$, is estimated by solving the following soft STMs optimization problem
\begin{equation}
\label{eq:hrstm1}
\begin{split}
& \min_{\tensor{W}, b, \boldsymbol\xi} \; \frac{1}{2}\langle \tensor{W}, \tensor{W}\rangle + C \sum_{m=1}^M \xi_m,\\
\mbox{s.t.} \quad & y_m \left(\langle\tensor{W}, \tensor{X}_m\rangle + b \right)   \geq 1- \xi_m, \quad \xi_m \geq 0,  \quad m=1, \ldots, M,
\end{split}
\end{equation}
where $b$ is the bias term, $\boldsymbol \xi = [\xi_1,\ldots, \xi_m]^{\text{T}}$ is the vector of slack variables and $C$ is the term which controls the relative importance of penalizing the training errors.  The training is performed in such a way that the margin of the support tensors is maximized while the upper bound on the misclassification errors in the training set is minimized.

Consider, as an example,  the case where the weight parameter, $\tensor{W}$, is  represented via a CP model, to give
\begin{equation}
\tensor{W} = \sum_{r=1}^R \mathbf{u}_r^{(1)}\circ \mathbf{u}_r^{(2)}\circ\cdots\circ\mathbf{u}_r^{(N)},
\end{equation}
where $\mathbf{u}_r^{(n)} \in\mathbb{R}^{I_n}, \;\; n =1,2,\ldots, N$. Then, the $n$-th mode matricization of $\tensor{W}$ can be written as
\begin{equation}
\mathbf{W}_{(n)} = \mathbf{U}^{(n)} (\mathbf{U}^{(N)} \odot \cdots \odot \mathbf{U}^{(n+1)}\odot \mathbf{U}^{(n-1)}\odot\cdots\odot \mathbf{U}^{(1)}) ^{\text{T}} = \mathbf{U}^{(n)}(\mathbf{U}^{(-n)})^{\text{T}},
\end{equation}
where $\odot$ denotes the Khatri-Rao product. Note that the inner product in (\ref{eq:hrstm1}) can be computed very efficiently e.g., in the form
\begin{equation}
\langle \tensor{W}, \tensor{W} \rangle = \tr[\mathbf{W}_{(n)} \mathbf{W}_{(n)}^{\text{T}}] = \text{vec}(\mathbf{W}_{(n)})^{\text{T}} \text{vec}(\mathbf{W}_{(n)}).
\end{equation}

The above optimization problem, however, is not convex with respect to all sets of parameters in $\tensor{W}$. For this reason, we  adopt an iterative scheme in which at each iteration we solve only for the parameters that are associated with the $n$-th mode of the parameter tensor $\tensor{W}$, while keeping all the other parameters fixed, similarly to the ALS algorithm. More specifically,  for the $n$-th mode, at each iteration we  solve the following optimization problem
\begin{equation}
\begin{split}
& \min_{\mathbf{W}_{(n)}, b, \boldsymbol\xi} \, \frac{1}{2} \tr \left(\mathbf{W}_{(n)}\mathbf{W}_{(n)}^{\text{T}}\right) + C\sum_{m=1}^M \xi_m,\\
\mbox{s.t.} \quad & y_m\left(\tr\left(\mathbf{W}_{(n)}\mathbf{X}_{m(n)}^{\text{T}}\right) + b\right) \geq 1-\xi_m, \quad \xi_m \geq 0, \\
& \quad m=1, \ldots, M, \quad n=1,\ldots,N,
\end{split}
\end{equation}

\noindent Under the assumption that the tensor $\tensor{W}$ is represented as a sum of rank one tensors, we can replace the above matrices by $\mathbf{W}_{(n)} = \mathbf{U}^{(n)} (\mathbf{U}^{(-n)})^{\text{T}}$, and thus the above equation becomes
\begin{equation}
\label{eq:hostm}
\begin{split}
 & \min_{\mathbf{U}^{(n)},b,\boldsymbol\xi} \, \frac{1}{2}\tr(\mathbf{U}^{(n)} (\mathbf{U}^{(-n)})^{\text{T}} (\mathbf{U}^{(-n)}) (\mathbf{U}^{(n)})^{\text{T}}) + C\sum_{m=1}^M \xi_m, \\
\mbox{s.t.} \quad & y_m \left(\tr\left(\mathbf{U}^{(n)}(\mathbf{U}^{(-n)})^{\text{T}} \mathbf{X}^{\text{T}}_{m(n)} \right)+ b \right) \geq 1- \xi_m, \quad \xi_m \geq 0, \\
& \quad m=1, \ldots, M, \quad n=1,\ldots,N.
\end{split}
\end{equation}

\noindent  At each iteration  the optimization problem is solved for  only one set of matrices $\mathbf{U}^{(n)}$ for the mode $n$,  while keeping the other matrices $\mathbf{U}^{(k)}$ for $k\neq n$, fixed.
Note that   the optimization problem defined in (\ref{eq:hostm}) can also be solved using a classic vector-based SVM implementation.

 Let us define $\mathbf{A} = \mathbf{U}^{(-n)^{\text{T}}} \mathbf{U}^{(-n)}$, which is a positive definite matrix, and  let $\widetilde{\mathbf{U}}^{(n)} = \mathbf{U}^{(n)}\mathbf{A}^{\frac{1}{2}}$. Then,
\begin{equation}
\begin{split}
\tr[\mathbf{U}^{(n)} (\mathbf{U}^{(-n)})^{\text{T}} (\mathbf{U}^{(-n)}) (\mathbf{U}^{(n)})^{\text{T}}] &= \tr[\widetilde{\mathbf{U}}^{(n)} (\widetilde{\mathbf{U}}^{(n)})^{\text{T}}  ] \\
&= \text{vec}(\widetilde{\mathbf{U}}^{(n)})^{\text{T}} \text{vec}(\widetilde{\mathbf{U}}^{(n)}).
\end{split}
\end{equation}
By letting $\widetilde{\mathbf{X}}_{m(n)} = \mathbf{X}_{m(n)} \mathbf{U}^{(-n)}\mathbf{A}^{-\frac{1}{2}}$, we have
\begin{equation}
\begin{split}
\tr[\mathbf{U}^{(n)} (\mathbf{U}^{(-n)})^{\text{T}} \mathbf{X}_{m(n)}^{\text{T}}] &= \tr[\widetilde{\mathbf{U}}^{(n)} \widetilde{\mathbf{X}}^{\text{T}}_{m(n)}] \\
&= \text{vec}(\widetilde{\mathbf{U}}^{(n)})^{\text{T}} \text{vec} (\widetilde{\mathbf{X}}_{m(n)}).
\end{split}
\end{equation}
Then, the optimization problem (\ref{eq:hostm}) can be simplified as
\begin{equation}
\label{eq:hostm2}
\begin{split}
& \min_{\mathbf{U}^{(n)}, b, \boldsymbol \xi} \, \frac{1}{2} \text{vec}(\widetilde{\mathbf{U}}^{(n)})^{\text{T}} \text{vec}(\widetilde{\mathbf{U}}^{(n)}) + C \sum_{m=1}^M \xi_m,\\
\mbox{s.t.} \quad & y_m \left(\text{vec}(\widetilde{\mathbf{U}}^{(n)})^{\text{T}} \text{vec}(\widetilde{\mathbf{X}}_{m(n)}) + b \right) \geq 1-\xi_m, \quad \xi_m \geq 0,  \\
& \quad m=1, \ldots, M, \quad n=1,\ldots,N.
\end{split}
\end{equation}
Observe that now the STM optimization problem for the $n$-th mode in (\ref{eq:hostm}) can be formulated as a standard vector SVM problem with respect to $\widetilde{\mathbf{U}}^{(n)}$.
In other words, such a procedure  leads to a straightforward implementation of the algorithm by solving (\ref{eq:hostm2}) with respect to $\widetilde{\mathbf{U}}^{(n)}$  using a standard SVM implementation, and then solving for $\mathbf{U}^{(n)}$ as $\mathbf{U}^{(n)} = \widetilde{\mathbf{U}}^{(n)} \mathbf{A}^{-\frac{1}{2}}$.

\section{Kernel Support Tensor Machines}

The class of support tensor machines has been recently  extended to the nonlinear case by using the kernel framework, and is referred to as  \textit{kernel support tensor regression} (KSTR) \citep{gao2012kernel}. After mapping  each row of every original tensor (or every tensor converted from original vector) onto a high-dimensional space, we obtain the associated points in a new high-dimensional feature space, and then compute the regression function. Suppose we are given a set of training samples $\{\tensor{X}_m, y_m\}, m=1,\ldots,M$, where each  training sample, $\tensor{X}_m$, is a data point in $\mathbb{R}^{I_1} \otimes \mathbb{R}^{I_2}$, where $\mathbb{R}^{I_1}$ and $\mathbb{R}^{I_2}$ are two vector spaces, and $y_m$ is the target value associated with $\tensor{X}_m$. Denote by $\bz_{mp}$  the $p$-th row of $\tensor{X}_m$,  we can then use a nonlinear mapping function $\varphi(\tensor{X}_m)$ to map $\tensor{X}_m$ onto a high-dimensional tensor feature space, and define a nonlinear mapping function for tensor $\tensor{X}_m$, given by
\begin{equation}
\phi (\tensor{X}_m) = \left[\begin{array}{c}
                        \varphi(\bz_{m1}) \\
                        \varphi(\bz_{m2})  \\
                        \vdots\\
                        \varphi(\bz_{mI_1})
                      \end{array}\right].
\end{equation}
In this way, we can obtain the new kernel function
\begin{equation}
\begin{split}
K(\tensor{X}_m,\tensor{X}_n) &= \phi (\tensor{X}_m) \phi (\tensor{X}_n)^{\text{T}} =
                     \left[\begin{array}{c}
                        \varphi(\bz_{m1}) \\
                        \varphi(\bz_{m2})  \\
                        \vdots\\
                        \varphi(\bz_{mI_1})
                      \end{array}\right]
                       \left[\begin{array}{c}
                        \varphi(\bz_{n1}) \\
                        \varphi(\bz_{n2})  \\
                        \vdots\\
                        \varphi(\bz_{nI_1})
                      \end{array}\right]^{\text{T}} \\
                      &  =
                      \left[
                        \begin{array}{ccc}
                          \varphi(\bz_{m1})\varphi(\bz_{n1})^{\text{T}} & \cdots & \varphi(\bz_{m1})\varphi(\bz_{nI_1})^{\text{T}} \\
                          \vdots & \ddots & \vdots \\
                          \varphi(\bz_{mI_1})\varphi(\bz_{n1})^{\text{T}} & \cdots & \varphi(\bz_{mI_1})\varphi(\bz_{nI_1})^{\text{T}} \\
                        \end{array}
                      \right].
\end{split}
\end{equation}
Note that such a kernel function is quite different from the function used in the SVR - the output of this kernel function is a matrix as opposed to a scalar in SVR.

For instance,  if we use an RBF kernel function within KSTR, then the $mn$-th element of the kernel matrix is expressed as
\begin{equation}
\varphi(\bz_{mp_1})\varphi(\bz_{np_2})^{\text{T}} = e^{-\gamma \|\bz_{mp_1} - \bz_{np_2} \|^2}.
\end{equation}
The support tensor regression with an $\varepsilon$-insensitive loss function is similar to standard support tensor regression, and the  regression function in  such a case can  be defined as
\begin{equation}
f(\tensor{X}) = \mathbf{u}^{\text{T}} \phi(\tensor{X})\mathbf{v} + b.
\end{equation}
This function can be estimated  by solving the following  quadratic programming problem
\begin{equation}
\label{eq:KSTM}
\begin{split}
& \min_{\mathbf{u,v}, b, \xi_m, \xi_m^*} \, \frac{1}{2}\left\|\mathbf{uv}^{\text{T}}\right\|^2 + C\sum_{m=1}^M (\xi_m+\xi_m^*)\\
\mbox{s.t.} \quad & \left\{ \begin{array}{ll}
             y_m - \mathbf{u}^{\text{T}}\phi(\tensor{X}_m)\mathbf{v}-b & \leq \varepsilon - \xi_m, \\
             \mathbf{u}^{\text{T}}\phi(\tensor{X}_m)\mathbf{v} + b - y_m  &\leq \varepsilon + \xi_m^*, \quad m =1, \ldots, M,\\
             \xi_m, \xi_m^* &\geq 0,
           \end{array}
           \right.
\end{split}
\end{equation}
where $C$ is a pre-specified constant, $\varepsilon$  is a user-defined scalar, and $\xi_m, \xi_m^*$ are slack variables which represent the upper and lower constraints on the outputs of the classification system.

Consider again (\ref{eq:KSTM}), and let $\mathbf{u}$ be a column vector of the same dimension as the row number of samples; we can then calculate the vector $\mathbf{v}$ using the Lagrange multiplier method, where the Lagrangian is constructed according to
\begin{equation}
\begin{split}
\max_{\alpha_m, \alpha_m^*, \eta_m, \eta_m^*} \min_{\mathbf{v},b,\xi_m,\xi_m^*} &
\left\{ \begin{array}{l}
L= \frac{1}{2}\|\mathbf{u}\mathbf{v}^{\text{T}}\|^2 + C\sum_{m=1}^M(\xi_m+\xi_m^*) \\-\sum_{m=1}^M(\eta_m\xi_m + \eta_m^*\xi_m^*) - \\
\sum_{m=1}^M\alpha_m(\varepsilon + \xi_m - y_m + \mathbf{u}^{\text{T}}\phi(\tensor{X}_m)\mathbf{v} + b)\\ - \sum_{m=1}^M \alpha_m^*(\varepsilon + \xi_m^* - y_m + \mathbf{u}^{\text{T}}\phi(\tensor{X}_m)\mathbf{v}+b)
\end{array} \right\}\\
\mbox{s.t.} \quad & \alpha_m, \alpha_m^*,\eta_m, \eta_m^* \geq 0, \quad m=1,\ldots,M,
\end{split}
\end{equation}
and $\alpha_m, \alpha_m^*, \eta_m,\eta_m^*$ are the Lagrange multipliers.
In next step, we  determine the Lagrange multipliers $\alpha_m, \alpha_m^*$ and $\mathbf{v}$, $\|\mathbf{v}\|^2$.

Alternatively, let $\mathbf{x}_n' = \mathbf{v}^{\text{T}}\phi(\tensor{X}_n) = \frac{1}{\|\mathbf{u}\|^4}\sum_{m=1}^M(\alpha_m - \alpha_m^*) K(\tensor{X}_m, \tensor{X}_n)$ be the new training samples. Then, we  construct another Lagrangian
\begin{equation}
\begin{split}
\max_{\alpha_m, \alpha_m^*, \eta_m, \eta_m^*} \min_{\mathbf{u},b,\xi_m,\xi_m^*}
& \left\{\begin{array}{l}
         L=\frac{1}{2}\|\mathbf{uv}^{\text{T}}\|^2 + C\sum_{m=1}^M (\xi_m + \xi_m^*) \\- \sum_{m=1}^M(\eta_m\xi_m + \eta_m^*\xi_m^*) -\\
         \sum_{m=1}^M \alpha_m (\varepsilon + \xi_m - y_m + \mathbf{x}'_m \mathbf{u} + b) \\- \sum_{m=1}^M (\alpha_m^* (\varepsilon + \xi_m^* + y_m - \mathbf{x}'_m \mathbf{u} - b))
       \end{array}
\right\}\\
\mbox{s.t.} \quad & \alpha_m, \alpha_m^*, \eta_m, \eta_m^* \geq 0, \quad m=1,\ldots,M,
\end{split}
\end{equation}
and obtain $\mathbf{u}$ and $b$. These two steps are performed iteratively to compute $\mathbf{v, u},b$.

The KSTR method has the following  advantages over the STR method: i) it has a strong ability to learn and a superior generalization ability; ii)  the  KSTR  is able to solve nonlinearly separable problems more efficiently. A disadvantage of KSTR is that the computational load of KSTR is much higher than  that of STR.

\section{Tensor Fisher Discriminant Analysis (FDA)}
\sectionmark{Tensor Fisher Discriminant Analysis}

Fisher discriminant analysis (FDA) has been widely applied for classification. It aims to find  a direction which separates the class means well while minimizing the variance of the total training measurements. This procedure is equivalent to maximizing the symmetric Kullback-Leibler divergence (KLD) between positive and negative measurements with identical covariances, so that the positive measurements are separated from the negative measurements. To extend the FDA to tensors, that is, to introduce Tensor FDA (TFDA), consider $M$ training measurements, $\tensor{X}_m\in\mathbb{R}^{I_1\times I_2\times \cdots\times I_N},  \;\; m=1, \ldots, M$, associated with the class labels $y_m\in\{+1,-1\}$. The mean of  positive training measurements is $\tensor{L}_{+} = (1/M_+)\sum_{m=1}^M ({I}(y_m=+1)\tensor{X}_m)$, the mean of the  negative training measurements is  $\tensor{L}_{-} = (1/M_-)\sum_{m=1}^M ({I}(y_m=-1)\tensor{X}_m)$, the mean of all training measurements is $\tensor{L} = (1/M)\sum_{m=1}^M \tensor{X}_m$, and $M_+(M_-)$ are the numbers of the positive (negative) measurements. The decision function is a multilinear function $y(\tensor{X}) = \text{sign}(\tensor{X}\bar\times_1 \mathbf{w}_1 \cdots \bar\times_N \mathbf{w}_N + b)$, while  the projection vectors, $\mathbf{w}_n\in\mathbb{R}^{I_n}, \;  n=1, \ldots, N$, and the bias $b$ in TFDA are obtained from
\begin{equation}
\max_{\mathbf{w}_n|_{n=1}^N} \, J(\mathbf{w}_n) = \frac{\|(\tensor{L}_+ -\tensor{L}_-)\bar\times_1 \mathbf{w}_1 \cdots \bar\times \mathbf{w}_N \|^2}{\sum_{m=1}^M \|(\tensor{X}_m-\tensor{L})\bar\times_1 \mathbf{w}_1 \cdots\bar\times \mathbf{w}_N\|^2}.
\end{equation}
Unfortunately, there is no closed-form solution for TFDA, however,  alternating projection algorithms can be applied  to iteratively obtain the desired  solution using the TFDA framework.
For the projection vectors $\mathbf{w}_n|_{n=1}^N$, we can obtain the bias $b$ using
\begin{equation}
b = \frac{M_- -M_+ -(M_+ \tensor{L}_+ + M_-\tensor{L}_-)\bar\times_1\mathbf{w}_1 \cdots\bar\times_N\mathbf{w}_N}{M_- + M_+}.
\end{equation}

The regularized Multiway Fisher Discriminant Analysis (MFDA) \citep{lechuga2015discriminant} aims to impose the structural constraints in such a way that the weight vector $\mathbf w$ will be decomposed  as $\mathbf w = \mathbf w_N \otimes \cdots\otimes \mathbf w_1$. Let $\mathbf X_{(1)} \in \mathbb R^{M\times I_1I_2\cdots I_N}$ be an unfolded matrix of observed samples and  $\mathbf Y\in \mathbb R^{M\times C}$  the matrix of dummy variables indicating the group memberships. The optimization problem of MFDA can then be formulated as
\begin{equation}
\begin{split}
\arg & \max_{\mathbf w} \frac{\mathbf w^{\text{T}} \mathbf S_B \mathbf w}{\mathbf w^{\text{T}} \mathbf S_T \mathbf w + \lambda \mathbf w^{\text{T}} \mathbf R \mathbf w}\\
\end{split}
\end{equation}
where  $\mathbf w$ denotes $\mathbf w_N \otimes \cdots\otimes \mathbf w_1$, $\mathbf S_B = \mathbf X_{(1)}^{\text{T}}\mathbf Y (\mathbf Y^{\text{T}} \mathbf Y)^{-1}\mathbf Y^{\text{T}}\mathbf X_{(1)}$ is the between covariance matrix, $\mathbf S_T =\mathbf X_{(1)}^{\text{T}}\mathbf X_{(1)}$ is the total covariance matrix, and $\mathbf R$ is usually an identity matrix. The regularization term $\lambda \mathbf w^{\text{T}}\mathbf R\mathbf w$ is added to improve the numerical stability when computing the inverse of $\mathbf S_T$ in high-dimensional settings $(M\ll I_1I_2\cdots I_N)$. This optimization problem can be solved in an alternate fashion by fixing all $\mathbf w_n$ except one.

The Linear Discriminant Analysis (LDA) can also be extended to tensor data, which is referred to as Higher Order Discriminant Analysis (HODA) \citep{phan2010tensor} and  Multilinear Discriminant Analysis (MDA) \citep{li2014multilinear}. Suppose that we have $M$ tensor samples $\tensor X_m\in\mathbb R^{I_1\times I_2\times\cdots\times I_N}, m=1,\ldots, M$ belonging to one of the $C$ classes, and $M_c$ is the number of samples in class $c$ such that $M=\sum_{c=1}^C M_c$. $y_m$ denotes the associated class label of $\tensor X_m$. The class mean tensor for class $c$ is computed by $\tensor L_c = \frac{1}{M_c}\sum_{m=1}^{M_c} \tensor X_m I(y_m=c)$, and the total mean tensor is $\tensor L = \frac{1}{M} \sum_{m=1}^M \tensor X_m$. The goal is to find the set of optimal projection matrices $\mathbf W_1, \mathbf W_2,\ldots,\mathbf W_N$ ($\mathbf W_n\in\mathbb R^{I_n\times J_n}$ for $n=1,\ldots,N$) that lead to the most accurate classification in the projected tensor subspace where
\begin{equation}
\tensor Z_m = \tensor X_m \times_1 \mathbf W_1^{\text{T}} \cdots \times_N \mathbf W_N^{\text{T}} \in\mathbb R^{J_1\times J_2\times \cdots \times J_N}.
\end{equation}
The mean tensor of class $c$ in the projected subspace is given by
\begin{equation}
\widetilde{\tensor{L}}_c = \frac{1}{M_c}\sum_{m=1}^{M_c} \tensor Z_m I(y_m=c) = \tensor L_c \times_1 \mathbf W_1^{\text{T}}\cdots \times_N \mathbf W_N^{\text{T}},
\end{equation}
which is simply the projection of $\tensor L_c$. Similarly, the total mean tensor of projected samples is $\tensor L\times_1 \mathbf W_1^{\text{T}}\cdots \times_N \mathbf W_N^{\text{T}}$.

The mode-$n$ between-class scatter matrix in the projected tensor subspace can be derived,
 as follows
\begin{equation}
\begin{split}
\mathbf B_n = & \sum_{c=1}^C M_c\left[ \left( \tensor L_c\prod_{n=1}^N \times_n \mathbf W_n^{\text{T}}\right)_{(n)} -   \left( \tensor L\prod_{n=1}^N \times_n \mathbf W_n^{\text{T}}\right)_{(n)}       \right]\\
& \left[ \left( \tensor L_c\prod_{n=1}^N \times_n \mathbf W_n^{\text{T}}\right)_{(n)} -   \left( \tensor L\prod_{n=1}^N \times_n \mathbf W_n^{\text{T}}\right)_{(n)}       \right]^{\text{T}}\\
= & \mathbf W_n^{\text{T}}\left \{ \sum_{c=1}^C M_c \left[ (\tensor L_c - \tensor L )\prod_{k=1,k\neq n}^N \times_k \mathbf W_n^{\text{T}}\right]_{(n)}  \right.\\
& \left. \left[ (\tensor L_c - \tensor L )\prod_{k=1,k\neq n}^N \times_k \mathbf W_n^{\text{T}}\right]_{(n)}^{\text{T}}  \right\}\mathbf W_n \\
=& \mathbf W_n^{\text{T}} \mathbf B_{-n}\mathbf W_n.
\end{split}
\end{equation}
Note that $\tensor L \prod_{n=1}^N \times_n \bW_n^{\text{T}} = \tensor L\times_1 \bW_1^{\text{T}}\times_2 \bW_2^{\text{T}}\cdots\times_N \bW_N^{\text{T}}$.
Here, $ \mathbf B_{-n}$ denotes the mode-$n$ between-class scatter matrix in the partially projected tensor subspace (by all tensor modes except for mode $n$). The mode-$n$ between class scatter matrix characterizes the separation between $C$ classes in terms of mode-$n$ unfolding of the tensor samples.

Similarly, the mode-$n$ within-class scatter matrix is
\begin{equation}
\begin{split}
\mathbf S_n = & \mathbf W_n^{\text{T}} \left\{ \sum_{c=1}^C\sum_{m=1}^{M_c} \left[ (\tensor X_m^c - \tensor L_c) \prod_{k=1,k\neq n}^N \times_k \mathbf W_k^{\text{T}} \right]_{(n)} \right. \\
& \left.  \left[ (\tensor X_m^c - \tensor L_c) \prod_{k=1,k\neq n}^N \times_k \mathbf W_k^{\text{T}} \right]_{(n)}^{\text{T}}                      \right\}\mathbf W_n\\
= & \mathbf W_n^{\text{T}} \mathbf S_{-n} \mathbf W_n,
\end{split}
\end{equation}
where $\mathbf S_{-n}$  represents the mode-$n$ within-class scatter matrix in the partially projected tensor subspace (in all tensor modes except for $n$).
Finally, the MDA objective function can be described as
\begin{equation}
\begin{split}
J(\mathbf W_n) =& \frac{\sum_{c=1}^C M_c \| (\tensor L_c - \tensor L)\prod_{n=1}^N\times_n \mathbf W_n^{\text{T}} \|_F^2}{\sum_{c=1}^C\sum_{m=1}^{M_c}  \| (\tensor X_m^c - L_c)\prod_{n=1}^N \times_n \mathbf W_n^{\text{T}}\|_F^2 }\\
= & \frac{\text{Tr}(\mathbf W_n^{\text{T}} \mathbf B_{-n} \mathbf W_n)}{\text{Tr}(\mathbf W_n^{\text{T}} \mathbf S_{-n} \mathbf W_n)},
\end{split}
\end{equation}
and the set of optimal projection matrices should maximize $J(\mathbf W_n)$ for $n=1,\ldots, N$ simultaneously, to best preserve the given class structure.

In this section, we have discussed  SVM and STM for  medium and large scale problems.
For big data classification problems  quite perspective and promising  is a quantum computing approach
\citep{Lloyd14SVM,schuld2015introduction,Li2015experimental,chatterjee2016generalized,biamonte2016quantum} and a
tensor network approach discussed in the  next Chapter \ref{Chapter3}.


\chapter{Tensor Train Networks for Selected Huge-Scale Optimization Problems}
\chaptermark{TT Networks for Optimization Problems}
\label{Chapter3}

\vspace{-0.5cm}

For extremely large-scale, multidimensional datasets, due to curse of dimensionality, most standard numerical methods for computation and optimization problems are intractable.
This chapter introduces feasible solutions for several generic huge-scale dimensionality reduction and related optimization problems, whereby the involved
optimized cost functions are approximated by suitable low-rank TT networks. In this way, a very large-scale optimization problem
can be converted into a set of much smaller optimization sub-problems of the same kind
\citep{schollwock11-DMRG,Schollwock13,Holtz-TT-2012,dolgovEIG2013,KressnerEIG2014,Cichocki2014optim,Lee-TTPI},
which can be solved using standard methods.

The related  optimization problems  often involve structured matrices and vectors with  over a billion entries
 (see \citep{Grasedyck-rev,Dolgovth,Garreis2016,Stoudenmire_optim,Hubig16} and references therein).
In particular, we focus on Symmetric Eigenvalue Decomposition (EVD/PCA) and Generalized Eigenvalue Decomposition (GEVD) \citep{dolgovEIG2013,KressnerEIG2014,Hubig15,KressnerEVD2016,Huckle2012Lanczos,zhang2016subspace},
 SVD  \citep{Lee-SIMAX-SVD},
solutions of overdetermined and undetermined  systems of linear algebraic equations  \citep{Oseledets-Dolgov-lin-syst12,DolgovAMEN2014},
 the Moore--Penrose pseudo-inverse of structured matrices \citep{Lee-TTPI}, and LASSO regression problems \citep{Lee-Lasso}.
 Tensor networks for extremely large-scale multi-block (multi-view) data are also discussed, especially TN models for orthogonal Canonical Correlation
  Analysis (CCA) and related Higher-Order Partial Least Squares (HOPLS) problems  \citep{NIPSQibin,Qibin-HOPLS,hou_qibin2016aaai,hou2016phd}.
  For convenience, all these problems are reformulated as constrained optimization problems which are then, by virtue  of low-rank tensor
  networks, reduced to manageable lower-scale optimization sub-problems. The enhanced tractability and scalability is achieved through tensor network contractions and other tensor network transformations.
  \markright{3.1.\quad TT Splitting}
%

Prior to introducing solutions to several fundamental optimization problems
for very large-scale data, we shall
describe basic strategies for optimization with cost functions in TT formats.


\section{Tensor Train (TT/MPS) Splitting and Extraction of Cores}
\sectionmark{TT Splitting}
\label{sect:orthog}

\subsection{Extraction of a Single Core and a Single Slice for ALS Algorithms}
\subsectionmark{Extraction of a Single Core}

\begin{figure}[t!]
(a)
\begin{center}
\includegraphics[width=8.99cm]{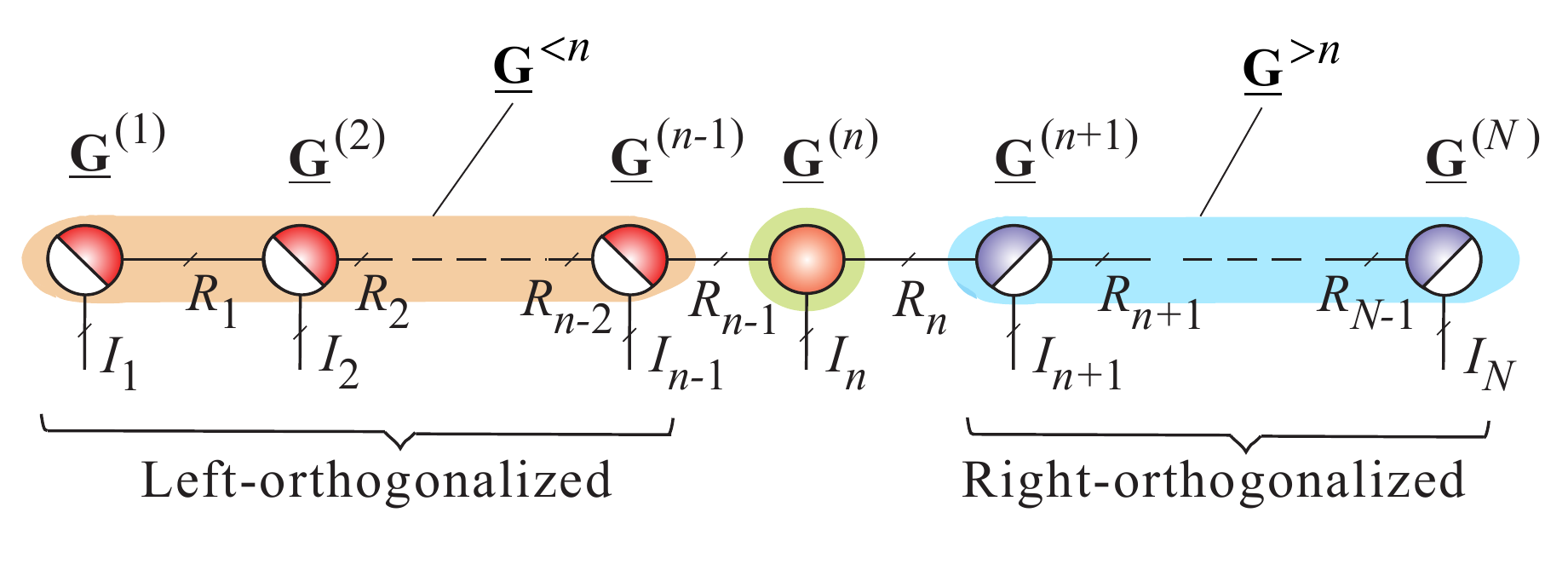}\\
\end{center}
(b)
\begin{center}
\includegraphics[width=10.5cm]{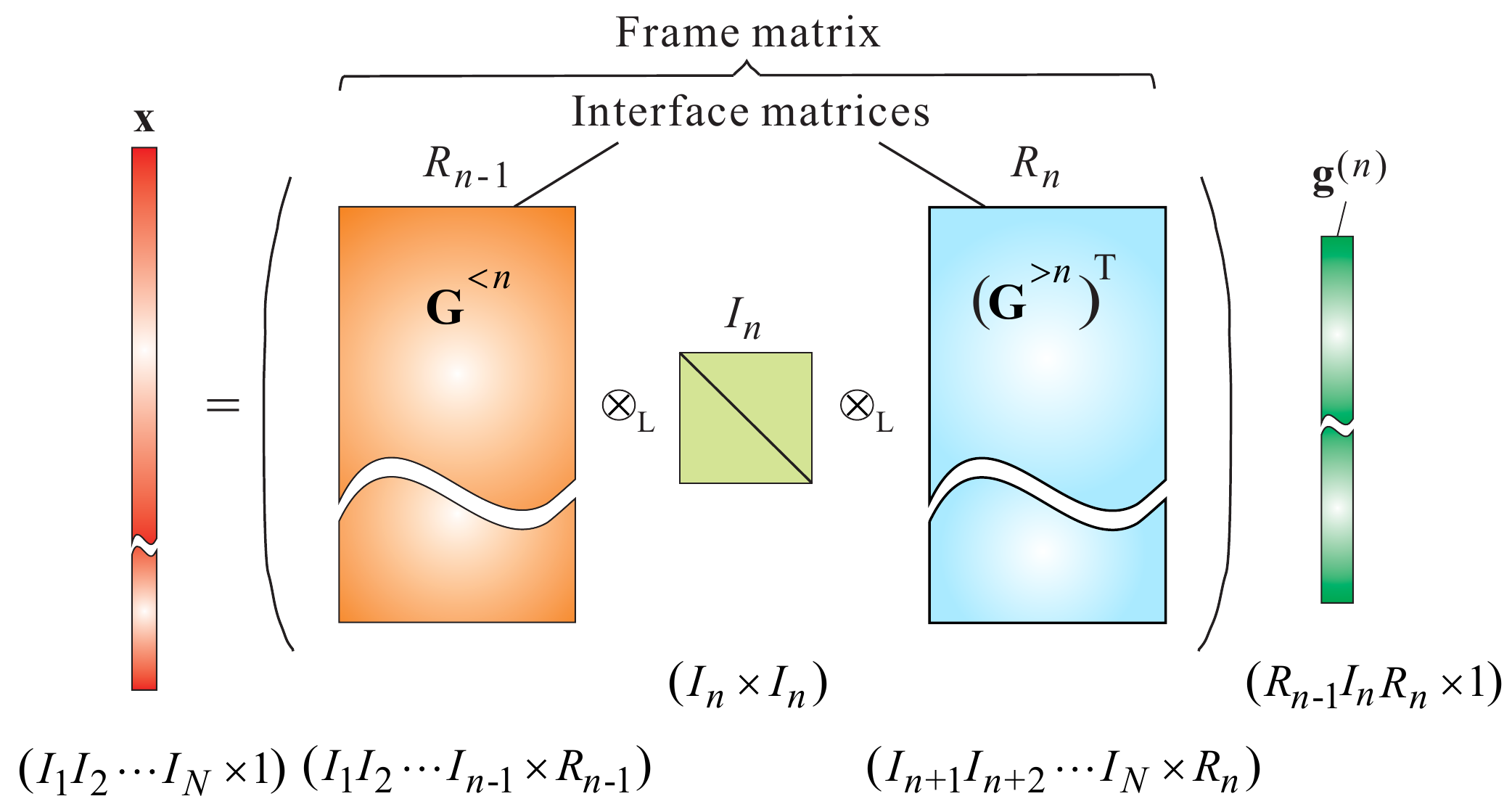}
\end{center}
\caption{Extraction of a single core from a TT network. (a) Representation of a tensor train by the left- and right-orthogonalized sub-trains, with respect to the $n$th core. (b) Graphical illustration of a linear frame equation expressed via a frame matrix or, equivalently, via interface matrices (see Eqs. (\ref{frame-eq}) and (\ref{frame-matrices})).}
\label{Fig:TTsplitt}
\end{figure}

For an efficient implementation of  ALS optimization algorithms, it is convenient to first divide a TT network which represents a tensor,\linebreak $\underline \bX = \llangle \underline \bG^{(1)}, \underline \bG^{(2)}, \ldots, \underline  \bG^{(N)} \rrangle \in \Real^{I_1 \times I_2 \times \cdots \times I_N}$, into sub-trains, as illustrated in Figure~\ref{Fig:TTsplitt}(a). In this way, a large-scale task is replaced by easier-to-handle sub-tasks,
whereby the aim is to extract a specific TT core or its  slices from the whole TT network.
For this purpose, the TT sub-trains can be defined  as follows \citep{Holtz-TT-2012,KressnerEIG2014}
\be
\label{TT-splitt1}
\underline \bG^{<n} &=&\llangle \underline \bG^{(1)}, \underline \bG^{(2)}, \ldots, \underline \bG^{(n-1)} \rrangle \in \Real^{I_1 \times I_2 \times \cdots \times I_{n-1} \times R_{n-1}} \\
\underline \bG^{>n} &=& \llangle \underline \bG^{(n+1)}, \underline \bG^{(n+2)}, \ldots, \underline \bG^{(N)}\rrangle \in \Real^{R_n \times I_{n+1}  \times \cdots \times I_{N}}
\label{TT-splitt2}
\ee
while their corresponding unfolding matrices, also called interface matrices, are given by:
\be
&&\bG^{< n}= [\tG^{< n}]_{<n-1>} \in \Real^{I_1 I_2  \cdots  I_{n-1} \times R_{n-1}}, \notag  \\
&&\bG^{> n} = [\tG^{> n}]_{(1)}  \in \Real^{R_n \times I_{n+1}   \cdots I_{N}}.
\label{interface-matrices}
\ee
For convenience, the left and right unfoldings of the cores are defined as
\be
\bG_L^{(n)}= [\tG^{(n)}]_{<2>} \in \Real^{R_{n-1} I_n \times R_n} \;\; \mbox{and} \;\; \bG_R^{(n)}= [\tG^{(n)}]_{<1>} \in \Real^{R_{n-1} \times  I_n R_n} \notag
\ee

It is important to mention here, that the orthogonalization of core tensors is an essential procedure in many algorithms for the TT formats (in order to  reduce computational complexity of contraction of core tensors and improve robustness of the algorithms) \citep{OseledetsTT11,dolgovEIG2013,Dolgovth,KressnerEIG2014,steinlechner15,Steinlechner_phd2016}.\\


 When considering the $n$th TT core, it is usually assumed that all cores to its left are left-orthogonalized, and all cores to its right are right-orthogonalized.

  Notice that if a TT network, $\underline \bX$, is $n$-orthogonal then the interface matrices are also orthogonal, i.e.,
\be
(\bG^{<n})^{\text{T}} \; \bG^{< n} & =& \bI_{R_{n-1}}, \notag\\
\bG^{> n} \; (\bG^{>n})^{\text{T}} & = & \bI_{R_{n}}. \notag
\ee
%

\noindent Through basic multilinear algebra, we can  construct a set of linear  equations, referred to as the frame equation, given by
\be
\bx = \bG_{\neq\,n} \; \bg^{(n)}, \qquad n=1,2,\ldots, N,
\label{frame-eq}
\ee
where $\bx= \mbox{vec}(\underline \bX) \in \Real^{I_1 I_2 \cdots I_N}$, $\;\;\bg^{(n)} = \mbox{vec}(\underline \bG^{(n)}) \in \Real^{R_{n-1} I_n R_n}$, while a tall-and-skinny matrix, $\bG_{\neq n}$, called the frame matrix, is formulated as
\begin{equation}
 \bG_{\neq \, n} = \bG^{<n} \; \otimes_L \; \bI_{I_n} \; \otimes_L \;(\bG^{>n})^{\text{T}} \in \Real^{I_1 I_2  \cdots  I_{N} \times R_{n-1} I_n R_n}.   \\
\label{frame-matrices}
\end{equation}
{\bf Remark.} From Eqs. (\ref{frame-eq}) and (\ref{frame-matrices}), observe that the frame and interface matrices indicate an important property of the TT format -- its linearity with respect to each core $\underline \bG^{(n)}$, when expressed in the vectorized form (see  Figure~\ref{Fig:TTsplitt}(b)).

\begin{figure}[t]
(a)\\
\begin{center}
\includegraphics[width=9.99cm]{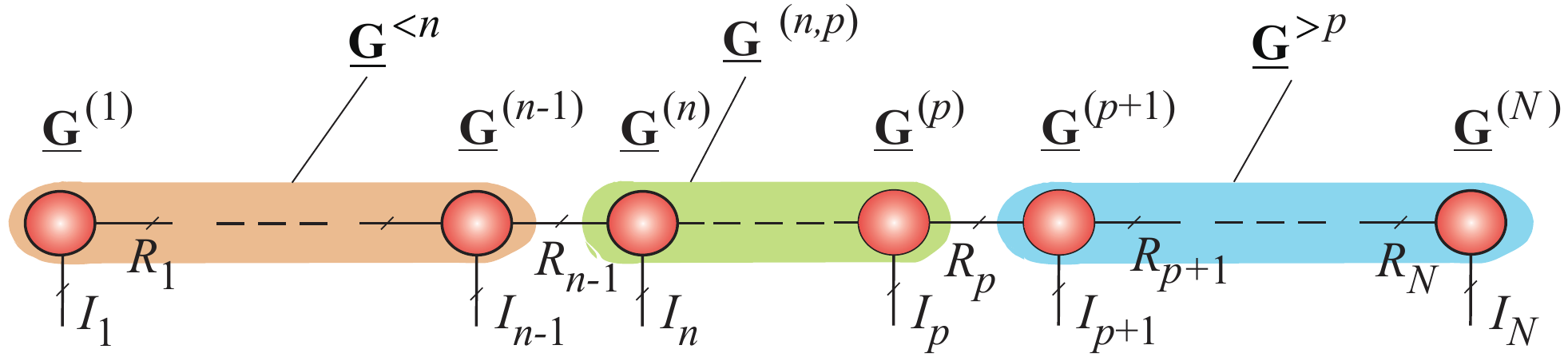}
\end{center}
(b)
\begin{center}
\includegraphics[width=11.5cm]{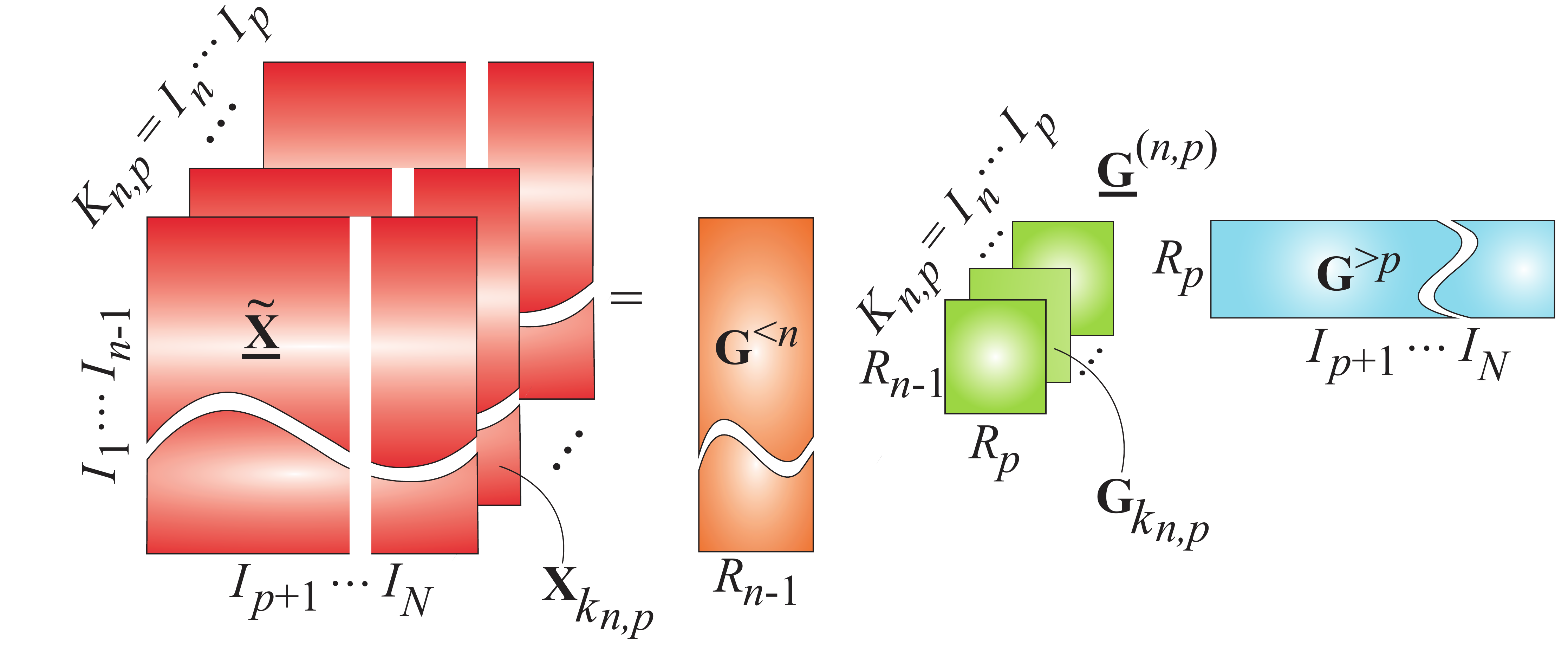}
\end{center}
\caption{Procedure for the extraction of an arbitrary group of cores from a TT network. (a) Splitting a TT network into three TT sub-networks.
(b)  Reshaping an $N$th-order data tensor to a 3rd-order  tensor  and its decomposition  using the Tucker-2/PVD model. This procedure
 can be implemented for any $n=1,2,\ldots, N-1$ and $ p>n$.
}
\label{Fig:TTsplitt3}
\end{figure}

Another important advantage of this approach is that by splitting the TT cores, we can express the data tensor in the following matrix form (see Figure~\ref{Fig:TTsplitt3} for $n=p$)
\be
\widetilde \bX_{k_n} =  \bG^{< n} \; \bG_{k_n}^{(n)} \; \bG^{>n} \in \Real^{I_1 \cdots I_{n-1} \times I_{n+1} \cdots I_N},
\label{eq:Grasedyck}
\ee
 where $\widetilde \bX_{k_n}$ are lateral slices of a 3rd-order reshaped  raw tensor,
 $\tilde {\underline \bX} \in \Real^{I_1 \cdots I_{n-1} \; \times \; I_n \; \times \; I_{n+1} \cdots I_N }$, and $n=1,2, \ldots,N, \;\; i_n=1,2, \ldots,I_n$,
 (obtained by  reshaping  an  $N$th-order  data tensor $\underline \bX \in \Real^{I_1 \; \times I_2 \; \times \; \cdots \; \times I_N }$).
Assuming that the columns of $\bG^{< n}$ and the rows of $\bG^{> n}$ are orthonormalized, the lateral slices of a core tensor $\underline \bG^{(n)}$ can be expressed as
\be
\bG^{(n)}_{i_n} =  (\bG^{< n})^{\text{T}}  \; \widetilde \bX_{i_n} \; (\bG^{>n})^{\text{T}} \in \Real^{R_{n-1} \; \times \; R_{n}}.
\label{eq:Grasedyck2}
\ee

\subsection{Extraction of Two Neighboring Cores  for  Modified ALS (MALS)}
\label{sect:two-cores}

The Modified ALS algorithm, also called the two-site Density Matrix Renormalization Group (DMRG2)
algorithm\footnote{The DMRG algorithm was first proposed by \citet{White1992}.
At that time, people did not know the relation between tensor network and
the DMRG. In  \citep{Ostlund95} pointed out that the wave function
generated by the DMRG iteration is a matrix product state.
The objective of the DMRG was  to compute the ground states (minimum
eigenpairs) of spin systems.} requires the extraction of two neighboring cores \citep{Holtz-TT-2012}.

\begin{figure}[t!]
(a)
\begin{center}
\includegraphics[width=9.99cm]{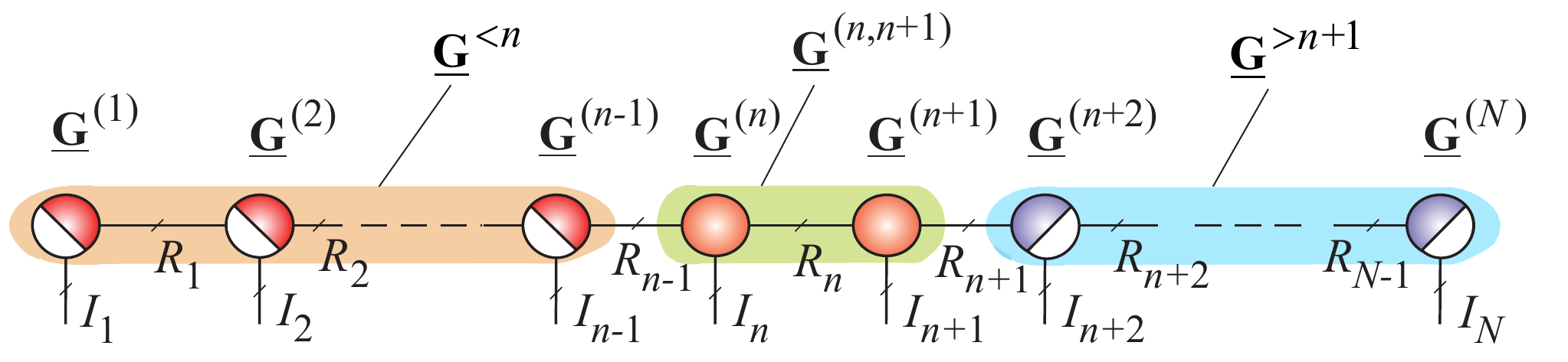}
\end{center}
(b)
\begin{center}
\includegraphics[width=11.5cm]{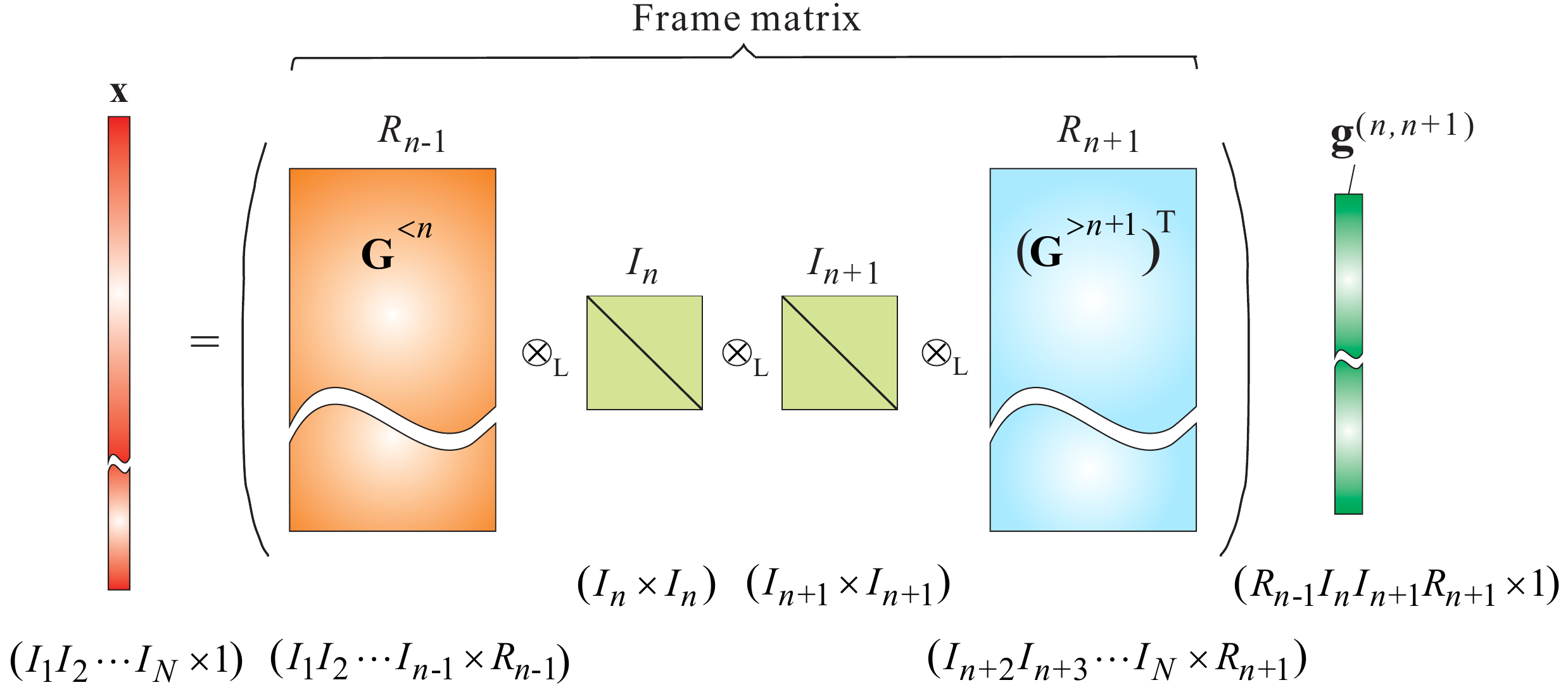}
\end{center}
\caption{Extraction of two neighboring cores. (a) Graphical representation of a tensor train and the left- and right-orthogonalized sub-trains.
(b) Graphical illustration of the  linear frame equation (see  Eqs. (\ref{frame-eq2}) and (\ref{frame-matrices2})).}
\label{Fig:TTsplitt2}
\end{figure}

Similarly to the previous section, the extraction of a block of two neighboring  TT cores is based on the following linear equation
\begin{equation}
\bx =\bG_{\neq \, n,n+1} \; \bg^{(n,n+1)}, \qquad n=1,2,\ldots, N-1,
\\
\label{frame-eq2}
\end{equation}
where the frame (tall-and-skinny) matrix is formulated as
\be
\bG_{\neq \, n,n+1} &=& \bG^{<n} \; \otimes_L \; \bI_{I_n} \; \otimes_L \; \bI_{I_{n+1}} \; \otimes_L \; (\bG^{>n+1})^{\text{T}}  \notag \\
&\in& \Real^{I_1 I_2  \cdots  I_{N} \times R_{n-1} I_n  I_{n+1} R_{n+1}}
\label{frame-matrices2}
\ee
and  $\bg^{(n,n+1)} =\mbox{vec}(\bG^{(n)\;\text{T}}_{L} \bG^{(n+1)}_{R})=
\mbox{vec}(\underline \bG^{(n,n+1)}) \in \Real^{R_{n-1} I_n I_{n+1} R_{n+1}}$, for $n=1,2,\ldots,N-1$
(see Figure~\ref{Fig:TTsplitt2}).

\begin{table}
\centering
\caption{Basic recursive formulas for the TT/MPS decomposition
of an $N$th-order tensor  $\underline \bX = \llangle \underline \bG^{(1)},\underline \bG^{(2)}, \ldots,\underline \bG^{(N)} \rrangle \in \Real^{I_1\times\cdots\times I_N}$, where $\;\underline \bG^{(n)} \in \Real^{R_{n-1} \times I_n \times R_n}$, $\; \bG^{(n)}_{(1)} \in \Real^{R_{n-1} \times I_n  R_n}$, $\;  \bG^{(n)}_{(2)} \in \Real^{I_n \times R_{n-1}  R_n}$, $\; \bG^{(n)}_{<2>} =[\bG^{(n)}_{(3)}]^{\text{T}} \in \Real^{R_{n-1} I_n \times R_n}$, and
$\;\;\underline \bX = \underline \bG^{\leq N} = \underline \bG^{\geq 1}$, $\;\; \underline \bG^{\leq 0} = \underline \bG^{\geq N+1}=1$.}
\vspace{.5pc}
{\footnotesize
\shadingbox{
\begin{tabular}{ll}
\hline\\[-0.8pc]
\multicolumn{1}{c}{TT (global)}   &   \multicolumn{1}{c}{TT (recursive)}\\
\hline
\multicolumn{2}{c}{Multilinear products}
\\
	$\displaystyle
	\underline \bX = 	\underline \bG^{(1)}\times^1 \underline \bG^{(2)}  \times^1\cdots \times^1 \underline \bG^{(N)}
	$
	&
	$\displaystyle
	\underline \bG^{\leq n} = \underline \bG^{\leq n-1} \times^1 \underline \bG^{(n)}
	$ 
\\[0.5ex]
	&
	$\displaystyle
	\underline \bG^{\geq n} =  \underline \bG^{(n)} \times^1 \underline \bG^{\geq n+1}
	$ 
\\[0.5ex]
\hline
\multicolumn{2}{c}{Vectorizations} \\
	$\displaystyle
	\text{vec}(\underline \bX)
	= \prod_{n=N}^1
	\left( \bG^{(n)\;\text{T}}_{(1)}  \otimes \bI_{I_1 I_2\cdots I_{n-1}} \right)
	$
	&
	\hspace{-1.2cm}
	$\displaystyle
	\text{vec} (\underline \bG^{\leq n} ) =
	$
\\[0.5ex]
	&
	\hspace{-2.1cm}	
		$=\displaystyle
		\left( \bG^{(n)\;\text{T}}_{(1)} \otimes  \bI_{I_1I_2\cdots I_{n-1}}  \right)
		\text{vec} \left(\underline \bG^{\leq n-1}\right)
		$
\\[0.5ex]
	$\displaystyle
	\text{vec}(\underline \bX)
	= \prod_{n=1}^N \left( \bI_{I_{n+1}\cdots I_N} \otimes  \bG^{(n)}_{<2>}  \right)
	$
	&
	\hspace{-1.2cm} 
	$\displaystyle
	\text{vec}\left(\underline \bG^{\geq n}\right) =
	$
\\[0.5ex]
	&	
	\hspace{-2.5cm} 
		$\displaystyle
		=\left( \bI_{I_{n+1}I_{n+2}\cdots I_N} \otimes  \bG^{(n)}_{<2>}  \right)
		\text{vec}\left(\underline \bG^{\geq n+1}\right)
		$
\\[0.5ex]
	{
	$\displaystyle
	\text{vec}(\underline \bX) =
	\left((\bG^{>n}_{(1)} )^{\text{T}} \otimes \bI_{I_n}\otimes (\bG^{<n}_{(n)} )^{\text{T}}  \right)
		\text{vec} \left(\underline \bG^{(n)} \right)	
	$
	}
	&
\\[1.4pc]
\hline
\multicolumn{2}{c}{Matricizations} \\
	$\displaystyle
	\bX_{(n)} =
	\bG^{(n)}_{(2)}
	\left(\bG^{>n}_{(1)} \otimes \bG^{< n}_{(n)}\right)
	$
	&
	\hspace{-1.3cm} $\displaystyle
	\bG^{<n}_{(n)} = \bG^{(n-1)}_{(3)}
	\left(\bI_{I_{n-1}} \otimes  \bG^{<n-1}_{(n-1)} \right)
	$
\\[0.5ex]
	$\displaystyle
	\bX_{<n>} = \bG^{\leq n}_{<n>} \bG^{>n}_{(1)}
	= \left( \bG^{< n+1}_{(n+1)} \right)^{\text{T}} \bG^{>n}_{(1)}
	$
	&
	\hspace{-1.3cm} $\displaystyle
	\bG^{\geq n}_{(1)} = \bG^{(n)}_{(1)}
	\left(\bG^{\geq n+1}_{(1)}  \otimes  \bI_{I_{n}}\right)
	$
\\[0.5ex]
$\displaystyle
	\bX_{<n>} = \left(\bI_{I_n}  \otimes  \bG^{\leq n-1}_{<n-1>}\right) \bG^{(n)}_{<2>} \bG^{\geq n+1}_{(1)}
	$
	&
\\[1pc]
$\displaystyle
	\bX_{<n-1>} = \bG^{\leq n-1}_{<n-1>} \bG^{(n)}_{(1)} \left(\bG^{\geq n+1}_{(1)} \otimes  \bI_{I_n}\right)
	$
	&
\\[1pc]
\hline
\multicolumn{2}{c}{Slice products} \\
	$\displaystyle
	x_{i_1,i_2,\ldots,i_N} =
	\bG^{(1)}_{i_1}
	\bG^{(2)}_{i_2}
	\cdots 	\bG^{(N)}_{i_N}
	$
	&
	\hspace{-0.7cm}
	$\displaystyle
	\bG^{\leq n}_{i_1,\ldots,i_n} =
	\bG^{\leq n-1}_{i_1,\ldots,i_{n-1}}
	\bG^{(n)}_{i_n}
	$
\\[0.5ex]
	&
	\hspace{-0.76cm}
	$\displaystyle
	\bG^{\geq n}_{i_n,\ldots,i_N} =
	\bG^{(n)}_{i_n}
	\bG^{\geq n+1}_{i_{n+1},\ldots,i_N}
	$
\\[0.5ex]
\hline
\end{tabular}
}}
\label{tab:TT_recursive}
\vspace{30pt}
\end{table}

Simple matrix manipulations now yield the following useful recursive formulas (see also Table  \ref{tab:TT_recursive})
\be
\bX_{(n)} &=& \bG^{(n)}_{(2)} \; ( \bG^{<n}_{(n)} \otimes_L \bG_{(1)}^{>n}) ,\\
\bG_{\neq \;n} &=& \bG_{\neq \;n,n+1} \;  (\bI_{R_{n-1}I_n} \otimes_L (\bG^{(n+1)}_{(1)})^{\text{T}}) ,\\
\bG_{\neq \;n+1} &=& \bG_{\neq \;n,n+1} \; ((\bG^{(n)}_{(3)})^{\text{T}} \; \otimes_L \; \bI_{I_{n+1} R_{n+1}}).
\label{frame-matrices12}
\ee
If the cores are normalized in a such way that all cores to the left of the currently considered (optimized) TT core $\underline \bG^{(n)}$ are left-orthogonal{\footnote{We recall here, that $\bG_{(3)}^{(n)} \in \Real^{R_n \times R_{n-1} I_n}$ and $\bG_{(1)}^{(n)} \in \Real^{R_{n-1} \times I_n R_{n}}$ are the mode-3 and mode-1 matricization of the TT core tensor $\underline \bG^{(n)} \in \Real^{R_{n-1} \times I_n \times R_n}$, respectively.}}
 \be \bG_{(3)}^{(k)}\;\bG_{(3)}^{(k)\;\text{T}} = \bG_{<2>}^{(k)\;\text{T}}\;\bG_{<2>}^{(k)}=\bI_{R_k}, \quad k<n,
 \ee
 and all the cores to the right of  $\underline \bG^{(n)}$ are right-orthogonal
\be
\bG_{(1)}^{(p)} \; \bG_{(1)}^{(p)\; \text{T}}=\bG_{<1>}^{(p)} \; \bG_{<1>}^{(p)\;\text{T}} =\bI_{R_{p-1}}, \quad p>n,
\ee
 then the frame matrices have orthogonal columns \citep{NRG-DMRG12,dolgovEIG2013,KressnerEIG2014}, that is
\be
\bG_{\neq\, n}^{\text{T}} \; \bG_{\neq \; n}&=&\bI_{R_{n-1} \,I_n \,R_n}, \\
\bG_{\neq \, n,n+1}^{\text{T}} \; \bG_{\neq\,n,n+1}&=&\bI_{R_{n-1} \, I_n \,I_{n+1}\, R_{n+1}}.
\label{orth-frames}
\ee

 Figure~\ref{Fig:TTsplitt3} illustrates that the  operation of splitting the TT cores into  three sub-networks can be expressed
 in the following matrix form (for $p=n+1$)
\be
\widetilde\bX_{k_{n,n+1}} =  \bG^{< n}  \; \bG_{k_{n,n+1}}^{(n,n+1)} \; \bG^{>n+2} \in \Real^{I_1 \cdots I_{n-1} \; \times \; I_{n+2} \cdots I_N},
\label{eq:Grasedyck3}
\ee
for $n=1,2, \ldots,N-1, \;\; k_{n,n+1}=1,2, \ldots,I_n I_{n+1}$, where $\widetilde \bX_{k_{n,n+1}}$
are the frontal slices of a 3rd-order reshaped data tensor
 $\tilde {\underline \bX} \in \Real^{I_1 \cdots I_{n-1} \; \times \; I_{n} I_{n+1} \; \times \; I_{n+2} \cdots I_N}$.

Assuming that the columns of $\bG^{< n}$ and rows of $\bG^{> n+2}$ are orthonormalized, the frontal slices of
a super-core tensor, $\underline \bG^{(n,n+1)} \in \Real^{R_{n-1} \times I_n I_{n+1} \times R_{n+1}}$, can be expressed in the following simple form
\be
\bG^{(n,n+1)}_{k_{n,n+1}} =  (\bG^{< n})^{\text{T}}  \; \widetilde \bX_{k_{n,n+1}} \; (\bG^{>n+2})^{\text{T}} \in \Real^{R_{n-1} \times R_{n+1}}.
\label{eq:Grasedyck4}
\ee

\section{Alternating Least Squares (ALS) and Modified ALS (MALS)}
\sectionmark{ALS and MALS for TT}
\label{sect:ALS-TT}

Consider the minimization of a scalar cost (energy) function, $J(\underline \bX) $, of an $N$th-order tensor variable expressed in the TT format as
$\underline \bX \cong \llangle \underline \bX^{(1)},\underline \bX^{(2)},\ldots, \underline \bX^{(N)}\rrangle$.
The solution is sought in the form of a tensor train; however, the  simultaneous minimization over all cores $\underline \bX^{(n)}$
is usually  too complex and  nonlinear. For feasibility, the procedure is replaced by a sequence of optimizations carried out over
one core at a time, that is
$\underline \bX^{(n)} = \argmin J (\underline \bX^{(1)}, \ldots, \underline \bX^{(n)}, \ldots, \underline \bX^{(N)})$.

\begin{figure}[p]
\begin{center}
\hspace{0.9cm} \includegraphics[width=9.0cm]{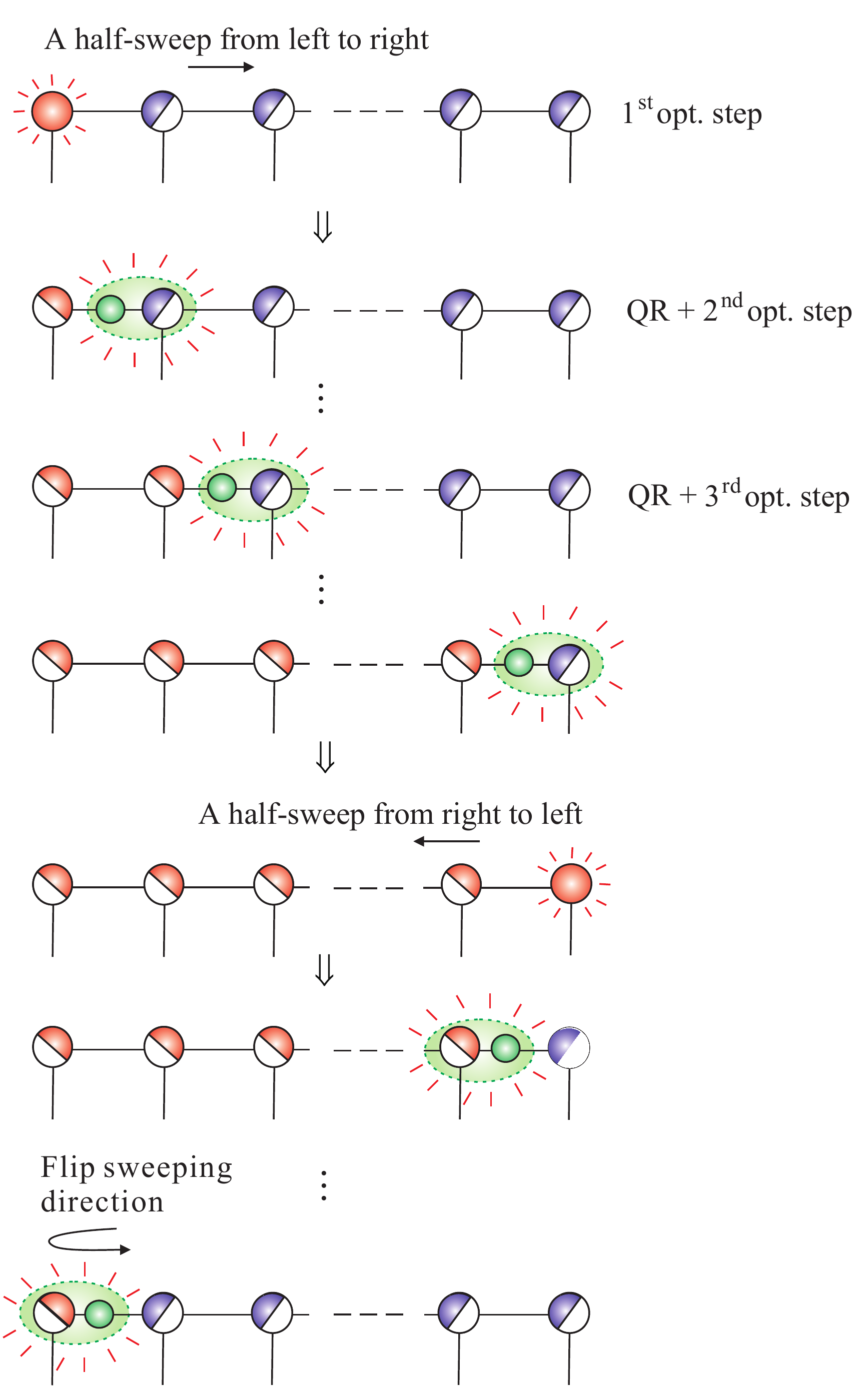}
\end{center}
\caption{The ALS algorithm for TT decomposition  corresponding to the DMRG1. The  idea is to optimize only one core tensor at a time (by the minimization of a suitable cost function), while keeping the other cores  fixed. Optimization of each core tensor is followed by an orthogonalization step via the QR/LQ  decompositions. Factor matrices $\bR$ are absorbed (incorporated) into the subsequent core. The small green circle denotes a triangular matrix $\bR$ or $\bL$ which is merged with a neighbouring TT core (as indicted by the green shaded ellipsoid).}
\label{Fig:TTALS}
\end{figure}

The idea behind the Alternating Least Squares (ALS) optimization{\footnote{The ALS algorithms  can be considered as the block nonlinear Gauss-Seidel iterations.}  (also known as the Alternating Linear Scheme or one-site DMRG (DMRG1))} is that in each
local optimization (called the micro-iteration step), only one core tensor, $\underline \bX^{(n)}$, is updated, while all other cores are kept fixed.
Starting from some initial guess for all cores, the method  first updates, say,  core $\underline \bX^{(1)}$,
while  cores $\underline \bX^{(2)}, \ldots \underline \bX^{(N)}$,
are fixed, then it  updates $\underline \bX^{(2)}$ with
$\underline \bX^{(1)}, \underline \bX^{(3)}, \ldots, \underline \bX^{(n)}$
fixed and so on, until $\underline \bX^{(N)}$ is optimized. After completing this forward half-sweep, the algorithm proceeds backwards along the sequence of cores $N,N-1, \ldots,1$, within the so-called backward half-sweep. The sequence of forward and backward iterations complete one full sweep\footnote{Note that this sweeping process operates in a similar fashion as the  self-consistent recursive loops, where the solution is improved iteratively and gradually.}, which corresponds to one (global-) iteration. These iterations are  repeated until some stopping criterion is satisfied, as illustrated in Algorithm~\ref{alg:ALS-TT} and  Figure~\ref{Fig:TTALS}.

\begin{algorithm}[t]
{\small
\caption{\textbf{Alternating Least Squares (ALS)}}
\label{alg:ALS-TT}
\begin{algorithmic}[1]
\REQUIRE Cost (energy) function $ J (\underline \bX)$ and an initial guess for an  \\
$N$th-order tensor in the TT format
$\underline \bX = \llangle \underline \bX^{(1)}, \underline \bX^{(2)}, \ldots , \underline \bX^{(N)}\rrangle$
\ENSURE A tensor $\tX$ in the TT format minimizes the cost\\
 function $J(\tX)$
\WHILE {Stopping condition not fulfilled}
\FOR{$n=1$ to $N$} %
\STATE \hspace{-0.3cm} Find:
$\hat{\underline \bX}^{(n)} = \displaystyle\mathop{\rm arg\,min}_{\underline \bX^{(n)}_*}
J(\hat{\underline \bX}^{(1)}, \ldots, \hat{\underline \bX}^{(n-1)}, \underline \bX_*^{(n)}, \underline \bX^{(n+1)},\ldots, \underline \bX^{(N)})$
\ENDFOR
\ENDWHILE
\end{algorithmic}
}
\end{algorithm}

As in any iterative optimization based on gradient descent, the cost function can only decrease, however, there is no guarantee that a global
minimum would be reached. Also, since the TT rank for the desired solution is unknown, the standard ALS relies on an initial guess for the TT rank;
 for some initial conditions the iteration process can therefore be very slow. To alleviate these problems, the modified ALS (MALS/DMRG2)
 scheme
  aims to merge two neighboring TT cores (blocks), optimize the resulting ``super-node''  also called ``super-core'' or ``super-block'', and
split again the result into separate factors by low-rank matrix factorizations, usually using truncated SVD {\footnote{Although DMRG2 and MALS are equivalent, the SVD was not used in original DMRG2
algorithm. As a matter of fact, quantum physicists use rather the time-evolving block decimation (TEBD) algorithm in order to optimize matrix product states (TT/MPS) (for more detail see \citet{Vidal03,OrusVidal2008,Orus2013}.)}}.

\begin{figure}[t!]
\begin{center}
\includegraphics[width=10.5cm]{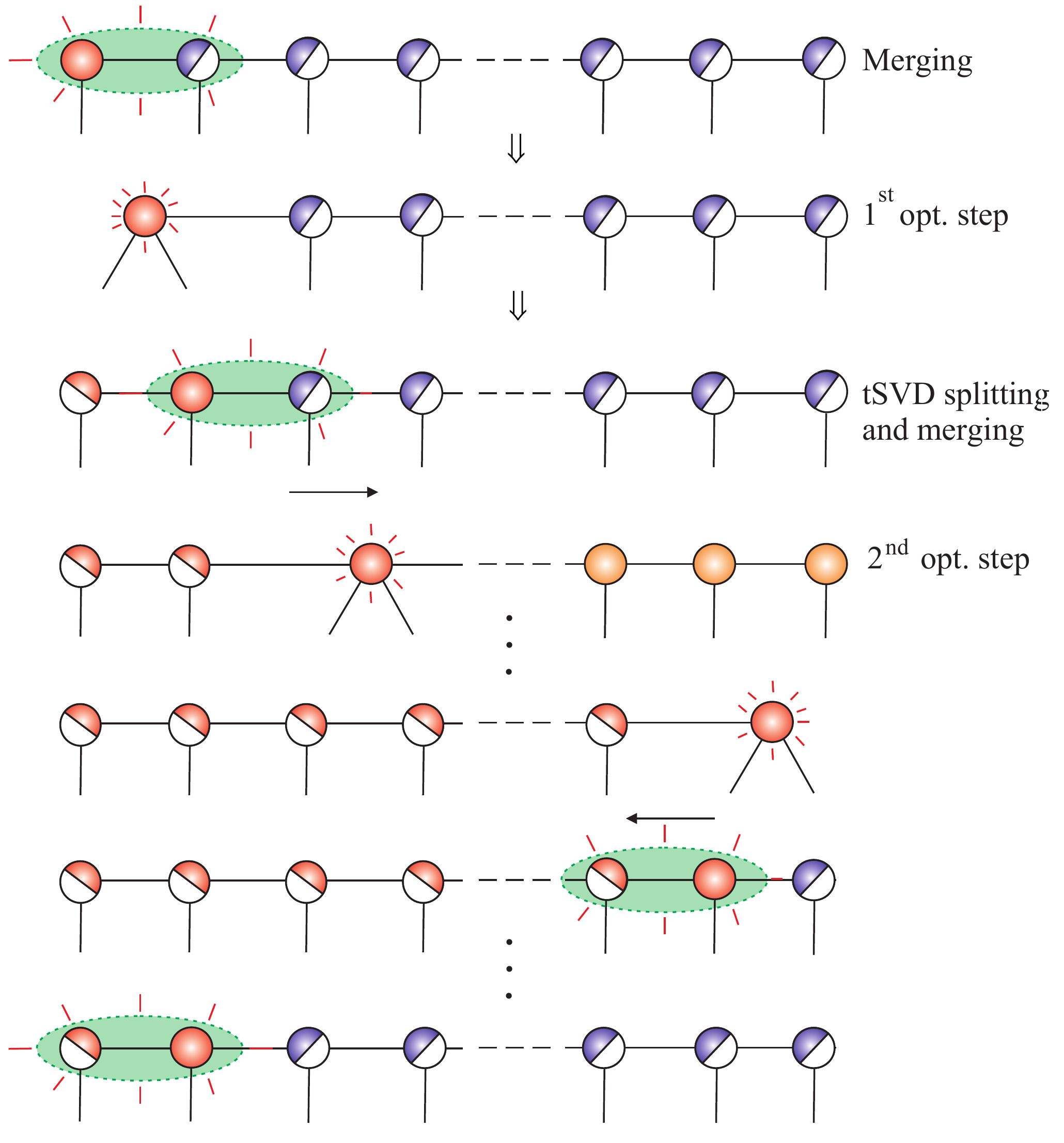}
\vspace{-12pt}
\end{center}
\caption{ The Modified ALS (MALS)  algorithm corresponding to the DMRG2. In every  optimization step, two neighboring  cores are merged and the optimization is performed over the  so obtained ``super-core''. In the next step, truncated SVD or other low-rank matrix factorizations (LRMF) are applied in order to split the optimized super-core into two cores.  Note that each optimization sub-problem is computationally more expensive than  the standard ALS. During the iteration process, convergence speed may increase considerably as compared to the standard ALS algorithm, and TT ranks can be adaptively estimated.} 
\label{Fig:MALS}

\end{figure}

\noindent{\bf Some remarks:}
\begin{itemize}
\item  If the global cost (loss) function $J(\underline \bX)$ is quadratic, so too are the local cost functions  $J_n(\underline \bX^{(n)})$, and the problem reduces to solving a small-scale optimization at each micro-iteration. In addition, local problems in the micro-iterations are usually
  better conditioned, due to the orthogonalization of TT cores.
\item In most optimization problems considered in this monograph the TT decomposition is assumed to be in an orthogonalized form, where all cores are left-- or right--orthogonal with respect to the core $\underline \bX^{(n)}$, which is being
optimized (see also Figure~\ref{Fig:TTALS} and Figure~\ref{Fig:MALS})  \citep{Holtz-TT-2012,QTT-Tucker}.
\end{itemize}

\section{Tensor Completion for Large-Scale Structured Data}
\sectionmark{Tensor Completion}
\label{sect:completion}

The objective of tensor completion is to reconstruct a high-dimensional structured multiway array for which a large proportion of entries is missing or noisy  \citep{Grasedyck15TTcompl,Steinlechner_phd2016,steinlechner15,zhao2015bayesian,BayesCPD-TNNLS,Yokota-SmCA}.

 With the assumption that a good low-rank TN approximation for an incomplete data tensor
$\underline \bY \in \Real^{I_1 \times I_2 \times \cdots \times I_N}$ does exist, the tensor completion problem can be formulated as the following  optimization problem
\be
\min_{\underline \bX} \; J(\underline \bX) = \| P_{\Omega} (\underline \bY) - P_{\Omega} (\underline \bX) \|_F^2
\label{TT-compl}
\ee
subject to  $\underline \bX \in {\cal{M}}_R :=\{\underline \bX \in \Real^{I_1 \times I_2 \times \cdots \times I_N} \; | \; \mbox{rank}_{TT} (\underline \bX) = \brr_{TT}\}$, where $\underline \bX $ is the reconstruction tensor in TT format and
$\brr_{TT}=\{R_1,R_2, \ldots, R_{N-1}\}$. The symbol $P_{\Omega}$ denotes the projection onto the sampling set $\Omega$ which corresponds to the indices of the known entries of $\underline \bY$, i.e., $ P_{\Omega} ( \underline \bY(i_1,i_2,\ldots, i_N)) =  \underline \bY(i_1,i_2,\ldots, i_N)$ if
 $(i_1,i_2, \ldots,i_N) \in \Omega$, otherwise zero.

For a large-scale high-order reconstruction tensor represented in the TT format, that is,
$\underline \bX= \llangle \underline \bX^{(1)}, \underline \bX^{2)} , \ldots, \underline \bX^{(N)}\rrangle$, the above optimization problem can be converted by the ALS approach to a set of \textit{smaller} optimization problems. These can be represented in a scalar form through slices of  core tensors, as follows 
\[
 \widehat{\underline \bX}^{(n)}
= \mbox{arg}\min_{\underline \bX^{(n)}} \, \sum_{\bi \in \Omega} \biggl(\underline \bY(i_1,\ldots,i_N)
 - \bX^{(1)}_{i_1} \cdots  \bX^{(n)}_{i_n} \bX^{(n+1)}_{i_n} \cdots \bX^{(N)}_{i_N} \biggr)^{\!2}. 
\]
For convergence, it is necessary to iterate through all cores, $\underline \bX^{(n)}$, $n=1,2,\ldots,N$, and over several sweeps.

 For efficiency, all the cores are orthogonalized after each sweep, and the cost function is minimized in its vectorized form by sequentially minimizing each individual slice, $\bX^{(n)}_{i_n}$, of the cores,
   $\underline \bX^{(n)} \in \Real^{R_{n-1}\times I_n \times R_n}$  (see also formulas (\ref{eq:Grasedyck}) and (\ref{eq:Grasedyck2})), that is
  \be
 \widehat{\bX}^{(n)}_{i_n} &=& \mbox{arg}\min_{\bX^{(n)}_{i_n}} \, \| \mbox{vec}(\widetilde{\underline \bY}_{i_n})
 - \mbox{vec} ( \bX^{< n} \; \bX^{(n)}_{i_n} \; \bX^{> n})\|_\Omega^2 \\
& =& \mbox{arg}\min_{\bX^{(n)}_{i_n}} \, \| \mbox{vec}(\widetilde{\underline \bY}_{i_n}
 - [ \bX^{< n} \; \otimes_L \; (\bX^{> n})^{\text{T}}]  \mbox{vec}(\bX^{(n)}_{i_n})\|_\Omega^2, \notag
\label{TT-compl4}
\ee
for all $i_n=1,2,\ldots,I_n$ and  $n=1,2,\ldots,N$, where $\| \bv \|_\Omega^2 = \sum_{i\in \Omega} v_i^2$.

The above optimization problem for $\mbox{vec}(\bX^{(n)}_{i_n})$
can be considered as a standard linear least squares (LS) problem,
 which can be efficiently solved even for huge-scale datasets.

In general, TT ranks are not known beforehand and must be estimated during the optimization procedure.
These can be estimated by, for example, starting with a maximum rank $R_{\max}$ and reducing gradually ranks by TT rounding.
Alternatively, starting from a minimum rank, $R_{\min}$, the ranks could be gradually increased.
 The rank increasing procedure can start with $\brr_{TT} =\{1,1,\ldots,1\}$, and the so obtained result is used for another run (sweep), but with $\brr_{TT} = \{1, 2,\ldots,1\}$. In other words, the rank $R_2$ between the second and third core is increased, and so on until either the prescribed residual tolerance $\varepsilon$ or $R_{\max}$ is reached. As the rank value is a key factor which determines the complexity of the ALS, the second approach is often more efficient and is computationally comparable with rank-adaptive approaches \citep{Grasedyck15TTcompl,Steinlechner_phd2016,steinlechner15}.

\section{Computing a Few Extremal Eigenvalues and Eigenvectors}


\subsection{TT Network for Computing the Single Smallest Eigenvalue and the
Corresponding  Eigenvector}
\subsectionmark{Computing Smallest Eigenvalue}
\label{sect:single-EVD}

Machine learning applications often require the computation of extreme (minimum or maximum) eigenvalues
and the corresponding eigenvectors of large-scale structured symmetric matrices. This problem can  be formulated as the standard symmetric
eigenvalue decomposition (EVD), in the form
\be
\bA \, \bx_k = \lambda_k \bx_k,   \qquad  k=1,2,\ldots,K,
\label{eq:eigp} \ee
where $\bx_k \in \Real^I$ are the orthonormal eigenvectors and $\lambda_k$ the corresponding eigenvalues of a symmetric matrix  $\bA \in \Real^{I \times I}$.
%

Iterative algorithms for  extreme eigenvalue problems often
exploit the Rayleigh Quotient (RQ) of a symmetric  matrix
 as the following  cost function
%
\be
J(\bx) = R(\bx, \bA) = \frac{\bx^{\text{T}} \bA \bx}{\bx^{\text{T}} \bx}= \frac{\langle \bA \bx, \bx\rangle}{\langle\bx,\bx\rangle}, \qquad \bx \neq 0.
\ee
Based on the Rayleigh Quotient  (RQ), the largest,  $\lambda_{max}$, and smallest, $\lambda_{min}$, eigenvalue of the matrix $\bA$ can be computed as
\be
\lambda_{max} = \max_{\bx} R( \bx, \bA ),\quad
 \lambda_{min} = \min_{\bx}  R( \bx, \bA ),
\ee
while the critical points and critical values of
$R(\bx, \bA)$ are the corresponding eigenvectors and eigenvalues of $\bA$.

The  traditional methods for solving eigenvalue problems for a symmetric matrix, $\bA \in \Real^{I \times I}$, are prohibitive for
 very large values of $I$, say $I=10^{15}$ or higher.
This computational bottleneck can be very efficiently dealt with through low-rank tensor approximations, and the last 10 years have witnessed the development of such techniques for several classes of optimization problems, including EVD/PCA and SVD \citep{dolgovEIG2013,KressnerEIG2014,Lee-SIMAX-SVD}.
The principle is to represent the cost function in a tensor format; under certain conditions, such as that tensors can be often quite well approximated in a low-rank TT  format, thus allowing for low-dimensional parametrization. 

In other words, if a structured symmetric  matrix $\bA$  and its eigenvector $\bx$ admit low-rank TT approximations,
 a large-scale eigenvalue problem can be converted into a set of smaller optimization problems
by representing the eigenvector, $\bx \in \Real^I $, and the matrix, $\bA \in \Real^{I \times I} $, in TT (MPS/MPO) formats, i.e.,
a TT/MPS for an $N$th-order tensor $\underline \bX \cong \llangle \underline \bX^{(1)}, \ldots, \underline \bX^{(N)} \rrangle \in \Real^{I_1 \times \cdots \times I_N}$
and a matrix TT/MPO for a $2N$th-order tensor  $\underline \bA \cong \llangle \underline \bA^{(1)}, \ldots, \underline \bA^{(N)} \rrangle \in \Real^{I_1 \times I_1 \times \cdots \times I_N \times I_N}$, where $I=I_1 I_2 \cdots I_N$.

\begin{figure}[ht!]
\begin{center}
\includegraphics[width=0.99\textwidth]{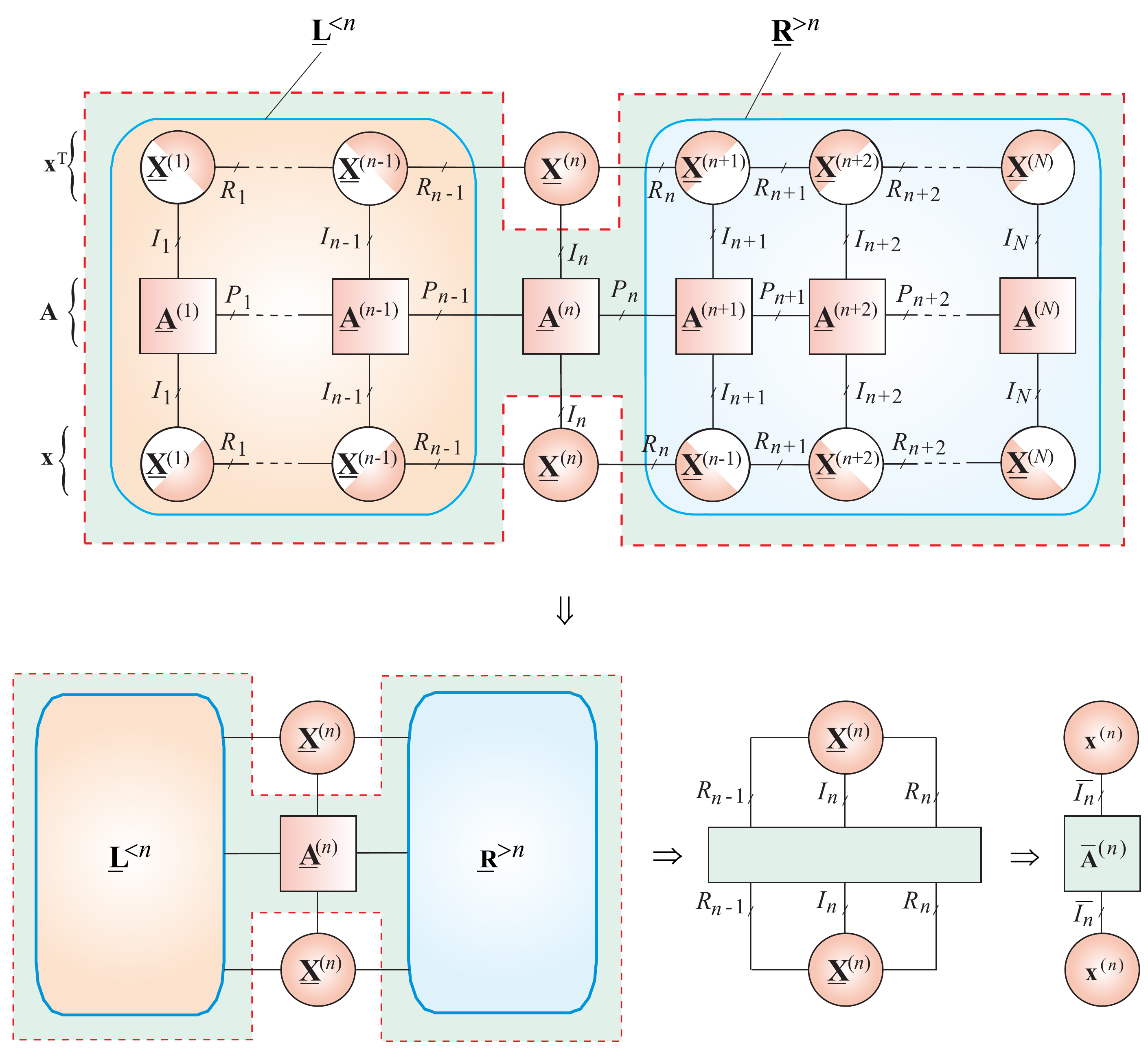}
\end{center}
\vspace{-12pt}
\caption{A conceptual TT network for the computation of a single extreme eigenvalue, $\lambda$, and the corresponding eigenvector,
$\bx \in \Real^I$,  for a symmetric  matrix $\bA \in \Real^{I \times I}$. The frame matrix  maps the TT core  into a large vector.  The  tensor network corresponds to the cost function (quadratic form), $\bx^{\text{T}} \bA \bx$, where the matrix $\bA$ and vectors $\bx \in \Real^I$ are given in the tensor train format with distributed  indices $I=I_1 I_2 \cdots I_N$. The cores in the  shaded areas form the matrix $\overline \bA^{(n)}$ (the effective Hamiltonian), which can be computed by a sequentially optimized contraction of the TT cores.}
\vspace{-12pt}
\label{Fig:EIG1}
\end{figure}

Figure~\ref{Fig:EIG1} illustrates this conversion into a  set of much smaller sub-problems by employing  tensor contractions and the frame equations (see section \ref{sect:orthog})
\be
\bx =\bX_{\neq n} \, \bx^{(n)}, \qquad n=1,2,\ldots,N,
\ee
with the frame  matrices
\be
\bX_{\neq n} = \bX^{<n} \otimes_L \bI_{I_n} \otimes_L (\bX^{>n})^{\text{T}} \in \Real^{I_1 I_2  \cdots  I_{N} \times R_{n-1} I_n R_n}.
\label{frame-matricesX}
\ee
Since the cores $\underline \bX^{(m)}$ for $m\neq n$ are constrained to be left- or right-orthogonal, the RQ can be minimized (or maximized), as follows
\be
\min_{\bx} \, J(\bx) &=& \min_{\bx^{(n)}} \, J(\bX_{\neq n} \bx^{(n)}) \\
&=& \min_{\bx^{(n)}} \frac{  \bx^{(n)\;\text{T}} \; \overline \bA^{(n)} \; \bx^{(n)} } {\langle\bx^{(n)},\bx^{(n)}\rangle}, \quad n=1,2,\ldots, N, \quad  \nonumber
\label{RQframe}
\ee
where $\bx^{(n)} = \mbox{vec}(\underline\bX^{(n)}) \in \Real^{R_{n-1} I_n R_n}$ and the matrix $\overline \bA^{(n)}$
(called  the effective Hamiltonian) can be expressed as
\be
\overline \bA^{(n)} = (\bX_{\neq n})^{\text{T}} \bA \bX_{\neq n} \in \Real^{R_{n-1}I_n R_n \times R_{n-1}I_n R_n}
\label{frameAxn}
\ee
for $n=1,2,\ldots,N$, under the condition $(\bX_{\neq n})^{\text{T}} \bX_{\neq n} = \bI$.

\begin{figure}[ht!]
(a)
\begin{center}
 \includegraphics[width=7.2cm]{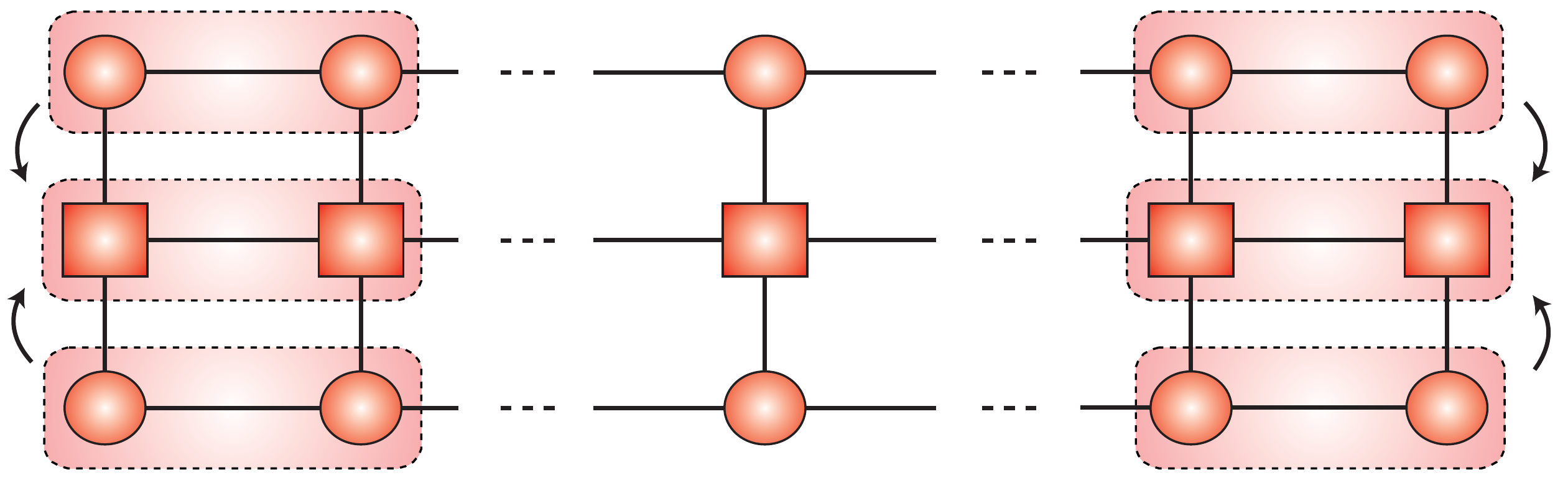}\\
 \end{center}
(b)
\begin{center}
  \includegraphics[width=7.2cm]{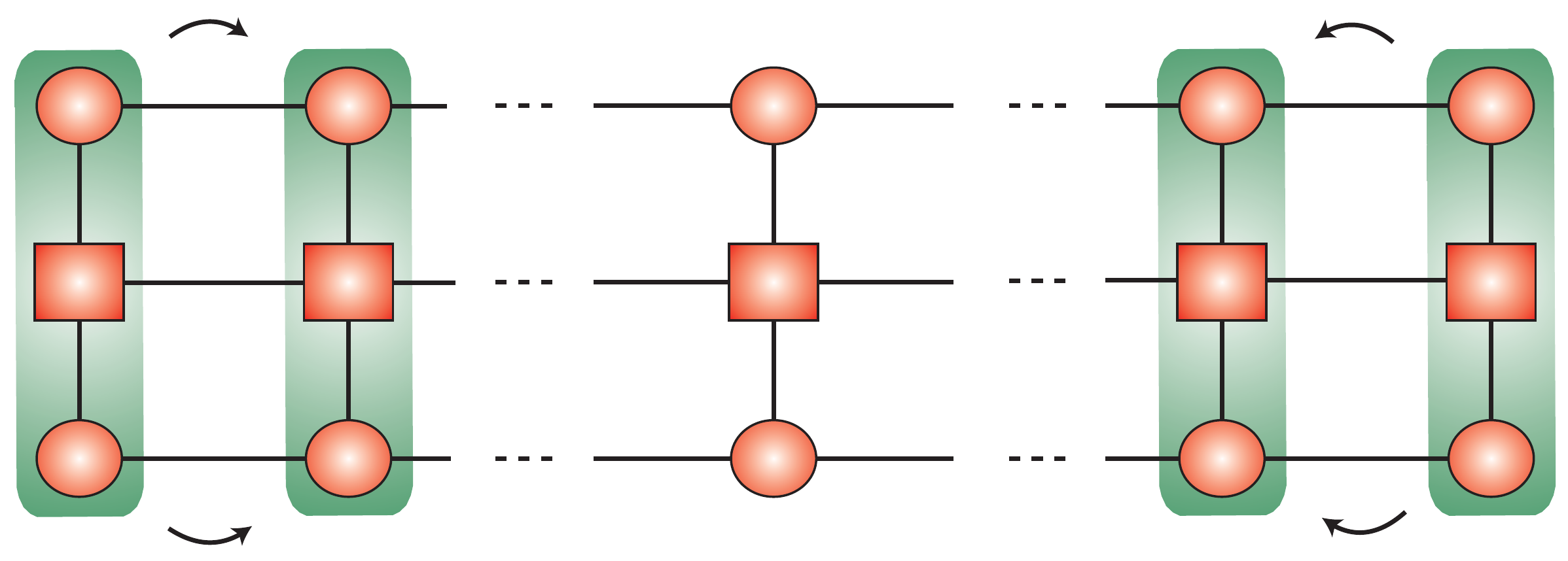}
  \end{center}
  (c)
\begin{center}
  \includegraphics[width=7.1cm]{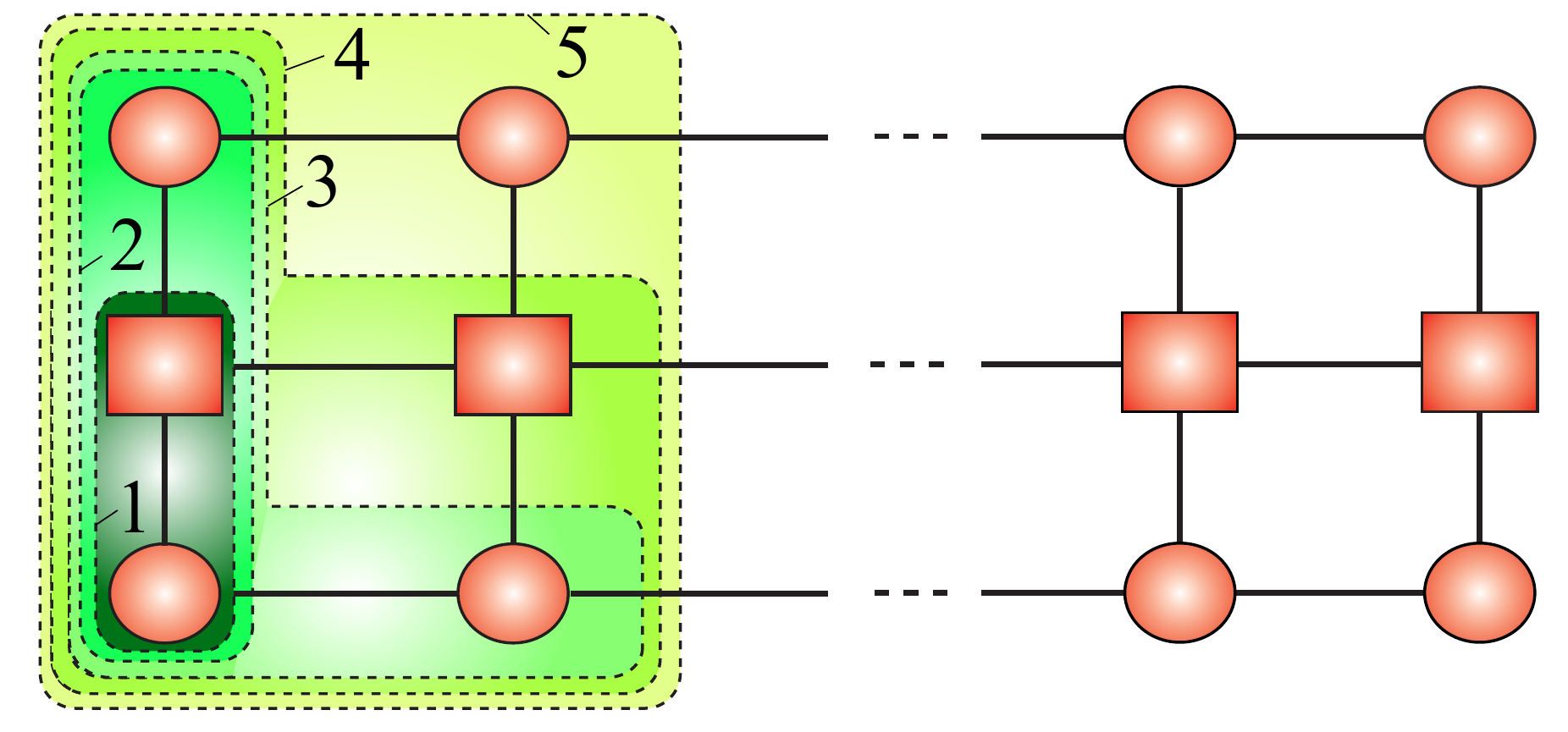}
  \end{center}
\caption{Sequential contraction of the TT  network. (a) Sub-optimal (inefficient) and (b) Optimal (efficient) contraction of the
TT (MPS/MPO) network. (c) Details of optimal sequential contraction (the numbers $1,2,3,\ldots$ indicate the order of tensor contractions).}
\label{Fig:TTcontrac}
\end{figure}

For relatively small TT ranks, the matrices $\overline \bA^{(n)}$ are usually much smaller than the original matrix $\bA$,
so that a large-scale EVD problem can be converted into a
set of  much smaller EVD sub-problems, which requires solving the  set of equations
\be
\overline \bA^{(n)} \, \bx^{(n)} = \lambda \bx^{(n)}, \quad n=1,2,\ldots,N.
\label{eq:Ax_lx}
\ee
In practice, the matrices $\overline \bA^{(n)}$ are never computed directly by (\ref{frameAxn}). Instead, we compute matrix-by-vector multiplication{\footnote{Such local matrix-by-vectors multiplications can be incorporated  to standard iterative methods, such as Arnoldi, Lanczos, Jacobi-Davidson and LOBPCG.}}, $ \bA^{(n)} \, \bx^{(n)}$,
 iteratively via optimized  contraction of TT cores within a tensor network. The concept of optimized contraction is  illustrated in Figure~\ref{Fig:TTcontrac}.

It should be noted that the contraction of the matrix $\bA$ in TT/MPO format, with corresponding TT cores of the vectors $\bx$ and $\bx^{\text{T}}$ in TT/MPS
format, leads to the left- and right
contraction tensors $\underline \bL^{<n}$ and $\underline \bR^{>n}$, as illustrated in Figure~\ref{Fig:EIG1}.
In other words, efficient solution of the matrix equation (\ref{eq:Ax_lx}) requires computation of the blocks $\underline \bL^{<n}$ and
$\underline \bR^{>n}$; 
these can be built iteratively so as to best  reuse available information,
which involves an optimal arrangement of a tensor  network contraction.
In a practical implementation, the full
network contraction is never carried out globally, but through iterative sweeps from right to left or vice-versa, to build up $\underline \bL^{<n}$
and  $\underline \bR^{>n}$ from the previous steps. 
The left and right orthogonalization of the cores can also be exploited to simplify the tensor
contraction process \citep{schollwock11-DMRG,dolgovEIG2013,KressnerEIG2014,Lee-SIMAX-SVD,Lee-TTPI}.

{\bf Remark.} Before we construct the tensor network shown in the Figure~\ref{Fig:EIG1}, we need to construct  approximate  distributed representation of the matrix $\bA$ in the TT/MPO format.
For ways to construct efficient representations of huge-scale matrix
in TT/MPO  format  while  simultaneously performing compression see
\citep{Hubig16,Huckle16} and also Part 1 \citep{Part1}.

\subsection{TT Network for Computing Several Extremal Eigenvalues and
Corresponding Eigenvectors for Symmetric EVD}
\subsectionmark{TT Networks for the EVD}

In a more general case, several (say, $K$) extremal eigenvalues  and the corresponding  $K$ eigenvectors of a symmetric structured matrix $\bA$ can be computed
through the following trace
minimization (or maximization) with orthogonality constraints
\be
\min_{\bX} \, \tr(\bX^{\text{T}} \bA \bX),\qquad \mbox{s.t.} \;\; \bX^{\text{T}} \bX =\bI_K,
\label{Trace1}
\ee
where $\bX=[\bx_1,\bx_2,\ldots,\bx_K] \in \Real^{I_1 I_2 \cdots I_N \; \times \; K}$ is a matrix composed of eigenvectors $\bx_k$, and
 $\bA \in \Real^{I \times I}$ is a huge-scale symmetric matrix.

Note that the global minimum of this cost function is equivalent to a sum of $K$ smallest eigenvalues, $\lambda_1 +\lambda_2 + \cdots + \lambda_K$.

 For a very large-scale  symmetric matrix $\bA \in \Real^{I \times I}$  with say, $I=10^{15}$ or higher, the optimization problem in
   (\ref{Trace1}) cannot be solved directly.  A feasible TT solution would be to first   tensorize the matrices $\bA \in \Real^{I \times I}$
   and $\bX \in \Real^{I \times K}$  into the respective  tensors $\underline \bA \in \Real^{I_1 \times I_1 \times  \cdots \times I_N  \times I_N}$
   and $\underline \bX \in \Real^{I_1 \times I_2 \times \cdots \times I_N \times K}$, where $I=I_1 I_2 \cdots I_N$ and the value of every  $I_n$ is much smaller than $I$; then, a tensor network structure can be exploited  to allow a good low-rank TN approximation of the underlying cost function.

\begin{figure}
\begin{center}
\includegraphics[width=9.5cm]{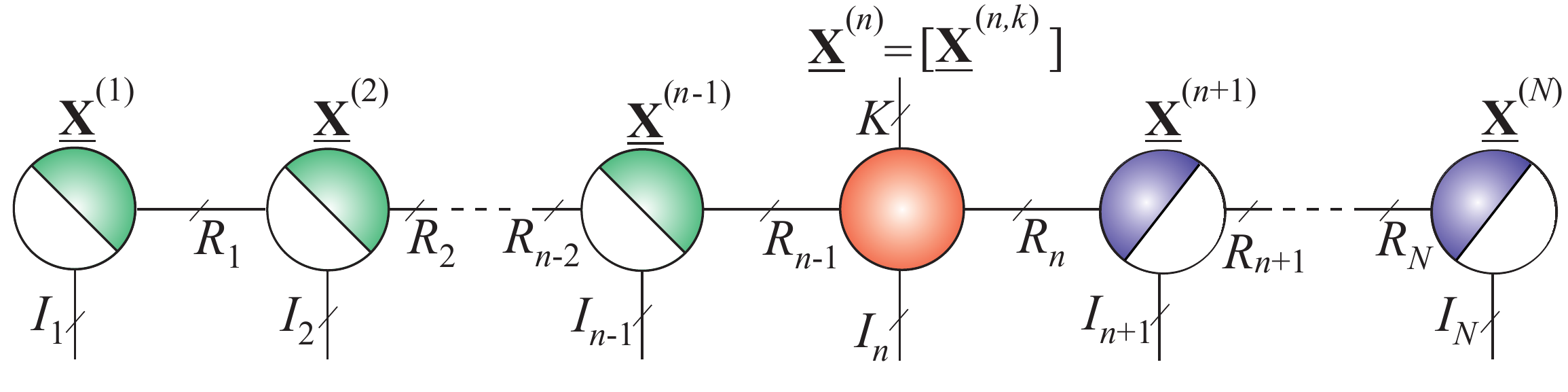}
\end{center}
\caption{Block-$n$  TT with left- and right-orthogonal cores, for  an $(N+1)$th-order tensor $\underline \bX \in \Real^{I_1 \; \times \; I_2 \; \times \cdots \times  \; I_N \; \times \; K}$.}
\label{Fig:EIG1KBTT}
\end{figure}

\begin{figure}[ht!]
\begin{center}
\includegraphics[width=11.99cm]{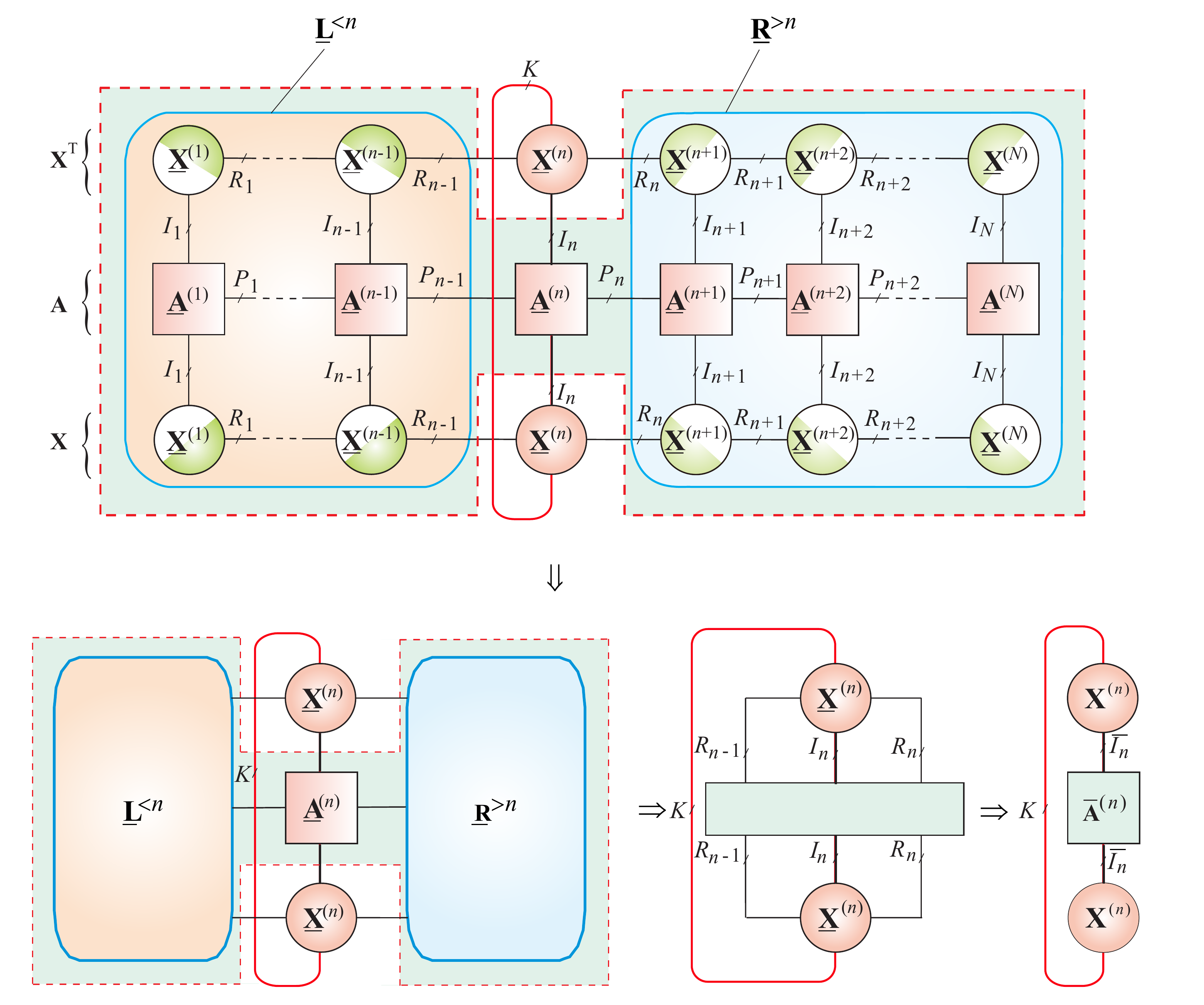}
\end{center}
\caption{Tensor network for the optimization problem in (\ref{Trace1})  and its contraction to reduce the scale of the problem. The network represents a cost function,  the minimization of which  allows us to compute in the TT format the $K$  eigenvectors corresponding to $K$ extreme eigenvalues.
The extreme eigenvalues are computed as $\mbi \Lambda = \bX^{(n) \;\text{T}} \; \overline \bA^{(n)} \, \bX^{(n)}$.}
\label{Fig:EIG1K}
\end{figure}

 Instead of representing each eigenvector $\bx_k$ individually in the TT format, we can parameterize them jointly in a block TT format (see Figure~\ref{Fig:EIG1KBTT})
  \citep{NRG-DMRG12,dolgovEIG2013}.  In general, the block TT format is suitable for approximate representation of tall and skinny  structured matrices using almost the same number of data samples as for representation a single vector/matrix.
%
Within the block-$n$ TT, all but one cores are 3rd-order tensors, the remaining core is a  4th-order tensor and has an  additional physical index $K$ which represents the number of eigenvectors, as shown
  in Figure~\ref{Fig:EIG1K}.
 The  position of this 4th-order core  $\underline \bX^{(n)}$
which carries the index $K$  is not fixed; during  sequential optimization, it is moved   backward and forward
 from position 1 to $N$ \citep{Holtz-TT-2012,dolgovEIG2013}.
 Observe that this 4th-order  core $\underline \bX^{(n)} =\left[\bX_k^{(n)}\right]^K_{k=1} = \left[\bX^{(n,k)}\right] \in \Real^{R_{n-1} \times I_n \times R_n \times K}$ comprises the $K$ different eigenvectors, while all remaining
 cores retain their original structure.
During the optimization process, the neighbouring
core tensors need to be appropriately reshaped, as explained in detail in Algorithm~\ref{alg:blockTTL} and Algorithm~\ref{alg:blockTTR}. 

\begin{figure}[t]
\begin{center}
\includegraphics[width=10.4cm]{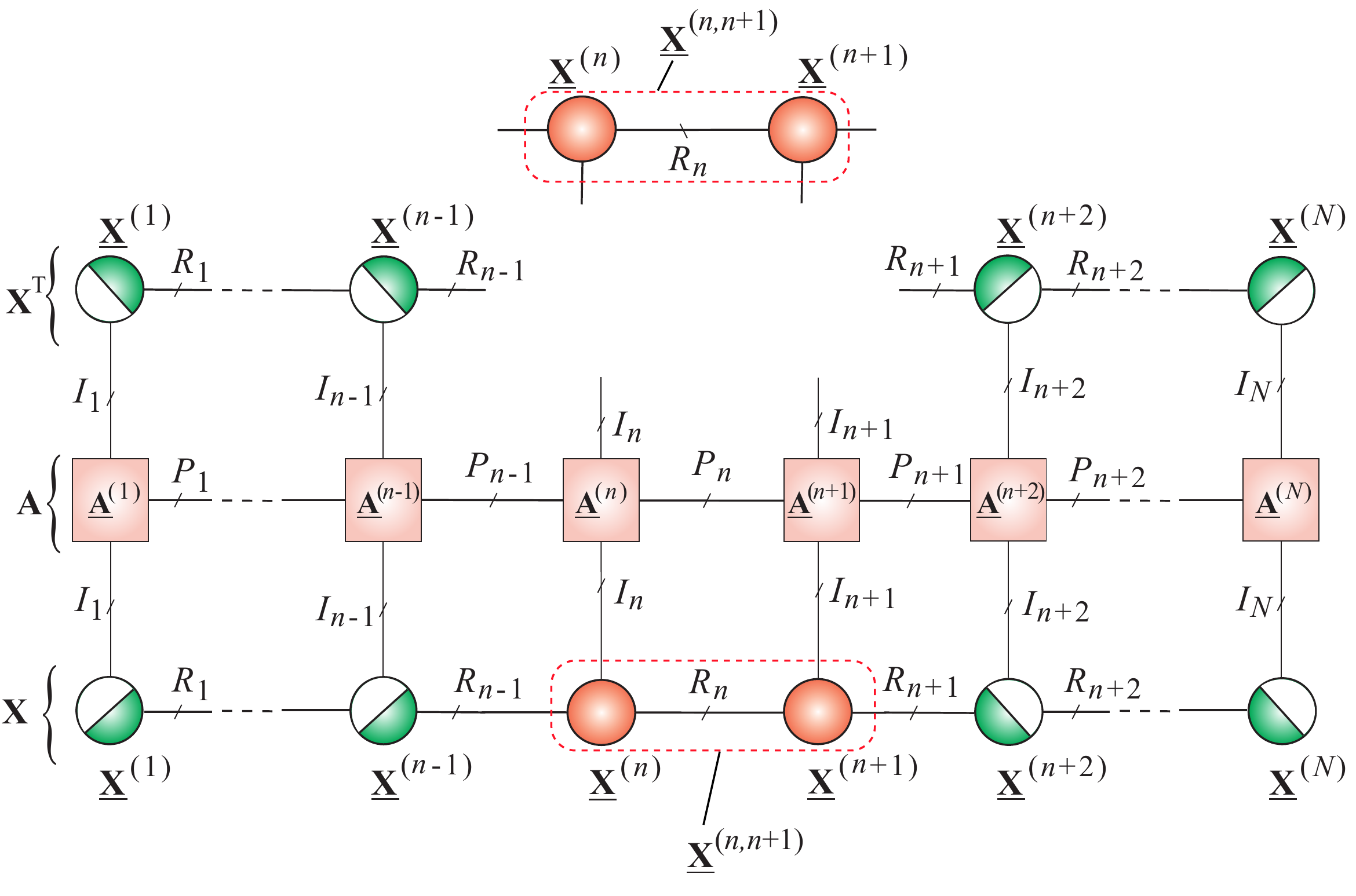}
\end{center}
\caption{Tensor network for the optimization problem in (\ref{Trace1})  and its contraction using the MALS (two-site DMRG2) approach.  The network  is the same and represents the same cost function as for the ALS but the optimization procedure is slightly different.  In this approach, we  merge two neighboring cores, optimize them and then split them back to the original cores via the SVD and estimate the optimal rank $R_n$ for a prescribed approximation accuracy.}
\label{Fig:MALS-EVD}
\end{figure}

 When using the block-$n$ TT model to represent the  $K$ mutually orthogonal vectors,
 the matrix frame equation takes a slightly modified form
\be
\bX= \bX_{\neq n} \; \bX^{(n)} \in \Real^{I \times K}, \quad n=1,2,\ldots,N,
\ee
where  $\bX^{(n)}= \bX^{(n)}_{<3>} = \bX^{(n)\;\text{T}}_{(4)} \in \Real^{R_{n-1} I_n R_n \times K}$.

Hence, we can express the trace in (\ref{Trace1}), as follows
\be
\tr(\bX^{\text{T}} \bA \bX) &=& \tr((\bX_{\neq n} \bX^{(n)})^{\text{T}} \bA \bX_{\neq n} \bX^{(n)}) \notag \\
&=&\tr(( \bX^{(n)})^{\text{T}} [\bX_{\neq n}^{\text{T}} \bA \bX_{\neq n}] \bX^{(n)})\notag \\
&=& \tr(\bX^{(n)\,\text{T}} \; \overline \bA^{(n)} \, \bX^{(n)}),
\label{eq:EVDc}
\ee
where $\overline \bA^{(n)} = \bX_{\neq n}^{\text{T}} \bA \bX_{\neq n} $.

\begin{algorithm} 
{\small
\caption{\textbf{Transformation of the block-$n$ TT  into the block-$(n+1)$ TT  \citep{dolgovEIG2013,KressnerEIG2014}}}
\label{alg:blockTTL}
\begin{algorithmic}[1]
\REQUIRE $N$th-order tensor, $\underline \bX$, in block-$n$ TT format with $n<N$
\ENSURE The tensor $\underline \bX$ in block-$(n+1)$  TT format,  with new  cores $\underline \bX^{(n)} \in \Real^{R_{n-1} \; \times I_n  \; \times \; R_n} $
and $\underline \bX^{(n+1)} \in \Real^{R_{n} \; \times \; I_{n+1} \; \times \; R_{n+1} \; \times K}$
%
    \STATE Reshape  the core $\underline \bX^{(n)} = \left[\underline \bX_k^{(n)}\right]^K_{k=1} = \left[\underline \bX^{(n,k)}\right] \in \Real^{R_{n-1} \times I_n \times R_n \times K}$ \\
    into  a 3rd-order tensor
      $\underline \bX_L^{(n)} \in \Real^{R_{n-1}  \; I_n \; \times \;  R_n \;\times K}$
     \STATE Perform a minimum rank decomposition  (using, e.g., QR or SVD)\\
     $\underline \bX_L^{(n)}(:,:,k) = \bQ \; \bP_k \in \Real^{R_{n-1} \, I_n  \; \times  \; R_{n}}, \;\; k=1,\ldots, K $  \\
      where
     $\bQ \in\Real^{R_{n-1} I_n \times R}, \qquad \bP_k \in \Real^{R \; \times \; R_{n}}$
     \STATE  Update the rank $R_n \leftarrow R$
    \STATE  Update new cores \\ $\bX^{(n)}_L \leftarrow \bQ, \;\;\; \underline \bX_{R}^{(n+1)}(:,:,k)  \leftarrow \bP_k \bX_{R}^{(n+1)} \in \Real^{R_n \times I_{n+1} R_{n+1}}, \;\; \forall k$
\RETURN  $\underline \bX= \llangle \underline \bX^{(1)}, \ldots , \underline \bX^{(n)}, \underline \bX^{(n+1)}, \ldots, \underline \bX^{(N)} \rrangle $
\end{algorithmic}
}
\end{algorithm}

\begin{algorithm}
{\small
\caption{\textbf{Transformation of the block-$n$ TT into the block-$(n-1)$ TT \citep{KressnerEIG2014}}}
\label{alg:blockTTR}
\begin{algorithmic}[1]
\REQUIRE $N$th-order tensor, $\underline \bX$, in the block-$n$ TT format with $n>1$
\ENSURE The tensor $\underline \bX$ in block-$(n-1)$  TT format,  with new cores $\underline \bX^{(n-1)}
\in \Real^{R_{n-2} \; \times \; I_{n-1} \; \times R_{n-1} \; \times \; K }$ and $\underline \bX^{(n)} \in \Real^{R_{n-1} \; \times \; I_n \times \; R_n}$
    \STATE Reshape  the core $\underline \bX^{(n)}  = \left[\underline \bX^{(n,k)}\right] \in \Real^{R_{n-1} \; \times \; I_n  \; \times \; R_n \; \times \; K}$ \\
    into a 3rd-order tensor
     for the TT core $\underline \bX_R^{(n)} \in \Real^{R_{n-1} \; \times  \; I_n R_n \times K}$
     \STATE Perform a minimum rank decomposition (using, e.g., RQ/LQ) \\
      $\underline \bX_R^{(n)}(:,:,k) = \bQ_k \bP \in \Real^{R_{n-1}  \; \times \; I_n R_{n}}, \;\; k=1,\ldots, K$  \\
      where $\bQ_k \in\Real^{R_{n-1} \; \times \; R}, \qquad \bP \in \Real^{R \; \times \; I_n R_n}$
     \STATE  Update the rank $R_{n-1} \leftarrow R$
    \STATE  Update new cores \\
    $\bX^{(n)}_{R} \leftarrow \bP, \;\;\; \underline \bX_{L}^{(n-1)}(:,:,k)  \leftarrow \bX_{L}^{(n-1)} \bQ_k \in \Real^{R_{n-2} I_{n-1} \;  \times \; R_{n-1}}, \;\; \forall k$
\RETURN  $\underline \bX= \llangle \underline \bX^{(1)}, \ldots , \underline \bX^{(n-1)}, \underline \bX^{(n)}, \ldots, \underline \bX^{(N)} \rrangle $
\end{algorithmic}
}
\end{algorithm}

\begin{algorithm}[t]
{\small
\caption{\textbf{One full sweep of the ALS algorithm for symmetric EVD  \citep{DolgovAMEN2014}}}
\label{alg:EVDALS}
\begin{algorithmic}[1]
\REQUIRE A symmetric matrix $\bA \in \Real^{I \times I}$, and initial guesses for \\
$\bX \in \Real^{I_1 I_2 \cdots I_N \times K}$ in block-1 TT format and with right orthogonal \\cores $\underline \bX^{(2)}, \ldots, \underline \bX^{(N)} $
\ENSURE Approximative solution $\bX$ in the TT format
\FOR{$n=1$ to $N-1$}
    \STATE Perform  tensor network contractions (Figure~\ref{Fig:EIG1K}) and solve a \\reduced  EVD problem (\ref{eq:EVDc}),
     for the TT core $\underline \bX^{(n)}$
   \STATE Transform block-$n$ TT to block-$(n+1)$ TT by applying \\
   Algorithm~\ref{alg:blockTTL},  such that the updated core $\bX^{(n)}$ is left-orthogonal
\ENDFOR
\FOR{$n=N$ to $2$}
    \STATE Perform  tensors network contraction and solve a reduced EVD problem (\ref{eq:EVDc})
     for the TT core $\underline \bX^{(n)}$
   \STATE Transform block-$n$ TT to block-$(n-1$) TT by applying \\ Algorithm~\ref{alg:blockTTR}, so that the updated core $\bX^{(n)}$ is right-orthogonal
\ENDFOR
\RETURN  $\underline \bX= \llangle \underline \bX^{(1)}, \underline \bX^{(2)}, \ldots , \underline \bX^{(N)} \rrangle $
\end{algorithmic}
}
\end{algorithm}

Assuming that the frame matrices have orthogonal columns,
 the optimization  problem in (\ref{Trace1}) can be converted into a set of linked
 relatively small-scale  optimization problems
\be
\min_{\bX^{(n)}} \tr( \bX^{(n) \,\text{T}} \; \overline \bA^{(n)} \, \bX^{(n)}),\quad \mbox{s.t.} \;\; \bX^{(n) \;\text{T}} \bX^{(n)}=\bI_K
\ee
for $n=1,2,\ldots,N$, where $\overline \bA^{(n)}$ is computed
recursively by  tensor network contractions, as illustrated in Figure~\ref{Fig:EIG1K}.
In other words, the EVD problem can be solved iteratively via optimized recursive contraction of the TT  network, whereby the active block (i.e., currently optimized core $\underline \bX^{(n)}=  \left[\bX^{(n,k)}\right]$) is sequentially optimized  by sweeping from left to right and back from right to left, and so on, until convergence  as illustrated in Algorithm~\ref{alg:EVDALS} \citep{dolgovEIG2013,KressnerEIG2014}.

During the optimization,  the global orthogonality constraint, $\bX^{\text{T}} \bX =\bI_K$, is equivalent to the set of local orthogonality
constraints, $(\bX^{(n)})^{\text{T}} \bX^{(n)}=\bI_K,\; \forall n$, so that due to left- and right-orthogonality of the TT cores, we can write
\be
\bX^{\text{T}} \bX &=& \bX^{(n)\,T } \; \bX^{\text{T}}_{\neq n} \; \bX_{\neq n} \bX^{(n)} \\
&=& \bX^{(n)\,\text{T}} \bX^{(n)}, \;\; \forall n. \nonumber
\ee
which allows for a dramatic reduction in computational complexity.

 \subsection{Modified ALS (MALS/DMRG2) for Symmetric EVD}

In order to accelerate the convergence speed and improve performance of ALS/DMRG1 methods,  the MALS (two-site DMRG2) scheme optimizes, instead of only one, simultaneously two neighbouring TT-cores at each micro-iteration,  while the other TT-cores remain fixed (assumed to be known). This reduces a large-scale optimization problem
 to a set of  smaller-scale optimization
problems,  although  the  MALS   sub-problems are essentially of bigger size than those for the ALS.

After  local optimization (micro-iteration) on a  merged core  tensor, $\underline \bX^{(n,n+1)} =\underline \bX^{(n)} \times^1  \underline \bX^{(n+1)}$,
 it  is in the next step decomposed  (factorized) into two separate TT-cores via a $\delta$-truncated SVD  of the unfolded  matrix,
$\bX^{(n,n+1)}_{<2>}  \in \Real^{R_{n-1} I_n  \; \times \; I_{n+1} R_{n+1}}$, as
\be
[\bU, \bS, \bV] = \mbox{tSVD} \left(\bX^{(n,n+1)}_{<2>}\right),
\ee
which allows us to both update   the   cores, $\bX^{(n)}_{<2>} =\bU$ and  $\bX^{(n+1)}_{<1>} = \bS \bV^{\text{T}}$,
and   to  estimate  the  TT rank between the cores as $R_n = \min(\mbox{size} (\bU,2), R_{\max})$.

In this way, the TT-ranks can be adaptively determined during the iteration process, while the TT-cores
can be left- or right-orthogonalized, as desired. The truncation parameter, $\delta >0$, can be selected heuristically or adaptively
(see e.g. \citep{Oseledets-Dolgov-lin-syst12,Lee-SIMAX-SVD}).

The use of truncated  SVD  allows us to:
(1) Estimate the optimal rank, $R_n$, and to
(2) quantify the error
associated with this splitting.
In other words, the power of MALS stems from the splitting process which allows for a
 dynamic control of  TT ranks during the iteration process, as well as from the ability to control the level of error
 of  low-rank TN approximations. 

\subsection{Alternating Minimal Energy Method  for EVD (EVAMEn)}
\subsectionmark{Alternating Minimal Energy Method  for EVD}
\label{Sect:Amen-Eigen}

Recall that the standard ALS method cannot adaptively change the TT ranks during the iteration process, whereas the MALS method achieves this in intermediate steps, through the  merging and then  splitting of two neighboring TT cores, as shown in the previous section.
The  cost function $J(\underline \bX)$ is not convex with respect to all the elements of the core tensors of $\underline \bX$ and as consequence the  ALS is prone to getting stuck in local minima and suffers from a relatively slow convergence speed.
On the other hand, the intermediate steps of the MALS  help  to alleviate these problems and improve the convergence speed, but there is still not guarantee  that the algorithm will converge  to a global  minimum.

The Alternating Minimal Energy  (AMEn) algorithm aims to  avoid the  problem of convergence to non-global minimum  by  exploiting the information about the gradient of a cost function or information about the value of a current residual by ``enriching'' the TT cores with additional information during the iteration process \citep{DolgovAMEN2014}.
In that sense, the AMEn can be considered as a subspace correction technique.
Such an enrichment was  efficiently implemented first by  Dolgov and Savostyanov to solve symmetric as well as nonsymmetric linear systems
\citep{Dolgov2016fast} and was later extended to symmetric EVD in \citep{KressnerEIG2014} and \citep{Hubig15}.


\begin{figure}[p]
(a)
\vspace{-0.4cm}
\begin{center}
\includegraphics[width=10.3cm]{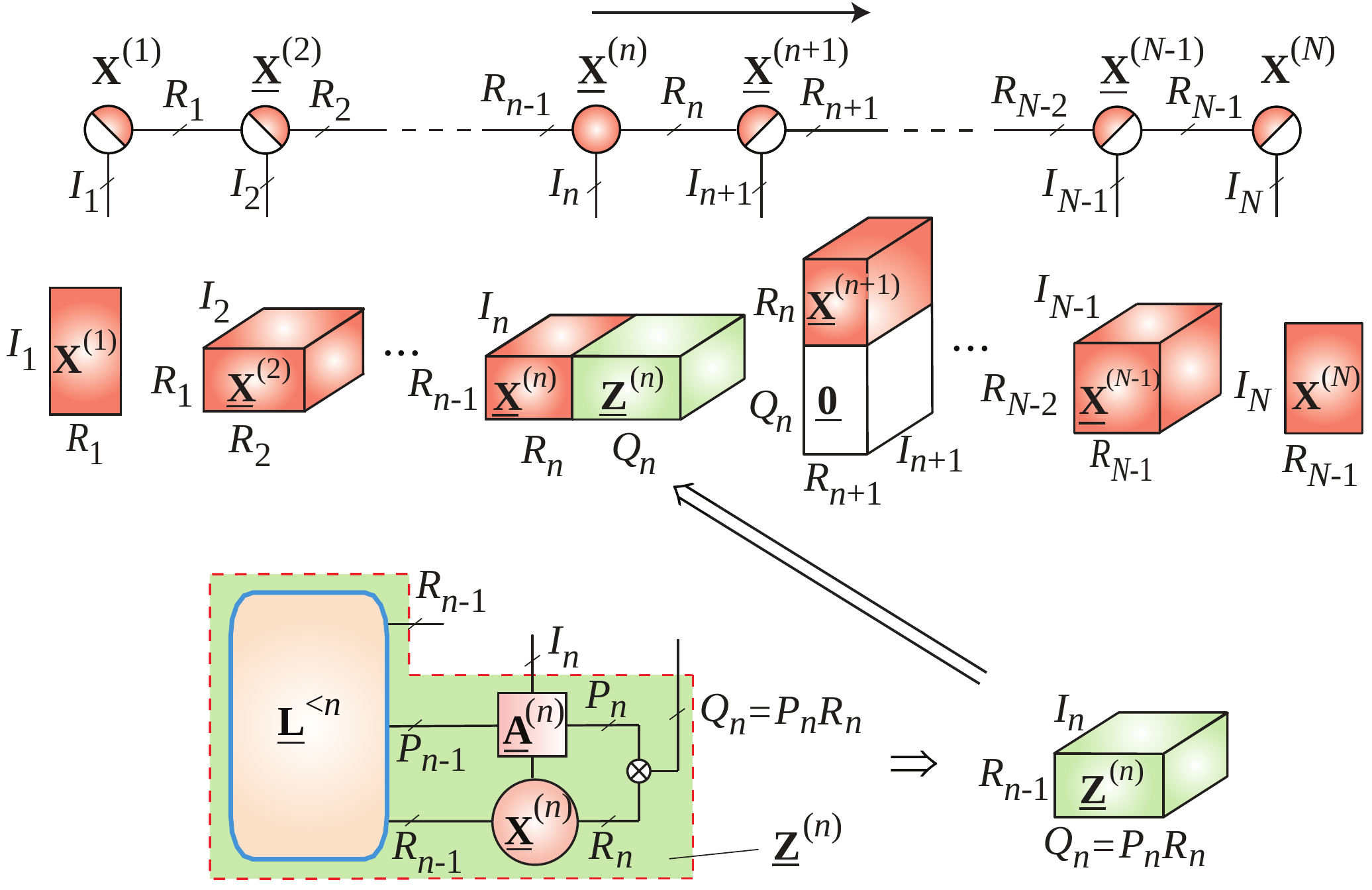}
\end{center}
(b)
\begin{center}
\includegraphics[width=10.3cm]{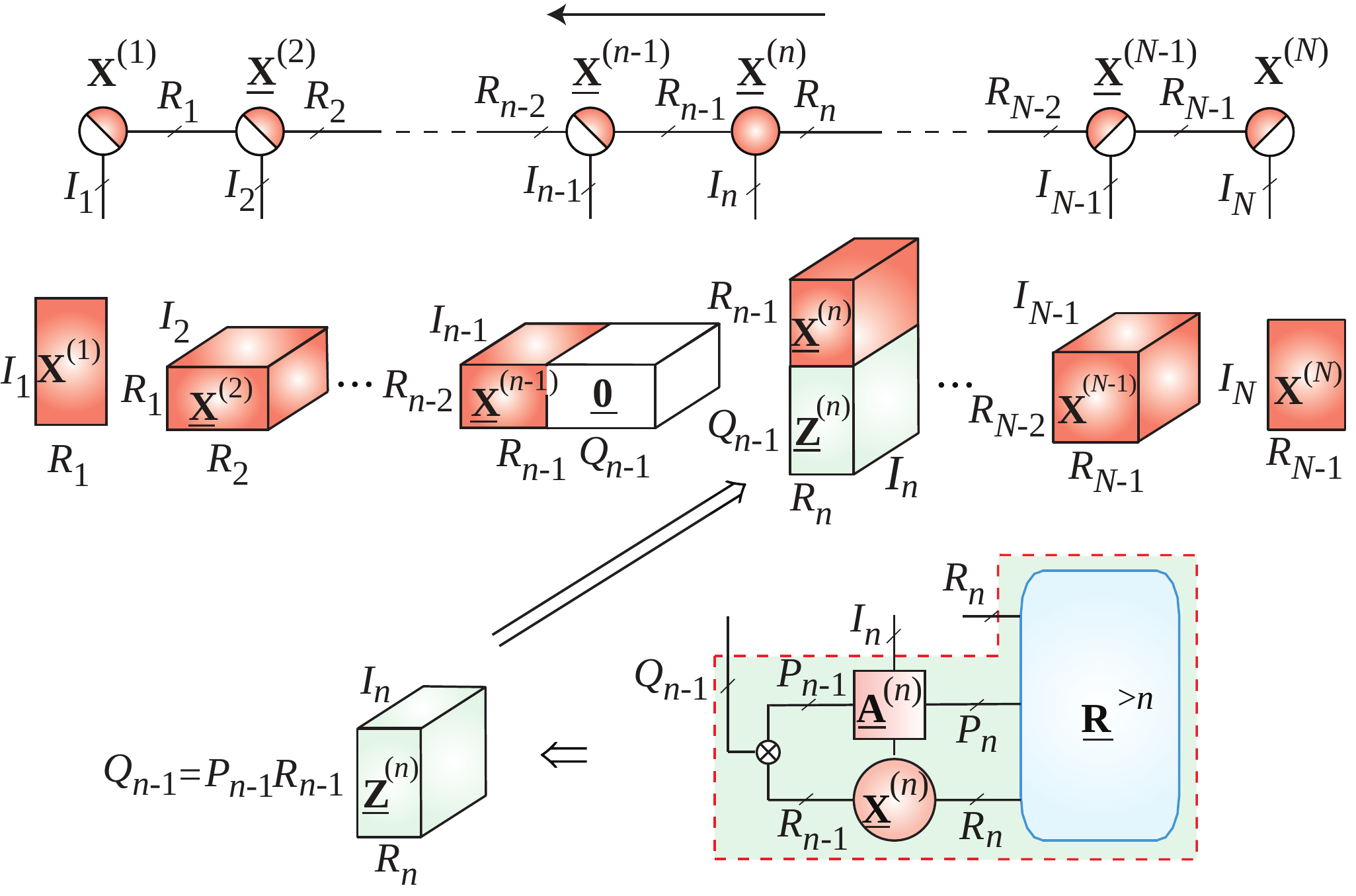}
\end{center}
\caption{Core enrichment step in the AMEn algorithm for the EVD problem, during: (a) Left-to-right
half-sweep ($n=1,2,\ldots,N-1$) and (b) Right-to-left half-sweep ($n=N,N-1,\ldots,2$) of the TT network (see also Figure~\ref{Fig:EIG1}). At the $n$th micro-iteration, only one core tensor, $\underline \bX^{(n)}$, is concatenated (enriched) with a core tensor $\underline \bZ^{(n)}$ obtained from the residual, while its neighboring core tensor is formally concatenated with a zero tensor. For simplicity, this is illustrated for $K=1$.}
\label{Fig:AMEn_Step_Enrichment}
\end{figure}

Similar to the ALS method, at each iteration step, the AMEn algorithm  updates a single core tensor (see Figure~\ref{Fig:AMEn_Step_Enrichment}). Then, it concatenates the updated TT core $\underline \bX^{(n)}$ with a core tensor  $\underline \bZ^{(n)}$ obtained from the residual vector.
At each micro-iteration, only one core tensor of the solution is \textit{enriched} by the core tensor computed based on the gradient vector.
By concatenating two core tensors $\underline \bX^{(n)}$ and $\underline \bZ^{(n)}$, the AMEn can achieve  global convergence,
while maintaining the computational and storage complexities  as low as those of the standard  ALS \citep{DolgovAMEN2014,Dolgov15CME}.

The concatenation step of the AMEn is also called the (core) enrichment, the basis expansion, or the local subspace correction. The concept was proposed by \citet{White2005correction} as a corrected one-side DMRG1 method, while
\citep{DolgovAMEN2014} proposed a significantly improved AMEn algorithm for solving large scale systems of linear equations together with theoretical convergence analysis. \citep{KressnerEIG2014} and \citep{Hubig15}  developed AMEn type methods for solving large scale eigenvalue problems.


In this case,  the goal is to compute  a single or only a few ($K \geq 1$) extreme eigenvalues and corresponding eigenvectors of a very large symmetric matrix $\bA \in \Real^{I\times I}$, with $I=I_1 I_2\cdots I_N$. To this end, we exploit a cost function in the form of the block Rayleigh quotient, $J(\bX) = \text{tr}( \bX^{\text{T}}\bA\bX) $, where the matrix $\bX = [\bx_1,\bx_2,\ldots,\bx_K] \in \Real^{I\times K}$ is assumed to be represented in a block TT format.

In the AMEn algorithm for the EVD, each core tensor, $\underline \bX^{(n)} \in \Real^{R_{n-1} \times I_n \times R_n \times K}$ of $\bX$, is updated at each micro-iteration, similarly to the ALS algorithm. 
Next, a core tensor $\underline \bZ^{(n)} \in \Real^{R_{n-1} \times I_n \times Q_n \times K}$ is constructed and  is concatenated
with the core $\underline \bX^{(n)}$ as
	\be
	\label{updated_evamen_xn}
	\widetilde{\underline \bX}^{(n)} \leftarrow \underline \bX^{(n)} \boxplus_3 \underline \bZ^{(n)}
	\in \Real^{R_{n-1} \times I_n \times (R_n + Q_n) \times K},
	\ee
or equivalently using  tensor slices
\be \widetilde {\bX}^{(n)}_{:,i_n,:,k} \leftarrow [\bX^{(n)}_{:,i_n,:,k},\;  \bZ^{(n)}_{:,i_n,:,k}] \in \Real^{R_{n-1} \times (R_n + Q_n)}, \;\;
\forall i_n, \forall k. \notag
\ee
For the consistency of sizes of TT cores, the neighboring core tensor is concatenated with a tensor consisting all zero entries, to yield
\be
\widetilde \bX^{(n+1)}_{:, i_{n+1}, :} \leftarrow \Bigl[ \begin{smallmatrix} \bX^{(n+1)}_{:, i_{n+1}, :}  \\  {\bf 0} \end{smallmatrix} \Bigr]
\in \Real^{(R_n+Q_n)  \times R_{n+1}}. \notag
\ee
Algorithm~\ref{alg:Eig_AMEN} describes the AMEn  for computing several extreme eigenvalues, $\bLambda = \text{diag}(\lambda_1,\ldots, \lambda_K)$,  and the corresponding eigenvectors, $\bX$.

\begin{algorithm}[htp]
\caption{\textbf{One full sweep of the AMEn algorithm for EVD (EVAMEn) \citep{KressnerEIG2014,Hubig15}}}
\label{alg:Eig_AMEN}
\begin{algorithmic}[1]
\REQUIRE A symmetric matrix $\bA \in \Real^{I_1 \cdots I_N \times I_1 \cdots I_N}$, and an initial \\  guess for $\bX \in \Real^{I_1 \cdots I_N \times K}$ in a block-1 TT format with  \\ right-orthogonal cores $\underline \bX^{(2)}, \ldots, \underline \bX^{(N)}$

\ENSURE Approximate solution $\bX$ in the TT format

\FOR{$n=1$ to $N-1$}
    \STATE Perform  tensor network contraction (Figure~\ref{Fig:EIG1K}) and solve\\ the reduced EVD problem (\ref{eq:EVDc}) for the TT core $\underline \bX^{(n)}$

    \STATE Compute the core $\underline \bZ^{(n)}$ defined in \eqref{eqn:enrichment_hubig} by a contraction of \\ the block $\underline \bL^{<n}$ and the cores $\underline \bA^{(n)}$, $\underline \bX^{(n)}$ (See Figure~\ref{Fig:AMEn_Step_Enrichment}(a))

    \STATE Augment the two cores $\underline \bX^{(n)}$ and  $\underline \bX^{(n+1)}$  with $\underline \bZ^{(n)}$ and $\underline {\bf 0}$ of\\ compatible sizes,  to give the slices \\
    $\bX^{(n)}_{:,i_n,:,k} \leftarrow \left[\bX^{(n)}_{:,i_n,:,k}, \bZ^{(n)}_{:,i_n,:,k} \right] $,
    $\qquad \bX^{(n+1)}_{:,i_{n+1},:} \leftarrow \left[ \begin{matrix} \bX_{:,i_{n+1},:}^{(n+1)} \\ {\bf 0}
    \end{matrix}\right]  $
   \STATE Transform block-$n$ TT into block-$(n+1)$ TT by applying \\ Algorithm~\ref{alg:blockTTL}, so that the updated core$\;\underline \bX^{(n)}$\,is left-orthogonal
\ENDFOR
\FOR{$n=N$ to $2$}

    \STATE Perform  tensors network contraction and solve the reduced EVD problem (\ref{eq:EVDc}) for the TT core $\underline \bX^{(n)}$.

    \STATE Compute the core $\underline \bZ^{(n)}$ by a contraction of the block $\underline \bR^{>n}$ \\
    and the cores $\underline \bA^{(n)}$, $\underline \bX^{(n)}$ (See Figure~\ref{Fig:AMEn_Step_Enrichment}(b)).

    \STATE Augment the two cores $\underline \bX^{(n-1)}$ and  $\underline \bX^{(n)}$  with $\underline{\bf 0}$ and $\underline \bZ^{(n)}$ of \\compatible sizes,   to give the slices \\
    $\bX^{(n-1)}_{:,i_{n-1},:} \leftarrow \left[\bX^{(n-1)}_{:,i_{n-1},:},
    {\bf 0} \right] $,
    $\qquad\bX^{(n)}_{:,i_n,:,k} \leftarrow \left[ \begin{matrix} \bX^{(n)}_{:,i_n,:,k} \\ \bZ^{(n)}_{:,i_n,:,k}
    \end{matrix}\right]  $

    \STATE Transform block-$n$ TT into block-$(n-1)$ TT by applying \\ Algorithm~\ref{alg:blockTTR} such that the updated core $\underline \bX^{(n)}$ is right-orthogonal.

\ENDFOR
\RETURN  $\underline \bX = \llangle \underline \bX^{(1)}, \underline \bX^{(2)}, \ldots , \underline \bX^{(N)} \rrangle $
\end{algorithmic}
\end{algorithm}
In principle, the TT cores $\underline \bZ^{(n)}$ are estimated from the
residual{\footnote{In order to improve convergence,  the residual $\bR$ is often defined as $\bR \cong \bM^{-1}(\bA \bX - \bX \bLambda)$, where $\bM$ is a suitably designed preconditioning matrix.}} as
	\be
	\bR \cong \bA \bX - \bX \bLambda
	\in \Real^{I \times K}
	\ee
with  $\bLambda = \bX^{\text{T}}\bA \bX \in \Real^{K \times K}.$

This corresponds to an orthogonal projection of the gradient of $J(\bX) = \text{tr}( \bX^{\text{T}}\bA\bX) $ onto the tangent space of the Stiefel manifold \citep{Absil-book}.
\citep{KressnerEIG2014} formulated the local residual of the MALS method for EVD as
	\begin{equation} \label{eqn:kressner_localcorrection}
	\begin{split}
	\bR_{n,n+1}
	&= \overline{\bA}^{(n,n+1)} \bX^{(n,n+1)} - \bX^{(n,n+1)} \bLambda  \\
	& =\bX_{\neq n,n+1}^{\text{T}} ( \bA\bX - \bX \bLambda)  \in \Real^{R_{n-1} I_n I_{n+1} R_{n+1} \times K} ,
	\end{split}
	\end{equation}	
where
$\overline{\bA}^{(n,n+1)} = \bX_{\neq n,n+1}^{\text{T}} \bA \bX_{\neq n,n+1}$
and $\bX^{(n,n+1)} (:,k) = \text{vec}(\underline \bX^{(n)}_{:,:,:,k} \times^1 \underline \bX^{(n+1)})$, $k=1,\ldots,K.$

Provided that the matrices $\bA$ and $\bX$ are given in a TT format, it is possible to build two core tensors $\underline \bZ^{(n)} \in \Real^{R_{n-1}\times I_n\times Q_n \times K}$ and $\underline \bZ^{(n+1)} \in \Real^{Q_n \times I_{n+1}\times R_{n+1}}$ separately, such that
	\be \label{eqn:kressner_local_coresZ}
 	\bR_{n,n+1}(:,k) \equiv  \text{vec}(\underline \bZ^{(n)}_{:,:,:,k} \times^1 \underline \bZ^{(n+1)}),
 	\ee
where the size $Q_n$ is equal to that of the residual, $Q_n = P_n R_n + R_n$, unless preconditioned.

\looseness=+1
Alternatively, \citep{Hubig15} developed a more efficient algorithm, whereby instead of using the exact residual $\bA\bX-\bX\bLambda$
to compute the enrichment term,  which is computationally expensive, they show that it is sufficient to exploit only the  $\bA\bX$ term.
Note that if $\bR = \bA\bX-\bX\bLambda$, then the column space is preserved as $\text{range}([\bX, \bR]) = \text{range}([\bX, \bA\bX])$. See also the discussion regarding the column space of $\bX_{\neq n+1}$ in Section~\ref{Sect:Amen-LinearSystem}.
Through this the simplified form, the expressions \eqref{eqn:kressner_localcorrection} and \eqref{eqn:kressner_local_coresZ} can be written as\\[-24pt]
	\begin{equation} \label{eqn:enrichment_hubig}
	\bR_{n,n+1} =   \bX_{\neq n,n+1}^{\text{T}} \bA\bX ,
	\qquad
	\bR_{n,n+1}(:,k) \equiv  \text{vec}(\underline \bZ^{(n)}_{:,:,:,k} \times^1 \underline \bZ^{(n+1)})	.
	\end{equation}
\looseness=+1
It is important to note that the TT cores $\underline \bZ^{(n)}$ and $\underline \bZ^{(n+1)}$ can be estimated independently without
 explicitly computing $\overline{\bA}^{(n,n+1)}$ or $\bX^{(n,n+1)}$.
Specifically, $\underline \bZ^{(n)}$ can be computed by a sequential contraction of the block $\underline \bL^{<n}$ and the cores $\underline \bA^{(n)}$ and $\underline \bX^{(n)}$. Figures \ref{Fig:AMEn_Step_Enrichment}(a) and (b) illustrate the core enrichment procedure in the AMEn algorithm during the left-to-right or the right-to-left half-sweep.  In the figure, for simplicity $\underline \bZ^{(n)}$ is computed based on the expression \eqref{eqn:enrichment_hubig} for eigenvalue problems, with $K=1$ eigenvalue  \citep{Hubig15}.
 %
Furthermore, it has been found that even a rough approximation to the residual / gradient for the enrichment $\underline \bZ^{(n)}$ will lead to
a faster convergence of the algorithm \citep{DolgovAMEN2014}.
See also Section~\ref{Sect:Amen-LinearSystem} for  another way of computing TT cores $\underline \bZ^{(n)}$.
%



In general, the core tensor $\underline \bZ^{(n)}$ is  computed in a way dependent on a specific optimization problem and an efficient approximation strategy.


\markright{3.5.\quad TT Networks for SVD}
\section{TT Networks for Tracking a Few Extreme Singular Values and Singular Vectors in SVD}\label{sect:SVD}\markright{\thesection.\quad TT Networks for SVD}


\begin{figure}
(a)
\vspace{-0.5cm}
\begin{center}
\includegraphics[width=11.3cm]{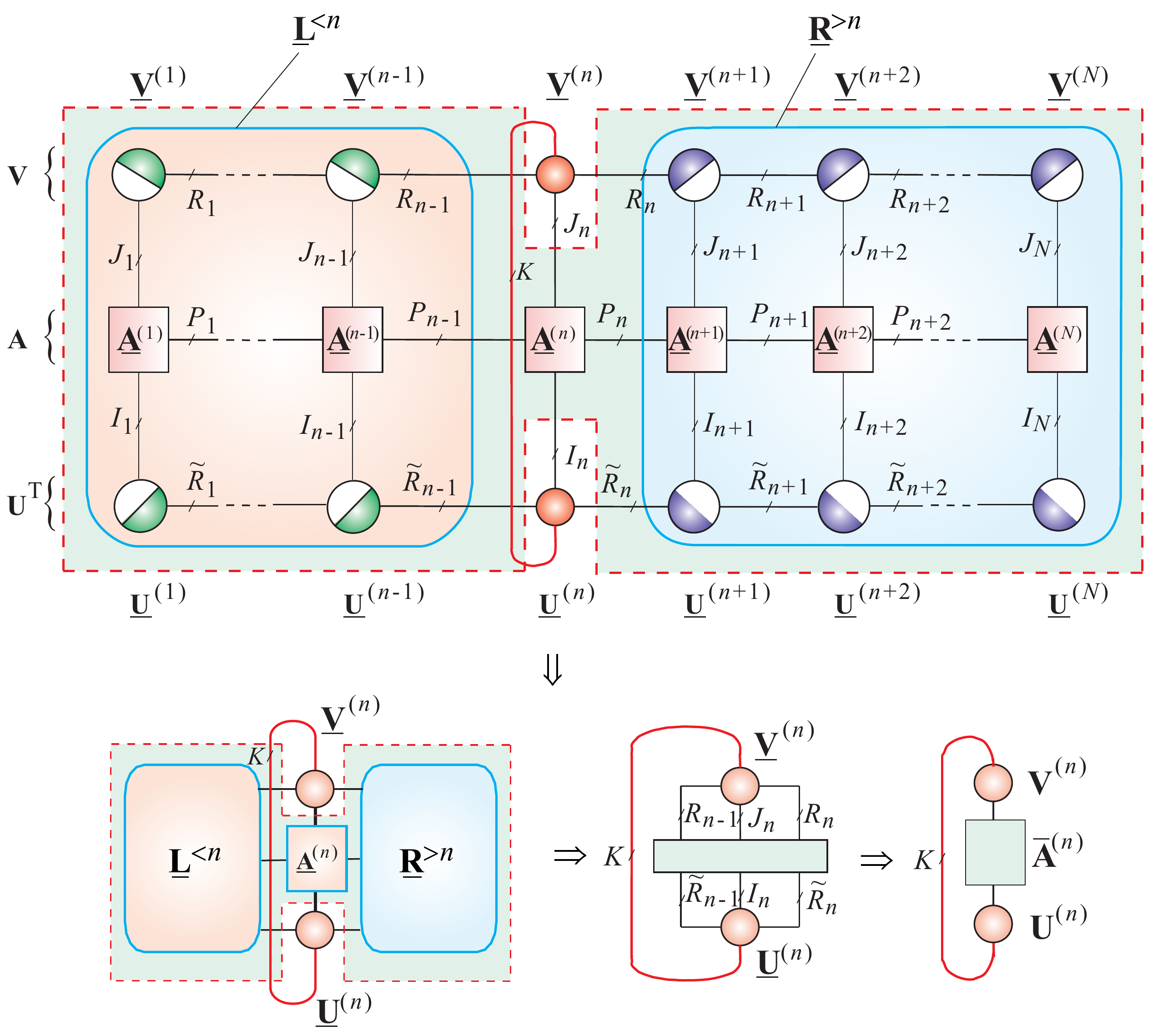}
\end{center}
(b)
\vspace{-0.6cm}
\begin{center}
\includegraphics[width=7.5cm]{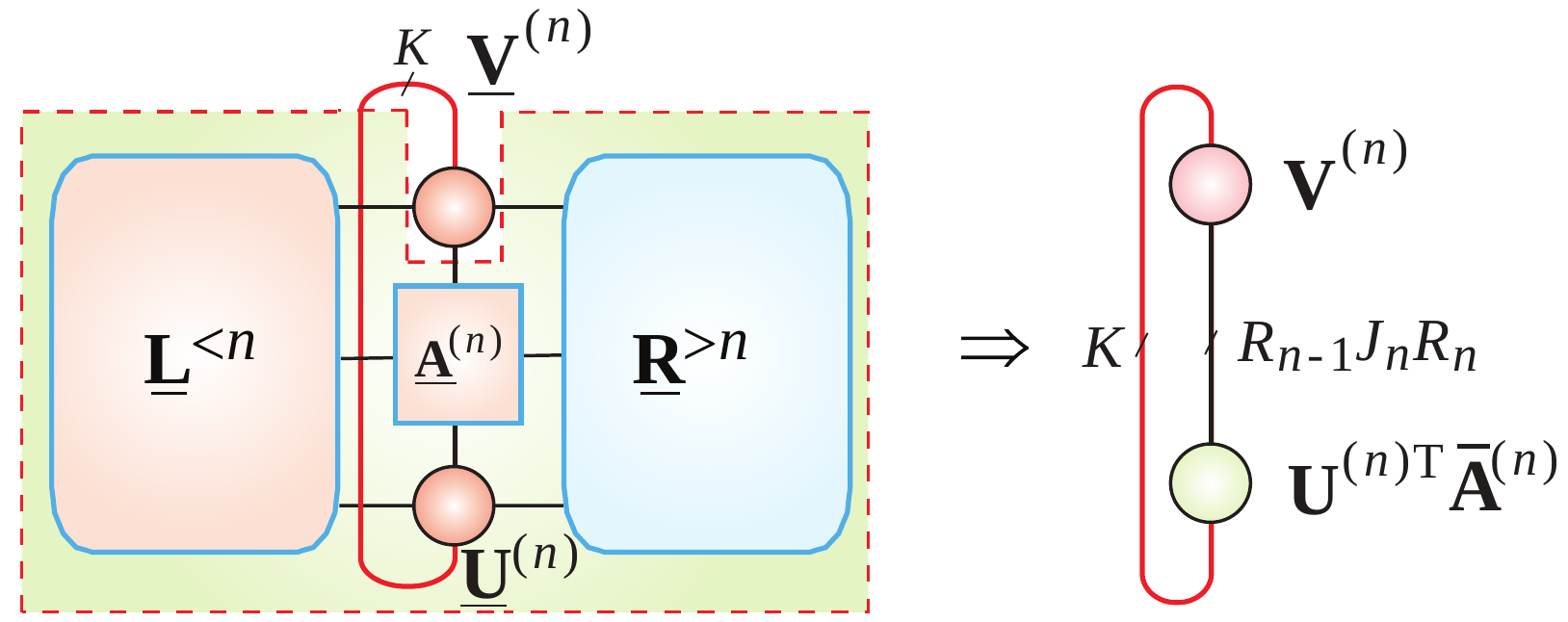}
\end{center}
\caption{Tensor network for computing the $K$  left- and right-singular vectors corresponding to the $K$ extremal singular values. (a) This is achieved via the maximization
of the trace, $\tr(\bU^{\text{T}} \bA \bV)$ or $\tr(\bV^{\text{T}} \bA^{\dagger} \bU)$, subject to the orthogonality constraints, $\bU^{\text{T}} \bU =\bI_k$ and $\bV^{\text{T}}\bV=\bI_K$. Note that the  singular values are computed as $\bS =\bU^{(n) \;\text{T}} \; \overline \bA^{(n)} \, \bV^{(n)}$. (b) In practice, we do not explicitly perform contraction, which leads to the matrix $\bar \bA^{(n)} \in \Real^{R_{n-1} I_n R_{n} \times R_{n-1} J_n R_{n}} $, but employ smaller matrix products $\bU^{(n)\;\text{T}} \bar \bA^{(n)} \in \Real^{K \times R_{n-1} J_n R_{n}}$.  }
\label{Fig:SVDK}
\end{figure}


Similarly to the symmetric EVD problem described in the previous section, the block TT concept can be employed
to compute  only  $K$ largest singular values  and the corresponding
singular vectors of a given matrix $\bA \in \Real^{I \times J}$, by performing the maximization of the following cost function
\citep{Lee-SIMAX-SVD,CichockiSISA,Cichocki2014optim}:
\be
J(\bU, \bV)= \tr(\bU^{\text{T}} \bA \bV), \quad \mbox{s.t.} \quad \bU^{\text{T}} \bU =\bI_K, \;\; \bV^{\text{T}} \bV =\bI_K,
\label{eq:maxKSVD}
\ee
where $\bU \in \Real^{I \times K}$ and  $\bV \in \Real^{J \times K}$.

Conversely, the computation of  $K$ smallest  singular values  and the corresponding left and right
singular vectors of a given matrix, $\bA \in \Real^{I \times J}$, can be  formulated as the following optimization problem:
\be
\max_{\bU,\bV} \tr(\bV^{\text{T}} \bA^{\dagger}  \bU), \quad \mbox{s.t.} \quad \bU^{\text{T}} \bU =\bI_K, \;\; \bV^{\text{T}} \bV =\bI_K,
\label{eq:minKSVD}
\ee
where $\bA^{\dagger} \in \Real^{J \times I}$ is the Moore-Penrose pseudo-inverse of the matrix $\bA$.
Now, after tensorizing the involved huge-scale matrices,  an asymmetric tensor network can be constructed
for the computation of $K$ extreme (minimal or maximal) singular values,
as illustrated in Figure~\ref{Fig:SVDK}
(see \citep{Lee-SIMAX-SVD} for detail and computer simulation experiments).

 \begin{figure}[t]
\begin{center}
\includegraphics[width=11.9cm]{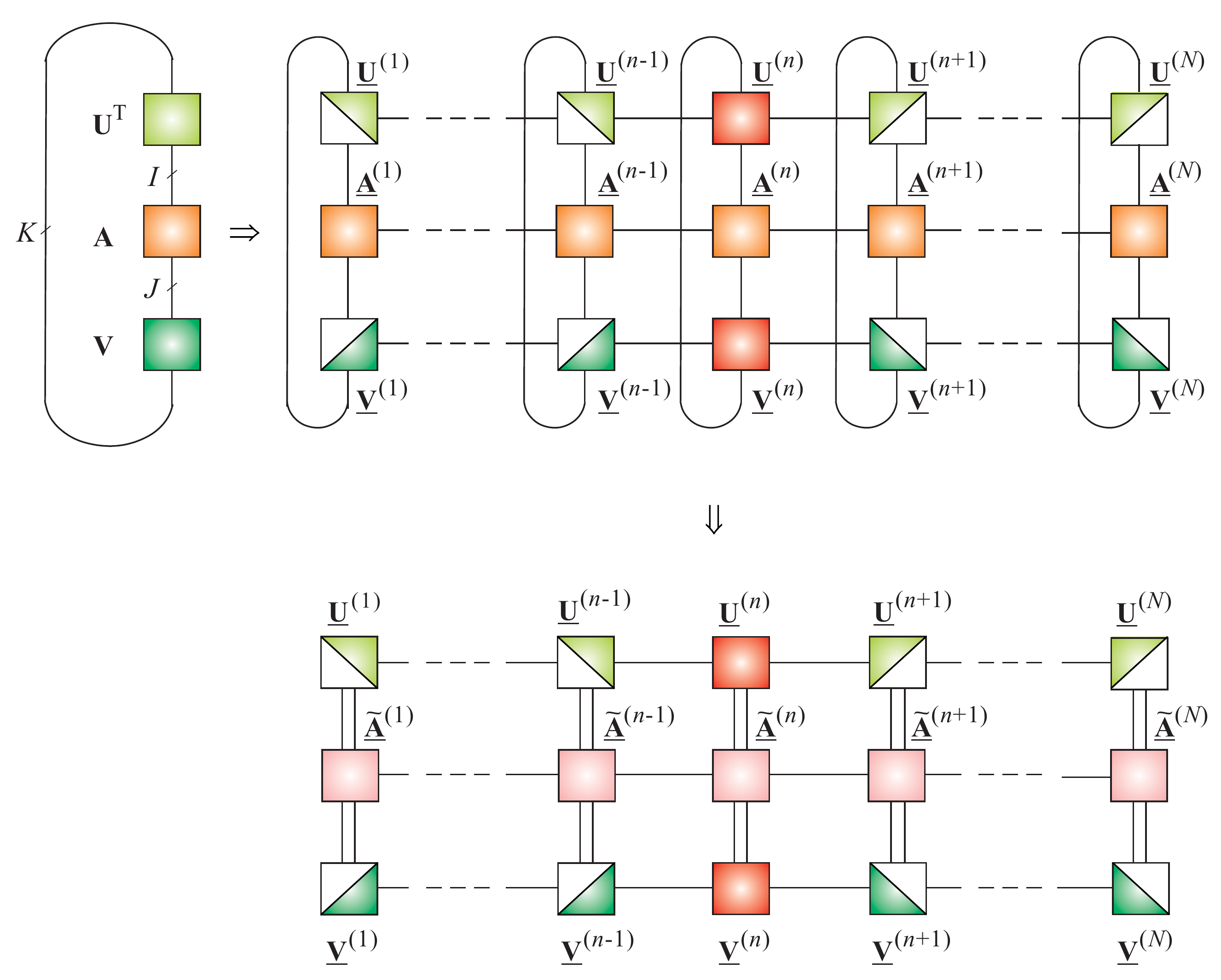}
\end{center}
\caption{MPS/MPO tensor network with the Periodic Boundary Conditions (PBC), for the  computation of a  large number  $K=K_1 K_2 \cdots K_N$
of extreme singular values and the corresponding left- and right-singular vectors (top).  Transformation into the standard  TT formats (MPS/MPO without OBC)  via vectorization, with the cores $\widetilde \bA^{(n)} = \bA^{(n)} \otimes \bI$ (bottom).}
\label{Fig:SVD-PBC}
\end{figure}%

The key idea  behind the solutions of (\ref{eq:maxKSVD}) and (\ref{eq:minKSVD}) is to perform TT core
contractions to reduce the  unfeasible huge-scale optimization problem to
relatively small-scale  optimization sub-problems, as follows:
\begin{itemize}
\item For the problem (\ref{eq:maxKSVD})
 \be
\max_{\bU^{(n)},\bV^{(n)}} \tr((\bU^{(n)})^{\text{T}} \; \overline \bA^{(n)}  \; \bV^{(n)}),
\ee
\item For the problem (\ref{eq:minKSVD})
\be
\max_{\bU^{(n)},\bV^{(n)}} \tr((\bV^{(n)})^{\text{T}}  \; (\overline \bA^{(n)})^{\dagger} \; \bU^{(n)}),
\ee
\end{itemize}
subject to the orthogonality constraints
\be
 (\bU^{(n)})^{\text{T}}  \bU^{(n)}= \bI_K, \quad (\bV^{(n)})^{\text{T}}  \bV^{(n)}= \bI_K, \;\; n=1,2,\ldots,N, \notag
 \ee
where $\bU^{(n)} \in \Real^{\tilde R_{n-1} I_n  \tilde R_n \times K}$ and  $\bV^{(n)} \in \Real^{R_{n-1} J_n R_n \times K}$, while the
  contracted matrices are formally expressed as
 \be
\overline \bA^{(n)}=\bU_{\neq n}^{\text{T}} \bA  \bV_{\neq n} \in \Real^{\tilde R_{n-1} \, I_n \, \tilde R_n \; \times  \;  R_{n-1} \, J_n \, R_n}\\
( \overline \bA^{(n)})^{\dagger}=\bV_{\neq n}^{\text{T}} \bA^{\dagger}  \bU_{\neq n} \in \Real^{ R_{n-1} \, J_n \, R_n  \; \times  \;  \tilde R_{n-1} \, I_n \, \tilde R_n}.
  \ee%
In this way, the  problem is reduced to the  computation of the largest or smallest singular values of a relatively small matrix $\overline \bA^{(n)}$, for which
 any efficient SVD algorithm can be applied.

Note that for a very large number of, say,  thousands or more singular values  $K =K_1K_2 \cdots K_N$,   the block-$n$ TT formats may  not be efficient, and alternative tensor networks should  be employed, as shown in Figure~\ref{Fig:SVD-PBC},  where each  TT core is a 4th-order tensor. One problem with such tensor networks with many
loops is the  computationally complex contraction of core tensors.
However,  all loops can be  easily  eliminated by vectorizing the orthogonal matrices
\be
\bU \in \Real^{I \times K} \rightarrow \bu \in \Real^{IK}, \quad \bV \in \Real^{J \times K} \rightarrow \bu \in \Real^{JK}
\ee
and by rewriting the cost function as (see bottom part of Figure~\ref{Fig:SVD-PBC})
\be
\nonumber
 J(\bU, \bV) &=& \tr(\bU^{\text{T}} \bA \bV) \rightarrow \\
 J(\bu, \bv)&=& \bu^{\text{T}} \widetilde \bA \bv = \bu^{\text{T}} (\bA \otimes \bI_K) \bv,
\ee
where $\widetilde \bA = \bA \otimes \bI_K \in \Real^{IK \times JK}$.

\section{GEVD and Related Problems using TT Networks}
\sectionmark{TT Networks for GEVD}

Table~\ref{Table:LPP} illustrates the rich scope of the Generalized Eigenvalue Decomposition (GEVD), a backbone
of many standard\footnote{In some applications, we need to compute only very few eigenvectors. For example,
in  spectral co-clustering, only one
  eigenvector (called Fiedler eigenvector) needs to be computed, which corresponds to the second smallest generalized eigenvalue
  \citep{Wu-Gleich2016}.}  methods for linear dimensionality reduction, classification, and clustering \citep{Kokiopoulou11,Cunningham2015linear,Benson2015-tensor,GleichPageRank2015,Wu-Gleich2016}.

\begin{table}[h!]
\caption{Cost functions and constraints used in classical feature extraction
(dimensionality reduction) methods that
can be formulated  as the symmetric (generalized) eigenvalue problem, $\min_{\bV} (\bV^{\text{T}} \bX \bC \bX^{\text{T}} \bV)$.
The objective is to find  an (orthogonal) matrix, $\bV$, assuming that data matrices, $\bX$, $\bW, \bD, \bH$, are known.
The symmetric matrix $\bC$ can take different forms, for example, $\bC= {\bf I}-\frac{1}{N} {\bf 11}^{\text{T}}, \; \bC=\bD-\bW,   \; \bC=(\bI - \bW^{\text{T}})(\bI -\bW), \; \bC= \bI-\bH,$ depending on  the method used (for more detail, see e.g. \citep{Kokiopoulou11,Cunningham2015linear}).}
\centering
\renewcommand{\arraystretch}{1.4}
{\footnotesize\tabcolsep 3pt \shadingbox{
\begin{tabular}{@{}l@{}|@{\hspace{.3ex}}c|c}
\hline
Method & Cost Function (min) & Constraints \\
\hline
\minitab[l]{Principal Component \\ Analysis/Multi-Dimensional \\ Scaling (PCA/MDS)} & {$\tr[-{\bf V}^{\text{T}}{\bf X}({\bf I}-\frac{1}{N} {\bf 11}^{\text{T}}){\bf X}^{\text{T}}{\bf V}]$} & ${\bf V}^{\text{T}}{\bf V}={\bf I}$ \\
 & & \\
\minitab[l]{Locally Preserving Projection \\(LPP)} & {$\tr[{\bf V}^{\text{T}}{\bf X}({\bf D}-{\bf W}){\bf X}^{\text{T}}{\bf V}]$} & ${\bf V}^{\text{T}}{\bf XDX}^{\text{T}}{\bf V}={\bf I}$ \\
 & & \\
Orthogonal LPP (OLPP) & {$\tr[{\bf V}^{\text{T}}{\bf X}({\bf {D-W}}){\bf X}^{\text{T}}{\bf V}]$} & ${\bf V}^{\text{T}}{\bf V}={\bf I}$ \\
 & & \\
\minitab[l]{Neighborhood Preserving \\ Projection (NPP)} & {$\tr[{\bf V}^{\text{T}}{\bf X}({\bf I}-{\bf W}^{\text{T}})({\bf I}-{\bf W}){\bf X}^{\text{T}}{\bf V}]$} & ${\bf V}^{\text{T}}{\bf XX}^{\text{T}}{\bf V}={\bf I}$ \\
 & & \\
Orthogonal NPP (ONPP) & {$\tr[{\bf V}^{\text{T}}{\bf X}({\bf I}-{\bf W}^{\text{T}})({\bf I}-{\bf W}){\bf X}^{\text{T}}{\bf V}]$} & ${\bf V}^{\text{T}}{\bf V}={\bf I}$ \\
 & & \\
\minitab[l]{Linear Discriminant Analysis \\(LDA)} & {$\tr[{\bf V}^{\text{T}}{\bf X}({\bf I}-{\bf H}){\bf X}^{\text{T}}{\bf V}]$} & ${\bf V}^{\text{T}}{\bf XX}^{\text{T}}{\bf V}={\bf I}$ \\
 & & \\
\minitab[l]{Spectral Clustering \\(Ratio Cut)} & {$\tr[{\bf V}^{\text{T}}({\bf D}-{\bf W}){\bf V}]$}
& ${\bf V}^{\text{T}}{\bf V}={\bf I}$ \\
& & \\
\minitab[l]{Spectral Clustering \\(Normalized Cut)} & {$\tr[{\bf V}^{\text{T}}({\bf D}-{\bf W})
{\bf V}]$} & ${\bf V}^{\text{T}}{\bf D V}={\bf I}$ \\
\hline
\end{tabular}
}}
\label{Table:LPP}
\vspace{36pt}
\end{table}

For example, the standard PCA,  PLS,  orthonormalized PLS (OPLS)  and CCA
can be formulated as the following EVD/GEVD problems
\be
&&\mbox{PCA :} \qquad \bC_x \bv = \lambda \bv  \\ \nonumber \\
&&\mbox{OPLS :}  \quad \;\; \bC_{xy} \bC^{\text{T}}_{xy} \bv = \lambda \bC_x \bv \\ \nonumber \\
&&\mbox{PLS :}  \qquad  \left[ \begin{matrix} \0 & \bC_{xy} \\ \bC^{\text{T}}_{xy} & \0 \end{matrix} \right]
\left[ \begin{matrix} \bv \\ \bu \end{matrix} \right]
=  \lambda  \left[ \begin{matrix} \bv \\ \bu \end{matrix} \right] \\ \nonumber \\
&&\mbox{CCA :}  \quad \;\; \left[ \begin{matrix} \0 & \bC_{xy} \\ \bC^{\text{T}}_{xy} & \0 \end{matrix} \right]
\left[ \begin{matrix} \bv \\ \bu \end{matrix} \right]
=  \lambda  \left[ \begin{matrix} \bC_x & \mbi \0 \\ \0 & \bC^{\text{T}}_{y} \end{matrix} \right]
\left[ \begin{matrix} \bv \\ \bu \end{matrix} \right],
\ee
where  for given data
matrices $\bX$ and $\bY$, $\bC_x =\bX^{\text{T}} \bX$ and $\bC_y =\bY^{\text{T}} \bY$  are the sample covariance matrices, while $\bC_{xy} =\bX^{\text{T}} \bY$ is the sample cross-covariance matrix.

In general, the  GEVD problem can be formulated as the  following trace optimization (maximization or minimization)  problem  
 \be
\max_{\bV \in \Real^{I \times K}} \tr(\bV^{\text{T}}  \bA  \bV), \quad \mbox{s.t.}  \quad \bV^{\text{T}} \bB \bV = \bI_K,
\label{GEVDK11}
 \ee
 where  the  symmetric matrix,  $\bA \in \Real^{I \times I}$, and symmetric positive-definite matrix,  $\bB \in \Real^{I \times I}$,
 are known, while the objective is to estimate $K$ largest (or the smallest in the case of minimization) eigenvalues and the
 corresponding eigenvectors represented by the orthonormal matrix $\bV \in \Real^{I \times K}$, with $I=I_1 I_2 \cdots I_N$ and $K=K_1 K_2 \cdots K_N$.
 For the estimation of the $K$ smallest eigenvalues, the
 minimization  of the cost function is performed instead of the maximization.

\begin{figure}[t!]
\begin{center}
\includegraphics[width=0.99\textwidth]{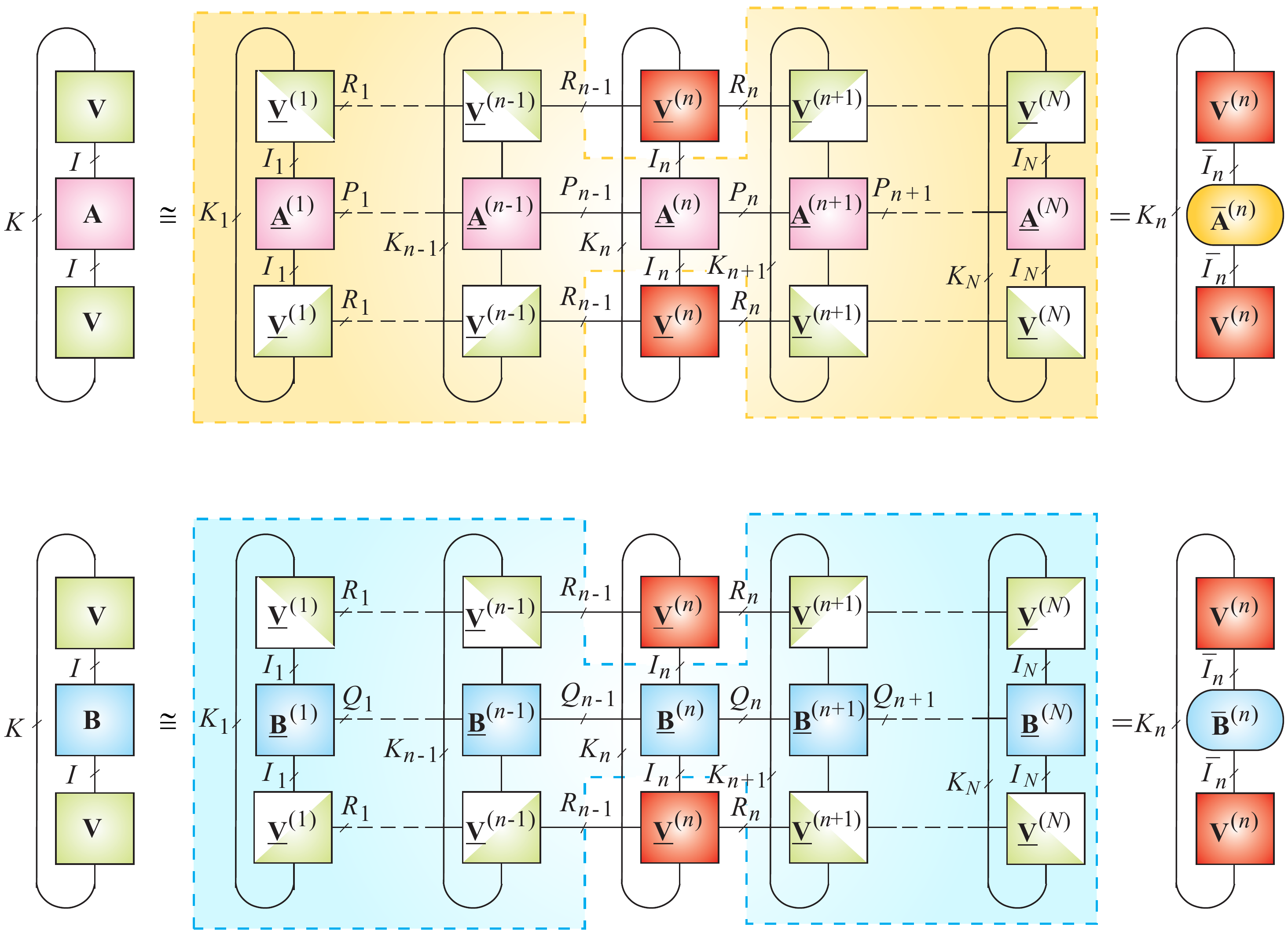}
\end{center}
\caption{Conceptual tensor networks for  trace ratio optimization  in TT formats. Note  that for a moderate  number $K$, all
 the loops  but one  may have a dimension (size) of unity and can thus be  removed.}
\label{Fig:GEVD}
\end{figure}

 It is sometimes more elegant to incorporate the constraint, $\bV^{\text{T}} \bB \bV = \bI_K$, into (\ref{GEVDK11}) into the trace operator, to give the following  optimization problem \citep{Absil-book}
   \be
\min_{\bV \in \Real^{I \times K}} \tr(\bV^{\text{T}}  \bA  \bV(\bV^{\text{T}} \bB \bV )^{-1}).
\label{eq:GEVD1}
 \ee
 or alternatively, to represent GEVD approximately as an unconstrained optimization problem
  \be
\min_{\bV \in \Real^{I \times K}}  \tr(\bV^{\text{T}}  \bA  \bV)+ \gamma \; \| \bV^{\text{T}} \bB \bV - \bI_K\|^2_F,
\label{GEVDK11u}
 \ee
where the parameter $\gamma > 0$ controls the orthogonality level of $\bB$.

 Also, by a change in  the variables,  $\bW=\bB^{1/2} \bV$,  the GEVD reduces to the standard symmetric EVD \citep{Cunningham2015linear}
  \be
\min_{\bW \in \Real^{I \times K}} \tr(\bW^{\text{T}} \bB^{-1/2}  \bA  \bB^{-1/2} \bW), \;\; \mbox{s.t.}  \;\; \bW^{\text{T}}  \bW = \bI_K. \nonumber
\label{EVDK12}
 \ee
 Finally, the GEVD  is closely related but not  equivalent, to the  trace ratio
 optimization problem, given by
   \be
\min_{\bV \in \Real^{I \times K}} \frac{\tr(\bV^{\text{T}}  \bA  \bV)}{\tr (\bV^{\text{T}} \bB \bV)}, \quad \mbox{s.t.} \quad \bV^{\text{T}} \bV = \bI_K.
\label{eq:GEVDRatio}
 \ee
Figure~\ref{Fig:GEVD}  illustrates the trace ratio optimization  for huge-scale matrices, whereby the two  TT networks  represent approximately
 $\tr(\bV^{\text{T}} \bA \bV)$ and  $ \tr(\bV^{\text{T}} \bB \bV)$.
  These sub-networks are simultaneously optimized in the sense that one is  being minimized while the other is being  maximized,  subject to orthogonality constraints.
 The objective is to estimate the  matrix $\bV \in \Real^{I \times K}$ in the  TT format, assuming that huge-scale structured
  matrices $\bA, \bB$ and $\bV$ admit good low-rank TT approximations.

Figure~\ref{Fig:GEVD} illustrates that a recursive contraction of TT cores reduces the trace ratio optimization problem (\ref{eq:GEVDRatio})
to a set of  following smaller scale trace ratio minimizations
 \be
\min_{\bV^{(n)}  \in \Real^{\overline I_n \times K_n}} \frac{\tr((\bV^{(n)})^{\text{T}} \; \overline \bA^{(n)}\; \bV^{(n)})}{\tr((\bV^{(n)})^{\text{T}} \; \overline \bB^{(n)}  \;\bV^{(n)})}, \quad \mbox{s.t.}  \quad  (\bV^{(n)})^{\text{T}} \bV^{(n)}= \bI_{K_n},\nonumber\\
\label{TTGEVDK}
 \ee
where the relatively small-size matrices, $\overline \bA^{(n)}$ and $\overline \bB^{(n)}$, are formally defined as
 \be
  \overline \bA^{(n)}= [\bV_{\neq n}^{\text{T}}  \bA  \bV_{\neq n}]\in \Real^{R_{n-1} I_n R_n \times R_{n-1} I_n R_n}
  \ee
  and
  \be
  \overline \bB^{(n)}= [\bV_{\neq n}^{\text{T}} \bB \bV_{\neq n}] \in \Real^{R_{n-1} I_n R_n \times R_{n-1} I_n R_n}.
   \ee

\subsection{TT Networks for  Canonical Correlation Analysis (CCA)}
\subsectionmark{TT Networks for CCA}

The Canonical Correlation Analysis (CCA) method, introduced by Hotelling in 1936, 
 can be considered as a generalization  of the PCA, and  is
 a classical method for determining the relationship between two sets of variables
 (for modern approaches to CCA, see \citep{Bach_PCCA2005,chu2013sparse,Cunningham2015linear} and references therein).

Given two  zero-mean (i.e., centered) datasets, $\bX \in \Real^{I \times J}$ and
$\bY \in \Real^{L \times J}$,   recorded from the same set of $J$ observations, the CCA seeks linear
combinations of the latent variables in $\bX$ and  $\bY$
that are maximally mutually correlated.
Formally, the classical CCA computes two projection vectors, $\bw_x =\bw^{(1)}_x  \in \Real^I$ and
$\bw_y =\bw^{(1)}_y \in \Real^L$, so as to maximize   the correlation coefficient
\be
\rho =\frac{\bw^{\text{T}}_x \bX \bY^{\text{T}}  \bw_y}{\sqrt{(\bw^{\text{T}}_x \bX \bX^{\text{T}} \bw_x)(\bw^{\text{T}}_y \bY \bY^{\text{T}} \bw_y)}}.
\ee
%
In a similar way,  kernel CCA can be formulated by replacing the inner product of matrices by the kernel matrices,
$\bK_x \in \Real^{J \times J}$ and $\bK_y \in \Real^{J \times J} $, to give
\be
\rho = \max_{\mbi \alpha_x, \mbi\alpha_y} \frac{\mbi \alpha^{\text{T}}_x \bK_x \bK_y \mbi \alpha_y}{\sqrt{(\mbi \alpha^{\text{T}}_x \bK_x \bK_x \mbi \alpha_x)(\mbi \alpha^{\text{T}}_y \bK_y \bK_y \mbi \alpha_y)}}.
\ee


Since the correlation coefficient $\rho$ is invariant to the scaling of the vectors $\bw_x$ and $\bw_y$, the standard CCA can be equivalently formulated as a
constrained optimization problem
\be
&&\max_{\bw_x,\bw_y} \bw_x^{\text{T}} \bX  \bY^{\text{T}} \bw_y  \label{CCA1} \\
 && \mbox{s.t.} \quad \bw_x^{\text{T}} \bX  \bX^{\text{T}} \bw_x=\bw_y^{\text{T}} \bY^{\text{T}}  \bY \bw_y=1.\nonumber
\label{CCA2}
\ee
The vectors $\bt_1=\bX^{\text{T}} \bw_x$ and $\bu_1=\bY^{\text{T}} \bw_y$ are referred to as the canonical variables.
%

%
%

Multiple canonical vectors for the classical CCA can be computed by reformulating the optimization problem in (\ref{CCA1}), as
follows
\be
\max_{\bW_x, \; \bW_y}\tr( \bW^{\text{T}}_x \; \bX \; \bY^{\text{T}} \; \bW_y)
\ee
\[
 \mbox{s.t.}  \quad \bW^{\text{T}}_x \; \bX  \bX^{\text{T}} \; \bW_x =\bI_K, \   \bW^{\text{T}}_y \; \bY  \bY^{\text{T}} \; \bW_y  =\bI_K, \  \bW^{\text{T}}_x \; \bX  \bY^{\text{T}} \; \bW_y =\mbi \Lambda,
\]
where $\bW_x =[\bw_x^{(1)},\bw_x^{(2)},\ldots,\bw_x^{(K)}] \in \Real^{I \times K}$ and
$\bW_y =[\bw_y^{(1)},\bw_y^{(2)},\ldots,\bw_y^{(K)}] \in \Real^{L \times K}$ and $\mbi \Lambda$ is any  diagonal matrix.
However, this may become   computationally very expensive due to orthogonality  constraints  for large-scale matrices.

An alternative approach is to use the orthogonal CCA  model,  which can be formulated as \citep{Cunningham2015linear}
\be
&&\max_{\bW_x,\bW_y} \left( \frac{\tr(\bW^{\text{T}}_x \bX \bY^{\text{T}}  \bW_y)}{\sqrt{\tr(\bW^{\text{T}}_x \bX \bX^{\text{T}} \bW_x) \tr(\bW^{\text{T}}_y \bY \bY^{\text{T}} \bW_y)}}\right) \\
&&\mbox{s.t.} \quad  \quad \bW^{\text{T}}_x  \bW_x = \bI_K, \quad  \bW^{\text{T}}_y \bW_y  =\bI_K. \notag
\label{eq:OCCA1}
\ee
\begin{figure}[t]
\begin{center}
\includegraphics[width=11.0cm]{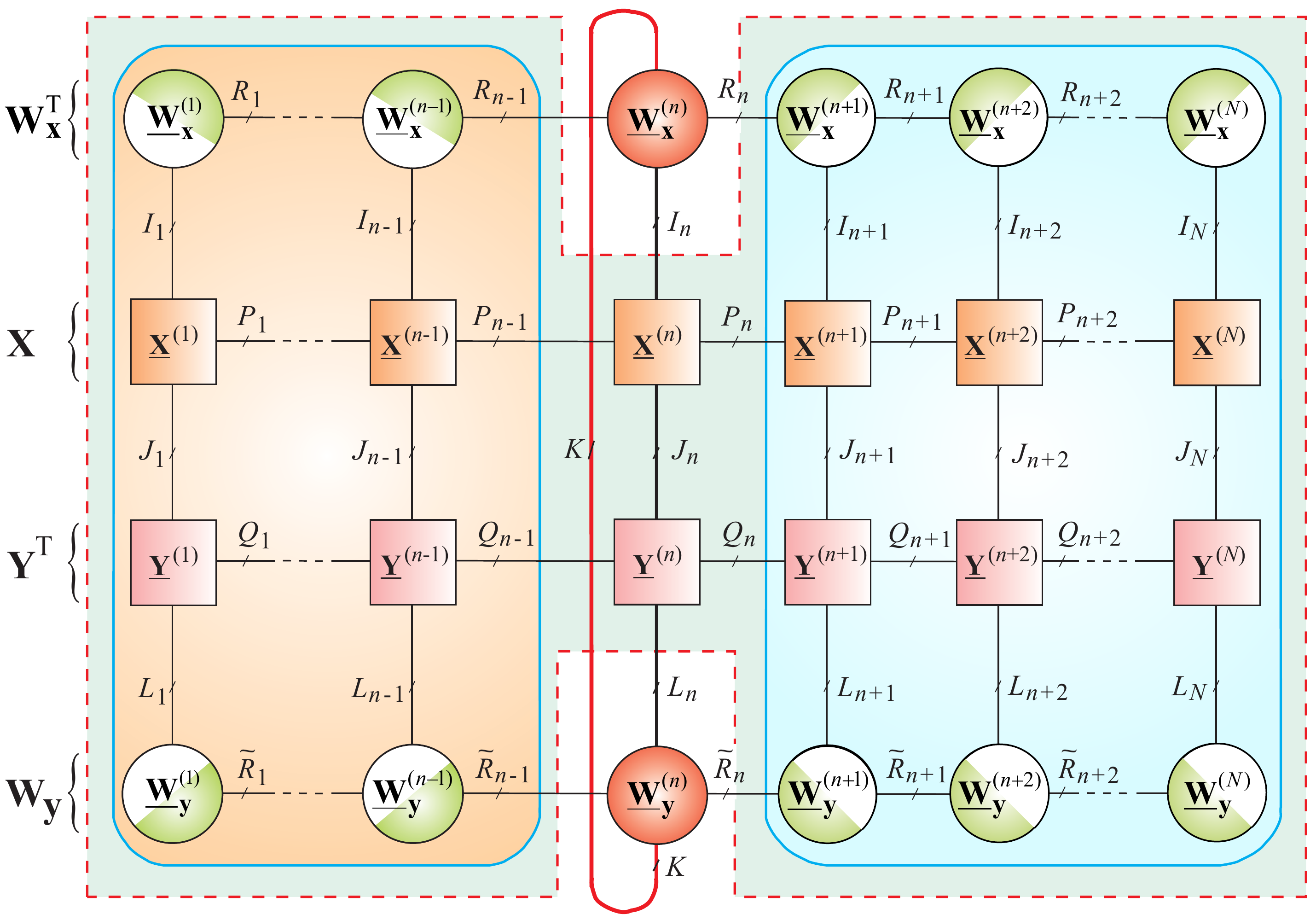}
\end{center}
\caption{A TT sub-network  for  multiple orthogonal CCA components, representing in the distributed form the term $\tr( \bW^{\text{T}}_x \bX \bY^{\text{T}} \bW_y)$. Similar TT sub-networks are constructed for the normalization terms $\tr( \bW^{\text{T}}_x \bX \bX^{\text{T}} \bW_x)$ and $\tr( \bW^{\text{T}}_y \bY \bY^{\text{T}} \bW_y)$, see (\ref{eq:OCCA1}). }
\label{Fig:CCAL}
\end{figure}
This model can be relatively easily implemented using the TT network approach, as illustrated in Figure~\ref{Fig:CCAL}.
By using contractions of TT sub-networks, the huge-scale optimization problem in (\ref{eq:OCCA1}) can be transformed into a set
of smaller optimization problems
\be
&&\max_{\bW^{(n)}_x,\bW^{(n)}_y} \left( \frac{\tr(\bW^{(n)\;\text{T}}_x \bC_{xy}^{(n)}  \bW^{(n)}_y)}{\sqrt{\tr(\bW^{(n)\;\text{T}}_x \bC^{(n)}_{xx} \bW_x) \tr(\bW^{(n)\;\text{T}}_y \bC^{(n)}_{yy} \bW^{(n)}_y)}}\right) \\
&&\mbox{s.t.} \quad  \quad \bW^{(n)\;\text{T}}_x  \bW^{(n)}_x = \bI_K, \quad  \bW^{(n)\;\text{T}}_y \bW^{(n)}_y  = \bI_K,\notag
\label{eq:OCCA2}
\ee
 for $n=1,2,\ldots,N$, where $\bC^{(n)}_{xy} = \bW^{\text{T}}_{x,\;\neq n}\bX \bY^{\text{T}} \bW_{y, \;\neq n}$, $ \;\bC^{(n)}_{xx} = \bW^{\text{T}}_{x,\;\neq n} \bX \bX^{\text{T}} \bW_{x, \;\neq n}$ and
 $\bC^{(n)}_{yy} = \bW^{\text{T}}_{y,\;\neq n}\bY \bY^{\text{T}} \bW_{y, \;\neq n}$.



For  computational tractability of huge-scale CCA problems and for  physical interpretability of latent variables,
it is often useful to impose   sparsity  constraints, to yield the sparse CCA.
Sparsity constraints can be  imposed, for example,  by applying the
$\ell_1$-norm penalty terms to the matrices $\bW_x$ and $\bW_y$, and assuming that the cross-product matrices $ \bX \bX^{\text{T}}$ and  $\bY \bY^{\text{T}}$
 can be roughly approximated by identity matrices \citep{SPCA-Witten,Witten-PHD}.
 Under such assumptions, a simplified optimization problem can be formulated as
\be
&&\min_{\bW_x,\bw_y} \left(- \tr(\bW_x^{\text{T}} \bX \bY^{\text{T}} \bW_y) +\gamma_1 \|\bW_x\|_1 +\gamma_2 \|\bW_y\|_1 \right) \\
 &&\mbox{s.t.} \quad (\bW_x)^{\text{T}} \; \bW_x = \bI_K \quad (\bW_y)^{\text{T}}  \; \bW_y = \bI_K. \notag
\label{PMD-CCA}
\ee
A basic optimization scheme  for an approximate huge sparse CCA in the TT format is illustrated in Figure~\ref{Fig:CCAL}, whereby the recursive
   contraction of TT networks reduces the  problem to a set of smaller-dimensional sparse CCA problems in the forms{\footnote{In quantum physics  TT/MPOs $\bX$ and $\bY^{\text{T}}$ are  often called transfer matrices and the
procedure to maximize the cost function $\bW_x^{\text{T}} \bX \bY^{\text{T}} \bW_y$
is called Transfer Matrix Renormalization Group (TMRG) algorithm \citep{Bursill1996,Wang1997,Xiang1999}.}}
  \be
 &&\min_{\bW^{(n)}_x,\bW^{(n)}_y} \left( -\tr( (\bW_x^{(n)})^{\text{T}}  \;  \bC^{(n)}_{xy} \; \bW^{(n)}_y)
 + \gamma_1\; \|\bW_{x}^{(n)}\|_1  + \gamma_2 \; \|\bW_{y}^{(n)}\|_1 \right)  \nonumber\\
 \\
&& \mbox{s.t.} \; \; (\bW_x^{(n)})^{\text{T}} \; \bW^{(n)}_x = \bI_K \quad (\bW_y^{(n)})^{\text{T}}  \; \bW^{(n)}_y = \bI_K,\nonumber
\ee
for $n=1,2,\ldots,N$, where $\bC^{(n)}_{xy} = \bW^{\text{T}}_{x,\;\neq n}\bX \bY^{\text{T}} \bW_{y, \;\neq n}$ are the contracted cross-covariance matrices.

Note that in the TT/MPO format, the sparsification penalty terms can be expressed as the $\ell_1$ norm  $\|\bW\|_1$,  e.g., as a sum of $\ell_1$ norms of fibers
of each TT-core, $\sum_{n=1}^N \sum_{r_{n-1},r_{n,k}}^{R_{N-1},R_N, K} \|\bw_{r_{n-1},r_n}^{(n,k)}\|_1$, or equivalently as a sum of the $\ell_1$ norms of slices,
 $\sum_{n=1}^N \|\bW_{(2)}^{(n)}\|_1$ (see Section \ref{sect:Lasso} for more detail).
In such cases, the orthogonalization of TT cores is replaced by a simple normalization of fibers to unit length.



\subsection{Tensor CCA for Multiview Data}

 The standard matrix CCA model has been generalized to tensor CCA,
  which in its simplest form can be formulated as the following optimization problem
 \be
 \label{eq:TCCA1}
&& \max_{\{\bw^{(n)}\}} (\underline {{\bC}} \; \bar \times_1 \; \bw^{(1)} \; \bar \times_2 \; \bw^{(2)}\;   \cdots \bar \times_N \; \bw^{(N)}) \\
&& \mbox{s.t.} \qquad  \bw^{(n)\;\text{T}} \;\bC_{nn} \; \bw^{(n)} = 1 \qquad (n=1,2,\ldots,N), \notag
 \ee
 where $\bw^{(n)} \in \Real^{I_n}$ are canonical vectors, $\bC_{nn} \in \Real^{I_n \times I_n}$ are covariance matrices and $\underline \bC \in \Real^{I_1 \times I_2 \times \cdots \times I_N}$ is an $N$th-order data tensor \citep{kim2007tensor,kim2009canonical,TCCA_Tao2015}.

\begin{figure}[t!]
(a)\\
\begin{center}
\includegraphics[width=8.1cm]{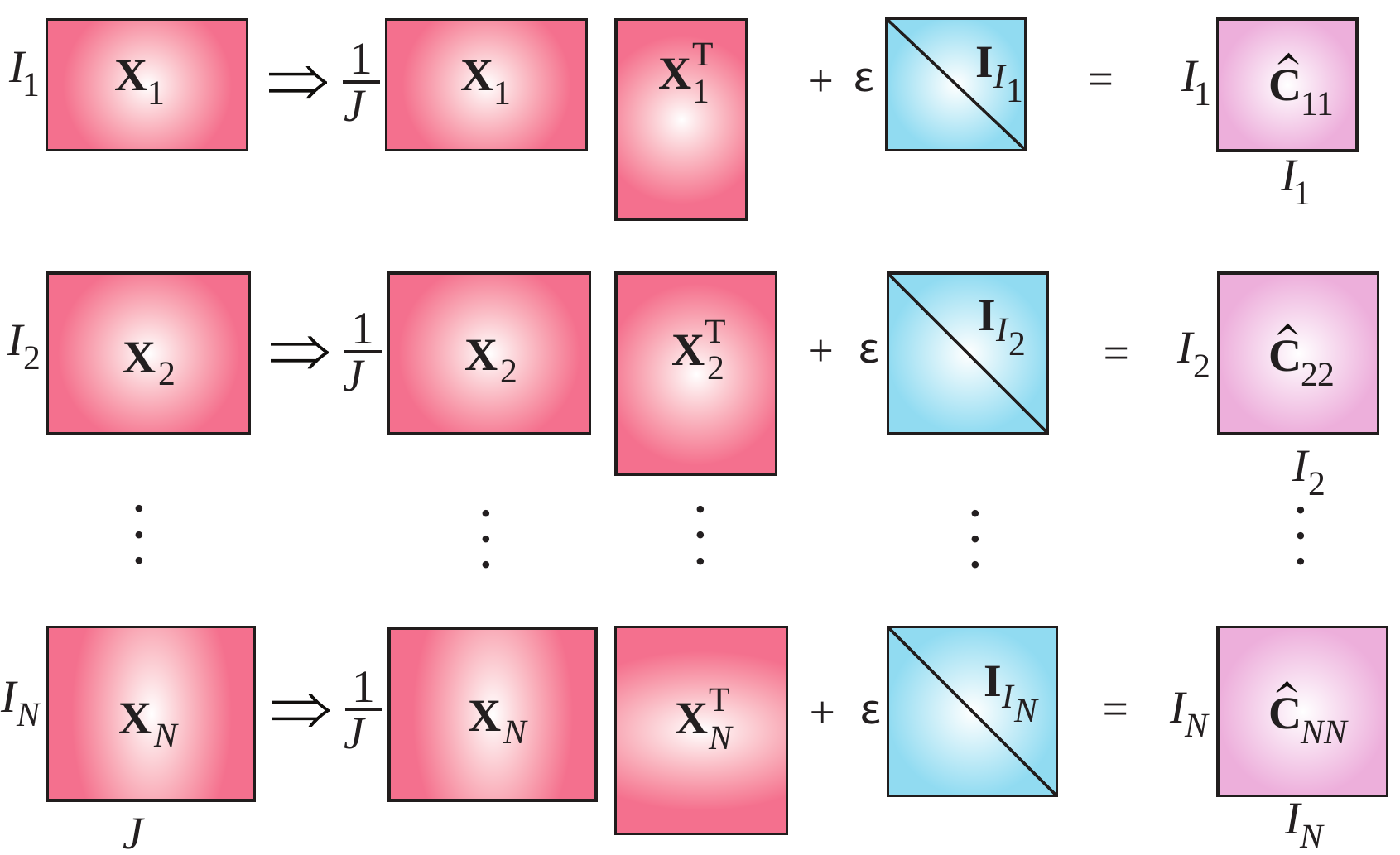}\\
\end{center}
(b)\\
\vspace{0.4cm}
\begin{center}
\includegraphics[width=9.5cm]{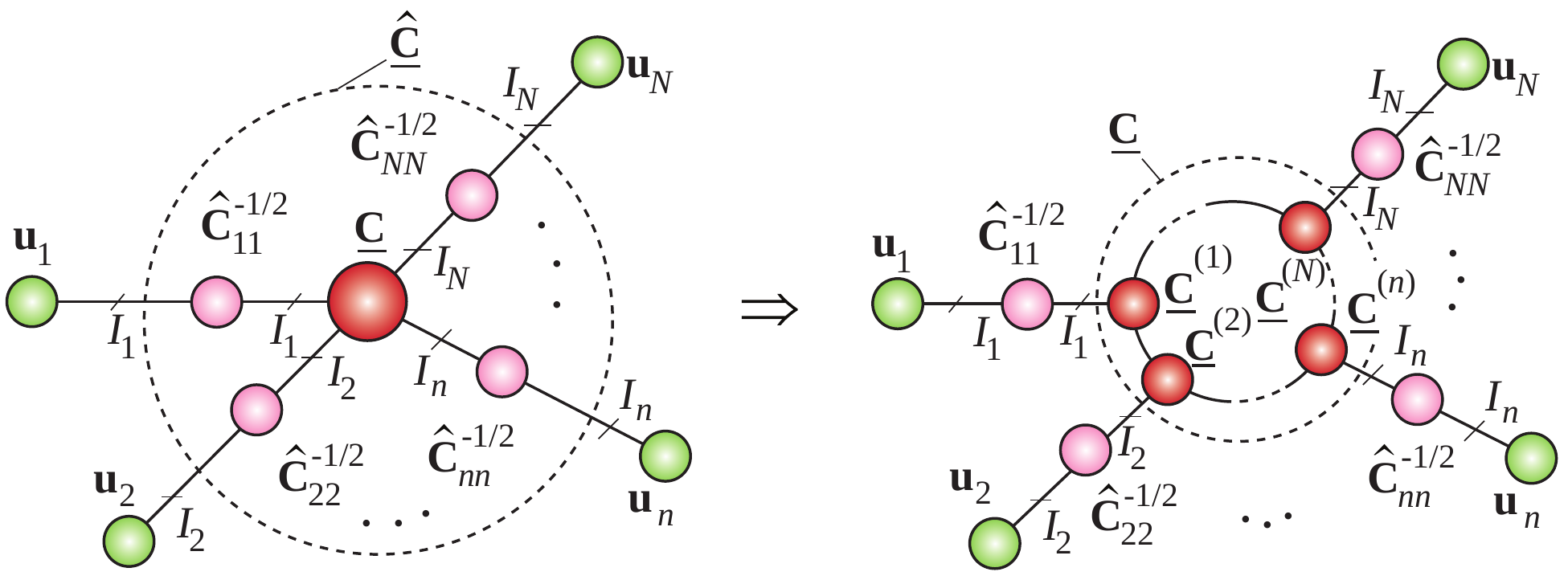}\\
\end{center}
\caption{Illustration of the computation of the tensor CCA for multi-view datasets. (a) Construction of regularized covariance matrices $\widehat{\bC}_{nn}
= (1/J)\bX_n \bX_n^{\text{T}} +\varepsilon \bI_{I_n}$. (b) Distributed representation of the cost function  (\ref{eq:TCCA2}) via tensor chain.}
\label{Fig:TCCA}
\end{figure}

 Similar to the standard CCA, the objective is to  find the canonical vectors which maximize the correlation function in (\ref{eq:TCCA1}).
Depending on the application, the matrices $\bC_{nn}$ and   tensor $\underline \bC$ can be constructed in many different ways. For example, for  multi-view data
  with $N$ different views $\{\bX_n\}_{n=1}^N$,  expressed as
$\bX_n =[\bx_{n1},\bx_{n2}, \ldots,\bx_{nJ}] \in \Real^{I_n \times J}$, the following  sampled covariance matrices and the tensor can be constructed  as \citep{TCCA_Tao2015}
\be
\bC_{nn} &=& \frac{1}{J}  \sum_{j=1}^J (\bx_{nj} \bx^{\text{T}}_{nj})= \frac{1}{J} \bX_n \bX_n^{\text{T}}, \\
\underline {\bC} &=& \frac{1}{J} \sum_{j=1}^J (\bx_{1j} \circ \bx_{2j} \circ \cdots \circ \bx_{Nj}).
\ee
Then the normalized data tensor, $\underline { \bC}$, can be expressed as
\be
\underline {\widehat{\bC}} = \underline {\bC} \; \times_1 \; \widehat{\bC}^{-1/2}_{11} \; \times_2 \; \widehat{\bC}^{-1/2}_{22} \; \cdots \times_N \; \widehat{\bC}^{-1/2}_{NN}
\ee
where the regularized covariance matrices are expressed as $\widehat{\bC}_ {nn}=\bC_ {nn} +\varepsilon \bI_{I_n}$, with small regularization parameters $\varepsilon >0$ (see Figure~\ref{Fig:TCCA}(a)).

In order to find canonical vectors, the  optimization problem (\ref{eq:TCCA1}) can be reformulated in a simplified form as
 \be
&& \max_{\{\bu^{(n)}\}}\ (\underline {\widehat{\bC}} \; \bar \times_1 \; \bu^{(1)} \; \bar \times_2 \; \bu^{(2)}\;   \cdots \bar \times_N \; \bu^{(N)}) \\
&& \mbox{s.t.} \qquad  \bu^{(n)\;\text{T}} \; \bu^{(n)} = 1 \qquad  (n=1,2,\ldots,N), \notag
\label{eq:TCCA2}
 \ee
 where $\bu^{(n)} = \widehat{\bC}^{-1/2}_{nn} \; \bw^{(n)}$, as illustrated in Figure~\ref{Fig:TCCA}(b).

 Note that the problem is equivalent to finding the rank-one approximation,
 for which many efficient algorithms exist.
However, if the  number of views is relatively large ($N>5$)  the core tensor becomes too large and we need to represent the core tensor in
 a distributed tensor network format (e.g., in TC format, as illustrated in Figure~\ref{Fig:TCCA}(b)).


\section{Two-Way Component Analysis in TT Formats}
\sectionmark{Large-Scale 2-way Component Analysis}
\label{sect:2CA-TT}

\begin{figure}[p!]
\begin{center}
\includegraphics[width=11.5cm]{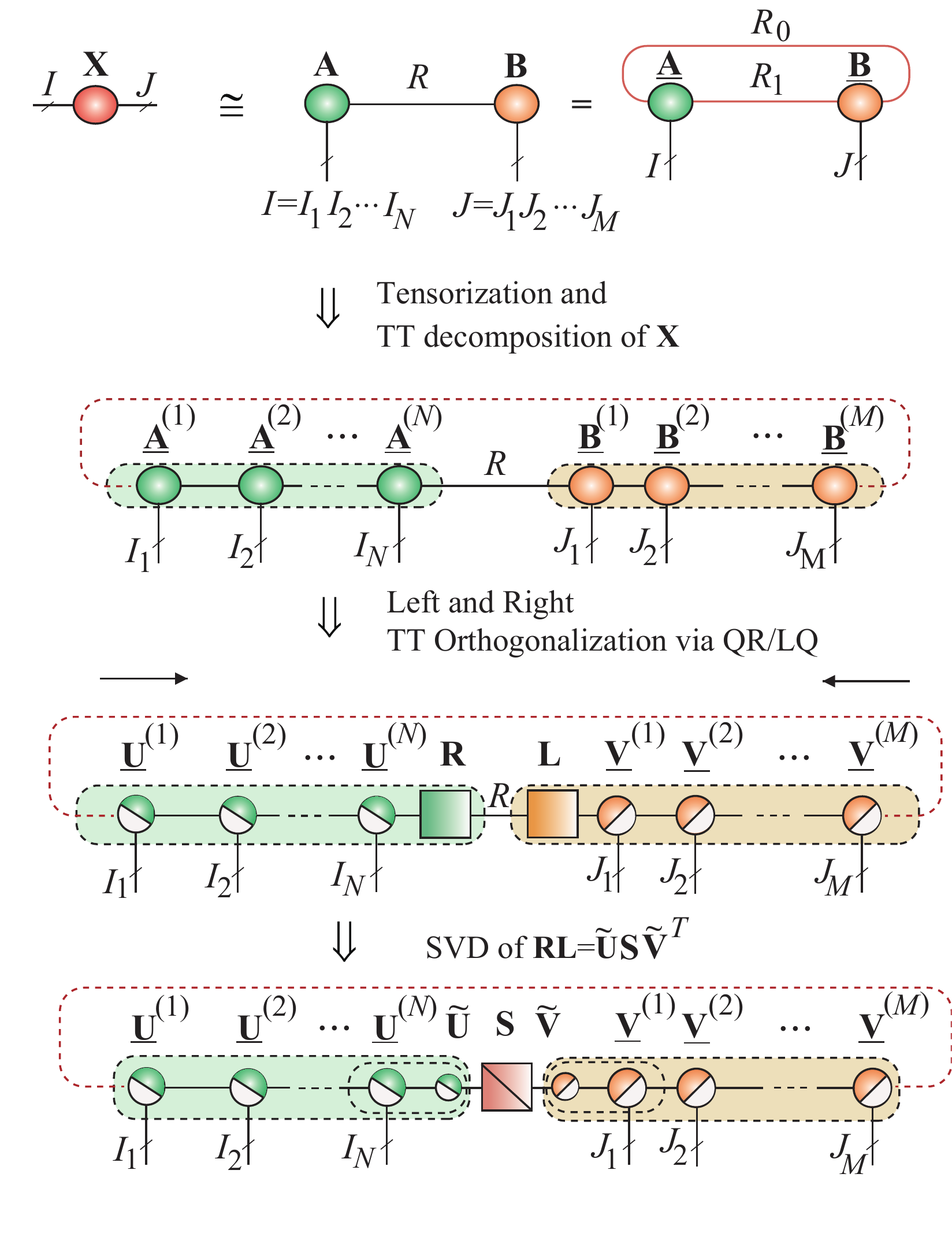}
\end{center}
\caption{Constrained low-rank matrix factorization (LRMF), $\bX \cong \bA \bB^{\text{T}}$, for a huge-scale  structured matrix, $\bX \in \Real^{I \times J}$, in the TT/TC format
(i.e., MPS with OBC or PBC).
By imposing suitable orthogonality constraints on  core tensors, large-scale  SVD/EVD can be performed in the TT/TC format.
Other constraints, such as nonnegativity and/or sparsity, are also conceptually possible.}
\label{Fig:CA-TT}
\end{figure}

Low-rank matrix factorizations  with specific constraints can be formulated in the following standard optimization setting \citep{cichocki1995multi,MultiNMF_Cichocki06,RALScichocki2007,NMF-book,CichockiJSICE2011}
\be
\min_{\bA,\bB} \, J(\bA,\bB)= \|\bX - \bA \bB^{\text{T}}\|_F^2,
\label{CAoptim}
\ee
 where a large-scale data matrix, $\bX \in \Real^{I \times J}$, is given and the objective is
to estimate the factor matrices, $\bA \in \Real^{I \times R}$ and $\bB \in \Real^{J \times R}$ (with the assumption that $R \ll \{I,J\}$),
subject to suitable constraints being imposed on one or both factor matrices.

If such a structured data matrix can be represented by a low-rank TT or TC network (e.g., through a cheap matrix cross-approximation the factor matrices
can be represented in a TT/TC  format, as illustrated in the top panel of Figure~\ref{Fig:CA-TT}),
then in the first step no constraint on the TT core tensors, $\underline \bA^{(n)}$ and $\underline \bB^{(m)}$, need to be imposed to represent $\bX$ in a distributed  form. However, such a factorization is not unique and the so obtained components (in TT format) do not have any physical meaning.

\markright{3.8.\quad Solving System of Linear Equations}

 With a matrix $\bX$ already in a compressed TT format,  in the second step, we may impose the desired constraints. For example, by imposing the left-orthogonality on TT cores of $\underline \bA \in \Real^{I_1 \times \cdots \times I_N \times R}$ and right-orthogonality on the  TT cores of $\underline \bB \in \Real^{R \times J_1 \times \cdots \times J_N}$, a
truncated SVD can be computed for the $R$ largest
singular values and the corresponding singular vectors, as illustrated in the bottom panel of Figure~\ref{Fig:CA-TT}.
In a similar way, we can compute the pseudo-inverse of a huge matrix.
%
%

%

\section{Solving Huge-Scale  Systems of Linear Equations}
\sectionmark{Solving System of Linear Equations}

Solving linear systems of large scale equations  arises throughout science  and engineering; e.g, convex optimization,  signal processing, finite elements and machine learning all rely partially on approximately solving linear equations, typically
using some additional criteria like sparseness or smoothes.

Consider  a huge system of linear algebraic equations
in a TT format (see also Part 1 \citep{Part1}), given by 
\be
\bA \bx \cong \bb
\ee
where $\bA \in \Real^{I \times J}$,  $\bb \in \Real^I$
 and  $\bx \in \Real^J$.
 The objective is to find an  approximative solution, $\bx \in \Real^J$, in a TT
 format, by imposing additional constraints (regularization terms) such as smoothness,
 sparseness and/or nonnegativity on the vector $\bx$. Several
innovative TT/HT network solutions to this problem do exist
(see \citep{Oseledets-Dolgov-lin-syst12,DolgovAMEN2014} and references therein), however, most
 focus on only symmetric square matrices. We shall  next consider several more general cases.

\subsection{\bf  Solutions for a Large-scale Linear Least Squares Problem using ALS}
\subsectionmark{Solutions for a  Linear Least Squares Problem using ALS}

The Least Squares (LS) solution  (often called the ridge regression)  minimizes the following  regularized cost function
\be
J(\bx) &=&  \|\bA \bx -\bb\|^2_2 +\gamma \|\bL \bx\|^2_2 \nonumber \\ 
&=& \bx^{\text{T}} \bA^{\text{T}} \bA \bx - 2 \bx^{\text{T}} \bA^{\text{T}} \bb +\bb^{\text{T}} \bb + \gamma \bx^{\text{T}} \bL^{\text{T}} \bL \bx,
\label{eq:AxbL2}
\ee
with the Tikhonov regularization term on the right hand side, while
 $\bL$  is the so-called  smoothing matrix, typically in the form of the discrete first-order or second-order  derivative matrix,
 and $\widetilde \bL=\bL^{\text{T}} \bL \in \Real^{J \times J}$.

\begin{figure}
\begin{center}
\includegraphics[width=0.99\textwidth]{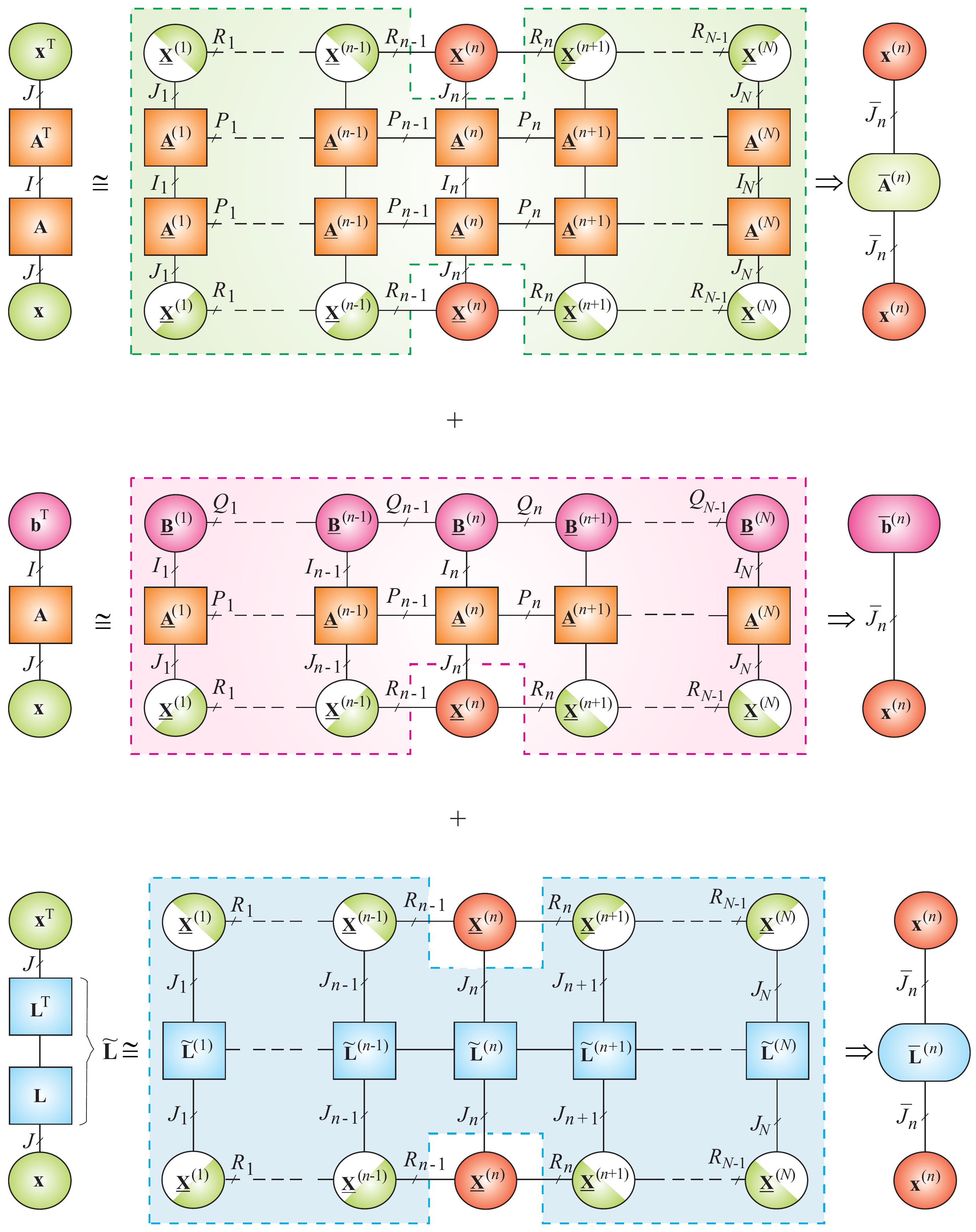}
\end{center}
\caption{A TT network, consisting of  three sub-networks, which solves (in the LS sense) a system of linear equations in the form $\bA \bx \approx \bb$,
where $\bA \in \Real^{I_1 \cdots I_N \; \times \; J_1 \cdots J_N}$ is a  huge non-symmetric matrix  and the Tikhonov regularized cost function is expressed  as $J(\bx)=\|\bA \bx -\bb\|_2^2 + \gamma \|\bL \bx\|^2_2 = \bx^{\text{T}} \bA^{\text{T}} \bA \bx - 2 \bb^{\text{T}} \bA \bx + \bb^{\text{T}} \bb + \gamma \bx^{\text{T}} \bL^{\text{T}} \bL \bx$.}
\label{Fig:Ax=b}
\end{figure}

  Upon neglecting the constant factor, $\bb^{\text{T}} \bb$, we arrive at a simplified form
\be
J(\bx) = \bx^{\text{T}} \bA^{\text{T}} \bA \bx - 2 \bx^{\text{T}} \bA^{\text{T}} \bb + \gamma \bx^{\text{T}} \widetilde \bL \bx.
\ee
The approximative TT representation  of the matrix $\bA$ and vectors $\bx$ and $\bb$ allows for the construction of three
tensor sub-networks, as shown in Figure~\ref{Fig:Ax=b}. Upon simultaneous recursive contractions of all the   three sub-networks
 for each node $n=1,2,\ldots,N$, a large-scale infeasible optimization problem in (\ref{eq:AxbL2}) can be converted into a set of much
 smaller optimization problems based on the  minimization of  cost functions
\be
J(\bx) &=& J(\bX_{\neq n} \bx^{(n)}) \\
&=&  (\bx^{(n)})^{\text{T}} \bX_{\neq n}^{\text{T}} \bA^{\text{T}} \bA \bX_{\neq n} \bx^{(n)}
- 2 (\bx^{(n)})^{\text{T}} \bX_{\neq n}^{\text{T}} \bA^{\text{T}}  \bb \notag \\
&&+ \gamma (\bx^{(n)})^{\text{T}} \bX_{\neq n}^{\text{T}} \bL^{\text{T}} \bL \bX_{\neq n} \bx^{(n)}, \quad n=1,2,\ldots,N. \notag
\ee
These can  be expressed in the following equivalent compact  forms
\be
J_n(\bx^{(n)}) =  (\bx^{(n)})^{\text{T}} \; \overline \bA^{(n)} \, \bx^{(n)}
- 2 (\bx^{(n)})^{\text{T}} \; \overline \bb^{(n)}  + \gamma (\bx^{(n)})^{\text{T}} \; \overline \bL^{(n)} \, \bx^{(n)}.\nonumber\\
\label{eq:AnBnLn}
\ee
where  $\bx^{(n)} \in \Real^{R_{n-1} J_n R_n}$,
$\overline \bA^{(n)} =\bX_{\neq n}^{\text{T}} \bA^{\text{T}} \bA \bX_{\neq n} \in \Real^{R_{n-1}  J_n  R_n \times R_{n-1}  J_n  R_n }$,
$ \overline \bb^{(n)} = \bX_{\neq n}^{\text{T}} \bA^{\text{T}} \bb   \in \Real^{R_{n-1} J_n R_{n}}$,
and $\overline \bL^{(n)} =  \bX_{\neq n}^{\text{T}} \bL^{\text{T}} \bL \bX_{\neq n} \in \Real^{R_{n-1}  J_n  R_n \times R_{n-1}  J_n  R_n }$.
%
%

In this way, a  large-scale system of linear equations is converted into a set of smaller-size  systems, which can be solved by any standard  method
\be
 \overline \bA^{(n)} \, \bx^{(n)} \cong \overline \bb^{(n)}, \qquad n=1,2,\ldots,N,
\label{eq:An}
\ee
It is important to note that if the TT cores are  left- and right-orthonormalized, the symmetric semi-positive definite matrix $\overline \bA^{(n)}$  is better conditioned than the original huge matrix $\bA$.

The  matrix multiplications  in (\ref{eq:AnBnLn}) is not performed explicitly, however, the TT format  alleviates the curse of dimensionality via a recursive contraction of cores, as shown in Figure~\ref{Fig:Ax=b}.
%

{\bf Remark.} The usual assumption that all data admits low-rank TT/QTT approximations is a key  condition for successfully employing  this approach.
However, for data with
a weak structure, TT approximations may have relatively large ranks which makes the calculation difficult or even impossible.
The way in which the TT ranks are  chosen and adapted is  therefore the key factor  in  efficiently solving  huge-scale
structured linear systems.

\subsection{Alternating Minimum Energy (AMEn)  Algorithm for Solving a Large-scale Linear Least Squares Problems}
\subsectionmark{AMEn  Algorithm for Solving a  Linear Least Squares Problems}
\label{Sect:Amen-LinearSystem}

The AMEn approach, introduced in Section~\ref{Sect:Amen-Eigen} for the EVD problem, has been developed first historically as an
efficient solution to large-scale least squares problems \citep{DolgovAMEN2014}.
The main difference from the EVAMEn method for solving eigenvalue problems is in the definition of the gradient of the cost function and in
the computation of enrichment cores $\underline \bZ^{(n)}$.

Algorithm~\ref{alg:AMEN} summarizes the AMEn algorithm for solving  huge systems of linear equations $\bA \bx \cong \bb$, for $\bA \in \Real^{I\times J}$ and $\bb \in \Real^{I}$ given in the TT format, with $I=I_1I_2\cdots I_N$ and $J=J_1J_2\cdots J_N$.
In the first preliminary step, we assume, that the core $\underline \bX^{(n)} \in \Real^{R_{n-1}\times J_n \times R_n}$ has been updated at the $n$th micro-iteration by solving the reduced linear system $\overline \bA^{(n)} \, \bx^{(n)} \cong \overline \bb^{(n)}$, similarly to the ALS method described in the previous section.
Note that solving this reduced linear system is equivalent to minimizing the cost function
$J(\bx)=\| \bA\bx - \bb \|^2$ with respect to the core $\bX^{(n)}$.
Similarly, for a residual vector defined as $\widetilde{\brr} = \bA^{\text{T}}(\bA \bx -\bb) \in\Real^J$ in the TT format
(i.e., $\underline {\widetilde{\bR}}= \llangle \underline {\widetilde{\bR}}^{(1)},  \underline {\widetilde{\bR}}^{(2)}, \ldots, \underline {\widetilde{\bR}}^{(N)}, \rrangle $,  each of its TT cores  $\underline {\widetilde{\bR}}^{(n)} \in \Real^{Q_{n-1}\times J_n \times Q_n}$ can be updated via the ALS scheme by $\widetilde{\brr}^{(n)} = \bR_{\neq n}^{\text{T}} \; \widetilde{\brr} = \arg\min_{\underline {\widetilde{\bR}}^{(n)}}\|\widetilde{\brr} - (\bA^{\text{T}} \bA\bx - \bA^{\text{T}} \bb)\|^2$.
In other words, the vector $\widetilde{\brr}$ approximates the gradient vector, $\nabla J(\bx) \propto \bA^{\text{T}} \bA \bx - \bA^{\text{T}} \bb$, and it can be efficiently updated via the contractions of core tensors.

\begin{algorithm}[t]
\caption{\textbf{AMEn for linear systems $\bA \bx \cong \bb$ \citep{DolgovAMEN2014}}}
\label{alg:AMEN}
\begin{algorithmic}[1]
\REQUIRE Matrix $\bA \in \Real^{I \times J}$, with $I\geq J$, vector $\bb \in \Real^{I}$, initial\\ guesses
for $\bx \in \Real^J$ and residual $\widehat{\brr}=\bA^{\text{T}} (\bA \bx-\bb) \in \Real^J$ in the\\ TT format
\ENSURE Approximate solution $\bx$ in the TT format  $\underline \bX = \llangle \underline \bX^{(1)}, \underline \bX^{(2)}, \ldots , \underline \bX^{(N)}\rrangle $
\WHILE {not converged or iteration limit is not reached}
\STATE Right-orthogonalize  cores $\underline \bX^{(n)}, \underline {\widehat{\bR}}^{(n)}$
 for $n=N,N-1,\ldots,2$

\FOR{$n=1$ to $N$}
    \STATE Update the  core $\underline \bX^{(n)}$ by solving $\overline \bA^{(n)} \, \bx^{(n)} \cong \overline \bb^{(n)}$
      \STATE  Update the core $\underline {\widehat{\bR}}^{(n)}$ by solving $\bR_{\neq n} \widehat{\brr}^{(n)} \cong \widehat{\brr}$
    \IF {$n<N$}
    \STATE Compute the core $\underline \bZ^{(n)}$ of the partially projected \\ gradient by solving
    		$(\bX^{< n} \otimes_L \bI_{J_n} \otimes_L (\bR^{>n})^{\text{T}})\; \bz^{(n)} \cong \widehat{\brr}$
    				
    \STATE Enrich the cores $\bX^{(n)}_{j_n} \leftarrow [\bX^{(n)}_{j_n}, \; \bZ^{(n)}_{j_n}]$ and $\bX^{(n+1)}_{j_{n+1}} \leftarrow \left[ \begin{matrix} \bX^{(n+1)}_{j_{n+1}}  \\  {\bf 0} \end{matrix} \right]$
    \STATE Left-orthogonalize cores $\underline \bX^{(n)}, \underline \bR^{(n)}$
\ENDIF
\ENDFOR
\ENDWHILE
\RETURN  $\underline \bX= \llangle \underline \bX^{(1)}, \underline \bX^{(2)}, \ldots , \underline \bX^{(N)} \rrangle $
\end{algorithmic}
\end{algorithm}

Next, for building an enrichment, $\underline \bZ^{(n)} \in \Real^{R_{n-1} \times J_n \times Q_{n}}$, \citep{DolgovAMEN2014} considered an approximation to the partially projected gradient
	\begin{multline}
	(\bX^{< n} \otimes_L \bI_{J_n\cdots J_N})^{\text{T}} (\bA^{\text{T}} \bA\bx - \bA^{\text{T}} \bb)   \\
	\cong \text{vec} \left( \llangle \underline \bZ^{(n)}, \underline \bR^{(n+1)}, \ldots, \underline \bR^{(N)} \rrangle \right) \in \Real^{R_{n-1} J_n\cdots J_N}
	\end{multline}
with respect to $\underline \bZ^{(n)}$ in the TT format. Since
	$
	\text{vec} ( \llangle \underline \bZ^{(n)}, \underline \bR^{(n+1)}, \ldots, \underline \bR^{(N)} \rrangle )
	= ( \bI_{R_{n-1}J_n} \otimes_L (\bR^{>n})^{\text{T}}) \; \bz^{(n)},
	$
the following closed form representation can be obtained
	\be
	\bz^{(n)} = (\bX^{< n} \otimes_L \bI_{J_n} \otimes_L (\bR^{>n})^{\text{T}})^{\text{T}} \; \brr
\ee
which can be efficiently computed via recursive contractions of core tensors,  similarly to the standard ALS method.

Note that due to the partial gradient projection, the sizes of the TT cores $\underline \bX^{(n)}$ and $\underline \bZ^{(n)}$ now become consistent, which allows us
 to perform the concatenation
	\be \label{updated_amen_xn}
	\widetilde{\underline \bX}^{(n)} \leftarrow \underline \bX^{(n)} \boxplus_{\;3} \underline \bZ^{(n)}
	\in \Real^{R_{n-1} \times J_n \times (R_n + Q_n)},
	\ee
or equivalently, $\widetilde {\bX}^{(n)}_{j_n} \leftarrow [\bX^{(n)}_{j_n},\; \bZ^{(n)}_{j_n}] \in \Real^{R_{n-1} \times (R_n + Q_n) }$.


After the enrichment and the subsequent orthogonalization of the TT cores, the column space (i.e., the subspace spanned by the columns) of the frame matrix, say $\bX_{\neq n+1}$, is expanded (see \citep{KressnerEIG2014} for more detail).

For rigor, it can be  shown that
	$$
	\text{range}\left( \widehat \bX_{\neq n+1} \right)
	\supset
	\text{range}\left( \bX_{\neq n+1} \right),
	$$
where   $\text{range}(\bA)$ denotes the column space of a matrix $\bA$, whereas  $\bX_{\neq n+1}$ and $\widehat \bX_{\neq n+1}$ are the frame matrices before and after the enrichment and orthogonalization.
It should be noted   that $\text{range}(\bX_{\neq n+1}) = \text{range}(\bX^{<n+1}) \otimes \Real^{J_{n+1}} \otimes \text{range}((\bX^{> n+1})^{\text{T}})$, and that
	$
	\bX^{<n+1}
	= \left( \bX^{<n}  \otimes_L \bI_{J_n} \right) \bX^{(n)}_{L}
	\in \Real^{J_1J_2\cdots J_n \times R_n} ,
	$
where $\bX^{(n)}_{L}$ is the left unfolding of $\underline \bX^{(n)}$.
     The left-orthogonalization of the enriched core $\widetilde{\underline{\bX}}^{(n)}$ in \eqref{updated_amen_xn} can be written as
	$$
	\widetilde \bX^{(n)}_{L} = \begin{bmatrix} \bX^{(n)}_{L} & \bZ^{(n)}_{L}
	\end{bmatrix}
	=
	\widehat \bX^{(n)}_{L} \bP
	\in\Real^{R_{n-1} J_n\times (R_n+Q_n)},
	$$
where $\widehat{\underline \bX}^{(n)} \in \Real^{R_{n-1}\times J_n \times \hat{R}_n}$ is the left-orthogonalized core and $\bP \in \Real^{\hat{R}_n \times (R_n + Q_n)}$ is a full row rank matrix.
From these expressions, it follows that $\text{range} (\widehat \bX^{<n+1} ) \supset \text{range} (\bX^{<n+1} )$.




\subsection{Multivariate Linear Regression and Regularized  Approximate Estimation of Moore-Penrose Pseudo-Inverse}
\subsectionmark{TT Networks for Linear Regression}

\begin{figure}[p]
\begin{center}
\includegraphics[width=0.95\textwidth]{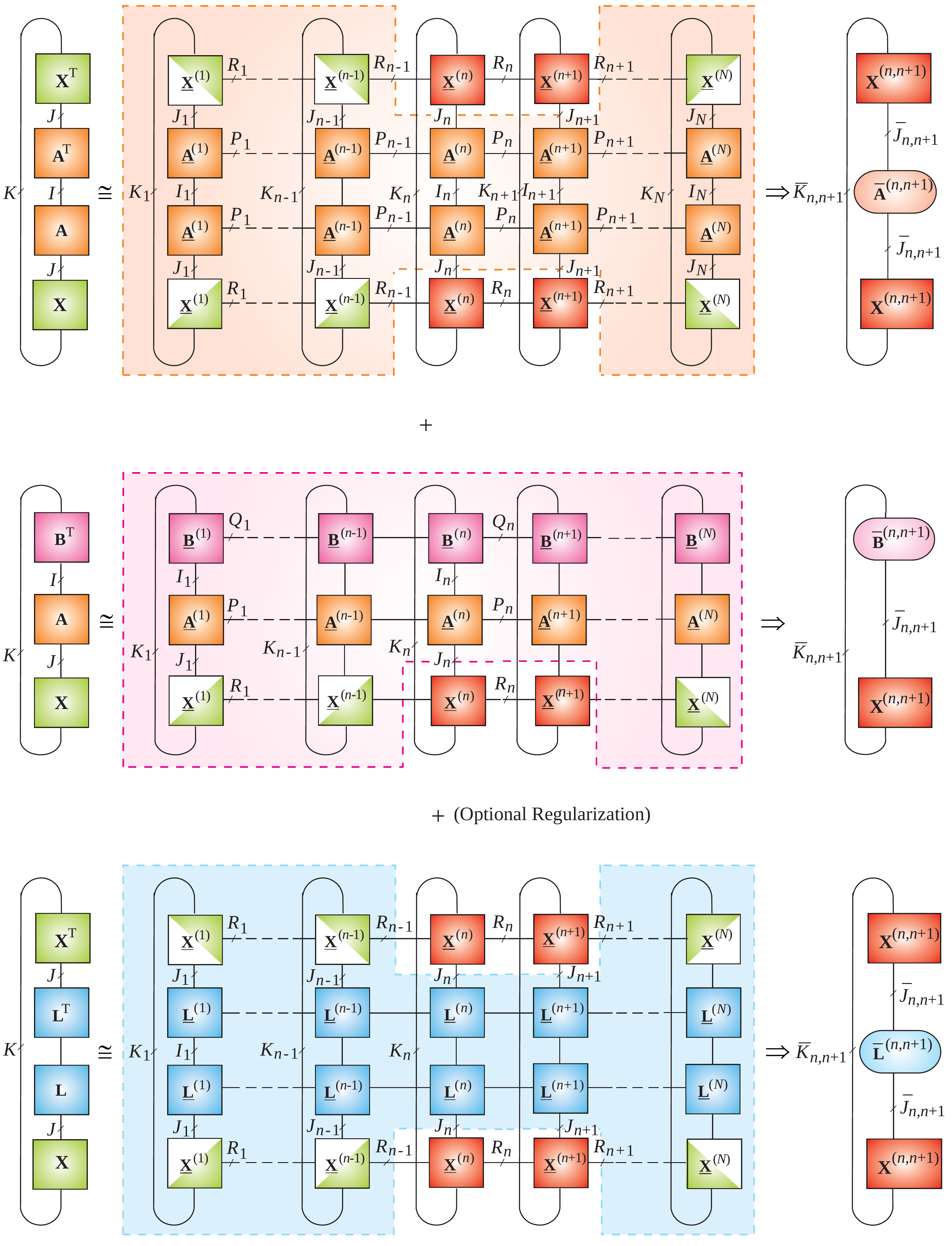}
\end{center}
\caption{TT sub-networks for the  minimization of a very large-scale regularized cost function $J(\bX)=\|\bA \bX -\bB\|_F^2 + \gamma \|\bL\bX\|^2_F =
\tr(\bX^{\text{T}} \bA^{\text{T}} \bA \bX) -2 \tr(\bB^{\text{T}} \bA \bX) +\gamma \tr(\bX^{\text{T}} \bL^{\text{T}} \bL \bX) + \tr(\bB^{\text{T}} \bB)$.
This is achieved by  contractions of core tensors  using the MALS (DMRG2) approach. These contractions convert the huge-scale cost function to relatively
smaller-scale cost functions. The top panel represents the cost function $\tr(\bX^{\text{T}} \bA^{\text{T}} \bA \bX)$, the middle panel describes the cost function $\tr(\bB^{\text{T}} \bA \bX)$, while the bottom panel represents regularization term $\tr(\bX^{\text{T}} \bL^{\text{T}} \bL \bX)$.}
\label{Fig:AXBL2}
\end{figure}

The approaches described in the two previous  sections can be straightforwardly  extended to regularized multivariate regression, which can be formulated in a  standard
way as the minimization of the cost function
\be
J(\bX) &=&  \|\bA \bX -\bB\|^2_F +\gamma \|\bL \bX\|^2_F  \\
&=& \tr(\bX^{\text{T}} \bA^{\text{T}} \bA \bX) - 2 \tr(\bB^{\text{T}} \bA \bX) + \tr(\bB^{\text{T}} \bB) + \gamma \tr( \bX^{\text{T}} \bL^{\text{T}} \bL \bX), \nonumber
\label{eq:AXBL22}
\ee
where  $\bA \in \Real^{I \times J}$  (with $I \geq J$),  $\bB \in \Real^{I \times K}$,
 $\bX \in \Real^{J \times K}$  and  $\bL \in \Real^{J \times J}$ (see Figure~\ref{Fig:AXBL2}).
 The objective is to find an approximative solution, $\bX \in \Real^{J \times K}$,  in a TT
 format, by imposing additional constraints  (e.g., via  Tikhonov regularization). In a special case
 when $\bB = \bI_I$, for $K=I$ the problem is equivalent to the computation of Moore-Penrose pseudo inverse in an appropriate TT format
 (see \citep{Lee-TTPI} for computer simulation experiments).

The approximative TT representation  of matrices $\bA$, $\bB$ and $\bX$ generates  TT sub-networks shown in Figure~\ref{Fig:AXBL2},
which can be optimized by the  ALS or MALS (DMRG2) approaches.

By minimizing  all three sub-networks simultaneously using recursive  contractions of TT cores and performing  ALS
sequentially for each node $n=1,2,\ldots,N$,  a large-scale infeasible optimization problem in (\ref{eq:AXBL22}) can be converted into a set of much smaller
optimization problems expressed by the set of (linked) cost functions
\be
\label{eq:583}
J_n(\bX^{(n)}) &=&  \tr\left((\bX^{(n)})^{\text{T}} \; \overline \bA^{(n)} \, \bX^{(n)}\right)
- 2 \tr\left((\overline \bB^{(n)})^{\text{T}} \, \bX^{(n)}\right) +  \notag \\
&& + \gamma \tr\left((\bX^{(n)})^{\text{T}} \; \overline \bL^{(n)} \, \bX^{(n)}\right), \;\;n=1,\ldots,N,
\ee
where
\be
\overline \bA^{(n)} &=&\bX_{\neq n}^{\text{T}} \bA^{\text{T}} \bA \bX_{\neq n} \in \Real^{\overline J_n  \times \overline J_n}, \notag \\
\overline \bB^{(n)} &=& \bX_{\neq n}^{\text{T}} \bA^{\text{T}} \bB   \in \Real^{\overline J_n  \times \overline K_n }, \notag \\
\overline \bL^{(n)} &=&  \bX_{\neq n}^{\text{T}} \bL^{\text{T}} \bL \bX_{\neq n}  \in \Real^{\overline J_n  \times \overline J_n} \notag
\ee
and  $\bX^{(n)} \in \Real^{\overline J_n  \times \overline K_n}$,  with $\overline J_n= R_{n-1} J_n R_n$ and $\overline K_n= R_{n-1} K_n R_n$.

The regularization term, $ \gamma \tr\left((\bX^{(n)})^{\text{T}} \; \overline \bL^{(n)} \, \bX^{(n)}\right)$, in (\ref{eq:583}) helps not only to alleviate the ill-posedness of the problem, but  also to
considerably reduce computational complexity and improve convergence by  avoiding  overestimation of the TT
ranks. 
In other words,   regularization terms, especially  those that impose  smoothness or sparsity, lead to  smaller TT ranks with  fewer parameters,
thus yielding a  better approximation.

Alternatively, the MALS (DMRG2) procedure produces somewhat larger-scale optimization sub-problems
\be
\label{eq:584}
&&J_n(\bX^{(n)}) =  \tr\left((\bX^{(n,n+1)})^{\text{T}}  \; \overline \bA^{(n,n+1)} \;  \bX^{(n,n+1)}\right) + \\
&&- 2 \tr\left((\overline \bB^{(n,n+1)})^{\text{T}} \; \bX^{(n,n+1)}\right)  + \gamma \tr\left((\bX^{(n,n+1)})^{\text{T}} \; \overline \bL^{(n,n+1)} \; \bX^{(n,n+1)}\right) \notag
\ee
where  $\bX^{(n,n+1)} \in \Real^{\overline J_{n,n+1}  \times \overline K_{n,n+1}}$ and
\be
\overline \bA^{(n,n+1)} &=&\bX_{\neq n,n+1}^{\text{T}} \bA^{\text{T}} \bA \bX_{\neq n,n+1} \in \Real^{\overline J_{n,n+1}  \times \overline J_{n,n+1} }, \notag \\
\overline \bB^{(n,n+1)} &=& \bX_{\neq n,n+1}^{\text{T}} \bA^{\text{T}} \bB   \in \Real^{\overline J_{n,n+1}  \times \overline K_{n,n+1} }, \\
\overline \bL^{(n,n+1)} &=&  \bX_{\neq n,n+1}^{\text{T}} \bL^{\text{T}} \bL \bX_{\neq n,n+1}  \in \Real^{\overline J_{n,n+1}  \times \overline J_{n,n+1} }
\ee
with $\overline J_{n,n+1}= R_{n-1} J_n J_{n+1} R_{n+1}$ and  $\overline K_{n,n+1}= R_{n-1} K_n K_{n+1} R_{n+1}$.

 One direction of  future work is to select an
optimal  permutation of data tensor  which gives  an optimal  tree for a TT approximation, given the available data samples.


\subsection{\bf Solving  Huge-scale LASSO and Related Problems}
\label{sect:Lasso}

One way to impose sparsity constraints on a vector or a matrix
is through the  LASSO approach, which  can be formulated as the following optimization problem \citep{Lasso1996,BoydFTML2011,TanLasso2015,kim2015new}:
 \be
 \min_{\bx} \{ \|\bA \bx -\bb\|_2^2 + \gamma \|\bx\|_1 \}
 \ee
 which is equivalent to
 \be
 \min_{\bx} \|\bA \bx -\bb\|_2^2, \quad \mbox{s.t.} \quad  \|\bx\|_1 \leq t
 \ee
 and is closely related to the Basis Pursuits (BP) and/or compressed sensing  (CS) problems, given by
  \be
&& \min \|\bx\|_1   \quad \mbox{s.t.} \quad \bA \bx =\bb, \nonumber \\
&& \min \|\bL \bx\|_1   \quad \mbox{s.t.} \quad \|\bA \bx -\bb\|_2^2 \leq \varepsilon. \nonumber
\ee
In many applications, it is possible to utilize  {\it a priori} information about a sparsity profile or sparsity pattern \citep{JordanL12016} in data.
For example, components of the vector $\bx$
may be clustered in groups or blocks, so-called group sparse components. In such a case, to model
the group sparsity, it is
convenient to replace the $\ell_1$-norm with the $\ell_{2,1}$-norm, given by
\be
\|\bx\|_{2,1} = \sum_k\|\bx_{g_k}\|_2, \quad k=1,2,\ldots,K,
\ee
where the symbol $\bx_{g_k}$ denotes the components in the $k$-th group and $K$ is the total  number of groups.

When the sparsity structures are overlapping, the problem can be reformulated as the  Group LASSO problem \citep{IRLS2014}
 \be
 \min_{\bx} \|\bA \bx -\bb\|_2^2 +\gamma \|\bG \mbi \Phi \bx\|_{2,1},
 \ee
where $\mbi \Phi$ is the optional sparse basis and $\bG$ a binary matrix  representing the group configuration.

 Similarly, a multivariate regression problem with sparsity constraints  can be formulated as
\be
J(\bX) &=&  \|\bA \bX -\bB\|^2_F +\gamma \|\bX\|_{{\cal {S}}_q} \\
&=& \tr(\bX^{\text{T}} \bA^{\text{T}} \bA \bX) - 2 \tr(\bB^{\text{T}} \bA \bX) + \tr(\bB^{\text{T}} \bB) + \gamma \tr( \bX^{\text{T}} \bX)^{q/2}, \nonumber
\label{eq:AXBL1}
\ee
where $\|\bX\|_{{\cal S}_q} = \tr (\bX^{\text{T}} \bX)^{q/2} 
= \sum_r \sigma_r^q (\bX)$
is the Schatten $q$-norm and $\sigma_r(\bX)$ is a singular value of the  matrix $\bX$.

Other  generalizations of the standard LASSO include the Block LASSO, Fused LASSO, Elastic Net and Bridge regression algorithms
\citep{OgutuLasso2014}. However,  these models have
not yet  been  deeply investigated or  experimentally tested  in the tensor network setting, where
 the challenge is to efficiently represent/optimize in a TT format  non-smooth
penalty terms in the form of the $\ell_1$-norm $\|\bx\|_1$, or more generally $\ell_q$-norm $\|\bL \bx\|^q_q$ and  $\|\bX\|_{{\cal S}_q}$ with  $ 0 < q \leq 1$.

\begin{figure}[t!]
\begin{center}
\includegraphics[width=0.99\textwidth]{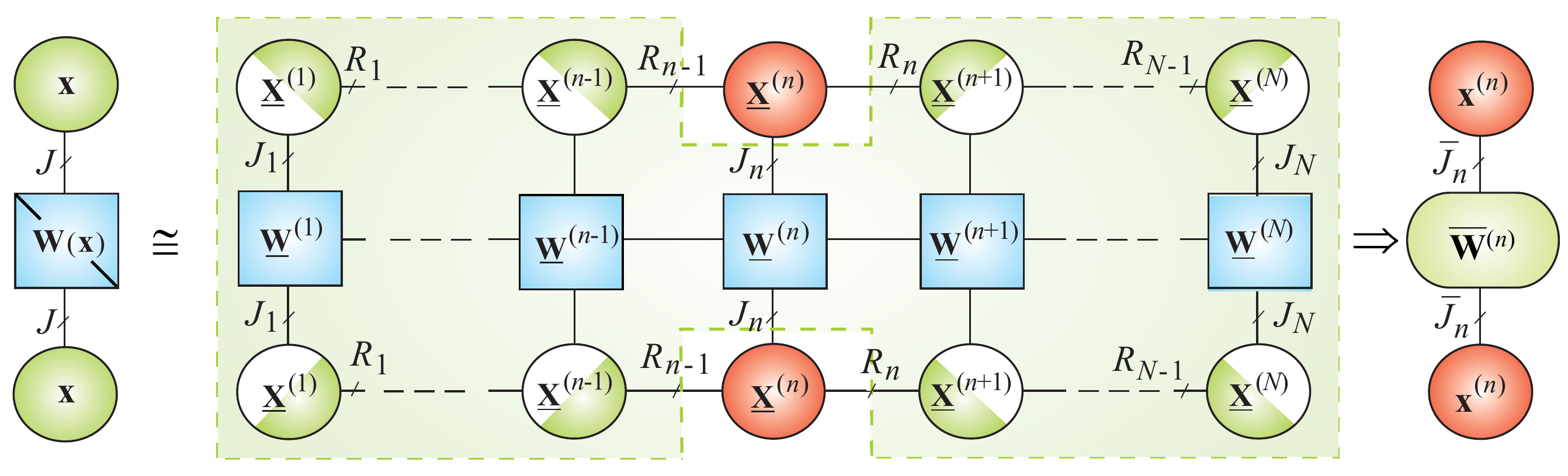}
\end{center}
\caption{A tensor network  to  compute, in the TT format,  the $\ell_q$-norm of a huge-scale vector $\bx \in \Real^J, \;\; J=J_1 J_2 \cdots J_N$,
using the  IRLS (weighted $\ell_2$-norm).}
\label{Fig:Lq-norm}
\end{figure}

A simple  approach in this direction would be  to apply  Iteratively Reweighted Least Squares (IRLS) methods{\footnote{
The IRLS approach for  the  $\ell_1$-norm
 is motivated by the variational representation of the
 norm $\|\bx\|_1 = \min_{w_j>0} (0.5)\sum_j (x^2_j/w_j + w_j)$.}} \citep{CandesIRLS2008},
whereby the $\ell_1$-norm is replaced by the reweighted $\ell_2$-norm (see Figure~\ref{Fig:Lq-norm}) \citep{Lee-Lasso}, that is
 \be
 \|\bx\|_1 = \bx^{\text{T}} \bW \bx =\sum_{j=1}^J w_j x_j^2,
 \ee
where $\bW (\bx)  =\diag(w_1,w_2, \ldots, w_J)$, with the diagonal elements $w_j = |x_j|^{-1}$.
 For a slightly more general scenario with the $\ell_q$-norm ($q \leq 1$), we have
\be
 \|\bx\|_q^q = \bx^{\text{T}} \bW(\bx) \bx,
 \ee
where the diagonal elements are $w_j  = [|x_j|^2  +\varepsilon^2]^{q/2-1}$, and
$\varepsilon >0$   is a very small number needed  to avoid divergence for a small $x_j$ \citep{CandesIRLS2008}.

Similarly, for the non-overlapping group LASSO \citep{IRLS2014}
\be
 \|\bx\|_{q,1}^q = \sum_k \|\bx_{g_k}\|^q_q = \bx^{\text{T}} \bW(\bx) \bx,
 \ee
 where $\bW$ is a block  diagonal matrix $\bW =\diag(\bW_1, \bW_2, \ldots, \bW_K)$, with every diagonal sub-matrix $\bW_k$, given by
 $[\bW_k]_{jj} =(|x_j^{(g_k)}|^2 +\varepsilon^2)^{q/2-1}$.
 In practice, it is  sufficient to approximate the diagonal matrix $\bW(\bx)$ in a tensorized form by a rank one tensor, that is, with the TT rank of unity.

\begin{figure}[t]
\begin{center}
\includegraphics[width=0.99\textwidth]{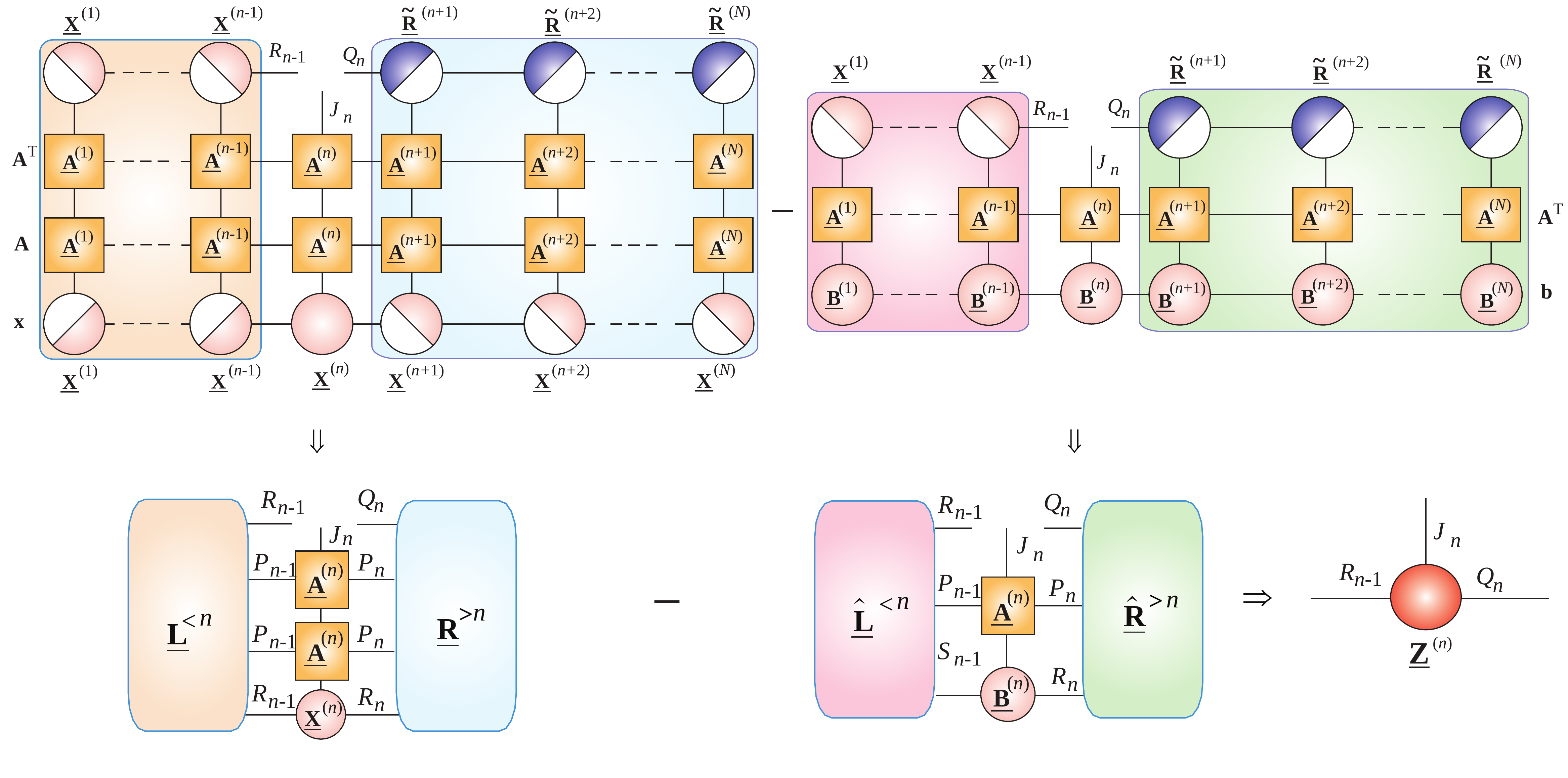}
\end{center}
\caption{Illustration of computation of the enrichment tensor $\underline \bZ^{(n)} \in \Real^{R_{n-1} \times J_n \times Q_n}$ in the AMEn algorithm for solving systems of huge linear equations. The bottom panel shows contraction of TT networks.}
\label{Fig:Ax=b2}
\end{figure}

In the simplest scenario, upon dividing the datasets into an appropriate number of sub-groups, the huge-scale standard LASSO problem  can be converted into a set of
 low-dimension
  group LASSO sub-problems, 
  given by (see also the first sub-networks in  Figure~\ref{Fig:Ax=b2})
\be
\min_{\bx^{(n)} \in \Real^{R_{n-1} J_n R_n}} \!\biggl(\!\bx^{(n)\;\text{T}} \overline \bA^{(n)} \bx^{(n)} -2 [\bx^{(n)}]^{\text{T}} \overline \bb^{(n)}
+\gamma \!\sum_{j_n=1}^{J_n} \|\underline \bX^{(n)}(:,j_n,:)\|_F \!\biggr)\nonumber\\
\ee
 which can be solved by any efficient algorithm for the group LASSO, e.g., by ADMM methods \citep{BoydFTML2011}.

\section[Truncated Optimization Approach]{Truncated Optimization Approach in TT Format}
\sectionmark{Truncated Optimization Approach}

A wide class of optimization problems can be solved by iterative algorithms which in general can be written as
\be
\bx_{k+1} =\bx_k + \Delta \bx_k=\bx_k + \bz_k,
\label{eq:gradstep}
\ee
where $\bx \in \Real^I$ is a vector of the  cost function $J(\bx)$ and vector
$\bz =\Delta \bx \in \Real^I$ is the update vector, which can take various forms, e.g.,

\begin{enumerate}

\item $\bz_k= - \eta_k  \nabla J (\bx_k)$,  for gradient descent methods, where $\eta_k >0$ is a learning rate, and $\nabla J(\bx_k)$ is the gradient vector of the cost function $J(\bx)$ at the current iteration;

\item $\bz_k= - \eta_k \bM^{-1}_k \nabla J (\bx_k)$, for preconditioned gradient descent, e.g., preconditioned conjugate gradient (CG), where $\bM \in \Real^{I \times I}$ is the preconditioned matrix;

\item $\bz_k= - \eta_k \bH^{-1}_k \nabla J (\bx_k)$ for quasi-Newton methods, where $\bH \in \Real^{I \times I}$ is an approximate Hessian;

\item  $\bz_k= - \eta_k  \bg_k$, for subgradient methods, where $\bg$ is subgradient.

\end{enumerate}

To date, the following  cost functions and corresponding gradients have been investigated in connection with tensor networks \citep{Lebedeva2011eigenTT,Ballani2013arnoldiHT,DolgovGMERES}

	\begin{itemize}

	\item $\nabla J(\bx_k) \propto \bA^{\text{T}}\bA\bx_k - \bA^{\text{T}}\bb$ for the squared residual $J(\bx) = \| \bA\bx - \bb \|^2$;

	\item $\nabla J(\bx_k) \propto \bA\bx_k - \bb$ for $J(\bx) = \bx^{\text{T}} \bA \bx - 2\bb^{\text{T}} \bx $. In this case the process  (\ref{eq:gradstep}) is often referred to as a preconditioned (non-stationary) Richardson iteration;

	\item $\nabla J(\bx_k) \propto \bA\bx_k - \bx_k \lambda_k$ with $\lambda_k = \bx_k^{\text{T}} \bA\bx_k $, for the Rayleigh quotient $J(\bx) = \bx^{\text{T}}\bA\bx$ with $\|\bx\|^2=1$, which is an orthogonal projection of the gradient of $J$ onto the tangent space of a sphere, or a Stiefel manifold in general \citep{Absil-book}.
	
\end{itemize}

\begin{figure}[t]
\begin{center}
\includegraphics[width=11cm]{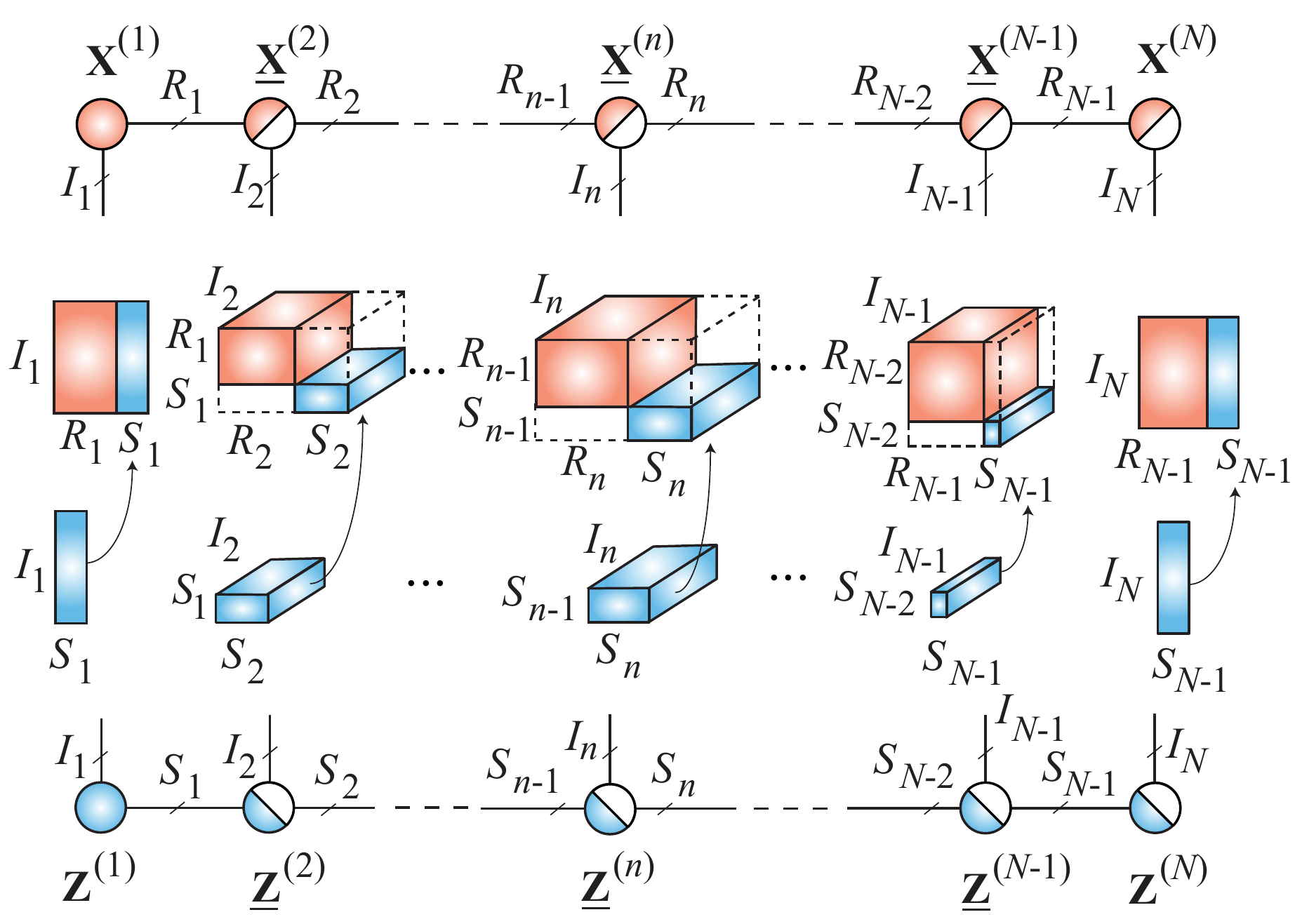}
\end{center}
\caption{A (preconditioned) gradient step within  a truncated iteration method in the TT format. The current solution vector $\bx_k$ and the (preconditioned) gradient vector $\bz_k = -\alpha_k \bM^{-1}_k \nabla J(\bx_k)$ are represented in the TT format, $\bx_k = \text{vec}( \llangle \bX^{(1)}, \underline \bX^{(2)}, \ldots, \bX^{(N)} \rrangle)$ and $\bz_k = \text{vec}( \llangle \bZ^{(1)}, \underline \bZ^{(2)}, \ldots, \bZ^{(N)} \rrangle)$. The addition of these two vectors can be performed by a partial (mode-2) direct sum between core tensors as $\underline \bX^{(n)} \oplus_2 \underline \bZ^{(n)} \in \Real^{(R_{n-1}+S_{n-1}) \times I_n \times (R_n + S_n)}$.}
\label{Fig:Truncated_Method_Gradient}
\end{figure}

Since huge scale optimization problems are intractable, we need to represent vectors and corresponding matrices in distributed tensor formats.  In a simplest scenario we can construct two TT networks, one representing vector $\bx$  and the other representing approximately gradient vector $\bz$, as illustrated in Figure~\ref{Fig:Truncated_Method_Gradient}.
Assuming that both huge vectors $\bx$ and $\bz$ are approximately represented in TT formats, the summation of two TT networks
leads to a new  TT network, for which the slices of the cores are given by
({\it cf.} Part 1 \citep{Part1})
\be
  \bX_{i_n}^{(n)} &\leftarrow& \left[ \begin{matrix}
 \bX_{i_n}^{(n)} & \0 \\ \0 & \bZ_{i_n}^{(n)} \end{matrix} \right],
 \quad (n=2,3\ldots,N-1)
 \ee
with the  border cores
 \be
\bX_{i_1}^{(1)} \leftarrow \left[\bX_{i_1}^{(1)},  \bZ_{i_1}^{(1)}\right], \qquad
 \bX_{i_N}^{(N)} &\leftarrow& \left[ \begin{matrix} \bX_{i_N}^{(N)}) \\ h \;  \bZ_{i_N}^{(N)} \end{matrix} \right],
 \ee
while for the right to left sweep we have
\be
\bX_{i_1}^{(1)} \leftarrow [\underline \bX_{i_1}^{(1)},  \bZ_{i_1}^{(1)}], \;\;  \bX_{i_n}^{(n)} &\leftarrow& \left[ \begin{matrix}
 \bX_{i_n}^{(n)} & 0 \\ 0 & \bZ_{i_n}^{(n)} \end{matrix} \right],
 \;\; (n=1,2,\ldots,N-1) \nonumber \\
 \bX_{i_N}^{(N)} &\leftarrow& \left[ \begin{matrix} \bX_{i_N}^{(N)}) \\ h \;  \bZ_{i_N}^{(N)} \end{matrix} \right].
 \ee
After the enrichment of the core $\underline \bX^{(n)}$, it  needs to be orthogonalized in order
to keep slice matrices orthogonal.
Figure~\ref{Fig:Truncated_Method_Gradient} illustrates the preconditioned gradient step (\ref{eq:gradstep}), where the solution vector $\bx_k$ and the preconditioned gradient $\bz_k = -\alpha_k \bM^{-1}_k \nabla J(\bx_k)$ are represented in TT formats. Note that the TT ranks $\{S_n\}$ of the preconditioned gradient $\bz_k$ are typically quite large due to the operations involved in the computation of the gradient.


Recall that the basic linear algebraic operations such as the matrix-by-vector multiplication can be performed efficiently, with a logarithmic computational complexity, via the TT representations for large matrices and vectors.
Since the algebraic operations on TT formats usually result in increased TT ranks, truncation (e.g., TT-rounding \citep{OseledetsTT11}) needs to be  subsequently performed.
 It should be emphasized that truncated iteration methods are not limited to the TT format but also apply to any low-rank tensor formats.

%


Let $\bT_\epsilon$ denote a truncation operator which truncates ranks of a tensor format with a tolerance $\epsilon>0$.
Truncated iteration methods can be described as a preconditioned gradient step combined with truncation, as
	$$
	\bx_{k+1} = \bT_\epsilon \left( \bx_k - \alpha_k \bT_{\eta} ( \bM^{-1}_k \nabla J (\bx_k) ) \right) , \quad \epsilon, \eta >0.
	$$
The rank truncation can be carried out by a hard thresholding scheme or a soft thresholding scheme \citep{Bachmayr2016}. In the iterative hard thresholding scheme, singular values below some threshold are set to zero as in the truncated SVD or TT rounding \citep{OseledetsTT11}.
In the iterative soft thresholding scheme, on the other hand, each of the singular values is shrunk to zero by the function $s_\kappa(x) = \text{sign}(x)\max(|x|-\kappa,0)$ for some $\kappa\geq 0$. The soft thresholding is equivalent to minimizing the nuclear norm of a matrix, which is a convex relaxation of rank. Both types of thresholding schemes adaptively change the ranks during iteration process.

The truncated iteration approach has been already investigated and extensively tested for several very large-scale data applications using various low-rank tensor formats, especially for:
	
\begin{itemize}
	
	\item Solving huge systems of linear equations and discretization of PDEs
			by iterative algorithms (e.g., Richardson, CG GMRES) and combined with CP format \citep{beylkin,Khoromskij2011pde}, CP and Tucker formats \citep{Billaud2014mres}, TT format \citep{Khoromskij2010quantics,DolgovGMERES}, HT format \citep{Kressner2011HT,Bachmayr2015ANR,Bachmayr2016soft}, and a subspace projection method combined with HT format \citep{Ballani2013arnoldiHT};
	
%
	
	\item A subspace projection method combined with HT format \citep{Ballani2013arnoldiHT};
			

	\item Computing  extreme eigenvalues and corresponding eigenvectors
		by Lanczos method combined with TT format \citep{Huckle2012Lanczos,Zhao_Lanczos},
		preconditioned inverse iteration with TT format \citep{Mach2013eig},
 block CG method (LOBPCG) combined with TT format \citep{Lebedeva2011eigenTT} and with HT format \citep{Kressner2011eigenHT};
		

\item Iterative hard thresholding for low-rank matrix/tensor completion \citep{Foucart:2013CS,Tanner2013matrixhard}.
	
\end{itemize}

An advantage of truncated iteration methods is that we do not need to construct tensor networks for a specific cost (loss) function. Moreover, theoretical global convergence properties hold under some restricted conditions (see \citep{Bachmayr2016} for more detail).
In addition, due to the analogy between low-rank truncation and  sparse signal estimation techniques, truncated iteration-type algorithms are also suitable for large-scale compressed sensing \citep{Blumensath2009CShard}.

\markright{3.10.\quad Riemannian Optimization}

However, truncated iteration methods have a few drawbacks. They require  the estimation of several  auxiliary  vectors.  In addition to the desired vector, $\bx$, right-hand vector, $\bz$, need to have  a low-rank  representation in the specific tensor format.
For example, in the GMRES method \citep{DolgovGMERES}, even if the solution vector and the right-hand side vector are well approximated by the TT format, the residuals and Krylov vectors involved in intermediate iterations usually have high TT ranks.
Moreover, all the core tensors and factor matrices of the solution vector have to be truncated at each iteration step, which often incurs high computational costs, especially when ranks are large.


\section{Riemannian Optimization for Low-Rank Tensor Manifolds}
\sectionmark{Riemannian Optimization}

Methods of Riemannian Optimization (RO) have been recently a subject of great interest in data analytics communities; see, for example, \citep{Absil-book,abg-trust-2007,bonnabel-stochriem-2013,CambierAbsil2016,bento2016iteration}.
Some optimization problems discussed in previous sections, can be naturally formulated on Riemannian manifolds, so as to directly benefit from the underlying geometric  structures that can be exploited to significantly reduce the cost of obtaining solutions.
Moreover, from a Riemannian geometry point of view, constrained optimization problems  can often be  viewed as unconstrained ones.
It is therefore natural to ask whether Riemannian
optimization can also help open up new research directions in conjunction with tensor networks; see for example \citep{IshtevaAHL11,kasai2015riemannian,Steinlechner_phd2016,zhou2016riemannian,sato2017riemannian}.

Riemannian optimization for tensors can be formulated in the following  generalized form
\be
\label{rieman:minim}
\min_{\mftensor{X} \in \Mr_{r}}\  J(\mftensor{X}),
\ee
where  $J(\underline \bX)$ is a scalar cost function of tensor variables and the search space, $\Mr_{r}$, is smooth, in the sense of a differentiable manifold structure. Note that to  perform optimization, the  cost function  needs to be defined for points on the manifold $\Mr_{r}$.
%
%
The optimization problem in the form \eqref{rieman:minim} is quite general, and most of the basic optimization problems can be cast into this form.

In our applications of RO, $\Mr_{r}$ is a certain low-rank tensor  manifold.
 One of the potential advantages of using
Riemannian optimization tools for tensors is therefore that the intrinsic geometric and algebraic properties of the manifold allow us to convert a constrained
optimization problems to an unconstrained one,
in other words, to perform unconstrained optimization on a constrained space.

{\bf Manifold structure for Tucker model.}
The Riemannian optimization (also called the differential geometric optimization)
has been successfully applied to the computation of the best Tucker approximation
\citep{Ishteva2009,IshtevaAHL11,Steinlechner_phd2016},
whereby the Tucker manifold structure can be formulated  through a set of Tucker tensors of fixed multilinear rank $r_{ML}=\{R_1,\ldots,R_N \}$
\be
\Mr_{r} := \{ \underline \bX \in \Real^{I_1 \times \cdots \times I_N} \; | \;
\mbox{rank}_{ML}(\underline \bX) = r_{ML}\}, \notag
\ee
which forms a smooth  embedded submanifold of $\Real^{I_1 \times \cdots \times I_N}$ of
dimension
\be
\dim \Mr_{r}= \prod_{n=1}^{N} R_n + \sum_{n=1}^{N}( R_n I_n -R^2_n).
\ee

{\bf Manifold structure for  TT model.}
\noindent The set of TT tensors of fixed TT-rank, $r_{TT}=\{R_1,\ldots,R_{N-1} \}$, given by
\be
\Mr_{r} := \{ \underline \bX \in \Real^{I_1 \times \cdots \times I_N} \; | \;
\mbox{rank}(\underline \bX)_{TT} = r_{TT}], \notag
\ee
forms a smooth  embedded
submanifold of $\Real^{I_1 \times \cdots \times I_N}$ of dimension
\be
\dim \Mr_{r}= \sum_{n=1}^{N} R_{n-1}I_n R_{n} - \sum_{n=1}^{N-1} R^2_n,
\ee
with $R_N=1$.
If the inner product of two tensors induces a Euclidean metric on the embedding space
$\Real^{I_1 \times \cdots \times I_N}$, then the above  submanifold  $\Mr_{r}$ is a Riemannian manifold \citep{Holtz-manifold-2012,Steinlechner_phd2016,Uschmajew-Vander2013}.

Similar results are also available for the more general hierarchical Tucker models.
For example, \citep{Uschmajew-Vander2013} developed  the manifold structure for the HT  tensors,
while \citet{Lubich-Schneider13} developed the concept of dynamical low-rank approximation
for both HT  and TT formats. Moreover, Riemannian optimization in the Tucker and TT/HT  formats
has been successfully applied to large-scale tensor completion problems
\citep{sh-htcompl-2013,kasai2015riemannian,zhou2016riemannian}.


It is important to emphasize that the condition of full TN rank is a pre-requisite for the sets $\Mr_{r}$ to be smooth manifolds. Otherwise,  a set of
tensors for which the ranks are only bounded from above by $R_{\max}$ would no longer be a smooth manifold, but an algebraic
variety, for which special optimization methods are needed, at least for theoretical justification.

The fundamental and basic algorithms for \emph{Riemannian optimization} \citep{Absil-book} are quite attractive,
but are still not widely used as they are more technically complicated than, e.g., standard
gradient descent or conjugate gradient (CG) used in the  classical optimization algorithms of
low-rank tensor manifolds. However, in many cases Riemannian manifolds
can be relatively easily implemented and often provide much better performance than the standard algorithms
\citep{KressnerSV2013,kasai2015riemannian}.

The MANOPT package \citep{boumal-manopt-2014} (for a Python version see \citep{kw-pymanopt-2016})
provides a useful interface for many standard matrix Riemannian optimization techniques, however,
extensions for tensors (especially in high dimensions) are not easy and should be performed very carefully.
The key idea is to work in the full space, but with the restrictions:
i) the solution lies on a non-linear smooth manifold, and
ii) the manifold condition is enforced implicitly
(i.e., the derivatives with respect to the parametrization are not explicitly used).

In order to apply optimization algorithms based on line search, we must consider a direction on
a manifold. The tangent space at $\underline \bX \in \Mr_r$ is the set of
all tangent vectors (directions) at  $\underline \bX$, denoted by $T_{\underline \bX} \Mr_r$.
The tangent space is a linear space, which, at a point on the manifold, provides us with a vector space of tangent vectors
that allow us to find  a search direction on the manifold.
A Riemannian metric allows us to compute the angle and length of directions (tangent vectors).
Roughly speaking, the optimization is performed on the tangent space, by either performing linear search or building a local
model, in order to estimate the tangent vector and perform  the next iteration on the manifold. The retraction
operator provides a method to map the tangent vector to the next iterate (see the next section).

The smoothness condition of the manifold is crucial, since this allows us to define a tangent space, $T_{\mftensor{X}} \Mr_{r}$,
for each point on the manifold, which locally approximates the manifold with second-order accuracy.
Such a linearization of the manifold allows us to compute
the next approximation restricted to the tangent space, which is a linear subspace.
We shall next describe the general setting in which Riemannian optimization resides.

\subsection{Simplest Riemannian Optimization Method: Riemannian Gradient Descent}
\subsectionmark{Riemannian Gradient Descent}

 Consider a smooth manifold, $\Mr_{r}$, in a vector space of tensors.
 The simplest optimization method is gradient descent which has the form
 \begin{equation*}
     \mftensor{X}_{k+1} = \mftensor{X}_k - \alpha_k \nabla J(\mftensor{X}_k),
 \end{equation*}
 where $\alpha_k > 0$ is the step size.
 For a sufficiently small $\alpha$, the value of the  cost function $J(\mftensor{X})$ will decrease.
 However,  this does not provide a restriction that $\mftensor{X}_{k+1}$ should lie on a manifold;
although $\mftensor{X}_k$ is on the manifold, the step
 in the direction of the negative gradient will leave the manifold.
 If the manifold is a linear subspace, then the solution is quite simple ---
 we just need to take a projected gradient step, i.e.,
 project $\nabla J(\mftensor{X}_k)$ to this subspace.
 For a general smooth manifold we can do the same,
 but the projection is computed onto the tangent space of the manifold at the current iteration point, that is
 \begin{equation}\label{riem:gradstep}
     \mftensor{X}_{k+1} = \mftensor{X}_k - \alpha_k P_{T_{\mftensor{X}_k} \Mr_{r}}( \nabla J(\mftensor{X}_k)).
 \end{equation}
 The continuous version of \eqref{riem:gradstep} will have the following projected gradient flow
 \begin{equation}\label{riem:gradflow}
     \frac{d \mftensor{X}}{dt} =  -  P_{T_{\mftensor{X}} \Mr_{r}} (\nabla J(\mftensor{X})),
 \end{equation}
 and it is quite easy to see that if $\mftensor{X}(0)$ (initial conditions) is on the manifold, then the whole trajectory $\mftensor{X}(t)$ will be on the manifold.
 This is, however, not guaranteed for the discretized version in \eqref{riem:gradstep},
 which can be viewed as a forward Euler scheme applied to \eqref{riem:gradflow}.
 This also has important computational consequences.

 An iterative optimization algorithm involves the computation of a search
direction and then performing a step in that direction.
Note that all  iterations on a manifold $\Mr_r$
should be performed by following geodesics{\footnote{A geodesic is the shortest path on the manifold, which  generalizes the concept of straight lines in Euclidean space.}}. However, in general geodesics  may be  either expensive to compute or  even not available in a closed form.
To this end, in practice, we relax the constraint of moving along geodesics by applying the concept of {\it Retraction}, which is any map
$R_{\bX} : T_{\underline \bX} \Mr_r  \rightarrow \Mr_r$ that locally approximates a geodesic, up to first order, thereby reducing the computational cost of the update and  ensuring convergence of iterative algorithms.

\begin{figure}[t]
\centering
\includegraphics[width=6.5cm]{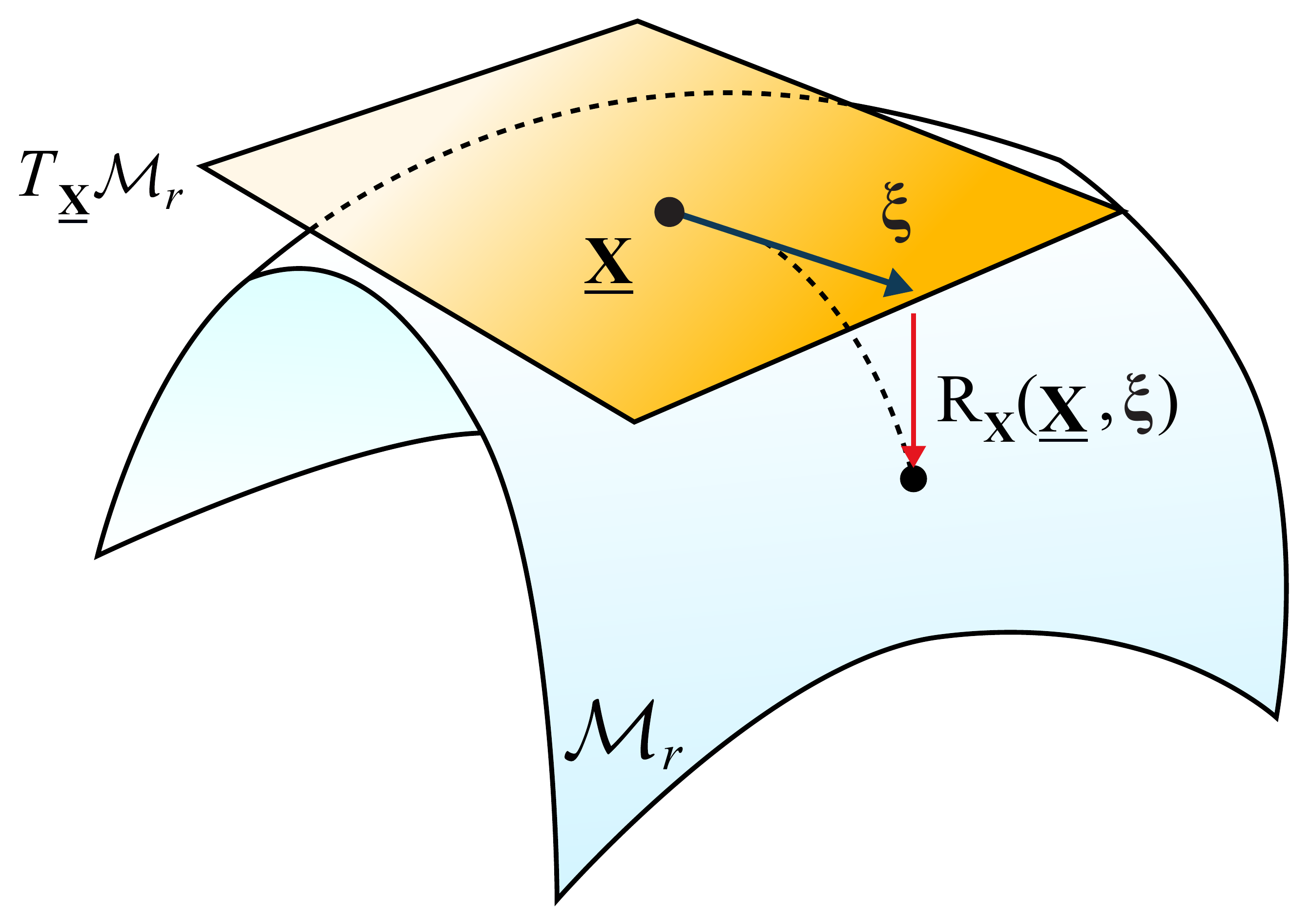}
\caption{The concept of retraction for a smooth submanifold. The iteration process for Riemannian gradient descent for optimization on a smooth submanifold $\Mr_{r}$ embedded in a Euclidean space consists of three-steps: 1) an update in the Euclidean space (which is not computed explicitly), 2) a linear projection of the gradient to the tangent space $T_{\mftensor{X}} \Mr_{r}$, and 3) a retraction to the manifold $\Mr_{r}$. See also Equation \eqref{riem:gradstepfinal} and \citep{Kressner2015}.}
\label{fig_riemannian_gradient}
\end{figure}

More precisely the  Retraction can be defined, as follows
\begin{definition} \citep{adm-riem-2002}
 A mapping $R_{\bX}$ is called the retraction, if
 \begin{enumerate}
     \item $R_{\bX}(\mftensor{X}, \mftensor{Z}) \in \Mr_{r}$ for $\mftensor{Z} \in T_{\mftensor{X}} \Mr_r$,
     \item $R_{\bX}$ is defined and smooth in a neighborhood of the zero section in $T_{\mftensor{X}} \Mr_r$, 
     \item  $R_{\bX}(\mftensor{X}, \mftensor{0}) = \mftensor{X}$ for all $\mftensor{X} \in \Mr_{r}$,
     \item $\frac{d}{dt} R_{\bX}(\mftensor{X}, t\mftensor{Z}) \mid_{t = 0} =
     \mftensor{Z}, \quad \mbox{for all }
     \mftensor{X} \in \Mr_r \quad \mbox{and } \;\; \mftensor{Z} \in T_{\mftensor{X}} \Mr_r$. 
 \end{enumerate}
 \end{definition}

 For low-rank matrices and tensors with fixed TT-ranks, the simplest retraction is provided by respective SVD and TT-SVD algorithms, however, there are many other types of retractions; for more detail we refer to the survey \citep{AbsilO2014retraction}.

 Finally, the Riemannian gradient descent method has the form
 \begin{equation}\label{riem:gradstepfinal}
     \mftensor{X}_{k+1} = R_{\bX}( \mftensor{X}_k , \; -\alpha_k P_{T_{\mftensor{X}_k} \Mr_{r}}  ( \nabla J(\mftensor{X}_k))).
 \end{equation}
Figure~\ref{fig_riemannian_gradient} illustrates the three-step procedure within the Riemannian gradient descent iteration given in \eqref{riem:gradstepfinal}.
From the implementation viewpoint, the main difficulty is to compute the projection of the gradient $P_{T_{\mftensor{X}_k} \Mr_{r}}(\nabla J(\mftensor{X}_k) )$.
For the optimization problems mentioned above, this projection can be computed without forming the full gradient.


\subsection{Advanced Riemannian Optimization Methods}

 {\bf Vector transport.} More advanced optimization techniques, such as conjugate-gradient (CG)-type methods
 and quasi-Newton methods,
 use directions from previous iteration steps. However, these directions lie in different tangent spaces which
 are different from the corresponding optimization in Euclidean spaces.
To this end, we need to employ the concept of \emph{vector transport}, which plays a crucial role in Riemannian optimization.
The idea of vector transport is quite simple: it represents a mapping from one tangent space $T_{\mftensor{X}}
\mathcal{M}_{r}$, to another tangent space, $T_{\mftensor{Y}} \mathcal{M}_{r}$.
 For low-rank matrix/tensor manifolds,
the orthogonal projection to a tangent space $P_{T_{\mftensor{Y}} \mathcal{M}_{r}}$ is an example of a vector transport.

\begin{figure}[t]
\centering
\includegraphics[width=8.5cm]{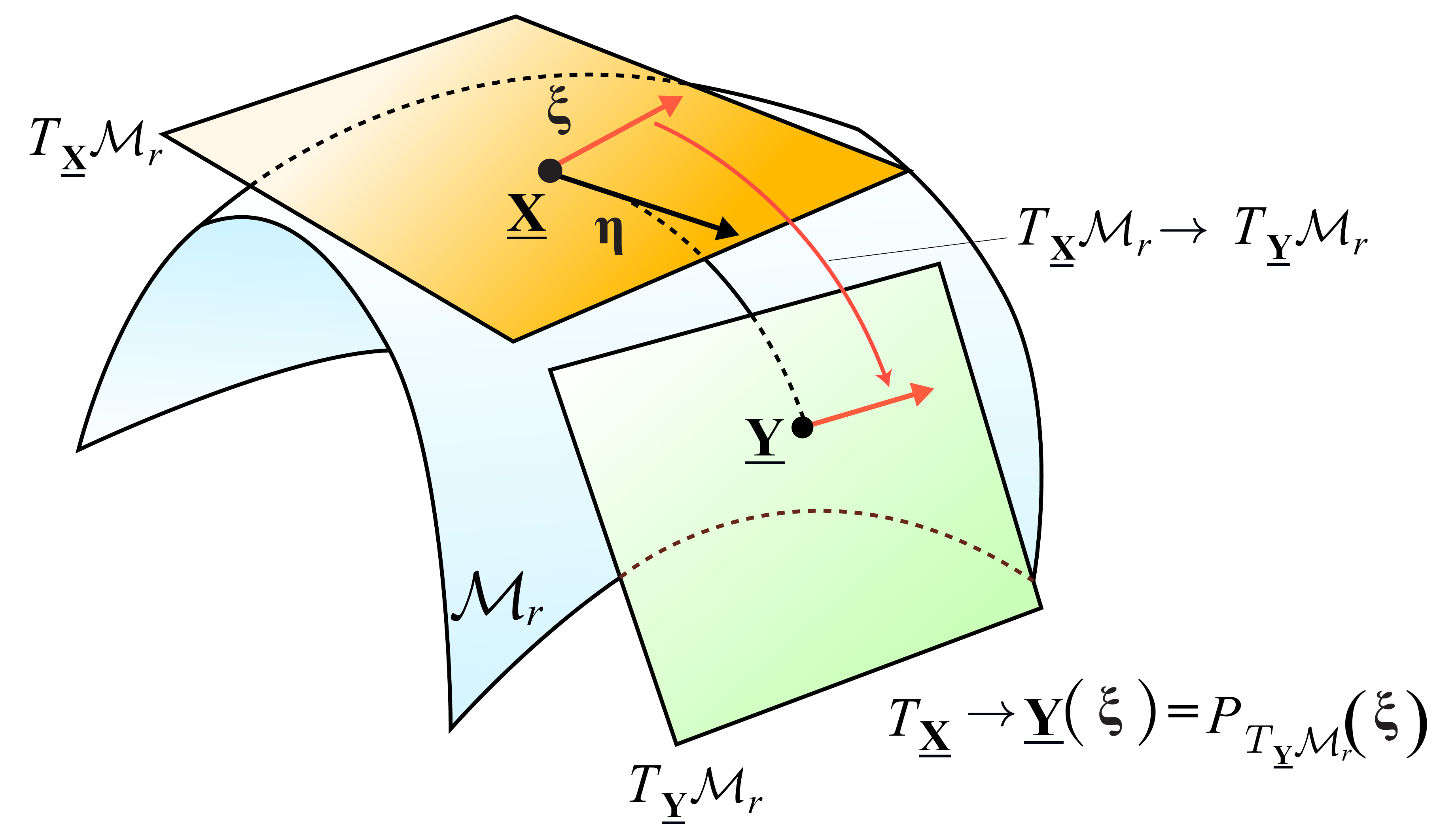}
\caption{Concept of the vector transport for the Riemannian CG method for the optimization on a manifold $\Mr_{r}$. The  update direction, $\mbi \eta$, at the previous iteration, which belongs to the tangent space $T_{\underline \bX} \Mr_{r}$, is mapped to a tangent vector in the tangent space $T_{\underline \bY} \Mr_{r}$, via a suitably defined mapping (i.e., the vector transport) denoted by $T_{\underline \bX \rightarrow \underline \bY}$.
See also Algorithm \eqref{riem:alg1}.}
\label{fig_riemannian_cg}
\end{figure}


{\bf Riemannian CG Method.}
Figure~\ref{fig_riemannian_cg} illustrates the vector transport during a Riemannian CG iteration, which
 transforms a tangent vector in $T_{\underline \bX} \Mr_{r}$ at a certain point $\underline \bX\in\Mr_{r}$ to another tangent space $T_{\underline \bY} \Mr_{r}$ at $\underline \bY\in\Mr_{r}$.
In the Riemannian CG iteration, the previous update direction $\mbi \eta_{k-1} \in T_{\underline \bX_{k-1}} \Mr_{r}$ is mapped to a tangent vector in $T_{\underline \bX_{k}} \Mr_{r}$, and
the transformed vector is combined with the Riemannian (projected) gradient vector to complete the CG iteration.
Now, by using the projection onto the tangent space, retraction, and vector transport, we can immediately implement a general
nonlinear CG method (Riemannian CG method) (see Algorithm~\ref{alg:RCG1}).

\begin{algorithm}[t!]
\caption{Riemannian CG method}\label{riem:alg1}
{\small
\begin{algorithmic}[1]
\STATE Initial guess $\mftensor{X}_0 \in \Mr_{r}$
\STATE $\mbi \eta_0 := - P_{T_{\mftensor{X}_0} \Mr_{r}}( \nabla J(\mftensor{X}_0) )$ \\
$\alpha_0 = \arg \min_{\alpha} J(\mftensor{X}_0 + \alpha \mbi \eta_0)$
\STATE  $\mftensor{X}_1 = R_{\bX}(\mftensor{X}_0, \alpha_0 \mbi \eta_0).$
\FOR{$k=1, \ldots$}
\STATE $\mbi \xi_k := P_{T_{\mftensor{X}_k} \Mr_{r}}( \nabla J(\mftensor{X}_k) )$
\STATE $\mbi \eta_k := -\mbi \xi_k + \beta_k T_{\mftensor{X}_{k-1} \rightarrow \mftensor{X}_k} \mbi \eta_{k-1}$
\STATE  $\alpha_k  = \arg \min_{\alpha} J(\mftensor{X}_k + \alpha \mbi \eta_k)$
\STATE Find the smallest integer $m \geq 0$ such that \\
 $J(\mftensor{X}_k) - J(R_{\bX}(\mftensor{X}_k ,2^{-m} \alpha_k \mbi \eta_k )) \geq \delta \; \langle \mbi \eta_k, 2^{-m} \alpha_k \mbi \eta_k \rangle$
 \STATE  $\mftensor{X}_{k+1} := R_{\bX}(\mftensor{X}_k, 2^{-m} \alpha_k \mbi \eta_k).$
\ENDFOR
\end{algorithmic}
}
\label{alg:RCG1}
\end{algorithm}

The main problem  in Riemannian CG methods is the computation of the parameter $\beta_k$  within the conjugate gradient direction.
To this end, as suggested in \citet{Kressner2015} the Polak-Ribiere update formula \citep{nw-book-2006} can be adapted to Riemannian optimization \citep{Absil-book}.

\subsection{Riemannian Newton Method}

The implementation of Riemannian Newton-type algorithms is much more complicated, but follows
the standard concept of Newton optimization methods. The local search direction, $\mbi \xi_k$, is determined from the
correction equation \citep{Absil-book} as
\begin{equation}\label{riem:newton}
 \bH_{\mftensor{X}_k} \mbi \xi_k = -P_{T_{\tX_k} \Mr_{r}} \nabla J(\mftensor{X}_k),
\end{equation}
where the Riemannian Hessian, $\bH_{\mftensor{X}_k} : T_{\mftensor{X}_k} \Mr_r \rightarrow T_{\mftensor{X}_k} \Mr_{r}$, is defined as \citep{amt-hessian-2013,Kressner2015}
\begin{equation}\label{riem:hessian}
 \bH_{\mftensor{X}_k} = P_{T_{\tX_k} \Mr_{r}}  \left(\nabla^2 J(\mftensor{X}_k) + P'_{\mftensor{X}_k} \nabla J(\mftensor{X}_k)\right) P_{T_{\tX_k} \Mr_{r}} ,
\end{equation}
and $P'_{\mftensor{X}_k}$ is the derivative of the projector $P_{T_{\tX} \Mr_{r}}$ with respect to $\mftensor{X}$.
The second term in the sum corresponds to the curvature of the manifold.
If the true minimizer is on the manifold, then $\nabla J(\mftensor{X}_*) = \mftensor{0}$, and
at least in the vicinity of the solution the second term plays no role.
However, this term can cause problems with stability if the solution is close to a singular point,
where the tangent space can exhibit considerable variation.
In many cases, only the first term in \eqref{riem:hessian} is used, leading to
the \emph{approximate Newton method}, also known as the Gauss-Newton method in constrained optimization, given by
\begin{equation}
\widehat{\bH}_{\mftensor{X}_k} = P_{T_{\tX_k} \Mr_{r}} \, \nabla^2 J \, P_{T_{\tX_k} \Mr_{r}}.
\end{equation}
 If $\nabla^2 J$ is positive definite (e.g., for strongly convex $J(\underline \bX)$) the matrix
$\widehat{\bH}_{\mftensor{X}}$  is also nonnegative definite.

The Gauss-Newton method has a  simple interpretation: we approximate the manifold by a tangent space,
and seek for the minimum of a specific cost function $J(\underline \bX)$ on the tangent space.
However, the theoretical justification of such methods, even for relatively simple low-rank matrix cases,
is still an open problem.
Recently, in order to achieve a better global convergence, an alternative version of  the Riemannian Newton method called the Trust-Region scheme was developed  by
\citet{abg-trust-2007,IshtevaAHL11,BoumalAbsil2011}.

Retraction and vector transport are  critical operations to the success of sophisticated  Riemannian optimization algorithms,
such as Riemannian CG,  Newton and/or quasi-Newton methods. Retraction is used to obtain the next iterate and
vector transport is used to compare tangent vectors in different tangent spaces and to transport
operators from one tangent space to another tangent space.

{\bf Parametrization.} When considering the tensor case, the whole concept above applies to the matrix-based tensor formats, such as
the Tucker format and TT format.
However, the concept does not apply to the CP decomposition, since the set
of all tensors with a given canonical rank does not form a manifold (CP decomposition is unique).

\subsection{Low-Rank TT-Manifold}

 Consider a set of $N$th-order
tensors in TT formats, with all TT-ranks equal to $(R_1,\ldots, R_{N-1})$. This set can also be thought of as an intersection
of $(N-1)$ low-rank matrix manifolds, since the $n$th TT-rank of a tensor $\underline \bX$ is
equal to the matrix rank of its $n$th mode unfolding $\bX_{<n>}$.

A tensor in the TT-format can be parameterized in the form
$$ \underline \bX(i_1, \ldots, i_N) = \bX^{(1)}_{i_1} \cdots \bX^{(N)}_{i_N},$$
where $\bX^{(n)}_{i_n} \in \Real^{R_{n-1} \times R_n}$ are slices of TT-cores $\tX^{(n)}$ .

%

{\bf Tangent space of TT.} The tangent space at a point $\underline \bX$ in the TT-manifold is defined as a set of tensors of the form
\begin{equation}\label{riem:tttang}
  \delta \underline \bX (i_1,\ldots,i_N) = \delta \bX^{(1)}_{i_1} \bX^{(2)}_{i_2} \cdots  \bX^{(N)}_{i_N} + \cdots + \bX^{(1)}_{i_1} \bX^{(2)}_{i_2}  \cdots  \delta \bX^{(N)}_{i_N},
\end{equation}
and is parameterized by the variations of the TT-cores $\delta \underline \bX^{(n)}$.


We denote the QR decompositions of unfoldings $\bX^{\leq n}$ and $\bX^{\geq n+1}$ as
$$ \bX^{\leq n} = \bQ_{\leq n} \bR_{\leq n}, \quad (\bX^{\geq n+1})^{\text{T}} = \bQ_{\geq n+1} \bR_{\geq n+1},$$
and the orthogonal projectors  onto the column spaces of these matrices can be introduced in the form
$$ \bP_{\leq n} = \bQ_{\leq n} \bQ^{\text{T}}_{\leq n}, \quad \bP_{\geq n} = \bQ_{\geq n} \bQ^{\text{T}}_{\geq n}.
$$
Using these projectors the orthogonal projector onto the tangent space, $P_{T_{\tX} \mathcal{M}_{r}}(\tZ)$, can be written as \citep{Holtz-manifold-2012,LubichTITT14}
\be
P_{T_{\tX} \mathcal{M}_{r}}(\tZ) = \sum_{n = 1}^{N}    P_{\tX}^{\leq n}( P_{\tX}^{\geq n+1} (\tZ) )\, ,
\ee
where
\be
P_{\tX}^{\leq n}(\tZ)  &=& \Ten(\bP_{\leq n}  \bZ_{<n>}) \, ,\\
P_{\tX}^{\geq n+1}(\tZ)  &=& \Ten( \bZ_{<n>}  \bP_{\geq n+1} ),
\ee
and $\Ten(\bZ)$ is an inverse operator of the matricization of the tensor $\tX$, i.e., creating a tensor of the same size as the tensor $\tX$.
%

It should be noted that any point
on a tangent space has TT-ranks bounded by $(2R_1, \ldots, 2R_{N-1})$,
thus the computation of the TT-rank $(R_1, \ldots, R_{N-1})$ approximation
is possible by running the rounding procedure (or the TT-SVD) at the
$\mathcal{O}(NIR^3)$ cost, where we
assumed that $I_n = I, \quad R_n = R.$
Therefore, many Riemannian optimization algorithms  for low rank matrices remain formally the same, making the Riemannian optimization
a promising  and powerful tool for optimization over TT-manifolds.

%

\subsection{Dynamical Low-Rank Matrix/Tensor Approximation}

Riemannian optimization is intricately connected with the concept of
\emph{dynamical low-rank approximation}, and although it was introduced
in mathematics  by \citet{lubich-koch-dynmat-2007}, its roots are in
 quantum mechanics and can be traced back to the so-called
Dirac-Frenkel principle. Indeed, some fundamental results can be found in the quantum molecular dynamics community which has used  the so-called
\emph{Multi Configurational Time-Dependent Hartree}  (MCTDH) 
method already since the 1990s \citep{mathe-mctdh-1992}.
As an example, consider a time-varying matrix $\bA(t)$ which we wish to approximate  by a low-rank matrix $\bX(t)$. Of course, this can be achieved
on a sample-by-sample basis, by computing SVD for each $t$, but this does not involve any ``memory'' about the previous steps.

A Dirac-Frenkel principle states that the optimal trajectory should not minimize the distance
$$ \Vert \bA(t) - \bY(t) \Vert$$
over all matrices of rank $R$, but instead the local dynamics given by
\be \Vert \frac{d\bA}{dt} - \frac{d\bY}{dt} \Vert.
\label{eq:dX}
\ee
For the Euclidean norm, the minimization of the cost function (\ref{eq:dX}) leads
the following differential equations on manifold $\Mr_{r}$
\begin{equation}\label{riem:dyntime}
  \frac{d\bY}{dt} = P_{T_\bY \mathcal{M}_{r}} \frac{d\bA}{dt},  \quad \bY(0) \in \Mr_r.
\end{equation}
In other words, the approximated trajectory $\bY(t)$ is such that
 the velocity $\frac{d\bY}{dt}$ always lies in the tangent space,
and thus the trajectory always stays on the manifold.
The dynamical approximation in \eqref{riem:dyntime} allows us to construct
(and actually strictly define) low-rank approximation to the dynamical systems.
If $\bA(t)$ admits the representation
\begin{equation}
  \frac{d\bA}{dt} = F(\bA, t),
\end{equation}
which is a typical case in quantum mechanics, where the idea was first proposed,
we obtain a dynamical system for the point on a low-rank manifold in the form
\begin{equation}\label{riem:dynsys2}
  \frac{d\bY}{dt} = P_{T_\bY \mathcal{M}_{r}} F(\bY, t).
\end{equation}
Such a system can be efficiently treated using only the parametrization of a point on a such manifold,
thereby greatly reducing the cost (and if the manifold approximates the solution well,with a sufficiently good accuracy).
Note that the continuous variant of the Riemannian gradient  in \eqref{riem:gradstep}
can be written exactly in the form \eqref{riem:dynsys2}, with
$$ F(\bY, t) = -\nabla J.$$
Retraction is needed only in the discrete case, since in the continuous case the exact solution
lies exactly on the manifold. Thus, if we have an efficient way to solve \eqref{riem:dynsys2} numerically,
that would give an efficient way of solving Riemannian optimization problems.

There are two straightforward ways of solving \eqref{riem:dynsys2}. The original approach of \citet{lubich-koch-dynmat-2007} (later generalized to the Tucker and TT models \citep{lubich-koch-dynten-2010})
is to write down ordinary differential equations for the parameters $\bU(t), \bS(t), \bV(t)$ of the SVD like decomposition in the form
$$ \bY(t) = \bU(t) \bS(t) \bV^{\text{T}}(t).$$
The second approach is also straightforward: apply any time integration scheme to the equation \eqref{riem:dynsys2}.
In this case, a standard method will yield the solution which is not on the manifold, and a retraction would be needed.
In \citet{LubichO2014splitting} a simple and efficient solution to this problem, referred to as the projector-splitting scheme,
was proposed, based on the special structure of the manifold. This  stable and explicit second-order scheme for the integration of \eqref{riem:dynsys2} examines  the projector onto the tangent space
$$ P_{T_\bY \Mr_{r}}  \bZ = \bU\bU ^{\text{T}} \bZ + \bZ \bV\bV^{\text{T}} - \bU\bU^{\text{T}} \bZ \bV\bV^{\text{T}} = \bP_U + \bP_V - \bP_{UV},$$
and represents it as a sum of three projectors, in order to apply the classical Strang-Marchuk splitting scheme.
For each projector, the equations of dynamical low-rank approximations can now be easily integrated.
Indeed, the equation
\begin{equation*}
  \frac{d\bY}{dt} =  \frac{d\bA}{dt} \bV\bV^{\text{T}},
\end{equation*}
has a simple solution; the $\bV$ component does not change and the new point is just
\begin{equation}\label{riem:projsplit}
  \bY_1 = \bY_0 + (\bA_1 - \bA_0) \bV,
\end{equation}
where  $\bA_1 = \bA(t + h), \quad \bA_0 = \bA(t).$
Similar formulas are valid for $\bP_V$ and $\bP_{UV}$. What is not trivial  is the order in which
the splitting steps are taken: we should make a step in $\bP_U$, then in $\bP_{UV}$, and then
in $\bP_V$. This leads to a much better convergence, and moreover an exactness result \citep{LubichO2014splitting} can
be proven
if $\bA_1$ and $\bA_0$ are on the manifold, and $\bY_0 = \bA_0$, then this scheme is exact.

%

Through the link with the Riemannian optimization, the projector-splitting scheme can be also used here,  which can be viewed as a second-order retraction \citep{AbsilO2014retraction}. A simple form
can be obtained for a given $\bX$ and a step $\bZ$, where for the matrix case the retraction is
$$ \bU_1  \bS_{\frac{1}{2}} = \QR((\bX + \bZ) \bV_0),
\quad
\bV_1 \bS^{\text{T}}_1 = \mathrm{QR}((\bX+\bZ)^{\text{T}} \bU_1).$$
This is just one step of a block power method, applied to $(\bX+\bZ)$. Although this is only a crude approximation
of the SVD of $(\bX+\bZ)$ which is the standard retraction, the projector-splitting retraction is also a second-order retraction.

The generalization to the TT  network was proposed by \citet{LubichTITT14},
and can be implemented within the framework of sweeping algorithms, allowing for the efficient TT-approximation of  dynamical systems
and solution of optimization problems with non-quadratic functionals.

\subsection{Convergence of Riemannian Optimization Algorithms}

 Convergence theory for the optimization over non-convex manifolds is much more complicated than the corresponding theory in Euclidean spaces,
 and far from complete. Local convergence results follow from the general theory \citep{Absil-book},
 however the important problem of the curvature and singular points is not yet fully addressed. One way
 forward is to look for the desingularization \citep{ko-desing-2016pre}, another technique is
 to employ the concept of tangent bundle. Even if the final point is on the manifold, it is not clear
 whether it is better to take a trajectory on the manifold. An attempt to study the global convergence was presented
 in \citep{ko-linconv-2016pre}, and even in this case,  convergence to the spurious local minima
 is possible in a carefully designed example. Also,  convergence should be considered together with
 low-rank approximation to the solution itself. When we are far from having converged, a rank-one approximation to the solution will be satisfactory, then the rank should gradually increase. This was observed
 experimentally in \citep{ko-linconv-2016pre} in the form of a ``staircase'' convergence. So, both the local
 convergence, in the case when the solution is only approximately of low-rank, and the global convergence of Riemannian optimization have to be studied, but numerous applications to  machine learning problems (like matrix/tensor completion)
 confirm its efficiency, while the applications to tensor network optimization have just started to appear.

\subsection{Exponential Machines for Learning Multiple Interactions}
\subsectionmark{Exponential Machines}

 Riemannian optimization methods have been recently demonstrated as a powerful approach for solving a wide variety of standard machine learning problems.
 In this section, we will briefly describe  one such promising application \citep{novikov2016exponential}.
Consider machine learning tasks with categorical data, where for
improved performances it is important to properly model complex interactions between
the multiple features. 
Consider the training data, $\{ (\bx_m, y_m) \}_{m=1}^M$, where $\bx_m = [x_{1,m}, \ldots, x_{N,m}]^{\text{T}}$ is an $N$-dimensional feature vector and $y_m $  the target value of the $m$-th object.

In the special case of $N=3$, the linear model which incorporates interactions between the features $x_1,x_2$ and $x_3$, can be described by
	\begin{equation*}
	\begin{split}
	\widehat{y}(\bx) & = w_{000} + w_{100} x_1 + w_{010} x_2 + w_{001} x_3 \\
				   & + w_{110} x_1 x_2 + w_{101} x_1 x_3 + w_{011} x_2 x_3 + w_{111} x_1 x_2 x_3,
	\end{split}
	\end{equation*}
where the parameters $w_{ijk}$ denote the strength of interactions between the features $x_i,x_j$ and $x_k$  (see also Chapter \ref{Chapter1}).

Note that in the general case of $N$ features, the number of pair-wise interactions is $\mathcal{O}(N^2)$, and the number of all possible interactions grows exponentially as $\mathcal{O}(2^N)$.
Such interactions occur, for example,  in  language/text analysis or in sentiment analysis where each individual word can interact with other words.
To deal with the exponentially large number of parameters,  \citet{novikov2016exponential,Stoudenmire_optim} proposed the so-called Exponential Machines (ExM), where a large tensor of parameters is represented compactly in the TT format in order to provide low-rank regularization of the model and a reduction of the problem to a manageable scale. The  so obtained linear model can be described in a general tensor form as
	$$
	\widehat{y}(\underline \bX) = \langle \underline \bW, \underline \bX \rangle
	= \sum_{i_1=0}^1\cdots \sum_{i_N=0}^1 w_{i_1,\ldots ,i_N} \prod_{n=1}^N x_n^{i_n},
	$$
where $\underline \bW \in \Real^{I_1 \times \cdots \times I_N}$ $(I_n=2\;, \forall n)$ is the tensor of parameters and $\tX= [1, x_1]^{\text{T}} \circ \cdots \circ [1, x_N]^{\text{T}}$
is the rank-1 tensor of all interaction terms.

The TT representation of the weight tensor is given by
\be
	\underline \bW &=& \underline \bW^{(1)} \times^1 \underline \bW^{(2)} \times^1 \cdots \times^1 \underline \bW^{(N)}  \notag \\
&=& \llangle \underline \bW^{(1)},\underline \bW^{(2)},\ldots,
	\underline \bW^{(N)} \rrangle, \notag
\ee
where $\underline \bW^{(n)} \in \Real^{R_{n-1} \times I_n \times R_n}$ are third-order core tensors for $n=1,\ldots,N$ with $R_0=R_N=1$.

Such a linear model can be applied to various machine learning problems for which the parameters are learned by minimizing the cost function
	\be
	J(\underline \bW) = \sum_{m=1}^M l \left( \left\langle \underline \bW, \underline \bX_m \right\rangle  , y_m \right) + \lambda \| \underline \bW \|^2_F,
	\label{eqn_penalized_loss_novikov}
\ee
where $l(\widehat{y}, y)$ is a squared loss function (in least squares), hinge loss (in support vector machines), and logistic loss (in logistic regression), and
$\| \underline \bW \|^2_F = \langle \underline \bW , \underline \bW \rangle$.

Riemannian gradient descent can be applied to solve  very efficiently the problem in \eqref{eqn_penalized_loss_novikov}, where the gradient of the cost function is computed as
	\be
	\nabla J(\underline \bW) = \sum_{m=1}^M \frac{\partial l}{\partial \widehat{y}}
\underline \bX_m + \lambda \underline \bW .
	\ee
Due to its linearity, the projection of the gradient to the tangent space is expressed as the sum
	\be
	P_{T_{\tW} \Mr_{r}} (\nabla J(\underline \bW)) = \sum_{m=1}^M
 \frac{\partial l}{\partial \widehat{y}} P_{T_{\tW} \Mr_{r}}(\underline \bX_m)
 + \lambda \underline \bW .
	\ee
%
Moreover, a stochastic version of the Riemannian gradient descent, summarized in Algorithm~\ref{alg_stochastic_riemannian_novikov}, can be derived by sampling randomly a mini-batch (subset) of data points from the total of $M$ training data, which makes the method more useful for big data applications. See \citep{novikov2016exponential} for more technical details, such as the initialization, step-size selection and dropout.

\begin{algorithm}[t]
\caption{Stochastic Riemannian optimization for learning interactions (ExM) \citep{novikov2016exponential}
}\label{alg_stochastic_riemannian_novikov}
{\small
\begin{algorithmic}[1]
\REQUIRE Training data $\{ (\bx_m, y_m ) \}_{m=1}^M$, desired TT rank\\
 $\{R_1,\ldots, R_{N-1}\}$, number of iterations $K$, mini-batch\\
  size $P$, $0<c_1<0.5$, $0<\rho<1$
\ENSURE Weight tensor $\underline \bW$ in TT format which approximately\\
 minimizes  the loss function \eqref{eqn_penalized_loss_novikov}
\STATE Initialize $\alpha_0=1$
\FOR {$k=1,\ldots,K$}
 \STATE    Sample $P$ data objects randomly, denote their indices by\\
  $h_1,\ldots,h_P \in \{1,\ldots,M\}$
    \STATE $\underline \bG_k = \sum_{p=1}^P \frac{\partial l}{\partial \widehat{y}}   P_{T_{\tW_{k-1}} \Mr_{r}}  ( \underline \bX_{h_p} ) + \lambda \underline \bW_{k-1} $
 \STATE Find the smallest integer $s \geq 0$ such that \\
 $J(\underline \bW_{k-1}) - J(R_{\bW}(\underline \bW_{k-1} , -\rho^{s} \alpha_k c_1 \underline \bG_k )) \geq \rho^{s} \alpha_k c_1  \langle \underline \bG_k,  \underline \bG_k \rangle$ \\
\STATE $\underline \bW_k := R_{\bW}(\underline \bW_{k-1}, - \rho^{s} \alpha_k \underline \bG_k).$
\ENDFOR
\end{algorithmic}
}
\end{algorithm}


This Riemannian optimization algorithm has been  successfully applied to analyze formally a tensor with $2^{160}$ entries (i.e. to analyze a text with 160 words, where interactions between  all words have been exploited). Furthermore, the Exponential Machine has been applied to a recommender system datasets MovieLens with $10^5$ users \citep{novikov2016exponential}.

 \subsection{Future Perspectives of Riemannian Optimization}
\subsectionmark{Future Applications of Riemannian Optimization}

 There has been a steady growth of the interest in Riemannian optimization for machine learning problems,
 especially in the context of low-rank matrix constraints \citep{BoumalAbsil2011,hs-riemnips-2015}.
  As an alternative, the nuclear norm regularization can be  used, since it often leads to convex problems with provable guarantees of convergence. The disadvantage
 is that without additional tricks, the optimization is performed with full matrices,
 whereas the Riemannian approach still works in the low-parametric representation,
 which makes it more efficient especially in the tensor case.
 \markright{\thesection.\quad Riemannian Optimization}%
 Attempts have been made to generalize the nuclear norm concept to the tensor case \citep{phhbd-ttnuc-2016}, but
 a recent negative result shows that it is not possible to provide a good convex surrogate for the TT-manifold
 \citep{rss-iterthreshold-2016}, thus making the Riemannian optimization the most promising tool for low-rank constrained
 optimization. For enhanced efficiency it is desirable to create instruments
 for the construction of such methods in the spirit of modern deep learning frameworks that are
 based on automatic differentiation techniques.
 The Pymanopt \citep{kw-pymanopt-2016} is the first step in this direction, but it is still quite far from
 big-data problems, since it works with full matrices even for a low-rank manifold,
 and in the tensor case that would be prohibitive.
 The stochastic variants of the Riemannian gradient can be readily derived \citep{shalit-riemonline-2010,bonnabel-stochriem-2013,novikov2016exponential} and such methods are applicable to
 a large number of data points (similar to stochastic gradient descent). Derivation of methods which are more complex than Riemannian stochastic gradient descent is still work in progress.

\markright{3.11.\quad Software and Simulations  for TNs and TDs}\section{Software and Computer Simulation Experiments with TNs for Optimization Problems}\markright{\thesection.\quad Software and Simulations  for TNs and TDs}
\sectionmark{Software and Simulations  for TNs and TDs}

Tensor decompositions and tensor network algorithms require sophisticated software libraries, which are  being rapidly developed.

The TT Toolbox,  developed by Oseledets and  coworkers, (\url{http://github.com/oseledets/TT-Toolbox}) for MATLAB
and (\url{http://github.com/oseledets/ttpy}) for PYTHON is currently the most complete  software for
 the TT (MPS/MPO) and QTT  networks \citep{oseledets2012tt}.
  The TT toolbox  supports advanced  applications, which rely on  solving sets of linear equations (including the AMEn algorithm),
  symmetric eigenvalue decomposition (EVD), and inverse/psudoinverse of huge matrices.

  Related and complementary  algorithms  implemented by \citet{KressnerEIG2014} are available within the MATLAB TTeMPS Toolbox (\url{http://anchp.epfl.ch/TTeMPS}).
   This MATLAB toolbox is designed to accommodate various algorithms in the Tensor Train (TT) / Matrix Product States (MPS) format, making use of the object-oriented programming techniques in current MATLAB versions. It also provides an implementation of the efficient AMEn algorithm
   for calculating multiple extremal eigenvalues  and eigenvectors of high-dimensional symmetric eigenvalue decompositions.

For standard TDs (CPD, Tucker models) the Tensor Toolbox for MATLAB,  originally developed by Kolda and Bader, provides several general-purpose functions and special facilities for handling sparse, dense, and structured TDs  \citep{KoldaTensorToolbox,Bader04matlabtensor}, while the $N$-Way Toolbox for MATLAB, by Andersson and Bro, was  developed mostly for Chemometrics applications \citep{Nwaytoolbox}.

The  Tensorlab toolbox developed by Vervliet, Debals, Sorber, Van Barel and De Lathauwer builds upon a  complex optimization framework and offers efficient numerical algorithms for computing the CP, BTD, and constrained Tucker decompositions.
The toolbox includes a library of many  constraints (e.g., nonnegativity, orthogonality) and offers the possibility to combine and jointly factorize dense, sparse and incomplete tensors \citep{Vervliet2016tensorlab}. The new release introduces a tensorization framework and offers enhanced support for handling large-scale
and structured multidimensional datasets.
 Furthermore, sophisticated tools  for the visualization of tensors of an arbitrary order are also available.

Our own developments include the TDALAB (\url{http://bsp.brain.riken.jp/TDALAB}) and TENSORBOX Matlab toolboxes (\url{http://www.bsp.brain.riken.jp/~phan}), which provide a
user-friendly interface and advanced algorithms
 for basic tensor decompositions (CP, Tucker, BTD) \citep{tdalab,tensorbox}.

%


  The Hierarchical Tucker toolbox  by Kressner and Tobler   (\url{http://www.sam.math.ethz.ch/NLAgroup/htucker_toolbox.html})
  focuses mostly on HT type  of tensor networks \citep{KressnerTobler14}, which  avoid  explicit computation of
the SVDs when truncating a tensor which  is already in a HT format \citep{Grasedyck-rev,KressnerTobler14,espigtensorcalculus}.


In  quantum physics and chemistry, a number of  software packages have been developed in the context
of DMRG techniques for simulating quantum networks. One example is the intelligent Tensor (iTensor) by  \citet{iTensor14},  an open source C++ library for rapid development of tensor network algorithms. The iTensor is competitive against other available software packages for  basic DMRG calculations, and   is especially well suited for developing next-generation tensor network algorithms.

The Universal Tensor Network Library  (Uni10),  developed in C++  by
 Ying-Jer Kao,
provides algorithms for   tensor network contraction and a convenient user interface (\url{http://uni10.org/about.html}). The library is
 aimed towards more complex tensor networks such as PEPS and MERA \citep{Uni10Kao2015}.
The Kawashima group in the University of Tokyo has developed a
parallel  $C++$ library for tensor calculations, ``mptensor''
 (\url{https://github.com/smorita/mptensor}).
This is the first library which can perform  parallel tensor
network computation on  supercomputers.

Tensor Networks for huge-scale optimization problems are still in their infancy, however, comprehensive computer experiments over a
number of diverse applications have validated the significant performance advantages, enabled by this technology.
 The following paradigms for  large-scale structured datasets represented by TNs have been considered in the recent publications listed below:

\begin{itemize}

\item Symmetric EVD, PCA, CCA, SVD  and related problems \citep{dolgovEIG2013,KressnerEIG2014,Kressner2014MC,Lee-SIMAX-SVD},

\item Solution of huge linear systems of equations, together with the inverse/pseudoinverse of matrices  and related problems
\citep{Oseledets-Dolgov-lin-syst12,DolgovGMERES,DolgovAMEN2014,Multigrid_Amen2016,Lee-TTPI},

\item  Tensor completion problems \citep{KressnerSV2013,steinlechner15,karlsson2015parallel,Grasedyck15TTcompl},

\item Low-rank tensor network methods for dynamical and parametric problems  \citep{LubichTITT14,TPA2015,cho2016numerical,Garreis2016},

\item TT/HT approximation of nonlinear functions and systems \citep{Khoromskij-TT,QTT-Tucker,KhoromMiao2014superfast},

\item  CP decomposition  and its applications \citep{Choi2014dfacto,Giga-tensor,SmolaNIPS15,Lath_CP2016,Tensor-dispalys1,kamal2016tensor,shin2017fully},

\item  Tucker decompositions and their applications \citep{UKang2016,billion-scale,DINTucker,oseledets2008tucker,BayesCPD-TNNLS,Caiafa2012-NC,caiafa2015stable},

\item TT/HT decompositions and their applications in numerical methods for scientific computing \citep{Litsarev16Integral,zhang-Osel15Anova,Zorin-TT,Benner16reduced,Khoromskij16efficient,Khoromskaia16_Fast,Dolgov2016fast,Khoromskij16efficient},

\item Classification and clustering problems using the TT format \citep{Stoudenmire_optim,Sun2016,NTucker}.

\end{itemize}

When it comes to the implementation, the recently introduced TensorFlow system is an open-source platform, developed by Google, which uses hardware tensor processing units and provides particularly strong support for deep neural networks applications. 



\section[Challenges in Applying TNs for Optimization]{Open Problems and Challenges in Applying TNs for Optimization and Related Problems}
\sectionmark{Challenges in Applying TNs for Optimization}
\label{sect:challenges}

In this chapter, we have demonstrated that tensor networks  are promising tools  for very large-scale optimization problems, numerical methods
(scientific computing) and dimensionality reduction. As we move towards real-world applications, it is important to highlight that the TN approach requires the following two assumptions to be satisfied:
\begin{enumerate}

\item   Datasets  admit sufficiently good low-rank TN approximations,

\item Approximate solutions are acceptable.

\end{enumerate}

 Challenging  issues that remain to be addressed in future research include the
extension of the low-rank tensor network  approach  to complex constrained  optimization problems
 with multiple  constraints, non-smooth penalty terms, and/or
sparsity and/or nonnegativity constraints.
It would be particularly interesting  to investigate  low-rank tensor networks  for the following
  optimization problems, for which efficient algorithms do exist for small- to medium-scale datasets:

\begin{enumerate}

\item Sparse Nonnegative Matrix Factorization (NMF) \citep{NMF-book}, given by
\be
\min_{\bA,\bB} \|\bX - \bA \bB^{\text{T}}\|^2_F + \gamma_1\|\bA\|_1 +\gamma_2 \|\bB\|_1, \;\; \mbox{s.t.} \; a_{ir} \geq 0,\;\;b_{jr} \geq 0,\nonumber \\
\ee
where $\bX \in \Real_+^{I \times J}$ is a given nonnegative data matrix and the $\ell_1$ norm $\|\cdot\|_1$ is defined element-wise,
as the sum of absolute values of the entries. The objective is to estimate huge nonnegative and sparse factor matrices
$\bA \in \Real_+^{I \times R}$, $\bB \in \Real_+^{J \times R}$.

\item Linear and Quadratic Programming (LP/QP) and related problems in the form
\be
\min \{ \frac{1}{2} \bx^{\text{T}} \bQ \bx +\bc^{\text{T}} \bx\}, \;\; \mbox{s.t.} \;\;\bA\bx =\bb, \;\ \bx \geq 0, \nonumber
\ee
where $\bx \in \Real^J$ and the matrix $\bQ \in \Real^{J \times J}$ is symmetric positive definite. When $\bQ=\mathbf{0}$ then the above QP problem reduces to the standard  form of Linear Program (LP).

If the nonnegativity constraints, $\bx \geq 0$, are replaced by   conic constraints, $\bx \in {\cal {K}}$, then the standard QP becomes a quadratic conic program, while
Semidefinite Programs (SDPs) are the generalization of linear programs to matrices. In the standard form,  an SDP minimizes a linear function of a matrix,
 subject to linear equality constraints and a positive definite matrix constraint
\be
\min_{\bX \succeq 0} \tr(\bC \bX) +\gamma \tr(\bX^q), \;\; \mbox{s.t.}\; \tr(\bA_i \bX)=b_i, \;i=1,\ldots, I,\nonumber
\ee
where $\bA_i, \bC$ and $\bX$ are $N \times N$ square matrices and $\bX$ must be a positive definite matrix.

Another  important generalization of LPs is Semi-Infinite Linear Programs (SILPs), which are linear programs with infinitely many constraints. In other words, SILPs
minimize a linear cost function subject to an infinite number of linear constraints.

\item Sparse Inverse Covariance Selection \citep{FriedmanSICov2008}, given by
\be
\min_{\bX \in {\cal S}_{++}^{N \times N}} \tr(\bC_x \bX) - \log \det\bX + \gamma \|\bX\|_1, \;\;\mbox{s.t.} \;\; \alpha \bI_N \preceq \bX \preceq \beta \bI_N, \nonumber \\
\ee
where $\bC_x$ is the given sample covariance matrix, ${\cal S}_{++}^{N \times N}$ denotes the set of $N \times N$ square strictly positive definite matrices and parameters
$\alpha,\beta>0$ impose bounds on the eigenvalues of the solution.

\item Regressor Selection, whereby
\be
\min \|\bA \bx -\bb\|^2_2, \quad \mbox{s.t.} \quad \|\bx\|_0 = \mbox{card}(\bx) \leq c,
\ee
and the operator card$(\bx)$ denotes the cardinality of $\bx$, that is, the number of non-zero components in a vector $\bx$.
The objective is to find the best fit to a given vector $\bb$ in the form of a linear combination of no more than $c$ columns of the matrix $\bA$.

%

\item Sparse Principal Component Analysis (SPCA) using the Penalized Matrix Decomposition (PMD)   \citep{SPCA-Witten,Witten-PHD,JordanSPCA2007}, given by 
\be
&&\max_{\bu,\bv} \{ \bu^{\text{T}} \bA  \bv \}, \\
&& \mbox{s.t.} \quad \|\bu\|_2^2 \leq 1,\;\;  \|\bv\|_2^2 \leq 1, \;\;  P_1(\bu) \leq c_1, \;\; P_2(\bv) \leq c_2, \notag
\label{PMD-PCA}
\ee
where the positive parameters, $c_1,c_2$, control the sparsity level, and the convex penalty functions $P_1,\; P_2$ can take a variety of forms.
Useful examples of   penalty functions which promote sparsity  and smoothness are \citep{Witten-PHD}
\be
\label{PvPMD1}
P(\bv)&=& \|\bv\|_1 =\sum_{i=1}^I |v_i|, \quad \mbox{(LASSO)} \\
P(\bv)&=& \|\bv\|_0 = \sum_{i=1}^I |\sign (v_i)|, \quad (\ell_0-\mbox{norm)} \\
P(\bv) &=& \sum_{i=1}^I |v_i| + \gamma \sum_{i=2}^I |v_i-v_{i-1}|.  \; \mbox{(Fused-LASSO)}
\label{PvPMD2}
\ee


\item Two-way functional PCA/SVD   \citep{Huang-2way} in the form 
\be
&& \max_{\bu,\bv} \{ \bu^{\text{T}} \bA  \bv- \frac{\gamma}{2} P_1(\bu) P_2(\bv) \}, \notag \\
 &&\mbox{s.t.} \quad \|\bu\|_2^2 \leq 1,\;\;  \|\bv\|_2^2 \leq 1.
\label{FPCA}
\ee

\item Sparse SVD  \citep{Lee-2010-bicluster}, given by 
\begin{equation}
\max_{\bu, \bv} \{ \bu^{\text{T}} \bA  \bv- \frac{1}{2} \bu^{\text{T}} \bu \bv^{\text{T}} \bv - \frac{\gamma_1}{2} P_1(\bu) - \frac{\gamma_2}{2} P_2(\bv) \}.\\
\label{SSVD}
\end{equation}


\item Generalized nonnegative SPCA   \citep{Allen-NGPCA} which solves
\be
&& \max_{\bu,\bv} \{ \bu^{\text{T}} \bA \bR \bv - \alpha \|\bv\|_1 \}, \notag\\
&& \mbox{s.t.} \quad \bu^{\text{T}} \bu  \leq 1,\;\;  \bv^{\text{T}} \bR \bv  \leq 1,\;\; \bv \geq 0.
\label{NPCA}
\ee

\end{enumerate}
%

The solution to these open  problems  would be important for a more wide-spread use of low-rank tensor approximations in practical applications.

Additional important and challenging issues that have to be addressed include:

\begin{itemize}

\item Development of new parallel algorithms for wider classes of  optimization problems, with the possibility to incorporate diverse cost functions and constraints. Most of  algorithms  described in this monograph are sequential;
parallelizing  such algorithms is important but not trivial.

\item Development of efficient algorithms for tensor contractions.

Tensors  in TT and HT formats can be relatively easily  optimized,
however, contracting a complex tensor network with cycles (loops) such as PEPS, PEPO, and MERA,  is
a complex  task which requires large computational resources. Indeed,
 for complex tensor  network structures, the time  needed to perform an exact
contraction may  in some cases increase exponentially with  the order of the tensor.
The work in this direction is important because, although approximate contractions may be sufficient for certain applications,
they may also brings difficulties in  controlling the approximation error.

\item Finding the ``best'' tensor network structures or tensor network formats for specific datasets or specific applications, so as to
provide low-rank TN approximations of a sufficiently low-order and with a  relatively low computational complexity.

The term ``best'' is used here in the sense of  a tensor network which provides linear or sub-linear compression of large-scale
 datasets, so as to
 reduce storage cost and computational complexity to affordable levels.

In other words, the challenge is to  develop more sophisticated and flexible tensor networks, whereby the paradigm shift is in a full integration of complex systems and optimization problems, e.g., a system  which simulates  the  structure of biological
 molecule \citep{Savostyanov2014exact}.


\item Current implementations of tensor train decompositions and tensor contractions still require a number of tuning parameters,
such as the approximation accuracy  and the estimation of TN ranks. A step forward would be the development of improved
   and semi-automatic  criteria for the control of approximation accuracy, and {\it a priori} errors bounds.
   In particular, the  unpredictable accumulation of the rounding error and the problem of TT-rank explosion should be carefully addressed.

\item Convergence analysis tools for TN algorithms are still being developed but are critical in order to better understand convergence properties
of such algorithms.



\item Theoretic and methodological approaches are needed to
determine the  kind of constraints  imposed on factor matrices/cores, so as to be able to extract only the  most significant features or only
the desired hidden (latent) variables
 with meaningful physical interpretations.


\item Investigations into the uniqueness of  various  TN decompositions  and their optimality properties (or lack thereof)
are a prerequisite for the development of faster and/or more reliable algorithms.

\item  Techniques   to visualize  huge-scale tensors in tensor network formats are still missing and are urgently needed.



\end{itemize}


\chapter{Tensor Networks for Deep Learning}
\chaptermark{Links Between Tensor Networks and Deep Learning}
\label{Chapter4}


Recent breakthroughs  in the fields of Artificial Intelligence (AI) and Machine Learning (ML) have been largely triggered by the emergence of the  class of deep convolutional neural networks (DCNNs), often   simply called CNNs. The CNNs have become a  vehicle for a large number of practical applications and commercial ventures in  computer vision, speech recognition, language processing, drug discovery, biomedical informatics, recommender systems, robotics, gaming, and artificial creativity, to mention just a few.

It is important to recall that widely used  linear classifiers, such as the SVM and STM   discussed in  Chapter {\ref{Chapter2}, can be considered as shallow feed-forward neural networks (NN).
They therefore inherit the limitations of shallow NNs, when compared to deep  neural networks{\footnote{A deep neural network comprises two or more processing (hidden) layers in between the input  and output layers, in contrast to a shallow NN, which has only one hidden layer.} (DNNs), which include: 		
\begin{itemize}

\item Inability to deal with highly nonlinear data, unless  sophisticated  kernel techniques  are employed;

\item  Cumbersome  and non-generic extensions to multi-class paradigms{\footnote{The multi-class version of the linear STM model is essentially either a one-versus-rest or  one-versus-one extension of the two-class version.}};

\item Incompatibility   of their shallow architectures to learn high-level features which can only be learnt  through a multiscale approach, that is, via DNNs.

\end{itemize}

The renaissance of deep learning  neural networks  \citep{schmidhuber2015deep,LeCun15Nature,Bengio2016book,schneider2017DLSpectrum}
 has both created an active frontier of research in  machine learning and has provided many advantages in  applications, to the extent that the performance of DNNs  in multi-class classification problems can be similar or even better than what is achievable by humans.

Deep learning is also highly interesting in very large-scale data analytics,  since:
\begin{enumerate}

\item Regarding the degree of nonlinearity and multi-level representation ability of features, deep neural networks often significantly outperform their shallow counterparts;

\item High-level representations learnt by  deep  NN models, that are easily interpretable by humans, can also help us to understand the complex information processing  mechanisms and  multi-dimensional interaction of neuronal populations in the human brain;

\item In big data analytics, deep learning is very promising  for mining structured data, e.g., for hierarchical multi-class classification of a huge number of images.

 \end{enumerate}


\noindent {\bf Remark.} It is well known that both shallow and deep  NNs are universal function approximators{\footnote{Universality refers to the ability of a neural network to  approximate  any function,  when no restrictions are imposed on its size  (width). On the other hand, depth efficiency refers to the phenomenon whereby a function realized by polynomially-sized deep neural network requires shallow neural networks to have super-polynomial (exponential) size for  the same accuracy of approximation (curse of dimensionality). This is often referred to as the {\it expressive power of depth}.} in the sense of their ability to approximate arbitrarily well any continuous
 function of $N$ variables on a compact domain; this is achieved under the condition that a shallow neural network has an unbounded width (i.e., the size of a hidden layer), that is, an unlimited number of parameters. In other words, a shallow NN may require a huge (intractable) number of parameters (curse of dimensionality),  while  DNNs can perform such approximations using a much smaller number of parameters.

Despite recent advances in the theory of DNNs, for a continuing success and future perspectives of DNNs (especially DCNNs), the following fundamental challenges need to be further addressed:


\begin{itemize}

\item Ability to generalize while avoiding overfitting in the learning process;

\item Fast  learning and the  avoidance of ``bad'' local and spurious minima, especially for highly nonlinear score (objective) functions;      

\item Rigorous investigation of the  conditions under which deep neural  networks are ``much better'' than the shallow networks (i.e., NNs with one hidden layer);    

\item Theoretical and practical bounds on  the expressive power of a specific architecture, i.e., quantification of the ability to approximate or learn wide classes of unknown nonlinear functions; 

\item Ways to reduce the number  of parameters without a dramatic reduction in performance.

\end{itemize}

The  aim of this chapter is to  discuss the many advantages of tensor networks in addressing the last two of the above challenges  and to build up  both intuitive and mathematical links  between DNNs and TNs.
Revealing  such  inherent connections will both cross-fertilize the areas of deep learning and tensor networks and provide new insights.
In addition to establishing the existing and developing new links, this will also help  to optimize  existing DNNs and/or  generate new, more powerful, architectures.

We  start with an intuitive account of DNNs  based on a simplified hierarchical
Tucker (HT) model,  followed by  alternative simple and  efficient architectures. Finally,   more
sophisticated  TNs,  such as MERA  and 2D tensor network models are briefly discussed in order to enable more flexibility,
 potentially improved performance, and/or higher expressive  power  of
 the next generation  DCNNs. \\

\section{A Perspective of Tensor Networks for Deep Learning}

Several research groups have recently investigated the application of tensor decompositions to simplify DNNs and to establish links between  deep learning and  low-rank tensor  networks
 \citep{lebedev2015fast,Novikov2015tensorizing,Poggio2016,DNN_TF2016,Shashua-HT,RBMTN17}.
For example, \citet{RBMTN17}  presented a general and constructive
connection between Restricted Boltzmann Machines (RBM) and TNs,
together with the correspondence between general  Boltzmann machines and  the TT/MPS model (see detail in the next section).

In a series of research papers \citep{Shashua-HT,Shashua_GTD2016,Cohen_IB,TMM_Sharir} the expressive power of a class of DCNNs  was analyzed using simplified Hierarchal Tucker (HT) models (see the next sections). Particulary, Convolutional Arithmetic Circuits (ConvAC), also known as Sum-Product Networks, and Convolutional Rectifier  Networks (CRN) have been theoretically analyzed as variants of the HT model. The authors  claim that a shallow (single hidden layer) network corresponds to the  CP \hbox{(rank-1)} tensor decomposition, whereas a deep network with $\log_2 N$ hidden layers realizes or corresponds to the HT  decomposition (see the next section).
Some authors also argued
that the ``unreasonable success'' of deep learning can
be explained by inherent laws of physics within the theory of TNs,
which often employ physical constraints of
 locality, symmetry, compositional hierarchical functions,  entropy, and polynomial
log-probability, that are imposed on the  measurements or  input training data \citep{Poggio2016,Tegmark16,RBMTN17}.
%

This all suggests that  a very wide range of  tensor networks can be potentially  used  to model and analyze some specific classes of DNNs, in order to obtain simpler and/or more efficient neural networks in the sense of enhanced expressive power or reduced complexity.
In other words,  consideration of tensor networks  in this context may give rise to new NN architectures which could  even be  superior to the existing ones; this topic has so far been overlooked by practitioners.
%

Furthermore, by virtue of redundancy present in both TNs and DNNs,
 methods used for the reduction or approximation of TNs
  may become a vehicle to achieve  more efficient DNNs, and with a reduced number of parameters.
%
Also, given that usually neither TNs nor  DNNs are  unique
(two NNs with different connection weights and biases may be modeling  the
same nonlinear function), the knowledge about redundancy in
TNs can help simplify  DNNs.

The general concept of optimization of DNNs via TNs is illustrated in Figure~\ref{Fig:DLTN2}.
For a specific DNN, we first  construct its equivalent or corresponding TN representation;
 then the  TN is transformed into a reduced   or canonical form by performing, e.g.,  truncated SVD. 
 This will reduce the
rank  value to the minimal requirement determined by a desired accuracy of approximation.
Finally,  the  reduced and optimized TN is mapped back
 to another DNN of similar universal approximation power to the original DNN, but with a considerably reduced number of parameters.  This follows from the fact that a rounded (approximated) TN has a smaller number of parameters.

  \begin{figure}[t]
\centering
  \includegraphics[width=0.941\textwidth]{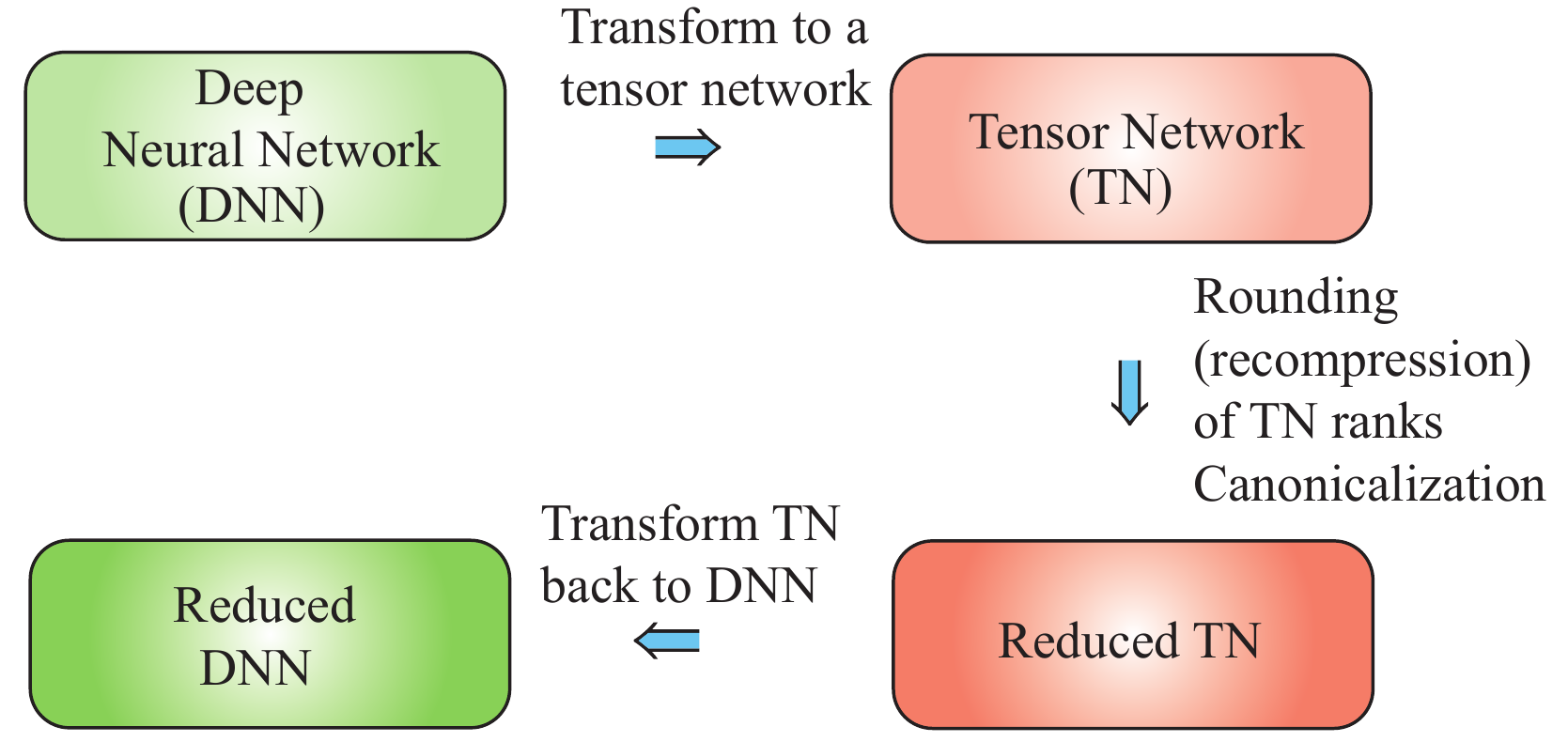}\\
  \caption{Optimization of  Deep Neural Networks (DNNs) using  the Tensor Network (TN) approach. In order to simplify a specific
   DNN and reduce the number of its parameters, we first  transform the DNN into  a specific TN, e.g., TT/MPS, then transform the so approximated (with reduced rank) TN back to a new, optimized DNN. Such a transformation can be performed in a layer by layer fashion, or globally for the whole DNN.  Optionally, we may choose to first construct and learn (e.g., via tensor completion) a tensor network and then transform it to an equivalent DNN.}
  \label{Fig:DLTN2}
\end{figure}

It should be noted that, in practice, the reduction of DNNs through low-rank TN  approximations  offers many potential advantages over a direct reduction of redundant DNNs,  owing  to the  availability of many  efficient optimization
methods to reduce the number of parameters
and achieve a pre-specified approximation error.
Moreover,  low-rank tensor networks  are capable of avoiding the curse of dimensionality through low-order sparsely interconnected core tensors.

On the other hand, in the past two decades, quantum physicists and computer scientists
 have developed solid theoretical understanding and efficient numerical techniques
for low-rank TN decompositions, especially the TT/MPS and TT/MPO.    
The entanglement  entropy{\footnote{Entanglement is a physical phenomenon that occurs when pairs or groups of particles,  such as photons, electrons, or qubits,
are generated or interact in such way that the quantum state of each particle cannot be described independently of the others,  so that
a quantum state must be described for the system as a whole. Entanglement entropy is therefore  a measure for the
amount of entanglement. Strictly speaking, entanglement entropy is a measure of how quantum information is stored in a quantum state and  is mathematically expressed as the von Neumann entropy of the reduced density matrix.}},  Renyi's entropy,
entanglement spectrum and long-range correlations are the most widely
used quantities (calculated from a spatial reduced density matrix) investigated
 in the theory of tensor networks.
The spatial reduced density matrix is determined by splitting
a TN  into two parts, say, regions A and B, where a density
matrix in region A is generated by integrating out all the degrees
of freedom in region B. The entanglement
spectra are determined by the eigenvalues of the reduced
density matrix \citep{Area_law94,eisert2010colloquium}.
Entanglement  entropy  also characterizes the
information content of a bipartition of a specific TN.  Furthermore,
the entanglement area law explains that the entanglement
entropy increases only proportionally to the boundary between
the two tensor sub-networks. Also, entanglement entropy characterizes
the information content of the distribution of  singular values of a matricized tensor,
and  can be viewed as a proxy for
the correlations between the two partitions; uncorrelated
data has zero entanglement entropy at any bipartition.

Note that  TNs are usually designed to efficiently represent
large systems which exhibit a relatively low entanglement
entropy. In practice, similar to deep  neural networks, we often need  to only  care about a small
fraction of the input measurements or training data among a huge number of possible
inputs.
This all suggest that certain guiding principles in DNNs correspond to the entanglement
area law used in the theory of tensor networks. These may then be used to
 quantify the  expressive power of a
wide class of DCNNs.  Note that long range correlations also typically increase with the entanglement.
We therefore conjecture that  realistic datasets in
most successful machine learning applications have relatively
low entanglement entropies \citep{calabrese2004entanglement}.
On the other hand, by exploiting the entanglement
entropy bound of TNs, we can rigorously quantify
the expressive power of a wide class of DNNs applied to complex and highly correlated datasets.

\section{Restricted Boltzmann Machines (RBM)}

Restricted Boltzmann machines (RBMs), illustrated in Figure~\ref{Fig:RBM}, are generative stochastic artificial neural networks, which  have found a wide range of applications in dimensionality reduction, feature extractions, and recommender systems, by virtue of their inherent modeling of the probability distributions of a variety of input data, including natural image  and speech signals.
%

The RBMs are fundamental basic building blocks in a class of Deep  (restricted) Boltzmann Machines (DBMs) and Deep Belief Nets (DBNs). In such deep neural networks, after training one RBM layer, the values of its hidden units can be treated as data for training  higher-level RBMs (see Figure~\ref{Fig:DBN}).

To create an RBM, we first perform a partition of variables into at least two sets: visible (observed) variables, $\mathbf v \in \mathbb{R}^{M}$, and  hidden variables, $\mathbf h \in \mathbb{R}^{K}$.
The goal of RBMs is
then to learn a higher-order latent representation, $\mathbf h \in\mathbb{R}^{K}$, typically for binary variables{\footnote{There are several extensions of RBMs for continuous  data values.}},
within a bipartite undirected graphical model encoding these two layers,
as illustrated in Figure~\ref{Fig:RBM}(a).

 \begin{figure}[htp]
\begin{center}
  \includegraphics[width=0.999\textwidth]{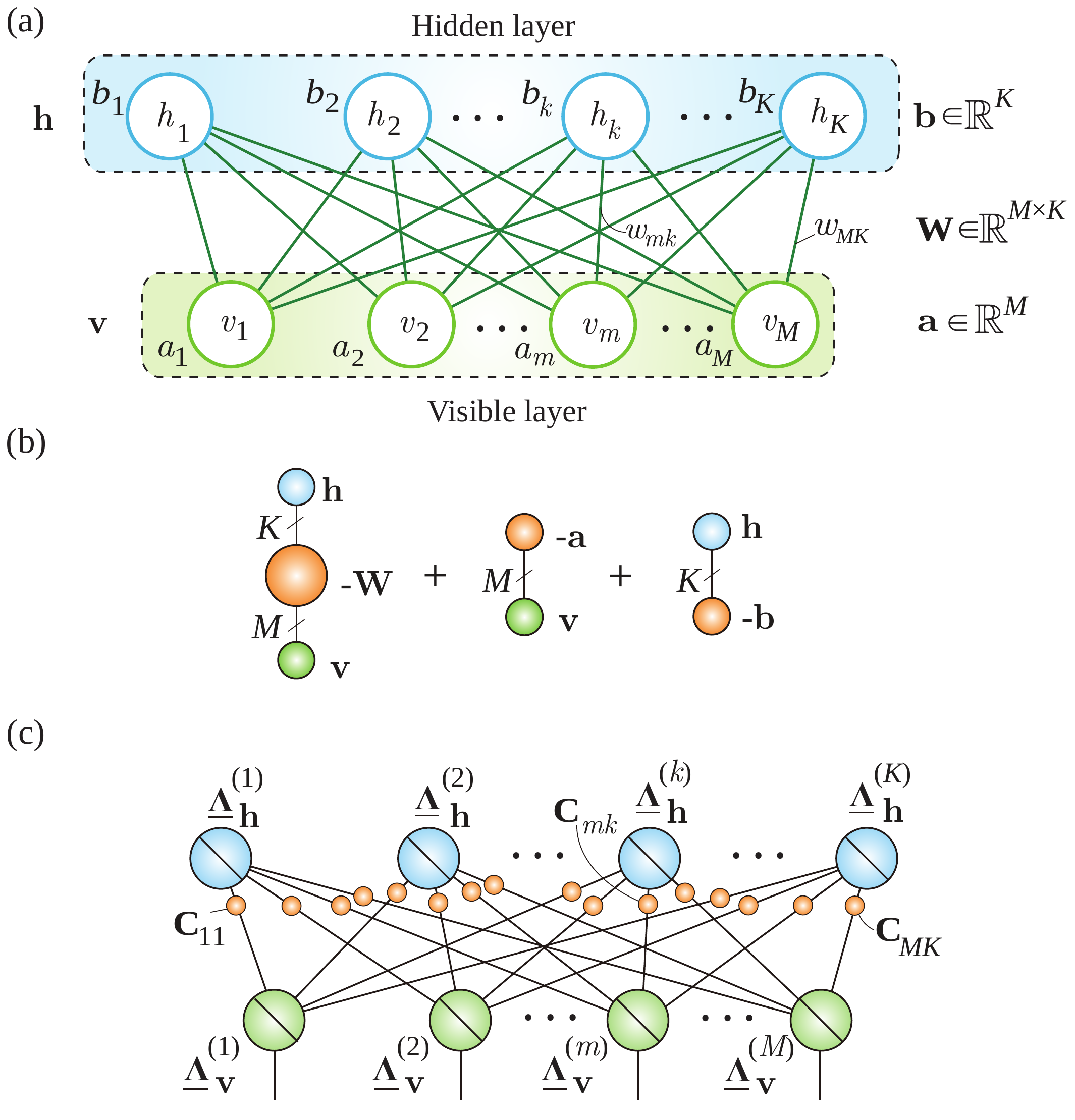}
  \end{center}
  \caption{Graphical  illustrations of the standard RBM. (a) Basic RBM model. The top (hidden) layer represents a vector  $\bh$ of stochastic binary variables, $h_k$ $(k=1,\ldots,K)$, and the bottom (visible) layer represents a vector,  $\bv$, of stochastic binary variables, $v_m$ $(m=1,\ldots,M)$. White round nodes represent observed (visible)  and  hidden variables while
   green lines show the links between the variables.   (b) Energy function, $ E(\bv,\bh; \theta)= -\bv^{\text{T}} \bW \bh - \ba^{\text{T}} \bv - \bb^{\text{T}} \bh $ for $\mathcal{F}(\bv)=0$,
  represented by TN diagrams (see (\ref{eq:Energy})).   (c) A TN diagram corresponding to an RBM   with diagonal core tensors,
 $\underline {\mbi {\Lambda}}_{\bv}^{(m)} =\mbox{diag}_{K}(1, e^{a_m})$ and
 $\underline {\mbi {\Lambda}}_{\bh}^{(k)} =\mbox{diag}_{M}(1, e^{b_k})$, and factor matrices,
  $\bC_{mk}=[1,1; 1, e^{w_{mk}}]$.
   The TN can be simplified through conversion to TT/MPS formats \citep{RBMTN17}. The tensor network represents a set of parameters $\theta= \{\bW,\ba, \bb\}$. Large scale matrices and vectors  can also be represented  by a TT  network, as explained in Chapter \ref{Chapter3}.}
  \label{Fig:RBM}
\end{figure}

 \begin{figure}[t!]
\begin{center}
  \includegraphics[width=0.999\textwidth]{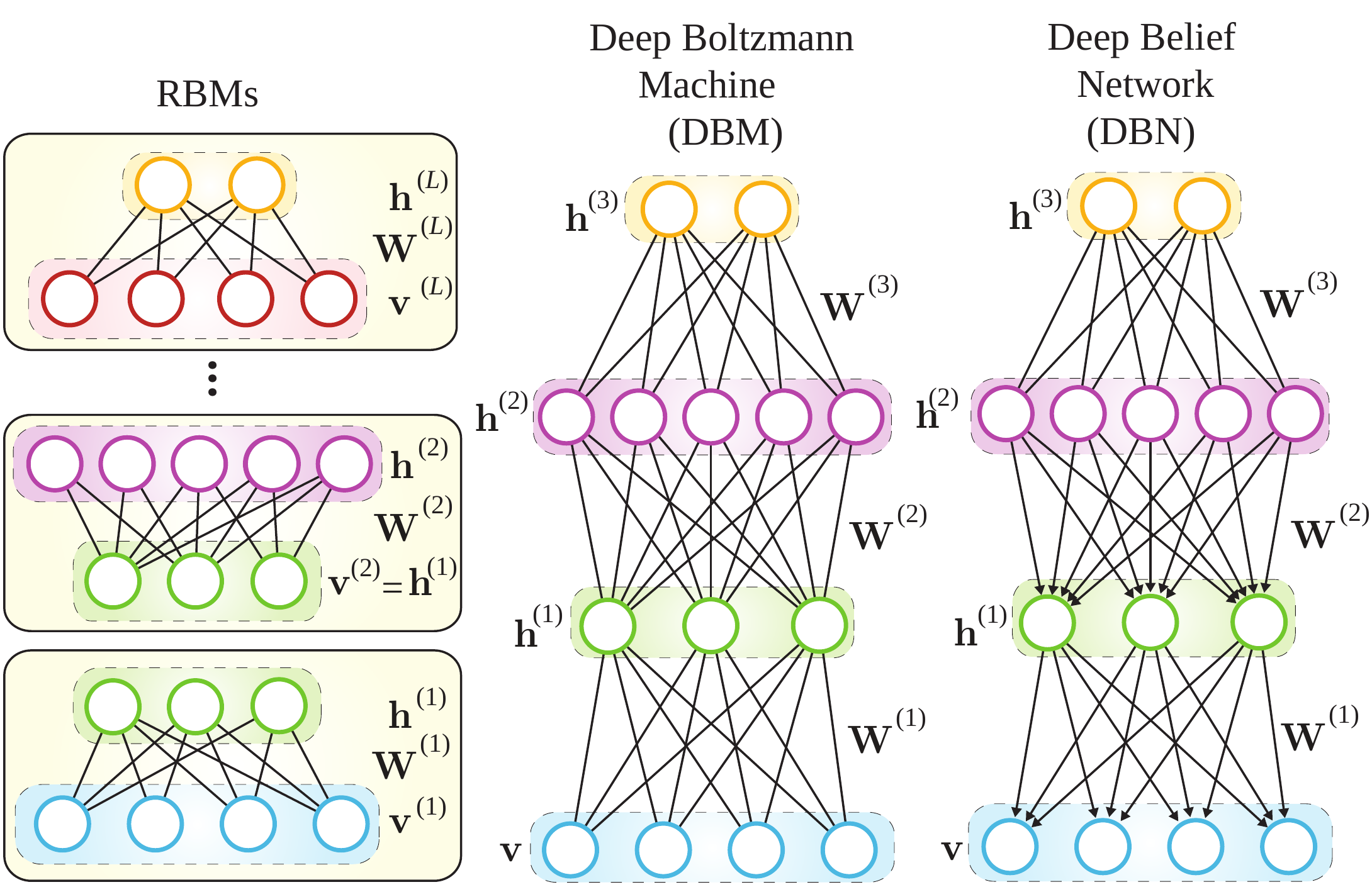}
  \end{center}
  \caption{Graphical illustrations of the construction of   Deep (restricted) Boltzmann Machines  (DBMs) and Standard Deep Belief Networks (DBNs).
Both  use initialization schemes based on greedy layer-wise training of Restricted Boltzmann Machines (RBMs) (left panel),
i.e., they  are  both probabilistic graphical models consisting of stacked layers of RBMs.
Although  DBNs and  DBMs look diagrammatically  quite similar, they are qualitatively very different.  The main
difference is in how the hidden  layers are connected. In a DBM, the connection between all layers is undirected,
thus each pair of layers forms an RBM, while a DBN  has undirected connections between
its top two layers and downward directed connections between all its lower layers (see \citet {Hinton09DBN} for detail).}
  \label{Fig:DBN}
\end{figure}


 The basic idea is that the hidden units of a trained RBM represent relevant features of observations. In other words, an RBM is a probabilistic graphical model based on a layer of hidden variables, which is used to build the probability distribution
over input variables. The standard RBM assigns energy to a joint configuration $(\mathbf v,\mathbf h)$ as (see Figure~\ref{Fig:RBM}(b))
\begin{equation}
E(\mathbf v, \mathbf h; \theta) = -(\mathcal{F}(\mathbf v) +
\mathbf v^{\text{T}}\mathbf W \mathbf h +
\mathbf a^{\text{T}} \mathbf v +
\mathbf b^{\text{T}}\mathbf h),
\label{eq:Energy}
\end{equation}
where $\mathbf a\in\mathbb R^M$ and $\mathbf b \in\mathbb R^K$ are the biases corresponding to the visible and hidden variables, $\mathbf W\in\mathbb R^{M\times K}$ is the (coupling) matrix of mapping parameters, $\theta = \{\mathbf a, \mathbf b, \mathbf W\}$ is  a set of the model parameters to be estimated,
and $\mathcal{F}(\mathbf v)$ is a type-specific function, e.g., $\mathcal{F}(\mathbf v)=0$ for
a binary input and $\mathcal{F}(\mathbf v)=-0.5\sum_m v_m^2$ for Gaussian variables.
Such an  RBM model obeys the Boltzmann distribution for binary variables, as follows
\be
\label{eq:rbm1}
p(\mathbf v, \mathbf h \;|\; \theta) &=& \frac{1}{Z(\theta)} \text{exp} \left( - E(\mathbf v, \mathbf h; \theta)\right),
\ee
where $Z(\theta)$ is the normalization constant.
%

%

%
The absence of connections between hidden binary variables within RBMs allows us to calculate relatively easily  the marginal distribution of the visible variables, as follows
\be
p(\bv \,| \, \theta) &=& \sum_{\bh} p(\bv,\bh \, | \, \theta) = \frac{1}{Z} \sum_{\bh} \exp\left({-E(\bv,\bh;\theta})\right)  \\
&=&  \frac{1}{Z} \prod_{m=1}^M \exp\left(a_m v_m \right) \prod_{k=1}^{K} \left(1+ \exp\left(b_k +  \sum_{m=1}^{M} w_{mk} v_m \right)\right), \notag
\ee
and the  distribution of the hidden variables
\be
\label{eq:phth}
p (\bh \;|\; \theta) &=& \sum_{\bv} p(\bv, \bh \,| \,\theta)= \frac{1}{Z} \sum_{\bv}
\exp( - E(\mathbf v, \mathbf h; \theta)) \\
&=&\frac{1}{Z} \exp \left( \sum_{k=1}^{K} b_k h_k \right) \prod_{m=1}^{M} \left(1+ \exp \left(a_m+\sum_{k=1}^{K} w_{mk} h_k\right)\right). \notag
\ee
The above formulas indicate that an RBM can be regarded as a product of ``experts'', in which a number of experts for individual observations are combined multiplicatively \citep{RBMintro}.

The RBM model can be converted into an equivalent (or corresponding) tensor network model.
 Figure~\ref{Fig:RBM}(c) \citep{RBMTN17} shows  that such a tensor network model
comprises  a set of interconnected   diagonal  core tensors:
$K$th-order  tensors $\underline {\mbi \Lambda}^{(m)}_{\bv} =\mbox{diag}_{K}(1, e^{a_m}) $ and $M$th-order tensors $\underline {\mbi \Lambda}^{(k)}_{\bh} =\mbox{diag}_{M}(1, e^{b_k}) $, all of sizes $2 \times 2 \times \cdots \times 2$. This means that the core tensors are extremely sparse since only the diagonal entries $\underline {\mbi \Lambda}^{(m)}_{\bv}(1,1,\ldots,1)=1$,  $\underline {\mbi \Lambda}^{(m)}_{\bv} (2,2,\ldots,2)=e^{a_m}$, and $\underline {\mbi \Lambda}^{(k)}_{\bh} (1,1,\ldots,1)=1$ ,  $\underline {\mbi \Lambda}^{(k)}_{\bh} (2,2,\ldots,2)=e^{b_k}$ are non-zero. The diagonal core tensors in the visible and hidden layer are interconnected via the $2 \times 2$ matrices $\bC_{mk} =[1,1;1,e^{w_{mk}}]$.
%

The bipartite structure of the RBM enables units in one layer to become conditionally independent of  the other layer. Thus,  the conditional distributions over the hidden and visible variables can be factorized as
\be
p(\mathbf v|\mathbf h; \theta) &= \prod_{m=1}^M p(v_m |\mathbf h), \quad
p(\mathbf h|\mathbf v; \theta) &= \prod_{k=1}^K p(h_k | \mathbf v),
\ee
%
and
\begin{equation}
\begin{split}
p(v_m=1|\mathbf h) &= \sigma\left(\sum_k w_{mk} h_k + a_m\right),\\
p(h_k =1| \mathbf v) & = \sigma\left(\sum_m w_{mk} v_m + b_k\right),
\end{split}
\end{equation}
where $\sigma(x) =1/(1+ \text{exp}(-x))$ is the sigmoid  function.

 \begin{figure}[t]
\begin{center}
  \includegraphics[width=0.75\textwidth]{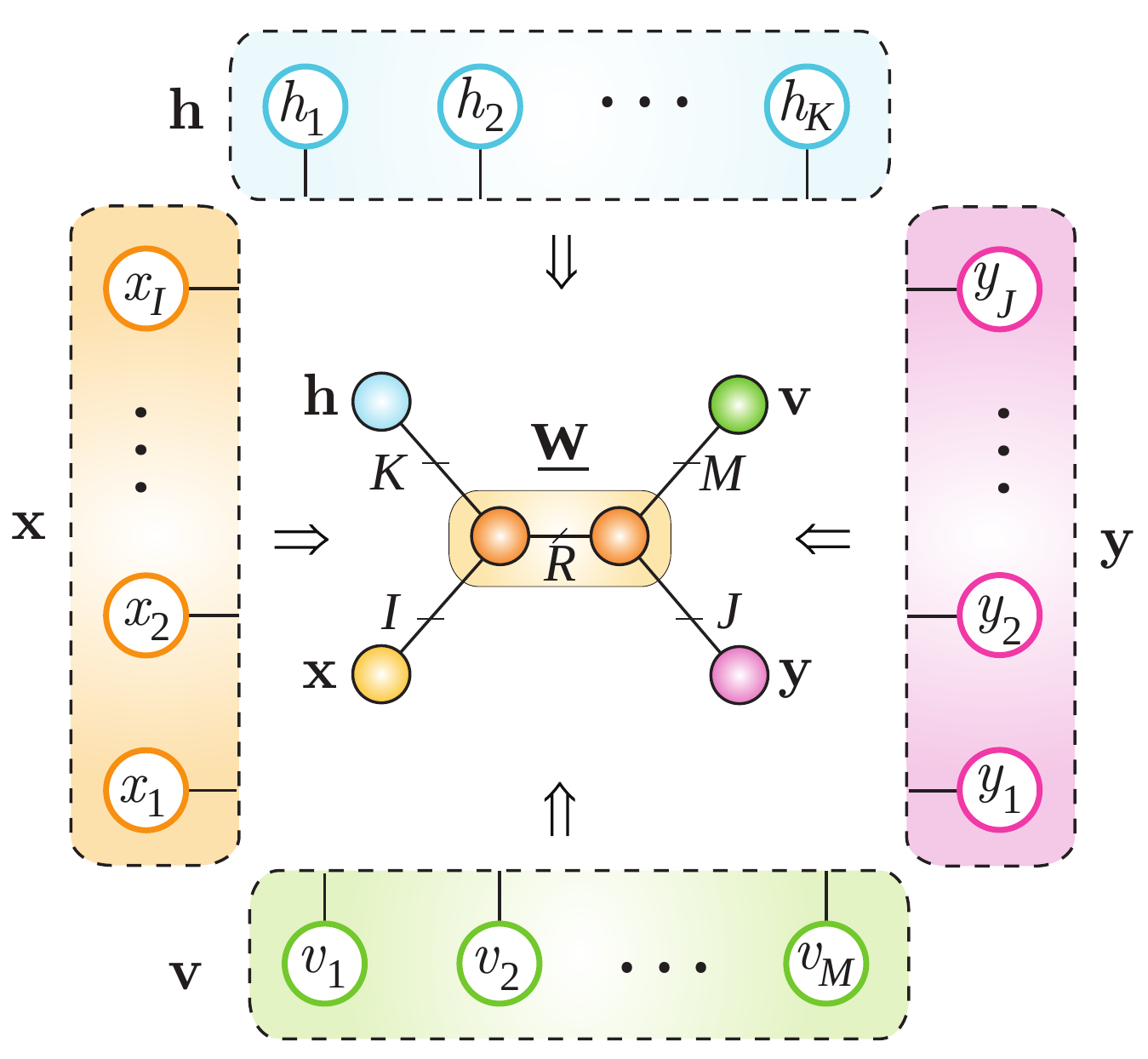}
  \end{center}
  \caption{Energy function of a higher-order (multi-way) RBM with four layers, given in (\ref{eq:way4}). In order to reduce the number of parameters, the 4th-order weight tensor $\underline \bW \in \Real^{I \times J \times M \times K}$ is represented in a distributed (tensor network)  form, in this particular case  as two 3rd-order cores. This allows us to reduce the number of parameters from $IJKM$ to only $IKR+JMR$, where typically $R \ll I,J,K,M$.}
  \label{Fig:RBMW4}
\end{figure}

{\bf Higher-order RBMs.} The modeling power of an RBM can be enhanced by increasing the number of hidden units, or by adding and splitting visible layers and/or  adding conditional and activation layers. The standard RBM can  be generalized and extended  in many different ways. The most natural way is to jointly express the energy  for more than two set of variables. For example, for four sets of binary variables,
$\bv, \bh, \bx, \by$, the energy function (with neglected biases) can be expressed in the following form (see Figure~\ref{Fig:RBMW4})
\be
E(\bv,\bh,\bx,\by) &=& - \sum_{i=1}^I\sum_{j=1}^J\sum_{m=1}^M\sum_{k=1}^K w_{ijmk} \; x_i \,y_j \, v_{m} \,  h_{k} \notag \\
&=& -\underline \bW \; \bar{\times}_1 \; \bx \; \bar{\times}_2 \; \by \; \bar{\times}_3 \; \bv \; \bar{\times}_4 \bh.
\label{eq:way4}
\ee
Unfortunately,  higher-order RBM models suffer from a large number of
parameters, which scales with the product of  visible, hidden and conditioning variables. To alleviate this problem, we can perform  dimensionality reduction
via low-rank tensor network approximations.

\subsection{Matrix Variate Restricted Boltzmann Machines}

The matrix-variate restricted Boltzmann machine (MVRBM) model was proposed as a generalization of the classic RBM to explicitly model matrix data \citep{qi2016matrix}, whereby  both input and hidden variables are in their matrix forms which are connected by bilinear transforms.
The MVRBM has much fewer model parameters than the standard RBM  and thus admits  faster training  while retaining a comparable performance to the classic RBM.

Let $\mathbf V \in \mathbb{R}^{J\times M}$ be a binary  matrix of visible variables, and
$\mathbf H \in \mathbb{R}^{K\times I}$ a  binary  matrix  of hidden variables.
 Given the parameters of a 4th-order tensor, $\tensor W\in\mathbb{R}^{I\times J\times M \times K}$, and bias matrices, $\bA \in \mathbb{R}^{J\times M}$ and $\bB \in\mathbb{R}^{K\times I}$, the energy function can be defined as (see Figure~\ref{Fig:RBMM})
\be
\label{eq:mvrbm1}
E(\mathbf V, \mathbf H; \theta) &=& - \sum_{i=1}^I\sum_{j=1}^J\sum_{m=1}^M \sum_{k=1}^K
v_{mj}w_{ijmk} h_{ki} \\
&&- \sum_{m=1}^M\sum_{j=1}^J v_{jm}a_{jm} - \sum_{k=1}^K \sum_{i=1}^I h_{ki}b_{ki}, \notag
\ee
where $\theta = \{\tensor W, \bA, \bB \}$ is the set of all model  parameters.

In order to reduce the number of free parameters, the higher-order weight tensor can be approximated by a suitable tensor network, e.g.  a TT or HT network. In the particular case of a four-way weight tensor, $\underline \bW \in \Real^{I \times J \times M \times K}$, it can be represented via truncated SVD by two \mbox{3rd-order} core tensors, $\underline \bW^{(1)} \in \Real^{M\times K \times R}$ and  $\underline \bW^{(2)} \in \Real^{I\times J \times R}$. For a   very crude  approximation with rank $R=1$, these core tensors simplify to matrices $\bW^{(1)}\in\mathbb{R}^{M\times K}$ and $\bW^{(2)}\in \mathbb{R}^{I\times J}$, while the energy function in (\ref{eq:mvrbm1}) simplifies to the following form
\begin{equation}
\label{eq:mvrbm2}
E(\mathbf V, \mathbf H; \theta) =-\text{tr}(\bW^{(1)} \mathbf H \bW^{(2)} \mathbf V) -\text{tr}(\mathbf V^{\text{T}}\mathbf A) -\text{tr}(\mathbf H^{\text{T}}\bB).
\end{equation}
Both matrices $\bW^{(1)}$ and $\bW^{(2)}$   approximate  the interaction between the input matrix $\mathbf V$ and the hidden matrix $\mathbf H$.
Note that in this way, the total number of free parameters is reduced from $IJMK$ to $IJ + KM + JM + IK$.  For the corresponding learning  algorithms see \citet{Hinton09DBN,RBMintro,qi2016matrix}.

 \begin{figure}[t]
\begin{center}
  \includegraphics[width=0.95\textwidth]{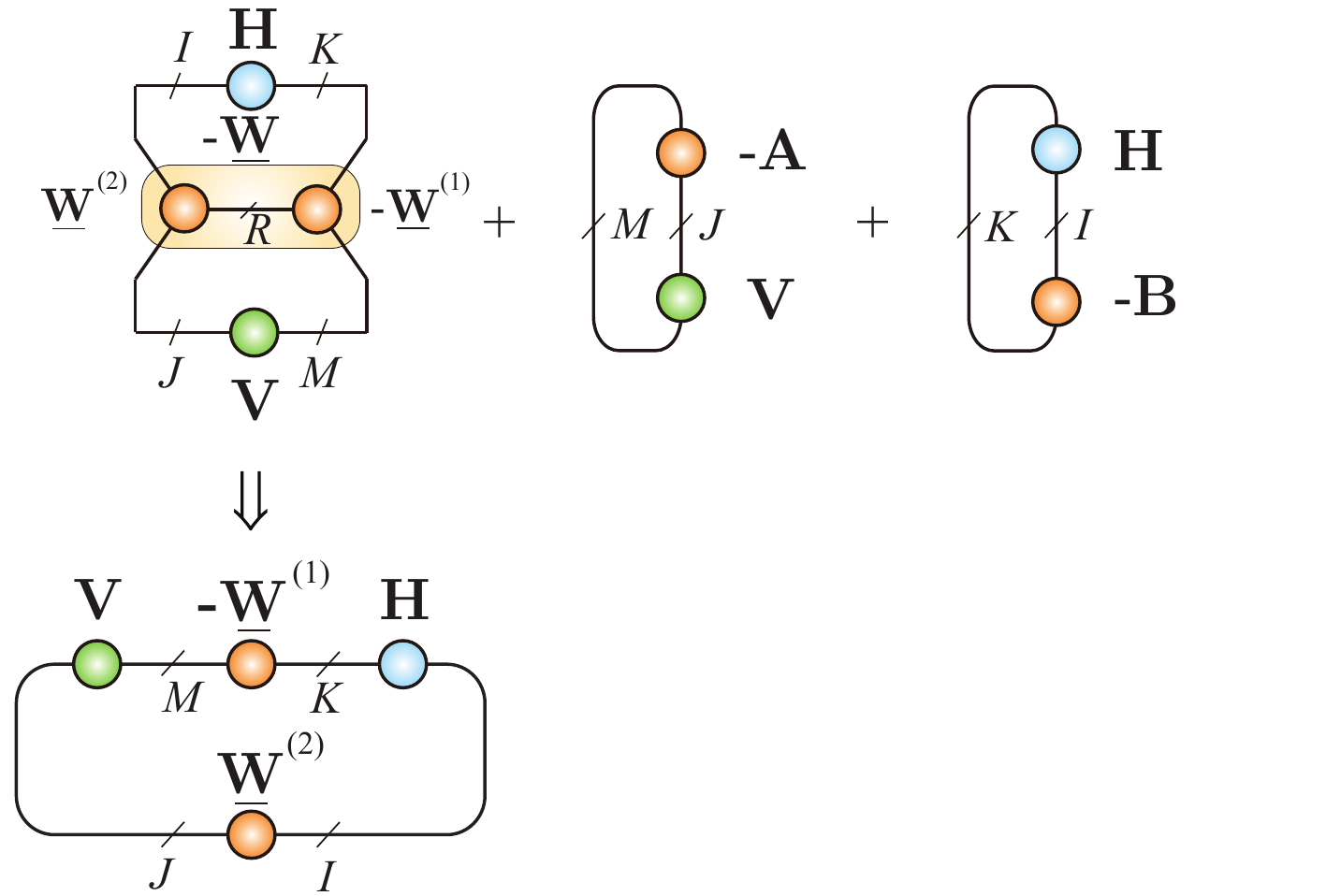}
  \end{center}
  \caption{Energy function  of the matrix-variate RBM represented by TN diagrams, given in (\ref{eq:mvrbm1}). In order to reduce the number of parameters, the 4th-order weight tensor can be represented by the outer product of two matrices (lower panel), $\underline \bW = \bW^{(1)} \circ \bW^{(2)}$. For a higher order weight tensor, we can use the
  low-rank  TT or HT representations.}
  \label{Fig:RBMM}
\end{figure}

\subsection{Tensor-Variate Restricted Boltzmann Machines}

  The tensor-variate RMBs (TvRBMs)  have  the ability to capture  multiplicative interactions between data modes and  latent variables (see e.g.  \citep{nguyen2015tensor}).
  The TvRBMs  are highly compact, in that the number of free parameters grows only linearly with the number of  modes, while their   multiway factoring of mapping parameters  link data modes and hidden units. Modes interact in a multiplicative fashion gated by hidden units, and TvRBM uses  tensor data and the hidden layer to construct an ($N+1$)-mode tensor.

In a TvRBM, the set of visible units is represented by an $N$th-order tensor, $\tensor V \in \mathbb{R}^{I_1\times I_2\times\cdots\times I_N}$, and the hidden units $\mathbf h$ are the same as in RBM.
 The goal is to model the joint distribution, $p(\tensor V, \mathbf h)$, by the optimization of
 the energy function
\begin{equation}
E(\tensor V, \mathbf h; \mbi \theta) = -\mathcal{F}(\tensor V) -\langle \tensor A, \tensor V\rangle - \mathbf b^{\text{T}}\mathbf h - \langle \tW , \tV \circ \bh \rangle,
\end{equation}
where $\mathcal{F}(\tensor V)$ is a type-specific function, $\tensor A \in\mathbb{R}^{I_1\times I_2\times \cdots\times I_N}$ are visible biases, and $\tensor W \in\mathbb{R}^{I_1\times I_2\times\cdots\times I_N\times K}$ are mapping parameters.

The hidden posterior is then given by
\begin{equation}
\label{eq:tvrbm}
p(h_k=1 | \tensor V; \theta) = \sigma \left( b_k + \langle\tensor V, \tW(:,\ldots,:,k)\rangle
 \right),
\end{equation}
where  $\sigma(x)$ is the sigmoid function.
The generative distribution, $p(\tensor V| \mathbf h)$, on the other hand  is type-specific.

The most popular cases are the binary and Gaussian inputs.
For  a  binary input, we have $\mathcal{F}(\tensor V) =0$, and the following generative distribution
\begin{equation}
p(v_{i_1 i_2 \ldots i_N} | \mathbf h) = \sigma \left(a_{i_1 i_2 \ldots i_N} +
\tW(i_1,\ldots,i_N,:)^{\text{T}} \; \bh\right).
\end{equation}
For a Gaussian input, assuming unit variance, i.e. $\mathcal{F}(\tensor V) =  -0.5\langle \tensor V, \tensor V\rangle$, the generative distribution becomes
\begin{equation}
p(v_{i_1 i_2\ldots i_N}| \mathbf h) = \mathcal{N} \left(\{ a_{i_1 i_2\ldots i_n} + \tW(i_1,\ldots,i_N,:)^{\text{T}} \; \bh\}; \tensor 1 \right),
\end{equation}
where $\tensor 1\in\mathbb{R}^{I_1\times I_2\times \cdots \times I_N} $ is the tensor with unit elements.


A major problem with the parameterizations  of $\tensor W$ is the excessively large number of free parameters, which scales with  the product of data modes and hidden dimensions.  In particular, the $(N+1)$th-order weight tensor $\tensor W$ has $K\prod_n I_n$ elements, and its size quickly reaches billions of entries, even when the mode dimensionalities $K, (I_1,\ldots, I_N)$  and $N$ are quite moderate.    
This makes learning in  the raw tensor format extremely difficult  or even impossible, since robust estimation of parameters would require a huge dataset.
To this end,  low-rank tensor network  decompositions can be employed to construct  approximate interactions between visible modes and hidden units.

In the simplest scenario, we can approximately  represent  $\tW$ by a rank-$R$  CP tensor
\be
\tensor W &=& \sum_{r=1}^R \lambda_r \cdot (\mathbf w_r^{(1)} \circ \cdots\circ \mathbf w_r^{(N)} \circ\mathbf w_r^{(N+1)}) \\
&=& \underline {\mbi \Lambda} \times_1 \bW^{(1)}  \cdots \times_N \bW^{(N)} \times_{N+1} \bW^{(N+1)}, \notag
\ee
where
the matrix $\bW^{(n)} = [\bw_1^{(n)}, \ldots,\bw^{(n)}_R] \in\mathbb{R}^{I_n \times R}$  for $n=1,\ldots,N$ represents the mode-factor weights, and $\bW^{(N+1)}=\bW^{(h)} \in\mathbb{R}^{K\times R}$ the hidden-factor matrix.
Such factoring allows for multiplicative interactions between modes to be moderated by the hidden units through the hidden-factor matrix $\mathbf W^{(h)}$.  Thus, the model captures the modes of variation through the new representation, enabled by $\mathbf h$. By employing the CP decomposition, the number of mapping parameters  may be   reduced down to $R(K+\sum_n I_n)$, which grows linearly rather than exponentially in $N$. Importantly, the conditional independence among intra-layer variables, i.e. $p(\mathbf h| \tensor V) = \prod_{k=1}^K p(h_k|\tensor V)$ and $p(\tensor V| \mathbf h) = \prod_n\prod_{i_n=1}^{I_n} p(v_{i_1 i_2\ldots i_N } | \mathbf h)$, is not affected. Therefore, the model preserves the fast sampling and interface properties of RBM.

\section{Basic Features of Deep Convolutional Neural Networks}
\sectionmark{Deep Convolutional Neural Networks}

Basic  DCNNs are  characterized by at least three features:
locality, weight sharing (optional) and pooling, as explained below:
\begin{itemize}

\item Locality refers to  the connection
of a (artificial) neuron to only  neighboring neurons
in the preceding layer, as opposed to being fed by the entire
layer (this is consistent  with biological NNs).

\item Weight sharing reflects the property  that different neurons in the same layer, connected
to different neighborhoods in the preceding layer, often share the same weights.
 Note that weight sharing, when combined with locality, gives rise to standard convolution{\footnote{Although weight sharing may  reduce the complexity of a deep neural network, this step is optional.  However, the locality at each layer is a key factor which gives DCNNs an exponential advantage over shallow NNs in terms of the representation ability \citep{anselmi2015deep,Poggio2016,PoggioDL2}.}}.

 \item Pooling is essentially an operator which gradually decimates (reduces)
layer sizes by replacing the local population of neural activations in a spatial window by a single value
(e.g., by taking their maxima, average values or their scaled products).
In the context of images, pooling induces invariance to translation, which often does not
affect semantic content, and  is interpreted as a way to create a hierarchy of abstractions in the patterns that neurons respond to \citep{anselmi2015deep,Shashua-HT}.

\end{itemize}

Usually,  DCNNs perform much better when dealing with compositional
function approximations{\footnote{A compositional function can take, for example, the following form $h_1(\cdots h_3(h_{21}(h_{11}(x_1,x_2)h_{12}(x_3,x_4)),h_{22}(h_{13}(x_5,x_6)h_{14}(x_7,x_8))\cdots)))$.}} and multi-class classification problems than shallow neural network architectures with one hidden layer.
In fact, DCNNs can even efficiently and  conveniently select a subset of features for multiple classes, while
for efficient learning a DCNN model can be pre-trained by first learning each DCNN layer, followed by fine tuning of the parameter of the entire model, using  e.g.,  stochastic gradient descent.
To summarize, deep learning neural  networks have the ability to  arbitrarily  well exploit and approximate the complexity of compositional hierarchical functions,  a wide class of functions to which shallow networks are blind.

\section{Score Functions for Deep Convolutional Neural Networks}
\sectionmark{Score Functions for DCNNs (ConvACs)}

Consider a multi-class classification task where the input training data,
also called local structures or instances, e.g., input patches in images,  are denoted by the set
$X= \{\bx_1,\ldots, \bx_N\}$, where $\bx_n \in \Real^S$, $\;(n=1,\ldots,N)$,
belong to one of distinct categories (classes) denoted by $c \in \{1,2, \ldots, C\}$.
Such  a representation is quite natural for many high-dimensional data -- in images,
 the local structures represent vectorization of  patches consisting of $S$ pixels, while in audio data  can be  represented through spectrograms.

For this kind of problems, DCNNs can be described  by  a set of  multivariate score functions
\be
y= h_{c}(\bx_1,\ldots , \bx_N) = \sum_{i_1=1}^{I_1} \cdots \sum_{i_N=1}^{I_N}
\underline \bW_{\,c}(i_1,\ldots,i_N)  \prod_{n=1}^N f_{\theta_{i_n}} (\bx_n), \notag \\
\label{eq:score1}
 \ee
 where  $\underline \bW_{\,c} \in \Real^{I_1 \times \cdots \times I_N}$   is an $N$th-order coefficient  tensor, with all dimensions $I_n=I, \; \forall n$, $N$ is the number of (overlapped) input patches $\{\bx_n\}$, $I_n$ is the size (dimension) of each mode $\underline \bW_{\,c}$,  and $f_{\theta_{1}}, \ldots, f_{\theta_{I_n}}$ are referred to as the representation functions (in the representation layer) selected from a parametric family of nonlinear functions{\footnote{Note that the representation layer can be considered as a tensorization of the input patches $\bx_n$.}}.

In general, the one-dimensional basis functions could be polynomials,  splines  or other sets
of basis functions.
 Natural choices for this family of nonlinear functions are also radial basis functions (Gaussian RBFs), wavelets{\footnote{Particularly interesting are Gabor wavelets, owing to their ability to induce features that resemble representations in the visual cortex of  human brain.}}, and  affine functions followed by point-wise activations.
 Note that the representation functions in standard (artificial) neurons have the form
 \be
 f_{\theta_{i}}(\bx)=\sigma(\tilde \bw_{i}^{\text{T}} \bx +b_{i}),
 \ee
for the set of parameters $\theta_{i}=\{\tilde \bw_{i}, b_{i}\}$, where $\sigma(\cdot)$ is a suitably chosen activation function.

\begin{figure}[t]
(a) \hspace{4cm} (b) \hspace{3cm} (c)\\ \\
 \includegraphics[width=0.99\textwidth]{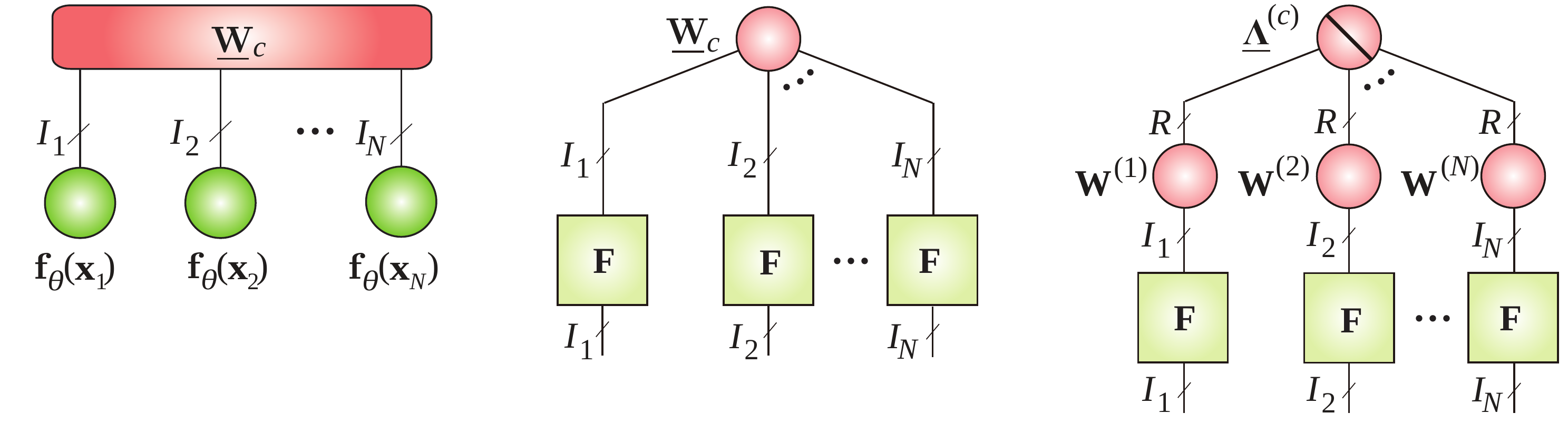}\\
  \caption{Various representations of the score function of a DCNN. (a)  Direct representation of the score function $h_{c}(\bx_1,\bx_2,\ldots,\bx_N)= \tW_c \bar \times_1  {\boldf}_{\mbi \theta}(\bx_1)  \bar \times_2  {\boldf}_{\mbi \theta}(\bx_2) \cdots \bar \times_N  {\boldf}_{\mbi \theta}(\bx_N)$. (b) Graphical illustration of  the  $N$th-order grid tensor of the score function $h_{c}$, which can be   considered as a special case of Tucker-$N$ model where the representation matrix $\bF \in \Real^{ I \times I}$ is built up of factor matrices; note that all the factor matrices are the same and $I_n =I,\; \forall n$.  (c)  CP decomposition of the coefficient tensor $\underline \bW_{\,c} = \underline {\mbi \Lambda}^{(c)} \times_1 \bW^{(1)} \times_2 \bW^{(2)} \cdots \times_N \bW^{(N)} =  \sum_{r=1}^{R} \lambda^{(c)}_r  \; (\bw^{(1)}_r \circ \bw^{(2)}_r \circ \cdots  \circ\bw^{(N)}_r)$, where $\bW^{(n)}=[\bw^{(n)}_1,\ldots, \bw^{(n)}_R] \in \Real^{I \times R}$. This CP model  corresponds to  a simple shallow neural network with one hidden layer  comprising weights $w^{(n)}_{ir}$, and the output layer comprising weights $\lambda^{(c)}_r $, $r=1,\ldots,R$. Note that the coefficient tensor $\underline \bW_c$ can be represented in a distributed form by any suitable tensor network.}
  \label{Fig:Grid}
\end{figure}

The representation layer play a key role to transform the inputs, by means of $I$ nonlinear functions, $f_{\theta_i}(\bx_n)$ $\;(i=1,2,\ldots,I)$, to template input patches, thereby creating $I $ feature maps \citep{Cohen_IB}.
Note that the representation layer can be described by a feature vector defined as
\be
{\boldf}={\boldf}_{\mbi \theta} (\bx_1) \otimes {\boldf}_{\mbi \theta}(\bx_2) \otimes \cdots  \otimes {\boldf}_{\mbi \theta}(\bx_N) \in \Real^{I_1 I_2 \cdots I_N},
\ee
where $ {\boldf}_{\mbi \theta}(\bx_n) =[f_{\theta_1}(\bx_n), f_{\theta_2}(\bx_n), \ldots, f_{\theta_{I_n}}(\bx_n)]^T \in \Real^{I_n}$ for $n=1,2,\ldots, N$ and $i_n=1,2,\ldots, I_n$.
Equivalently, the representation layer can be described as a rank one tensor
(see Figure~\ref{Fig:Grid}(a))
\be
\underline \bF =  {\boldf}_{\mbi \theta} (\bx_1) \circ {\boldf}_{\mbi \theta}(\bx_2) \circ \cdots  \circ {\boldf}_{\mbi \theta}(\bx_N)
\in \Real^{I_1 \times I_2 \times \cdots \times I_N}.
\ee
This allows us to represent the score function as an inner product of two tensors, as illustrated in Figure~\ref{Fig:Grid}(a)
\be
h_{c}(\bx_1,\ldots,\bx_N) = \langle \tW_c, \tF \rangle =
\tW_c \; \bar \times_1 \; {\boldf}_{\mbi \theta} (\bx_1) \; \bar \times_2 \; {\boldf}_{\mbi \theta} (\bx_2) \cdots \bar \times_N \; {\boldf}_{\mbi \theta} (\bx_N). \notag \\
\ee
To simplify the notations, assume that $I_n=I, \; \forall n$,   then we   can also construct a
 square matrix $\bF \in \Real^{I \times I}$ ,  with rows $\bF(i,:)=[f_{\theta_1}(\bx^{(i)}), f_{\theta_2}(\bx^{(i)}), \ldots, f_{\theta_{I}}(\bx^{(i)})]$ for $i=1,2,\ldots,I$ and entries  taken from values of nonlinear basis functions $\{ f_{\theta_1}, \ldots,  f_{\theta_I}\}$ on the selected templates
 $\{\bx^{(1)}, \bx^{(2)},\ldots, \bx^{(I)}\}$ \citep{Shashua_GTD2016}.

 For discrete data values, the score function can be represented by a grid tensor, as graphically illustrated in Figure~\ref{Fig:Grid}(b).
 The grid tensor of the nonlinear score function $h_c(\bx_1, ,..., \bx_N)$ determined
over all the templates $\bx^{(1)}, \bx^{(2)}, ..., \bx^{(I)}$  can be expressed, as follows
\be
\tW(h_c) = \tW_c \times_1 \bF  \times_2 \bF \cdots  \times_N \bF.
\ee
Of course, since the order $N$ of the coefficient (core) tensor is large, it cannot be
implemented, or even saved on a computer due to the curse of dimensionality.

To this end, we represent the weight tensor in some low-rank tensor network format with a smaller number of parameters. A simple model is the CP representation, which leads to a shallow network as illustrated in Figure~\ref{Fig:Grid}(c).  However, this approach is associated with   two problems:  (i) the rank $R$  of the coefficient tensor $\tW_c$ can be very large
(so compression ratio cannot  be very high), (ii) the  existing  CP decomposition algorithms are not very stable for very high-order tensors, and  so an alternative promising approach would be to apply tensor networks such as HT  that enable  us to avoid the curse of dimensionality.

Following the representation layer,  a  DCNN may consists of a cascade of  $L$ convolutional hidden layers with pooling in-between, where the number of layers $L$ should be  at least two. In other words, each hidden layer  performs 3D or 4D convolution followed by spatial window pooling, in order to reduce (decimate)  feature maps by e.g., taking a product of the entries in sub-windows. The output layer is a linear dense layer.

Classification can then be carried out in a standard way, through the maximization of a set of labeled score functions, $h_{c}$ for $C$ classes, that is,
 the predicted label for the input instants $X=(\bx_1,\ldots, \bx_N)$ will be the index $\hat y$ for which the score value attains a maximum, that is
\be
 \hat y = \displaystyle{\mbox{argmax}_{\{c\}}} h_{c}(\bx_1, \ldots , \bx_N).
\ee
Such score functions can be represented through their coefficient tensors which, in turn, can be approximated by low-rank tensor network decompositions. 

One restriction of the so formulated score functions (\ref{eq:score1}) is  that they allow  for a straightforward implementation of only a particular class of DCNNs, called  convolutional Arithmetic Circuit (ConvAC).
 However,  the score functions can be approximated indirectly and almost
 equivalently using more popular CNNs (see the next section).  For example, it was shown recently how NNs with a univariate rectified linear unit (ReLU) nonlinearity may perform  multivariate function approximation \citep{Poggio2016}.


 As  discussed  in Part 1 and in  Chapter \ref{Chapter1} the main idea is  to employ a low-rank tensor network representation to approximate and interpolate  a multivariate
function $h_c(\bx_1, \ldots, \bx_N)$ of $N$ variables
by a finite sum of separated products of  simpler functions  (i.e., via  sparsely interconnected
core tensors).

%


\noindent {\bf Tensorial Mixture Model (TMM).} Recently, the model in (\ref{eq:score1}) has been extended to a corresponding   probabilistic model, referred to as the Tensorial Mixture Model (TMM), given by  \citet{TMM_Sharir}
\be
 P(\bx_1,\bx_2, \ldots, &\bx_N&| \; y=c) \\
 &=& \sum_{i_1=1}^{I} \cdots \sum_{i_N=1}^{I}
\underline \bP_{y=c}(i_1,\ldots,i_N) \prod_{n=1}^N P(\bx_n | i_n; \theta_{i_n} ),\notag
 \label{eq:TMM}
 \ee
 where $\underline \bP_{(y=c)}(i_1,\ldots,i_N)$ (corresponding to $\underline \bW_{\,c}(i_1,\ldots,i_N)$) are jointly decomposed prior tensors and
 $\{P(\bx_n | i_n; \theta_{i_n} )\}$ (corresponding to $f_{\theta_{i_n}} (\bx_n)$) are mixing components  (marginal likelihoods) shared across all input samples $\bx_n$ $(n=1,\ldots,N)$ and classes $\{c\}$, $(c=1,2,\ldots,C)$. Note that for computational tractability,   suitable scaling and  nonnegativity constraints must be imposed  on the TMM model.

 The TMM can be considered as a generalization of  standard  mixture models widely used in machine learning, such as  the Gaussian Mixture Model (GMM).
However,  unlike the standard mixture models, in the TMM we cannot perform inference directly from (\ref{eq:TMM}), nor can we  even store the prior tensor,  $\underline \bP_{(y=c)}(i_1,\ldots,i_N) \in \Real^{I \times \cdots \times I}$,  given its exponential size of $I^N$ entries.
Therefore, the TMM presented by (\ref{eq:TMM}) is not tractable directly but only  in a distributed form,  whereby the coefficient tensor is approximated by low-rank tensor network decompositions with non-negative cores.
%
%
 The TMM  based on the HT decompositions  is promising for  multiple classification problems with partially missing training data, and potentially even  regardless of the distribution of missing data  \citep{TMM_Sharir}.

\section{Convolutional Arithmetic Circuits and HT  Networks}
\sectionmark{Convolutional Arithmetic Circuits}

Conceptually, the Convolutional Arithmetic Circuit (ConvAC)  can be divided into three parts:
(i) the first (input) layer   is the representation layer which transforms input vectors $(\bx_1,\ldots,\bx_N)$ into $N \, I$ real valued scalars $\{f_{\theta_i}(\bx_n)\}$ for $n=1,\ldots,N$ and $i=1,\ldots, I$. In other words, the representation functions, $f_{\theta_i}: \Real^S \rightarrow \Real$, $i=1,\ldots,I$, map each local patch $\bx_n$ into a feature space of dimension $I$;
(ii) the second, a key  part, is a  convolutional arithmetic circuit consisting of many hidden layers that takes the $N \, I$ measurements (training samples) generated by the representation layer;  (iii) the output layer, which can be represented by  a  matrix, $\tW^{(L)} \in \Real^{R \times C}$, which computes $C$ different score functions $h_{c}$ \citep{Shashua-HT}.

\begin{figure}[t!]
\begin{center}
 \includegraphics[width=0.95\textwidth]{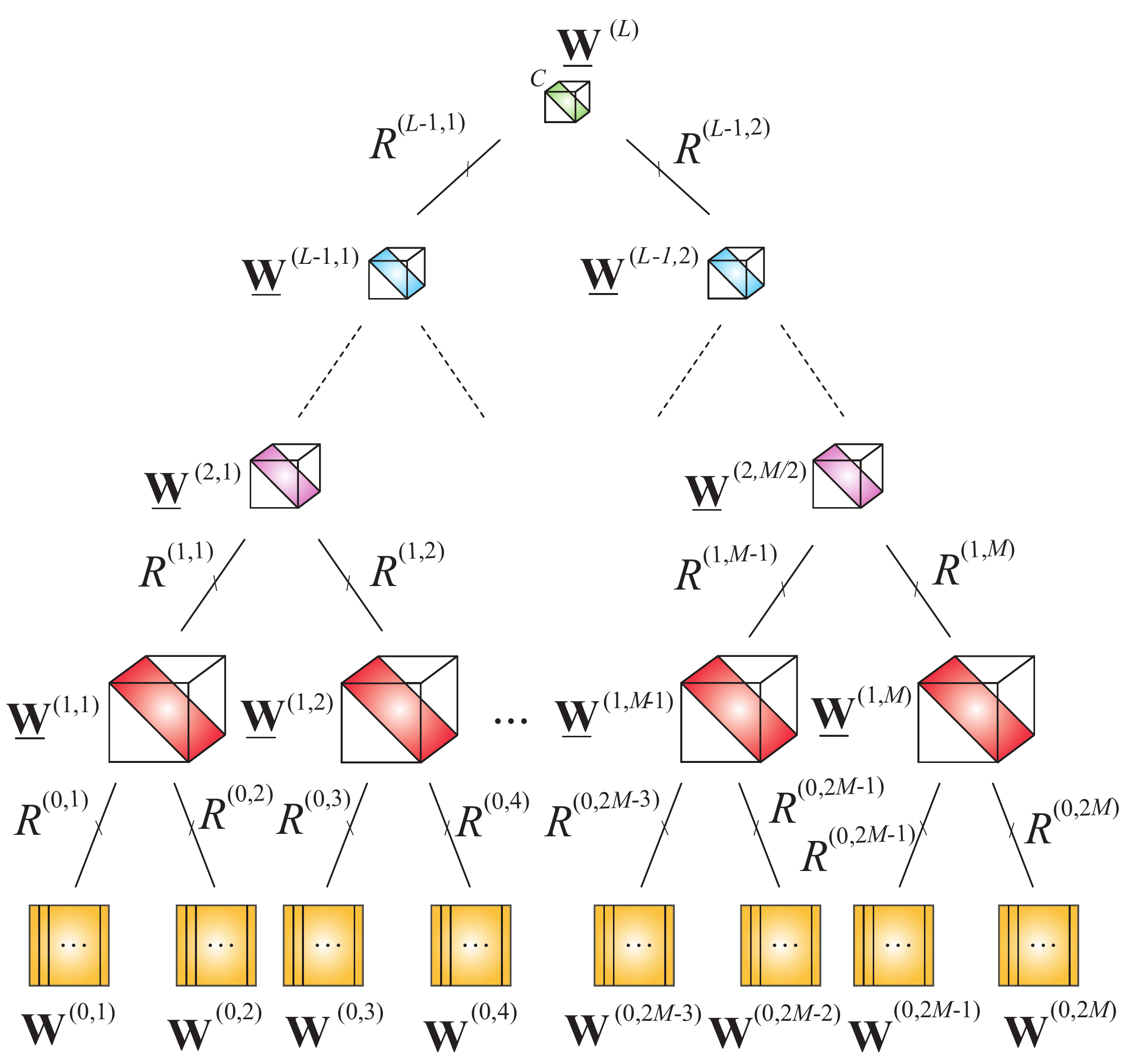}
\end{center}
  \caption{Architecture of a simplified  Hierarchical Tucker (HT) network consisting of  core tensors, $\tW^{(l,j)}$, with diagonal frontal slices, which approximates a grid tensor for a set of score functions $\{h_c\}$, for $c=1,\ldots,C$. This model corresponds to a ConvAC deep learning network. The HT tensor network consists of $L=\log_2(N)$ hidden layers with pooling-2 windows. For simplicity, we  assumed   $N=2M=2^L$ input patches, and equal HT-ranks $R^{(l,j)}$ in each layer, i.e., $R^{(l,1)}=R^{(l,2)}=\cdots =R^{(l)}$ for $l=0,1,\ldots, L-1$. Moreover, the cores have non-zero entries only on the diagonal of tensor cubes, as indicated by red rectangles. The representation layer is not  explicitly shown.}
  \label{Fig:DLHTC}
\end{figure}

 Once the set of score functions has been formulated as (\ref{eq:score1}), we can construct (design) a suitable multilayered  or distributed representation for DCNN representation. The objective is to estimate the parameters $\theta_1,\ldots, \theta_I$ and coefficient   tensors{\footnote{These tensors share the same entries, except for the  parameters in the output layer.}} $\underline \bW_{\,1},\ldots, \underline \bW_{\,C}$.
Since the tensors are of $N$th-order and each with $I^N$ entries, in order to avoid the curse of dimensionality  we employ low-rank tensor network representations. Note that a direct implementation of (\ref{eq:score1}) or (\ref{eq:TMM}) is intractable owing to a huge number of parameters.

The simplified HT tensor network, shown in Figure~\ref{Fig:DLHTC}, contains sparse 3rd-order core tensors $\underline \bW^{(l,j)}  \in \Real^{R^{(l-1,2j-1)} \times R^{(l-1,2j-1)} \times R^{(l,j)}}$  for $l=1,\ldots,L-1$  and  matrices $\bW^{(0,j)} = [\bw_1^{(0,j)}, \ldots, \bw_{R^{(0,j)}}^{(0,j)}] \in \Real^{I_j \times R^{(0,j)}}$ for $l=0$ and $j=1, \ldots, N/2^{\,l}$ (for simplicity, we assumed that $N=2M=2^L$).
In order to mimic basic features of the standard  ConvAC, we assumed that
$R^{(l,j)} = R^{(l)}, \; \forall j$ and  that frontal slices of the core tensors $\underline \bW^{(l,j)}$ are diagonal matrices,
 with entries $\tW^{(l,j)}(r^{(l-1)}, r^{(l-1)}, r^{(l)}) = w^{(l,j)}_{r^{(l-1)},r^{(l)}}$,
 as illustrated in Figure~\ref{Fig:DLHTC}.
  The top (output) layer is represented by a 3rd-order tensor $\tW^{(L)} \in \Real^{R^{(L-1)} \times R^{(L-1)} \times C}$ with diagonal frontal slices $\tW^{(L)}(:,:,c)= \mbox{diag}\{\lambda^{(c)}_1,\ldots,\lambda^{(c)}_{R^{(L-1)}}\}$, $c=1,\ldots,C$.
 Note that the sparse core tensors $\underline \bW^{(l,j)}$, with diagonal frontal slices, can be represented by  dense matrices defined as
 $\bW^{(l,j)} \in \Real^{R^{(l-1)} \times R^{(l)}}$ for $l=1,\ldots,L-1$, and the top tensor can be also represented by a dense matrix $\bW^{(L)}$ of size $R^{(L-1)} \times C$, in which each column corresponds to the diagonal matrix $\mbi \lambda^{(c)}= \mbox{diag}\{\lambda^{(c)}_1,\ldots,\lambda^{(c)}_{R^{(L-1)}}\}$.
	
The so simplified HT tensor network can be mathematically described in the following recursive form
\begin{align}
\bW^{(0, j )}  &= [\bw_1^{(0, j )}, \ldots, \bw^{(0, j )}_{R^{(0)}}] \in
\Real^{I_{j} \times R^{(0)}} \notag \\
	\bW^{( \leq 1, j )}_{r^{(1)}}  &= \sum_{r^{(0)}=1}^{R^{(0)}}
	w^{(1,j)}_{r^{(0)}, r^{(1)}}  \cdot
	\left( \bw^{( 0, 2j-1)}_{r^{(0)}} \circ \bw^{( 0, 2j)}_{r^{(0)}} \right) \in
\Real^{I_{2j-1} \times I_{2j}} \notag \displaybreak[3]\\
	&\cdots\label{eqn_convAC_vector}  \\
	 \underline \bW^{( \leq l, j )}_{\;r^{(l)}}  &= \sum_{r^{(l-1)}=1}^{R^{(l-1)}}
	w^{(l,j)}_{r^{(l-1)}, r^{(l)}}  \cdot
	\left( \underline \bW^{(\leq l-1, 2j-1)}_{\;r^{(l-1)}} \circ
\underline \bW^{(\leq l-1, 2j)}_{\;r^{(l-1)}}  \right) \notag\displaybreak[3] \\
 & \cdots  \notag \displaybreak[3] \\
	\underline \bW^{(\leq L-1, j )}_{r^{(L-1)}}  &= \sum_{r^{(L-2)}=1}^{R^{(L-2)}}
	w^{(L-1,j)}_{r^{(L-2)}, r^{(L-1)}}  \cdot
	\left(\underline \bW^{(\leq L-2, 2j-1)}_{\;r^{(L-2)}} \circ
 \underline \bW^{(\leq L-2, 2j)}_{\; r^{(L-2)}}  \right) \notag
 \end{align}
for $j=1,\ldots, 2^{L-l}$ and
 \begin{align}
	 \underline \bW_{c}
	&= \sum_{r^{(L-1)}=1}^{R^{(L-1)}} \lambda^{(c)}_{r^{(L-1)}} \cdot
	\left(\underline \bW^{(\leq L-1, 1)}_{\; r^{(L-1)}} \circ
  \underline \bW^{(\leq L-1, 2)}_{\; r^{(L-1)}} \right) \in \Real^{I_1 \times \cdots \times I_N},\notag\\
	\end{align}
where $\lambda^{(c)}_{r^{(L-1)}}= \tW^{(L)}(r^{(L-1)},r^{(L-1)},c)$.

In a  special case, when the weights in each layer are shared,  i.e.,
$\bW^{(l,j)}= \bW^{(l)}, \;\; \forall j$,
 the above equation can be considerably  simplified to
	\be
	\underline \bW^{(\leq l)}_{r^{(l)}} =
	\sum_{r^{(l-1)}=1}^{R^{(l-1)}}
	w^{(l)}_{r^{(l-1)}, r^{(l)}}  \;
	(\underline \bW^{(\leq l-1)}_{r^{(l-1)}}  \circ \underline \bW^{(\leq l-1)}_{r^{(l-1)}})
	\ee
for the layers $l= 1,\ldots, L$, 
where
$\underline \bW^{(\leq l)}_{r^{(l)}} = \underline \bW^{(\leq l)} (:,\ldots,:,r^{(l)}) \in \Real^{I\times\cdots \times I}$ are
sub-tensors of $\underline \bW^{(\leq l)}$, for each $r^{(l)} = 1,\ldots,R^{(l)}$,
and $w^{(l)}_{r^{(l-1)}, r^{(l)}}$ is the $(r^{(l-1)},r^{(l)})$th entry of the weight matrix $\bW^{(l)} \in \Real^{R^{(l-1)}\times R^{(l)}}$.\\

However, it should be noted that the simplified HT model shown in Figure~\ref{Fig:DLHTC} has a limited ability to approximate an arbitrary grid  tensor, $\underline \bW(h_{c})$, due to the strong constraints imposed on core tensors.
A more flexible and powerful HT model is shown in Figure~\ref{Fig:DLHT}, in which the constraints imposed on 3rd-order  cores have been completely removed.
In fact this is the standard HT tensor network, which can be mathematically expressed as
(with a slight abuse of notation), as follows
\be
\label{eq:DLTTF}
	\bw^{( \leq 0, j )}_{r}  &=& \bw^{(0,j)}_r  \notag \\ \notag  \\
	 \underline \bW^{( \leq l, j )}_{\;r}  &=& \sum_{r_1=1}^{R^{(l-1,2j-1)}} \sum_{r_2=1}^{R^{(l-1,2j)}}
	w^{(l,j)}_{r_1, r_2, r}  \cdot
	\left( \underline \bW^{(\leq l-1, 2j-1)}_{\;r_1} \circ
\underline \bW^{(\leq l-1, 2j)}_{\;r_2}  \right) \notag \\ \\
\underline \bW_{\,c} &=&    \underline \bW^{( \leq L, 1 )}, \notag
	\ee
for $l=1,\ldots,L$ and $\;\;j=1, \ldots, 2^{L-l}$. 

\begin{figure}[t!]
\begin{center}
 \includegraphics[width=0.99 \textwidth]{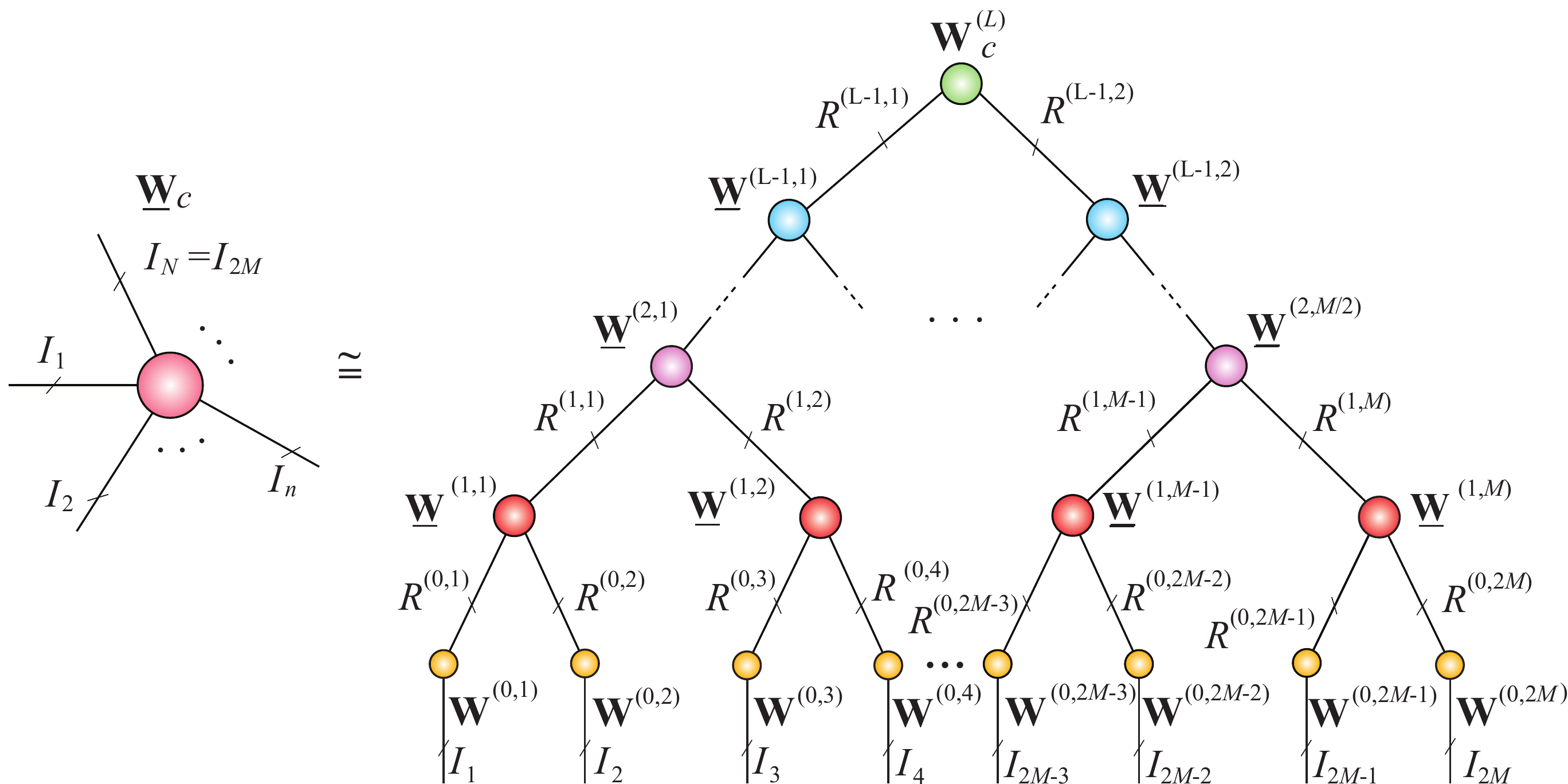}
\end{center}
  \caption{Hierarchical Tucker (HT) unconstrained  tensor network for the distributed representation and/or approximation of coefficient tensors, $\underline \bW_{\,c}$, of the score functions $h_{c}(\bx_1,\ldots,\bx_N)$. In general, no constraints are imposed on core tensors. The representation layer is not explicitly shown. For a grid tensor of the set of score functions, the matrices $\bW^{(L)}_c$ for $c=1,\ldots,C$ build  up a 3rd-order tensor corresponding to the output layer.}
  \label{Fig:DLHT}
\end{figure}

The HT model can be further extended by using more flexible and general  Tree Tensor Networks States (TTNS). As illustrated  in Figure~\ref{Fig:DLTTNS}, the use of the TTNS, instead HT tensor networks, allows for  more flexibility in the choice of the size of pooling-window, as the pooling size in each hidden layer can be adjusted  by applying core tensors with a suitable variable order in each layer.  For example,  if we use 5th-order (4rd-order) core tensors instead 3rd-order cores, then the pooling layer  will employ a size-4 pooling window (size-3 pooling) instead of only size-2 pooling window when using 3rd-order core tensors in HT  tensor networks. For more detail regrading HT networks and their generalizations to  TTNS, see Part 1 of our monograph \citep{Part1}.

\begin{figure}[t]
\begin{center}
 \includegraphics[width=0.799\textwidth]{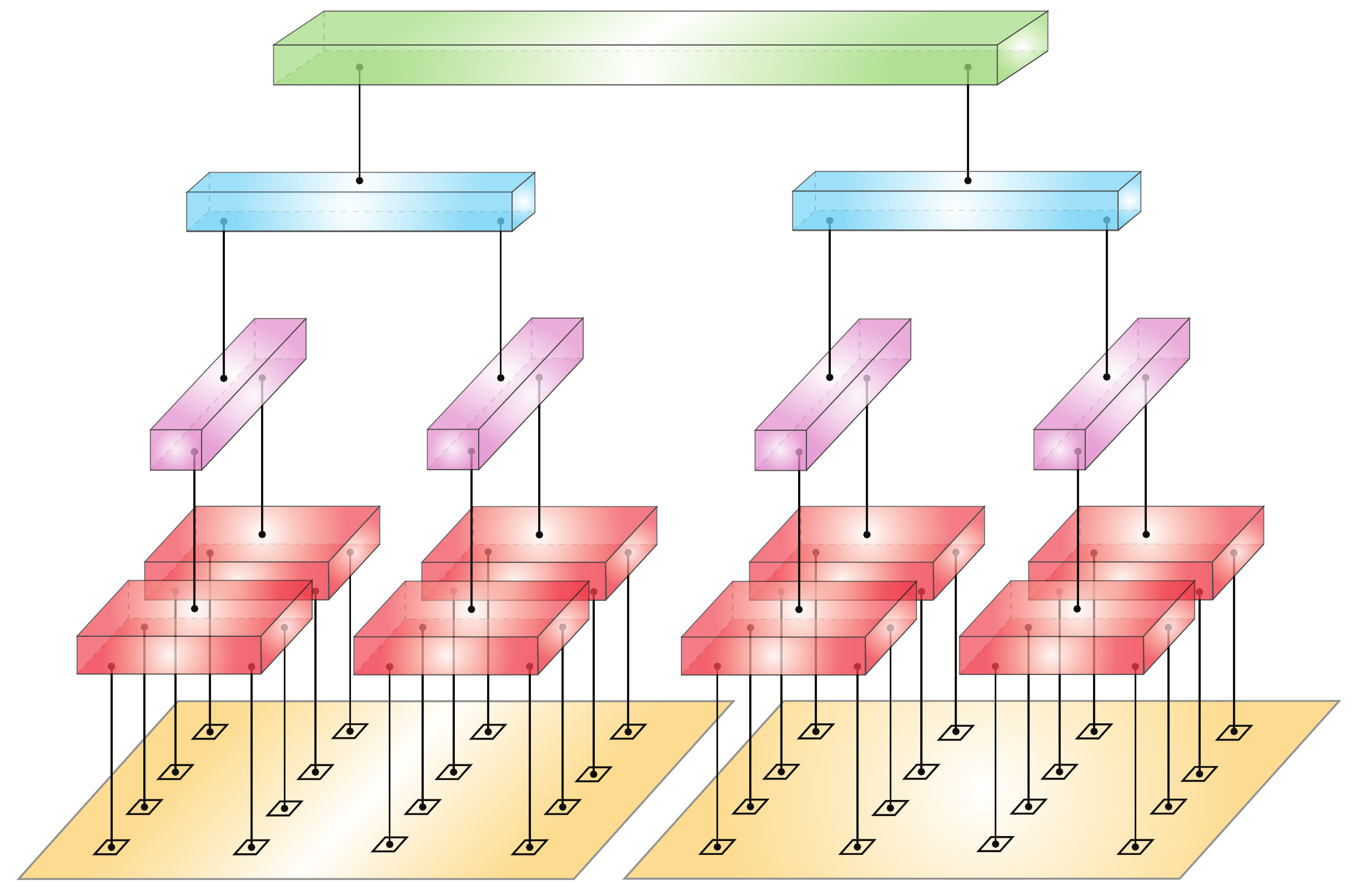}
\end{center}
  \caption{Tree Tensor Networks States (TTNS) with a variable order of core tensors. The rectangular prisms represent core tensors of orders 5 and 3, that allow pooling with the respective window sizes of  4 and 2. The first hidden layer comprises 5th-order cores, while the second and third layer consist of 3rd-order cores.}
  \label{Fig:DLTTNS}
\end{figure}

 \section{Convolutional Rectifier NN  Using  Nonlinear TNs}

 The convolutional arithmetic circuit (ConvACs) model employs the standard outer (tensor) products,  which for two tensors,
  $\underline \bA \in \Real^{I_1 \times \cdots \times I_N}$ and
  $\underline \bB \in \Real^{J_1 \times \cdots \times J_M}$, are defined as
\be
(\underline \bA \circ \underline \bB)_{i_1,\ldots, i_N,j_1,\ldots,j_M}= a_{i_1,\ldots, i_N} \cdot b_{j_1, \ldots, j_M}. \notag
\ee
In order to convert ConvAC tensor models to  widely used  convolutional rectifier networks, we need to employ the generalized (nonlinear) outer products, defined as \citep{Shashua_GTD2016}
 \be
(\underline \bA \circ_{\rho} \underline \bB)_{i_1,\ldots, i_N,j_1,\ldots,j_M}= \rho(a_{i_1,\ldots, i_N}, b_{j_1, \ldots, j_M}),
\ee
where the $\rho$ operator can take various forms, e.g.,
\be
\rho(a,b) = \rho_{\sigma,P} (a,b) = P[\sigma(a),\sigma(b)],
\ee
and is  referred to  as the activation-pooling  function{\footnote{The symbols $\sigma(\cdot)$ and $P(\cdot)$ are respectively the activation and pooling functions of the network.}}, which meets the  associativity and the commutativity requirements: $\rho(\rho(a,b),c)=\rho(a,\rho(b,c))$ and $\rho(a,b)=\rho(b,a),\;\; \forall a,b,c \in \Real$).
For   two vectors, $\ba \in \Real^J$ and  $\bb \in \Real^J$,  it is defined as  a matrix $\bC = \ba \circ_{\rho} \bb \in \Real^{I \times J}$, with entries $c_{ij}= \rho(a_i , b_j)$. Note that the nonlinear function $\rho$ can also take special form: $c_{ij}= \max\{a_i b_j,0\}$ or $c_{ij}= \max\{a_i,b_j,0\}$.

%
 For a particular case of the convolutional rectifier network with $\max$ pooling, we can use the following activation-pooling operator
 \be
 \rho_{\sigma,P} (a,b) = \max\{[a]_+,[b]_+\}=\max\{ a,b,0 \}.
\ee
In an analogous way, we can define the generalized Kronecker and the Khatri-Rao products.\\
\begin{example} Consider a generalized CP decomposition, which corresponds to a shallow  rectifier network in the form \citep{Shashua_GTD2016}
 \be
 \underline \bW_{\,c} = \sum_{r=1}^{R} \lambda^{(c)}_r  \cdot (\bw^{(1)}_r \circ_{\rho} \bw^{(2)}_r \circ_{\rho} \cdots  \circ_{\rho} \bw^{(N)}_r),
\label{eq:GCPD}
 \ee
 where the coefficients $\lambda^{(c)}_r $ represent weights of the output layer, vectors  $\bw^{(n)}_r  \in \Real^{I_n}$ are weights in the hidden layer and the operator $\circ_{\rho}$ denotes the nonlinear outer products defined above.
 %

The generalized  CP decomposition can be expressed in an equivalent vector  form, as follows
\be
\bw_c = \mbox{vec}(\tW_c)= [\bW^{(N)} \odot_{\rho} \bW^{(N-1)} \odot_{\rho} \cdots \odot_{\rho} \bW^{(1)}] \mbi \lambda^{(c)},
\ee
where $\odot_{\rho}$ is the generalized  Khatri-Rao product of two matrices.

 It should be noted that if we employ  weight sharing, then all vectors $\bw_r^{(n)}=\bw_r,\; \forall n$, and consequently the coefficient tensor, $ \underline \bW_{\,c}$, must be a symmetric tensor which further limits the ability of  this model to approximate a desired function.\\
\end{example}

\begin{example} Consider the simplified  HT model as shown Figure~\ref{Fig:DLHTC}, but
with  the generalized outer product defined above. Such nonlinear tensor networks can be rigorously  described for $I_1=I_2 =\cdots =I_N=I$, as follows
\be
\label{eqn_convAC_vector_rect}
\bw_{r^{(0)}}^{(\leq 0,j)} &=& \bw_{r^{(0)}}^{(0,j)} \notag \\
	 \underline \bW^{( \leq l, j )}_{\;r^{(l)}}  &=& \sum_{r^{(l-1)}=1}^{R^{(l-1)}}
	w^{(l,j)}_{r^{(l-1)}, r^{(l)}}  \cdot
	\left( \underline \bW^{(\leq l-1, 2j-1)}_{\;r^{(l-1)}} \circ_{\rho} \;
\underline \bW^{(\leq l-1, 2j)}_{\;r^{(l-1)}}  \right) \in \Real^{ I \times \cdots \times I} \notag \\
 \tW_{c} &=& \tW^{(\leq L, 1)}, \notag
	\ee
for $l=1,\ldots, L$ and $j=1,\ldots,2^{L-l}$, or in an equivalent matrix form as 
	\begin{equation}
\label{eqn_convAC_matrix}
	\begin{split}
	\bW^{(\leq 0, j)}  &= \bW^{(0, j)}  \in \Real^{I \times R^{(0)}}\\
	\bW^{(\leq l, j)}  &= \left(\bW^{(\leq l-1, 2j-1)} \odot_{\rho} \bW^{(\leq l-1, 2j)} \right) \, \bW^{(l,j)}   \in \Real^{I^{2^{\,l}} \times R^{\,(l)}} \\
	\text{vec} ( \underline \bW_{c} )
		& =	\left(\bW^{(\leq L-1, 1)} \odot_{\rho} \bW^{(\leq L-1, 2)} \right) \,
\mbi \lambda^{(c)} \in \Real^{I^N },
	\end{split}
	\end{equation}
for  $	 l=1,\ldots, L-1 $, where the core tensors are reduced to matrices $\bW^{(l,j)} 
\in \Real^{R^{(l-1)} \times R^{(l)}}$ with entries
$w^{(l,j)}_{r^{(l-1)},r^{(l-1)},r^{(l)}}$.
\end{example}

%

	 We should emphasize that the HT/TTNS architectures are not the only  suitable   TN  architectures
	  which can be used  to model  DCNNs, and the whole family of
	   powerful tensor networks can be employed for this purpose.
	Particularly attractive and simple are the TT/MPS, TT/MPO  and TC models for DNNs, for which efficient decomposition and tensor completion algorithms already exist.
The  TT/MPS, TT/MPO and TC networks provide not only simplicity in comparison to HT, but also very deep TN structures, that is, with $N$ hidden layers.
 Note that the  standard HT model generates architectures of DCNNs with $L=\log_2(N)$ hidden layers, while TT/TC networks may employ $N$ hidden layers.
 Taking into account the current  trend in deep leaning  to use a large number of hidden layers,
   it would be    quite attractive  to employ QTT tensor  networks with a relatively large number of hidden layers, $L=N \cdot \log_2(I)$.

 To summarize,  deep convolutional neural networks may be considered as a special case of hierarchical architectures,  which can be indirectly simulated and optimized  via
 relatively  simple  and well understood tensor networks, especially HT, TTNS,  TT/MPS, TT/MPO, and  TC (i.e., using unbalanced or balanced binary trees and graphical models). However,  more sophisticated tensor network  diagrams with loops, discussed in the next section may  provide potentially  better performance and the ability to generate novel architectures of DCNNs.

\section{MERA and 2D TNs for a Next Generation of DCNNs}

The  Multiscale Entanglement Renormalization Ansatz (MERA) tensor network was  first
introduced by \citet{Vidal2008MERA}, and  for this network  numerical algorithms
to minimize some specific cost functions or local Hamiltonians used in quantum physics already exist \citep{Evenbly-Vidal09alg}.
The MERA is a  relatively new tensor network, widely investigated in quantum
physics based on variational Ansatz, since it is capable of capturing many of the key complex physical properties
of strongly correlated ground states \citep{evenbly2015tensor}.
The MERA also shares many  relationships
with the AdS/CFT (gauge-gravity) correspondence, through its   complete holographic duality
with the tensor networks framework.
Furthermore, the MERA  can be regarded as a
TN  realization of an orthogonal wavelet
transform  \citep{Mera_wavelets1,Mera_wavelets2,matsueda2016MERA}.

\begin{figure}[htpp]
 (a) \hspace{4cm} (b)
 \vspace{-0.6cm}
\begin{center}
    \includegraphics[width=0.328\textwidth]{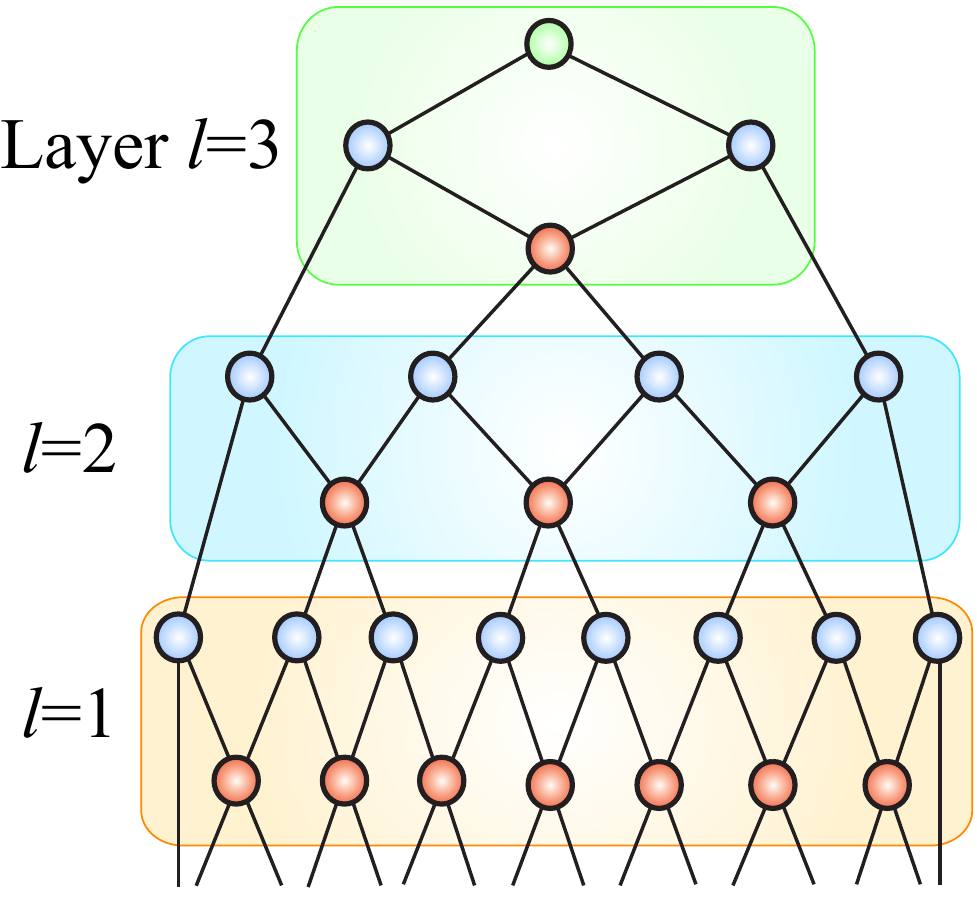} \hspace{0.4cm}
   \includegraphics[width=0.475\textwidth]{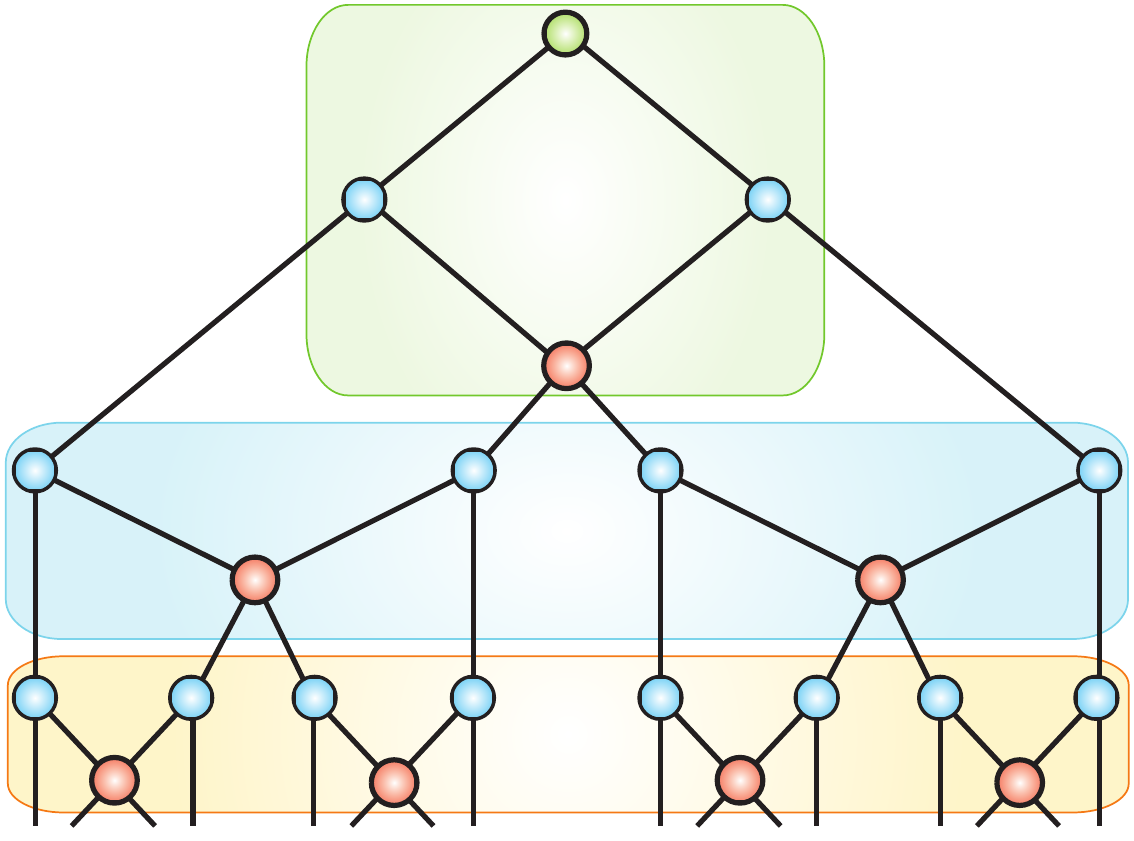} \\
   \end{center}
  (c)
\vspace{-0.6cm}
\begin{center} 
  \includegraphics[width=0.945\textwidth]{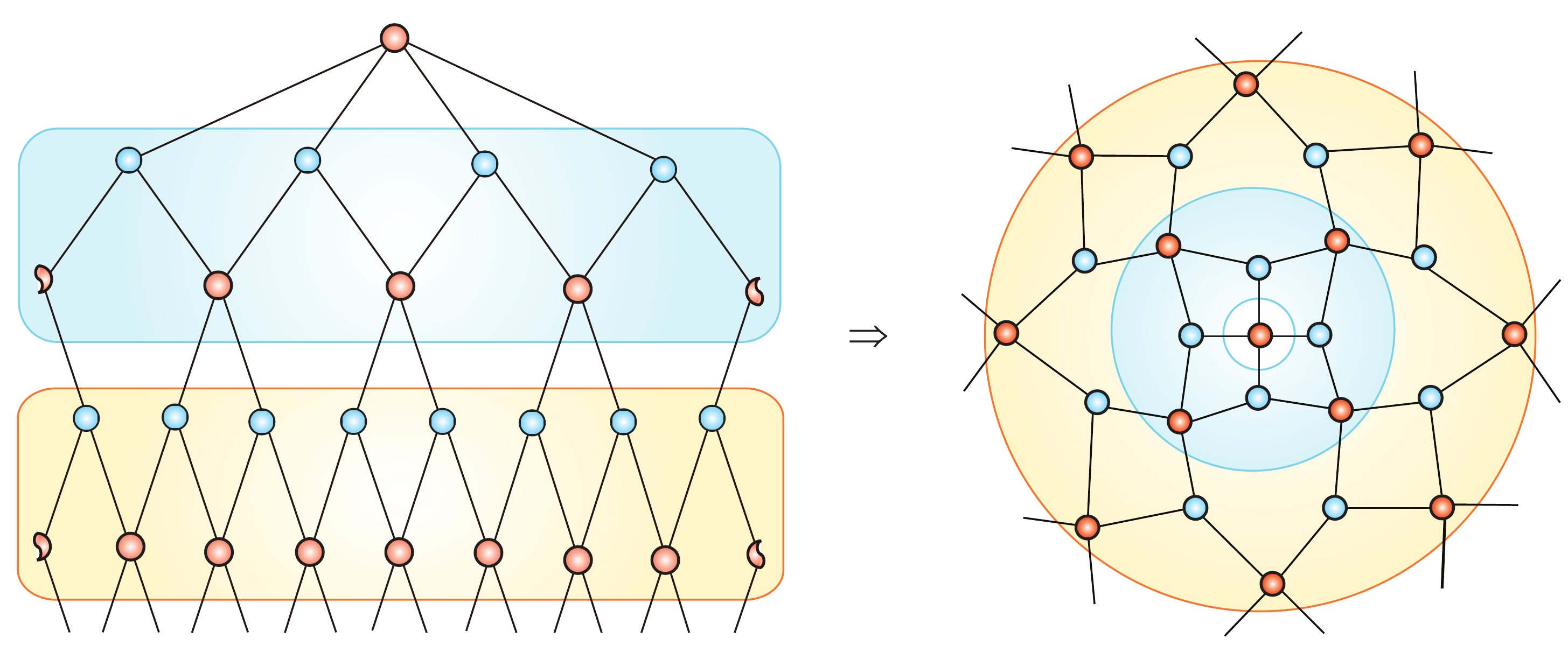}
   \end{center}
   (d)
   \vspace{-0.4cm}
\begin{center} 
  \includegraphics[width=0.945\textwidth]{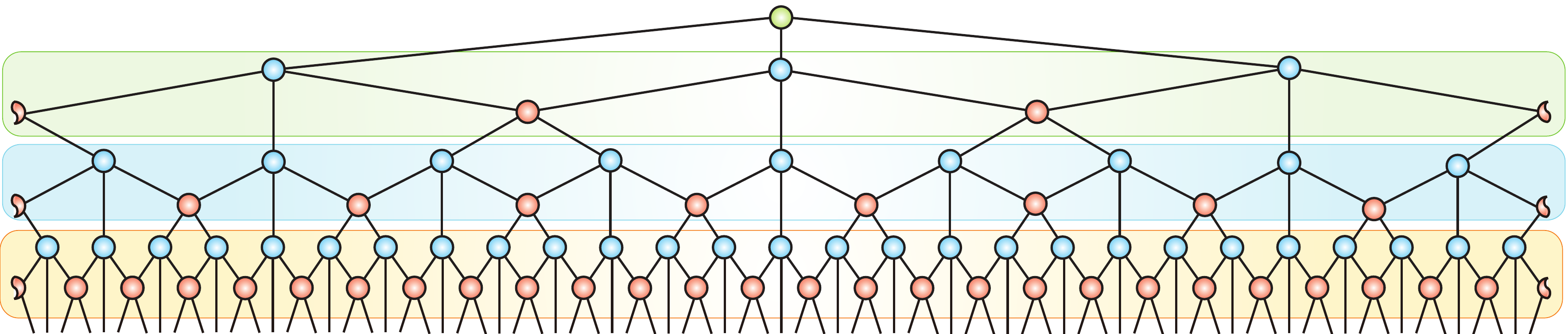}
   \end{center}
  \caption{Various architectures of MERA  tensor networks for the new generation of deep convolutional neural networks. (a) Basic binary MERA tensor network.  The
consecutive layers of disentangling and coarse-graining cores are indicated by different colors.  For the network shown in (a) the number of  cores after each such set of operations is approximately halved.
(b) Improved (lower complexity) MERA network. (c) MERA with periodic boundary conditions. Red half-circles mean that they are connected, i.e., they build up  one 4th-order core, as illustrated on the right panel of (c).
(d) Ternary MERA (with the Periodic Boundary Conditions (PBC)) in which coarse grainers are also 4th-order tensors, i.e., three sites (modes) are coarse-grained into one effective site (mode). }
  \label{Fig:DLMERA}
\end{figure}

For simplicity,  this section focuses on  1D
binary and ternary  MERA architectures (see Figure~\ref{Fig:DLMERA}(a) for basic binary MERA).
Instead of writing  complex mathematical formulas, it is more convenient to describe MERA tensor networks graphically, as illustrated in Figures \ref{Fig:DLMERA}(a), (b) and~(c).
Using  the terminology from quantum physics, the standard binary MERA architecture contains three classes of  core tensors: (i)  disentanglers -- 4th-order cores; (ii) isometries,  also called the coarse-grainers, which are typically 3rd-order cores  for binary MERA and 4th-order cores for ternary MERA;  and (iii) one output core  which is usually a matrix or a 4th-order core.
 Each MERA layer is constructed  of a row of disentanglers and a row of  coarse-grainers  or isometries. Disentanglers remove the short-scale entanglement
between the adjacent  modes,  while isometries renormalise each pair of modes  to a single mode.
Each  renormalisation layer performs these operations on a scale of   different length.

%
 From the  mapping perspective, the nodes (core tensors) can be considered as processing units, that is, the 4th-order cores map matrices to other matrices, while the coarse-grainers take  matrices and map them to  vectors.
 The key idea here is to realize that the ``compression'' capability arises from the hierarchy and the entanglement. As a matter of fact, the MERA network embodies the mutual information chain rule.
%

In other words, the main idea underlying MERA is that of disentangling the system
at scales of various lengths, following coarse graining Renormalization Group (RG) flow
in the system.
The MERA  is particularly effective for (scale invariant) critical points
of  physical systems.
The  key   properties of MERA can be summarized, as follows \citep{evenbly2015tensor}:

\begin{itemize}

\item MERA can capture scale-invariance in input data;

\item  It reproduces a polynomial decay of correlations between   inputs,
 in contrast to HT or TT networks which reproduce only exponential decay of correlations;

\item MERA  has the ability to  compress tensor data  much better than TT/HT tensor networks;

\item It reproduces a  logarithmic correction to the area law, therefore MERA is  more powerful tensor network than HT/TTNS or TT/TC networks;

\item MERA can be efficiently contracted due to  unitary constraints imposed on core tensors.

\end{itemize}

Motivated by these features, we are currently investigating  MERA tensor
networks as  powerful tools to model and analyze DCNNs.
The main objective  is to establish  a precise connection between MERA tensor networks and extended models of DCNNs. This connection may provide exciting new insights about deep learning while also allowing
for the construction of  improved families of DCNNs, with  potential application to more efficient data/image classification, clustering and prediction.
In other words,  we conjecture  that the
MERA and 2D TNs  (see Figure~\ref{Fig:DLHC}) will lead to useful new results, potentially
allowing  not only for better characterization of  expressive power of DCNNs, but also for new
practical implementations.
Conversely,  the links  and relations between TNs and DCNNs could
 lead to useful advances in the design of  novel deep neural networks.

 \begin{figure}[t]
  (a) \hspace{5cm} (b)
\vspace{-0.4cm}
\begin{center} 
  \includegraphics[width=0.255\textwidth]{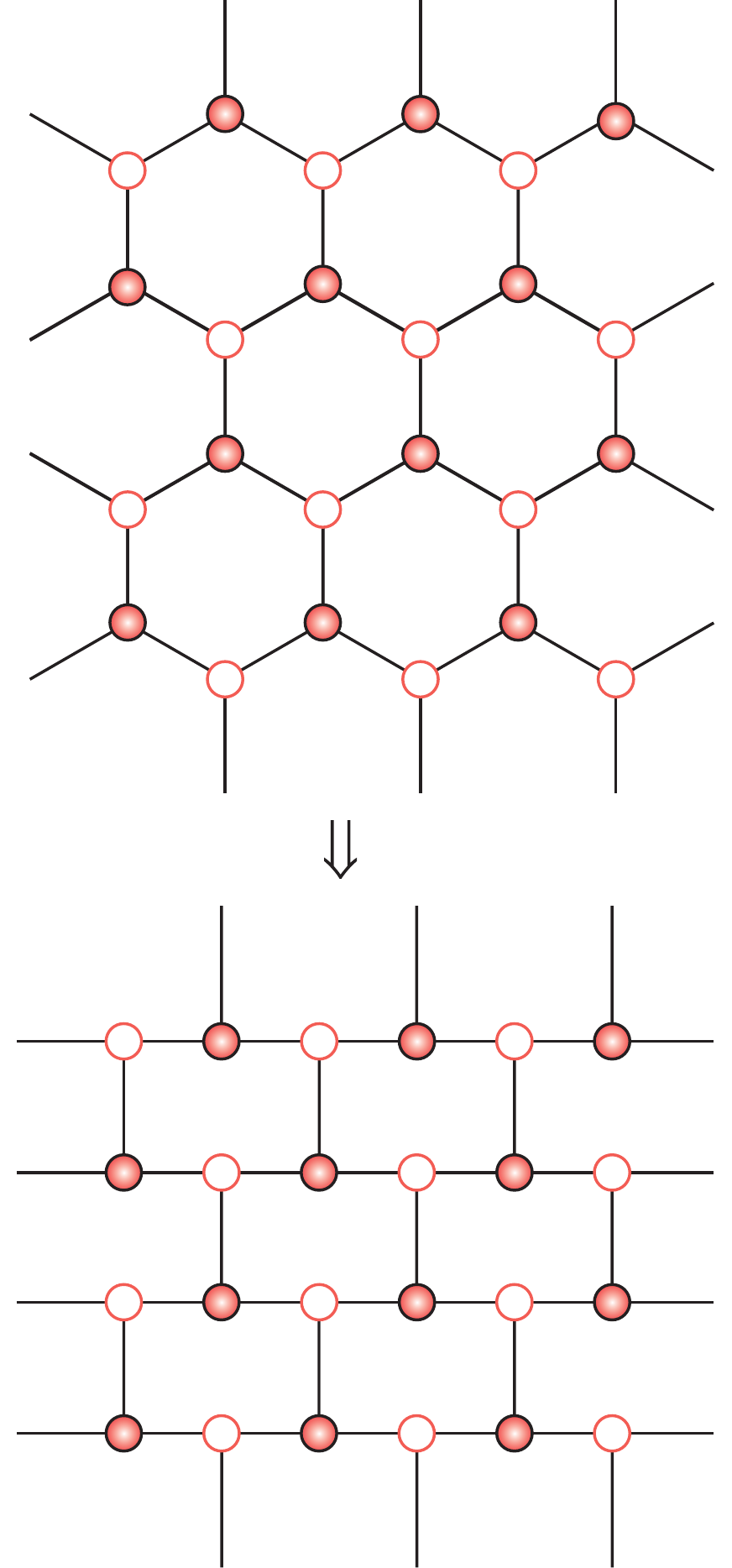}
  \hspace{0.9cm}
   \includegraphics[width=0.40\textwidth]{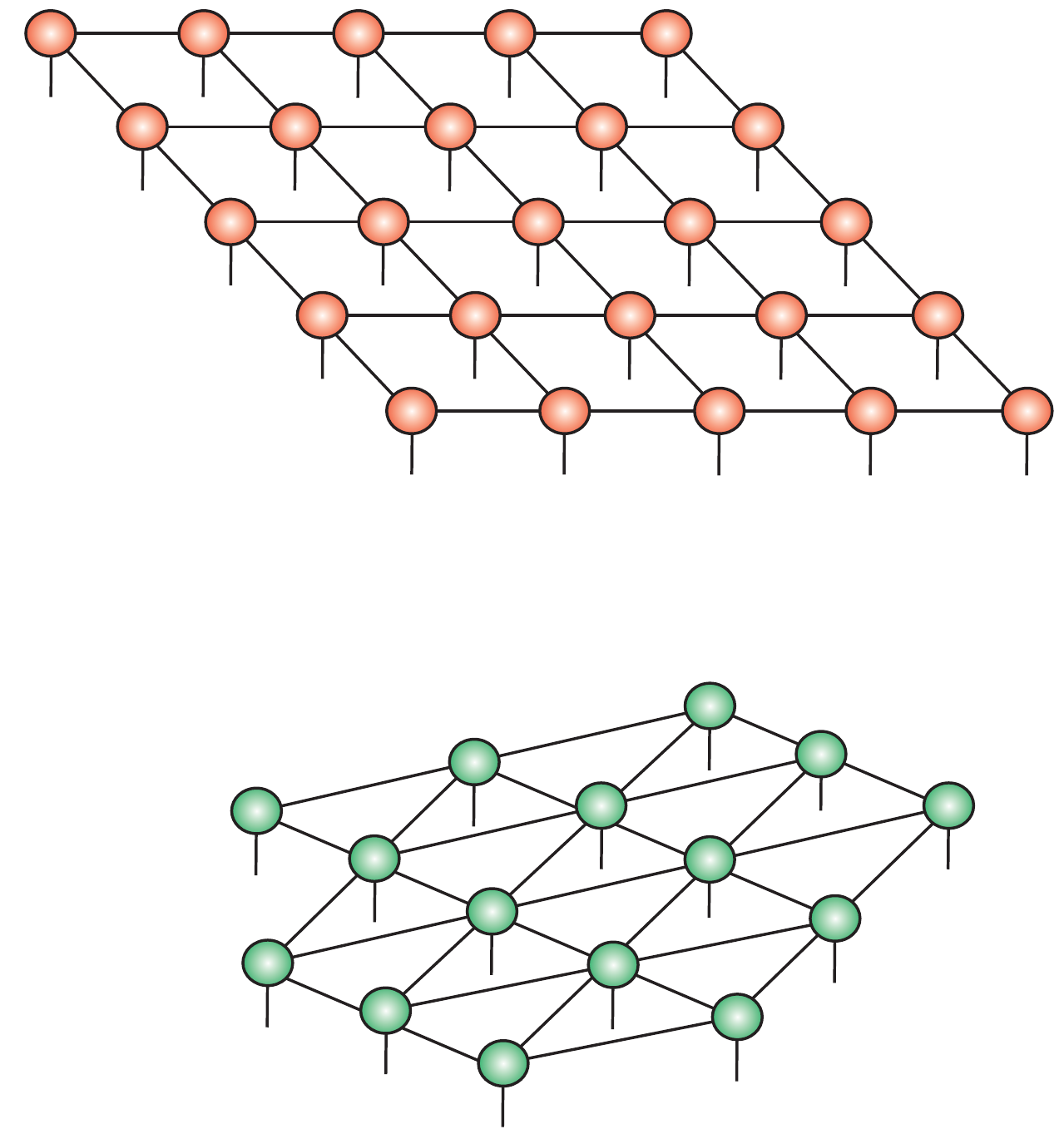}
   \end{center}
  \caption{Architectures of 2D tensor networks as a potential model for DNNs.
(a) A honey-comb lattice (top), which is equivalent to the brick-wall lattice (bottom), and (b) PEPS tensor  networks with 5th-order (top) and up to 7th-order cores (bottom).}
  \label{Fig:DLHC}
\end{figure}

The MERA tensor networks, shown in Figure~\ref{Fig:DLMERA}, may  provide a much higher expressive power of deep learning in comparison to networks corresponding to HT/TT architectures, since this class of tensor networks can model more complex long-term correlations between input instances. This follows from the facts that  for  HT and  TT tensor networks correlations between input variables  decay exponentially, and the
entanglement entropy saturates to a constant,  while  the more sophisticated MERA tensor networks
provide polynomially decaying correlations.

We firmly believe that the insights
into the  theory of tensor networks and quantum many-body physics
can provide better theoretical  understanding of  deep learning,
together with the guidance  for optimized  DNNs design.

To summarize, the  tensor network methodology   and architectures discussed briefly in this section may be extended to allow analytic construction of  new DCNNs. Moreover, systematic investigation of the correspondences between DNNs and  a wide spectrum of TNs
can provide a very fruitful perspective including verification of the existing conjectures and claims about
 operational similarities and correspondences between DNNs and TNs into  a more rigorous and
constructive framework.

\chapter{Discussion and Conclusions}
\label{Chapter5}

\vspace{-0.9cm}

Machine learning and data mining  algorithms are becoming increasingly important in the analysis of  large volume, multi-relational and multi-modal datasets, which are often conveniently represented as multiway arrays or tensors. The aim of this  monograph has therefore been to provide the ML and data analytics communities with both a state-of-the-art review and new directions in tensor decompositions and tensor network  approaches for fundamental problems of  large-scale optimization.  Our focus has been on supervised and unsupervised  learning with tensors, tensor regression, support tensor machines,  tensor canonical coloration analysis, higher-order and kernel partial least squares, and deep learning.

In order to demystify the concepts of tensor algebra and to provide a seamless transition from the flat view linear algebra to multi-view multilinear algebra, we have employed novel graphical representations of tensor networks  to interpret mathematical expressions directly and intuitively on  graphs rather than through tedious multiple-index relationships. Our  main focus has been on the Tucker, Tensor Train (TT)  and Hierarchical Tucker (HT) decompositions and their extensions or generalizations.
To make the material self-contained, we have also addressed the important concept of tensorization  and distributed  representation of structured lower-order data in tensor network formats, as a tool to provide efficient higher-order  representation of vectors, matrices and low-order tensors and the optimization of cost (loss) functions.

In order to combat the curse of dimensionality and possibly obtain  a linear or even sub-linear complexity of storage and computational complexities,  we have next addressed distributed representation of tensor functions through low-rank tensor networks. Finally, we have  demonstrated how such approximations can be used to solve a wide class of huge-scale linear/multilinear dimensionality reduction and related optimization problems
that are far from being tractable when using classical machine learning  methods.

The lynchpin of this work are low-rank tensor network approximations and the associated  tensor contraction algorithms, and we have elucidated how these can be used to convert otherwise intractable huge-scale optimization problems into a set of much smaller linked and/or distributed sub-problems of affordable sizes and complexity. In doing so, we have highlighted the ability of tensor networks to account for the couplings between the multiple variables, and for multimodal,
incomplete and/or noisy data.

The Part 2 finishes with a discussion of the potential of  tensor networks in the  design of improved and optimized deep learning neural network architectures. It is our firm conviction that the links between tensor  networks and  deep learning  provide both exciting new insights into  multi-layer  neural networks and a platform for the construction of  improved families of DNNs, with  potential applications in highly efficient data/image classifications, clustering and prediction.

\newpage
\appendix
\chapter*{Appendices}
\renewcommand{\thesection}{\arabic{section}}%
\def\theequation{A.\arabic{equation}}
\section{Estimation of Derivatives of GCF}\label{sec:apx_gcf}

The GCF of the observation and its derivatives are unknown, but can be estimated from the sample first GCF
\be
\hat{\phi}_{\bx}(\bu) = \frac{1}{T} \sum_{t = 1}^{T} \exp(\bu^{\text{T}} \bx_t) \, .
\ee
The $n$th-order partial derivatives of $\hat{\phi}_{\bx}(\bu)$ with respect to $\bu$ are $n$th-order tensors of size $I \times I \times \cdots \times I$, given by
\be
 \bpsi_{\bx}^{(n)}(\bu)  &=& \frac{\partial^{n} \hat{\phi}_{\bx}(\bu)}{\partial \bu^{n}}   = \frac{1}{T} \sum_{t = 1}^{T} \exp(\bu^{\text{T}} \bx_t) \underbrace{(\bx_t  \circ \bx_t \circ  \cdots \circ \bx_t)}_{\text{$n$ terms $\bx_t$}}\,  \notag \label{equ_dphi_x_n} \\
 &=& \diag_N(\frac{1}{T} \exp(\bu^{\text{T}} \bx_t) ) \times_1 \bX \times_2 \bX \cdots \times_N \bX  \notag .
\ee
The order-$N$ derivatives of the second GCF are then given by
\be
\tensor{\Psi}_{\bx}(\bu) = \sum_{k = 1}^{N}  \frac{(-1)^{k-1} \, (k-1)! }{\hat{\phi}_{\bx}(\bu)^k}\, \sum_{n_1, n_2, \ldots, n_k}  m_{n_1, n_2, \ldots, n_k} \, \calS(\bpsi_{n_1, n_2, \ldots, n_k}(\bu)) \, ,\;\notag\\ \label{equ_n_dgcf}
\ee
where $1\le n_1  \le n_2  \le \cdots \le  n_k \le N$, $n_1  + n_2  + \cdots +  n_k  = N$, $\bpsi_{\{n_1, n_2, \ldots, n_k\}}(\bu)$ are $N$th-order tensors constructed as an outer product of $k$ $n_j$th-order derivative tensors $\bpsi_{\bx}^{(n_j)}(\bu)$, that is
\be
\bpsi_{\{n_1, n_2, \ldots, n_k\}}(\bu) &=&  \bpsi_{\bx}^{(n_1)}(\bu) \circ \bpsi_{\bx}^{(n_2)}(\bu) \circ \cdots \circ \bpsi_{\bx}^{(n_k)}(\bu)   .
\ee
The operator $\calS(\tA)$ symmetrizes an $N$th-order tensor $\tA$, i.e., it yields a symmetric tensor, defined as
\be
	\calS(\tA) = \frac{1}{N!} \sum_{\pi_i}   \pi_i(\tA)   \, ,
\label{equ_sym_oper}
\ee
 where $\pi_i$ runs through the list of all $N!$ possible tensor permutations of the $N$th-order tensor $\tA$.
In the expression (\ref{equ_n_dgcf}), this symmetrization can be simplified to
operate only once on the resulting tensor $\tensor{\Psi}_{\bx}(\bu)$, instead on each tensor $\bpsi_{n_1, n_2, \ldots, n_k}(\bu)$.

The number $m_{n_1, n_2, \ldots, n_k}$ represents the total number of partitions of $\{1, 2, \ldots, N\}$ into $k$ distinct parts, each having $n_j$ entries, $j = 1, \ldots, k$. This number is given by
\be
m_{n_1, n_2, \ldots, n_k} = \frac{N!}{\prod_{j = 1}^{\bar{k}} (\bar{n}_j!)^{l_j} l_j!}  \, , \notag
\ee
where $\{1\le {\bar{n}}_1 < {\bar{n}}_2 < \cdots < {\bar{n}}_{\bar{k}} \le N \}$ denotes the set of distinct numbers of $\{n_1, n_2, \ldots, n_k\}$,  while $1\le l_j \le N$ represents the replication number of $\bar{n}_{j}$ in $\{n_1, n_2, \ldots, n_k\}$, for $j = 1, \ldots, \bar{k}$, i.e., $l_1 {\bar{n}}_1 + l_2  {\bar{n}}_2 +  \cdots + l_{\bar{k}}  {\bar{n}}_{\bar{k}} = N$.

The approximation of the derivative tensor $\tensor{\Psi}_{\bx}(\bu)$ for some low orders $N = 2, 3,\ldots, 7$, are given below
%
\begin{align}
\tensor{\Psi}_{\bx}^{(2)}(\bu) &= \frac{\bpsi^{(2)}(\bu)}{\hat{\phi}_{\bx}(\bu)}  - \frac{\bpsi^{(1)}(\bu) \circ \bpsi^{(1)}(\bu)}{\hat{\phi}_{\bx}^2(\bu)} \, , \notag \\
\tensor{\Psi}_{\bx}^{(3)}(\bu) &=  \frac{\bpsi^{(3)}{\hspace{-.3ex}}({\hspace{-.2ex}}\bu{\hspace{-.2ex}})}{\hat{\phi}_{\bx}(\bu)}  - \frac{3  \calS(\bpsi^{(1)}{\hspace{-.5ex}}({\hspace{-.2ex}}\bu{\hspace{-.2ex}}) \circ \bpsi^{(2)}{\hspace{-.5ex}}({\hspace{-.2ex}}\bu{\hspace{-.2ex}}))}{\hat{\phi}_{\bx}^2(\bu)}\notag\\
 &\quad + 2\,   \frac{\bpsi^{(1)}{\hspace{-.5ex}}(\bu) \circ \bpsi^{(1)}{\hspace{-.5ex}}(\bu) \circ \bpsi^{(1)}{\hspace{-.5ex}}(\bu)}{\hat{\phi}_{\bx}^3(\bu)} \, , \notag \\
%
\tensor{\Psi}_{\bx}^{(4)}(\bu) &=  \frac{\bpsi^{(4)}(\bu)}{\hat{\phi}_{\bx}(\bu)}  - \frac{4 \,\calS(\bpsi^{(1)}(\bu) \circ \bpsi^{(3)}(\bu)) + 3 \, \calS(\bpsi^{(2)}(\bu) \circ \bpsi^{(2)}(\bu))}{\hat{\phi}_{\bx}^2(\bu)}  \notag \\
& \quad + 2 \,
 \frac{ 6 \, \calS(\bpsi^{(1)}(\bu) \circ \bpsi^{(1)}(\bu)  \circ  \bpsi^{(2)}(\bu))}{\hat{\phi}_{\bx}^3(\bu)}  \notag \\
& \quad -6   \, \frac{\bpsi^{(1)}(\bu) \circ \bpsi^{(1)}(\bu) \circ \bpsi^{(1)}(\bu) \circ \bpsi^{(1)}(\bu)}{\hat{\phi}_{\bx}^4(\bu)} \,  , \notag \label{eq_Psi4}\displaybreak[3]\\
%
%
%
\tensor{\Psi}_{\bx}^{(5)}(\bu) &= \calS\left( \frac{\bpsi^{(5)}(\bu)}{\hat{\phi}_{\bx}(\bu)}  - \frac{5 \, \bpsi_{1,4}(\bu) + 10 \, \bpsi_{2,3}(\bu)}{\hat{\phi}_{\bx}^2(\bu)}   \right.\notag \\
& \quad   + 2 \,
 \frac{ 10 \bpsi_{1,1,3}(\bu)   +  15 \bpsi_{1,2,2}(\bu)}{\hat{\phi}_{\bx}^3(\bu)} 
- 6   \frac{15    \bpsi_{1,1,1,2}(\bu)}{\hat{\phi}_{\bx}^4(\bu)}\notag\\
  &\quad+ \left. 24   \frac{ \bpsi_{1,1,1,1,1}(\bu)}{\hat{\phi}_{\bx}^5(\bu)} \right)\, , \notag\displaybreak[3]\\
\tensor{\Psi}_{\bx}^{(6)}(\bu) &= \calS\left( \frac{\bpsi^{(6)}(\bu)}{\hat{\phi}_{\bx}(\bu)}  - \frac{6\,\bpsi_{1,5}(\bu) + 15 \, \bpsi_{2,4}(\bu)  + 10 \, \bpsi_{3,3}(\bu)}{\hat{\phi}_{\bx}^2(\bu)}   \right.\notag  \\
 &\quad  +2 \, \frac{15 \,\bpsi_{1,1,4}(\bu) + 60 \, \bpsi_{1,2,3}(\bu)  + 15 \, \bpsi_{2,2,2}(\bu)}{\hat{\phi}_{\bx}^3(\bu)}
  \notag  \\  &\quad
 - 6 \, \frac{20 \,\bpsi_{1,1,1,3}(\bu) + 45 \, \bpsi_{1,1,2,2}(\bu)}{\hat{\phi}_{\bx}^4(\bu)}  \notag  \\
  & \quad + 24 \, \frac{15 \,\bpsi_{1,1,1,1,2}(\bu)}{\hat{\phi}_{\bx}^5(\bu)}
\left.    - 120 \, \frac{ \bpsi_{1,1,1,1,1,1}}{\hat{\phi}_{\bx}^6(\bu)} \right)   \,,\notag \displaybreak[3]\\
%
\tensor{\Psi}_{\bx}^{(7)}(\bu) &= \calS\left( \frac{\bpsi^{(7)}(\bu)}{\hat{\phi}_{\bx}(\bu)}  - \frac{7\,\bpsi_{1,6}(\bu) + 21 \, \bpsi_{2,5}(\bu)  + 35 \, \bpsi_{3,4}(\bu)}{\hat{\phi}_{\bx}^2(\bu)}   \right.\notag  \\
 &\quad   + 2 \, \frac{21 \,\bpsi_{1,1,5}(\bu) + 105 \, \bpsi_{1,2,4}(\bu)  + 70 \, \bpsi_{1,3,3}(\bu)  + 105 \, \bpsi_{2,2,3}(\bu) }{\hat{\phi}_{\bx}^3(\bu)}
  \notag  \\  &
\quad - 6 \, \frac{35 \,\bpsi_{1,1,1,4}(\bu) + 210 \, \bpsi_{1,1,2,3}(\bu) + 105 \, \bpsi_{1,2,2,2}(\bu)}{\hat{\phi}_{\bx}^4(\bu)}  \notag  \\
&\quad  + 24 \, \frac{35 \,\bpsi_{1,1,1,1,3}(\bu) + 105 \, \bpsi_{1,1,1,2,2}(\bu)}{\hat{\phi}_{\bx}^5(\bu)}  \notag  \\
  &\quad  -120 \, \frac{21 \,\bpsi_{1,1,1,1,1,2}(\bu)}{\hat{\phi}_{\bx}^6(\bu)}
\left.    + 720 \, \frac{ \bpsi_{1,1,1,1,1,1,1}}{\hat{\phi}_{\bx}^7(\bu)} \right)  \notag  \,.
\end{align}
The symmetrization in derivation of $\tensor{\Psi}_{\bx}^{(4)}(\bu)$ can be performed only once, as in the derivatives of orders-5, 6 and 7.

\section{Higher Order Cumulants}\label{apx_cumulants}

When the mixtures are centered to have zero-mean, then the first derivative $\bpsi_{\bx}^{(1)}(\0) = \frac{1}{T} \sum_{t} \bx_t = \0$, which leads to $\bpsi_{\{n_1 = 1, n_2, \ldots, n_k\}}(\0)$ being all zero for arbitrary partitions $\{n_1 = 1, n_2, \ldots, n_k\}$. This simplifies the computation of the cumulant by ignoring the partitions $\{n_1 = 1, n_2, \ldots, n_k\}$. For example, cumulants of order $2, 3, \ldots, 7$ can be computed as
\begin{align*}
	\calK_{\bx}^{(2)} &= \bpsi^{(2)}(\0) = \frac{1}{T} \tI \times_1 \bX  \times_2 \bX  =   \frac{1}{T}  \bX \bX^{\text{T}}   \,,\notag \\
	\calK_{\bx}^{(3)} &= \bpsi^{(3)}(\0) = \frac{1}{T} \tI \times_1 \bX  \times_2 \bX \times_3 \bX \,,  \notag \displaybreak[3]\\
	\calK_{\bx}^{(4)} &= \bpsi^{(4)}(\0)  - 3 \calS(\bpsi^{(2)}(\0)  \circ \bpsi^{(2)}(\0) ) \,, \notag\displaybreak[3] \\
	\calK_{\bx}^{(5)} &= \bpsi^{(5)}(\0)  - 10 \calS(\bpsi^{(2)}(\0)  \circ \bpsi^{(3)}(\0)) \,, \notag\displaybreak[3] \\
	\calK_{\bx}^{(6)} &= \calS\left(\bpsi^{(6)}(\0)  - 15 \bpsi^{(2)}(\0)  \circ \bpsi^{(4)}(\0)  - 10 \bpsi^{(3)}(\0)  \circ \bpsi^{(3)}(\0)  \right. \notag \\ &\quad \left. + 30 \bpsi^{(2)}(\0) \circ \bpsi^{(2)}(\0) \circ \bpsi^{(2)}(\0) \right) \,, \notag \displaybreak[3] \\
	\calK_{\bx}^{(7)} &= \calS\left(\bpsi^{(7)}(\0)  - 21 \bpsi^{(2)}(\0)  \circ \bpsi^{(5)}(\0)  - 35 \bpsi^{(3)}(\0)  \circ \bpsi^{(4)}(\0)  \right. \notag \\ &\displaybreak[3] \left. + 210 \bpsi^{(2)}(\0) \circ \bpsi^{(2)}(\0) \circ \bpsi^{(3)}(\0) \right) \notag  \, .
\end{align*}

\section{Elementary Core Tensor for the Convolution Tensor}\label{sec:elementary_core_tensor}

The elementary core tensor $\tS$ of size  $N \times \underbrace{2 \times 2 \times \cdots \times 2}_{\text{$(N+1)$  dimensions}} \times N$ is constructed from $2N$ $N$th-order sparse tensors $\tS_{1}$, $\tS_{2}$, \ldots, $\tS_{2N}$ of size $2 \times 2 \times \cdots \times 2$, which take ones only at locations $(i_1, i_2, \ldots, i_{N})$ such that
\be
{\bar{i}}_1 + {\bar{i}}_2 + \cdots + {\bar{i}}_{N-1} + i_{N} = \begin{cases}
\max(N-n+1,0)	\quad &\text{odd} \; n,\\
\max(N-n+3,0)  \quad &\text{even} \; n,
\end{cases}
\ee
where $\bar{i}_n = 2 - i_n$.

The structure of the tensor $\tS$ based on $\tS_n$ is explained, as  follows
\be
\tS(1,:,\ldots,:,1,n)  =  \tS_{2n-1}\, , \quad
\tS(1,:,\ldots,:,2,n)  =  \tS_{2n} \, ,
\ee
for $n = 1, \ldots, N$.
The other sub-tensors $\tS(m,:, \ldots, :,n)$, $m = 2, \ldots, N$, $n = 1, \ldots, N$, are recursively defined as
\be
\tS(m,:,\ldots,:,1,n) &=& \tS(m-1,:,\ldots,:,2,n) ,    \\
\tS(m,:,\ldots,:,2,1) &=& \tS(m-1,:,\ldots,:,1, N) ,   \\
\tS(m,:,\ldots,:,2,n) &=& \tS(m-1,:,\ldots,:,1, n-1) .
\ee

According to definition of $\tS_n$, such tensors are orthogonal, i.e., their pair-wise inner products are zeros, and there are only $(N+1)$ such nonzero tensors, which are $\tS_1$,  \ldots, $\tS_N$ and $\tS_{N+2}$ for even $N$,
and $\tS_1$,  \ldots, $\tS_N$ and $\tS_{N+1}$ for odd $N$; the remaining  $(N-1)$  core tensors $\tS_n$ are zero tensors.

\chapter*{Acknowledgements}
\addcontentsline{toc}{chapter}{Acknowledgements}
The authors wish to express their sincere gratitude to the anonymous reviewers for their constructive and helpful comments. We also appreciate the helpful comments of Claudius Hubig (Ludwig-Maximilians- University, Munich)  and Victor Lempitsky (Skolkovo Institute of Science and Technology, Moscow), and help with the artwork of Zhe Sun (Riken BSI). We acknowledge the insightful comments and rigorous proofreading of selected chapters by Anthony Constantinides, Ilia Kisil, Giuseppe Calvi, Sithan Kanna, Wilhelm von Rosenberg, Alex Stott, Thayne Thiannithi, and Bruno Scalzo Dees (Imperial College London).

This work was partially supported by the Ministry of Education and Science of the Russian Federation (grant 14.756.31.0001), and by the EPSRC in the UK (grant EP/P008461).



\bibliographystyle{myplainnat}



\end{document}